\documentstyle[amstex,amssymb,12pt]{article}
  \input xypic
   \xyoption{curve}
 \newcommand{\lon}{\longrightarrow}
 
 \newcommand{\rar}{\rightarrow}
 
 \newcommand{\Proof}{\no{\bf Proof}.\, }
 
 \newcommand{\End}{{\mathrm End}}
 \newcommand{\Q}{{\mathbb Q}}

 \newcommand{\Z}{{\mathbb Z}}
 \newcommand{\K}{{\mathbb K}}
 \newcommand{\bS}{{\mathbb S}}
 \newcommand{\p}{{\partial}}

 \newcommand{\R}{{\mathbb R}}
  \newcommand{\N}{{\mathbb N}}
 \newcommand{\ot}{\otimes}
 
 \newcommand{\Id}{\mbox{Id}}
\newcommand{\Ass}{{\mathsf A\mathsf s\mathsf s}}
  \newcommand{\Liebi}{{\mathsf L\mathsf i\mathsf e^{\hspace{-0.2mm}_1}\hspace{-0.4mm}\mathsf B}}
 \newcommand{\Lie}{{\mathsf L\mathsf i\mathsf e}}
 
  \newcommand{\PROP}{{\mathsf P}}
\newcommand{\DefQ}{{\mathsf D\mathsf e\mathsf f \mathsf Q}}
\newcommand{\DefQh}{{\mathsf D\mathsf e\mathsf f \mathsf Q}^\hbar}

 \newcommand{\sT}{{\mathsf T}}
  \newcommand{\PP}{{\mathsf P^\circlearrowright}\hspace{-0.5mm}}
  \newcommand{\PPP}{{\mathsf P^+}\hspace{-0.5mm}}
   \newcommand{\LLL}{{\mathsf L\hspace{-0.6mm}^+}\hspace{-0.5mm}}
  
  \newcommand{\Graph}{\PROP^\circlearrowright}
  \newcommand{\CoLie}{{\mathsf C\mathsf o\mathsf L\mathsf i\mathsf e}}
  \newcommand{\EndM}{{{\sf End}\langle M\rangle}}
  
 \newcommand{\Prop}{\mathsf P\mathsf r\mathsf o\mathsf p}
 \newcommand{\Diop}{\mathsf D\mathsf i\mathsf o\mathsf p}
\newcommand{\no}{{\noindent}}
\newcommand{\sgn}{{\rm sgn}}
\newcommand{\sA}{{\mathsf A}}

\newcommand{\sC}{{\mathsf C}}
\newcommand{\sD}{{\mathsf D}}
\newcommand{\sE}{{\mathsf E}}
\newcommand{\sH}{{\mathsf H}}

\newcommand{\sF}{{\mathsf F}}
\newcommand{\sP}{{\mathsf P}}
\newcommand{\sQ}{{\mathsf Q}}

 \newcommand{\Beq}{\begin{equation}}
 \newcommand{\Eeq}{\end{equation}}
 \newcommand{\Beqr}{\begin{eqnarray}}
 \newcommand{\Eeqr}{\end{eqnarray}}
 \newcommand{\Beqrn}{\begin{eqnarray*}}
 \newcommand{\Eeqrn}{\end{eqnarray*}}
 \newcommand{\Ba}{\begin{array}}
 \newcommand{\Ea}{\end{array}}
 \newcommand{\Bi}{\begin{itemize}}
 \newcommand{\Ei}{\end{itemize}}
 \newcommand{\Bc}{\begin{center}}
 \newcommand{\Ec}{\end{center}}
 \newcommand{\fG}{{\mathfrak G}}
 \newcommand{\fGG}{{\mathfrak G^\circlearrowright\hspace{-0.5mm}}}
 \newcommand{\fGGG}{{\mathfrak G^+\hspace{-0.5mm}}}
 \newcommand{\f}{{\mathcal O}}

 \newcommand{\caD}{{\mathcal D}}
 \newcommand{\cE}{{\mathcal E}}
 \newcommand{\cF}{{\mathcal F}}
 \newcommand{\cG}{{\mathcal G}}
 
 \newcommand{\cI}{{\mathcal I}}

 \newcommand{\cM}{{\mathcal M}}

 \newcommand{\cU}{{\mathcal U}}

 \newcommand{\PBW}{{\mathcal P}{\mathcal B}{\mathcal W}}

 \newcommand{\al}{\alpha}
 \newcommand{\be}{\beta}
 \newcommand{\ga}{\gamma}
 
 \newcommand{\Ga}{\Gamma}
 
 \newcommand{\var}{\varepsilon}

 \newcommand{\Ker}{{\mathsf K\mathsf e\mathsf r}\, }

 \newcommand{\Hom}{{\mathrm H\mathrm o\mathrm m}}

 \newcommand{\cT}{{\cal T}}

 \newcommand{\sip}{\smallskip}
 \newcommand{\bip}{\bigskip}

\begin{document}

 \sloppy

 \title{ PROP profile of  deformation quantization and\\ graph complexes with
  loops and wheels\footnote{This work was
 partially supported by the G\"oran Gustafsson foundation.}}
 \author{ S.A.\ Merkulov}
 \date{}
 \maketitle

 %%%%%%%%%%%%%%%%%%%%%%%%%%%%%%%%%%%%%
\begin{center}
{\bf \S 1. Introduction}
\end{center}

\bip

\no
The first instances of graph complexes have been introduced in the
theory of operads and props which have found recently lots of
applications in algebra, topology and geometry. Another set of
examples has been introduced by Kontsevich \cite{Kon} as a way
to expose highly non-trivial interrelations between certain infinite
dimensional  Lie algebras and  topological objects, including moduli
spaces of curves, invariants of odd dimensional manifolds, and the
group of outer automorphisms of a free group.

\sip

Motivated by the problem of deformation quantization we introduce
and study
 directed graph complexes with oriented loops and wheels.
 We show that universal quantizations of Poisson structures can be
 understood as morphisms of dg  props,
 $$
Q: \DefQ \lon {\Liebi_\infty^\circlearrowright},
$$
 where

-- $\DefQ$ is the dg free prop
 whose representations in a graded vector space $V$ describe Maurer-Cartan elements in the
Hochschild dg Lie algebra, ${\mathcal D}_V$, of polydifferential operators
 on the ring, $\f_V:=\widehat{\odot^\bullet}V^*$, of smooth formal functions on $V$
 (see \S 2.7 for a precise definition); we call such Maurer-Cartan elements {\em star
products};  if $V$ is $\R^n$ concentrated in degree zero, then this notion coincides with the
ordinary notion of star product on smooth formal functions on $\R^n$.

\sip

-- $\Liebi_\infty^\circlearrowright$ is the {\em wheeled}\,
completion of the minimal resolution, $\Liebi_\infty$, of the prop,
$\Liebi$, of Lie 1-bialgebras; it is defined explicitly in \S 2.6 and is
proven to have the property that its representations in a {\em
finite-dimensional}\, graded vector space $V$ correspond to Maurer-Cartan elements in the Schouten Lie algebra,
$\wedge^\bullet \cT_V$, of polivector fields, where $\cT_V:=\mbox{Der} \f_V$; such Maurer-Cartan elements are called {\em Poisson structures}
\, on the formal graded manifold $V$; if $V$ is $\R^n$ concentrated in degree zero, then this notion coincides with the
ordinary notion of Poisson structure on $\R^n$.

% the symbol $\widehat{\Liebi^\circlearrowright_\infty}$ stands
%completion of $\Liebi^\circlearrowright_\infty$ with respect to the
%number of vertices.

\sip

 In the theory of props one is most interested in those directed
graph complexes which contain {\em no}\, loops and wheels. A major advance
in understanding the cohomology groups of such complexes was
recently accomplished in \cite{Ko2,MV, V} using key ideas of
$\frac{1}{2}$prop and Koszul duality. In particular, these authors
were able to compute cohomologies of directed versions (without
loops and wheels though)
%, $(\mathsf IB_\infty,\delta)$ and  $(\mathsf LieB_\infty,\delta)$,
of Kontsevich's ribbon graph complex and  ``commutative" graph
complex, and show that they both are acyclic almost everywhere. One of our purposes
in this
paper is  to extend some of  the results of
\cite{Ko2,MV,V} to a more difficult situation when the directed
graphs are allowed to contain loops and wheels (i.e.\ directed paths
which begin and end at the same vertex). In this case the answer
differs markedly from the unwheeled case: we prove, for example,
that while the cohomology of the wheeled extension,
$\Lie_\infty^\circlearrowright$, of the operad of
$\Lie_\infty$-algebras remains acyclic almost everywhere (see Theorem~4.1.1
for a precise formula for
$\sH^\bullet(\Lie_\infty^\circlearrowright)$), the cohomology of the
wheeled extension of the operad $\mathsf A\mathsf s\mathsf s_\infty$ gets more complicated.
Both these complexes describe irreducible summands of directed
``commutative" and, respectively, ribbon graph complexes
 with the
restriction on absence of wheels dropped.

\sip

 The wheeled complex
 $\Lie_\infty^\circlearrowright$ is a subcomplex of the above
 mentioned graph complex
 $\Liebi_\infty^\circlearrowright$ which describes Poisson structures\footnote{Strictly
speaking, a representation of $\Liebi_\infty$ in $V$ gives a degree $1$ element $\ga\in \wedge^\bullet \cT_V$
which satisfies not only the Maurer-Cartan condition, $[\ga,\ga]=0$, but also the vanishing condition,
$\ga|_{x=0}=0$, at $0\in V$ (see \S 2.6); given, however, $\ga\in \wedge^\bullet \cT_V$ with $\ga|_{x=0}\neq 0$, then, for
a formal parameter $h$ viewed as a coordinate on $\R$, the element $\bar{\ga}:=h\ga\in \wedge^\bullet
\cT_{\bar{V}}$, $\hat{V}:=V\times \R$, satisfies $\bar{\ga}|_{x=0}=0$ and hence comes from a representation
of the prop $\Liebi_\infty$; thus such representations encompass {\em arbitrary}\, formal Poisson structures.}.
Using Theorem 4.1.1
 on $\sH^\bullet(\Lie_\infty^\circlearrowright)$
 we show in \S 4.2 that a subcomplex of
   $\Liebi_\infty^\circlearrowright$ which is
 spanned by graphs with at most genus 1 wheeles is also acyclic almost
 everywhere. However this acyclicity breaks for graphs with higher
 genus wheels: we find an explicit cohomology class with 3 wheels
 in \S 4.2.4 which proves that the natural epimorphism, $\Liebi_\infty\rar
 \Liebi$, {\em fails}\, to stay quasi-isomorphism when extended to the
 wheeled completions, $\Liebi_\infty^\circlearrowright\rar
 \Liebi^\circlearrowright$.

\sip

Kontsevich's famous universal deformation quantization formulae
involve graphs with { wheels}
which encode {\em traces}\, of tensor powers of polyvector fields
and their partial derivatives and hence make, in general,
 sense only in finite-dimensions. We show in this paper that {\em every
 universal deformation quantization of  Poisson structures on graded manifolds
 {\em must}\, involve such
traces}, i.e.\ graphs with wheels are unavoidable. This fact is one of the motivation behind
our study of {\em wheeled}\, extensions of props, especially  the wheeled extensions, $\Liebi^\circlearrowright$ and
$\Liebi^\circlearrowright_\infty$, of the props related to Poisson
geometry.

\sip

We show that there exists a natural dg free
prop $[\Liebi^\circlearrowright]_\infty$ which extends the above prop,
 $\Liebi_\infty^\circlearrowright$, of polyvector fields and fits into a
commutative diagram
\[
 \xymatrix{
[\Liebi^\circlearrowright]_\infty \ar[dr]^{qis} \ar[r]^{\al} &
\Liebi^\circlearrowright_\infty
  \ar[d]%^{\pi^\circlearrowright}
  \\
 &
 \Liebi^\circlearrowright
 }
\]
where $\al$ is an epimorphism of nondifferential props, and
$qis$ a quasi-isomorphism of differential props. It is called a {\em
quasi-minimal prop resolution}\, of $\Liebi^\circlearrowright$.

\sip

 As we mentioned above, representations of $\Liebi_\infty^\circlearrowright$ in a finite-dimensional dg space
$V$ are precisely  Poisson structures on $V$.
%, i.e.\ Maurer-Cartan elements in the
%Lie algebra, $\wedge^\bullet \cT_V$, of polyvector fields.
This prompts us to call
representations of $[\Liebi^\circlearrowright]_\infty$ in a finite dimensional dg space $V$
{\em wheeled Poisson structures}\, on the formal graded manifold $V$.
Geometric meaning of such wheeled Poisson structures
 is not clear at present
--- the differential equations behind these new structures involve
not only Schouten brackets but also {\em traces}\, of derivatives of
tensor powers of polyvector fields\footnote{It is proven in \S 4.2.2 that this system  of differential equations can {\em not}\,
contain equations involving  just a {\em single}\, trace of some tensorial expression built from the Poisson tensor and its derivatives (which correspond to graphs
with genus {one} wheels).}
 (and hence makes sense only in
finite dimensions). Moreover, polyvector fields are only
part of the data --- there are other  tensors (including new ones in degree 0) in the content list of a
wheeled Poisson structure.
Remarkably, one can construct ``star products" out of
wheeled Poisson structures, i.e.\ they can be deformation quantized:

\sip

\no {\bf Main Theorem.} {\em There exists a morphism of dg props,
$$
\hat{Q}: \DefQ \rar {[\Liebi^\circlearrowright]}_\infty.
$$
%where $\widehat{[\Liebi^\circlearrowright]}_\infty$ is the
%completion of the prop ${[\Liebi^\circlearrowright]}_\infty$ with respect to
%the number of vertices.
Moreover, this morphism exists in the
category of props over the field, $\Q$, of rational numbers.}

\bip

\no {\bf Corollary.} {\em Every wheeled Poisson structure on a finite
dimensional formal manifold can be deformation quantized, i.e.\ there exists
an associated Maurer-Cartan element --- ``star product" ---  in the
Hochschild  dg Lie algebra ${\caD}_V[[\hbar]]$.}

\bip

\no The quantization morphism $\hat{Q}$ is very non-trivial: the proof of Theorem~\S
5.2 below implies that $\hat{Q}$ involves, e.g., {\em infinite}\, jets of
the polyvector fields constituent of a
wheeled Poisson structure.
There is a canonical monomorphism of dg props, $i: \Liebi^\circlearrowright_\infty\rar  {[\Liebi^\circlearrowright]}_\infty$,
so that those quantization morphisms $\hat{Q}$ which factor though $i$ provide us with quantizations of ordinary
Poisson structures; we shall discuss a purely propic construction of such quantization morphisms elsewhere.
%Every Kontsevich
%type quantization, $Q$, of ordinary Poisson structures
%$\Liebi_\infty^\circlearrowright$-structures gives canonically rise
%to the associated quantization, $\hat{Q}$, of
%${[\Liebi^\circlearrowright]}_\infty$-structures, but not vice
%versa.

\sip

One can get some intuition into the geometric meaning of
wheeled Poisson structures  from their simpler analogues
--- representations of quasi-minimal prop resolutions,
$[\Ass^\circlearrowright]_\infty$ and $[{\mathsf C\mathsf o\mathsf
m}^\circlearrowright]_\infty$, of the operads of associative and,
respectively, commutative algebras which have been computed in
\cite{MMS}.  We show some details in \S 4.5.2:  if an $\Ass_\infty^\circlearrowright$-structure in a finite-dimensional
dg space $V$ is just a homological vector field, $\eth$,  on $V[1]$ viewed as a non-commutative formal manifold,
then an   $[\Ass^\circlearrowright]_\infty$-structure  in $V$ is a pair $(\eth, f)$, involving a new
piece of data --- a cyclically invariant function $f$ on $V[1]$ --- all satisfying a system of differential equations
involving trace
of $\eth$.
%all the new
%fields in $[\Ass^\circlearrowright]_\infty$-structure in comparison with the well-known
%$\Ass_\infty$-structure as well as all the associated new (involving trace) equations
%can be explicitly described.
Unfortunately, such a detailed picture of wheeled Poisson structures
is out of reach at present: this problem appears to have complexity level comparable with that
of the problem
of computing homologies of the directed versions of famous Kontsevich's graph complexes \cite{Kon} (in fact, the graph complex behind
wheeled Poisson structures is closely related to the directed version of Kontsevich's ``commutative" graph complex) .

\sip

 A few words about our notations. The cardinality of a finite set
$I$ is denoted by $|I|$. The degree of a homogeneous element, $a$,
of a graded vector space is denoted by $|a|$ (this should never lead
to a confusion). $\bS_n$ stands for the group of all bijections,
$[n]\rar [n]$, where $[n]$ denotes (here and everywhere) the set
$\{1,2,\ldots,n\}$. The set of positive integers is denoted by
$\N^*$. If $V=\oplus_{i\in \Z} V^i$ is a graded vector space, then
$V[k]$ is a graded vector space with $V[k]^i:=V^{i+k}$. All our free props
are assumed to be completed with respect to the number of vertices (so that extra care
should be in place when composing prop morphisms).

\sip

We work throughout over the field $\K$ of characteristic 0.

\sip

The paper is organized as follows. In \S 2 we remind some
basic facts about props and graph complexes and describe a universal
construction which associates dg props to a class of sheaves of dg
Lie algebras on smooth formal manifolds; we then illustrate it with examples
which are relevant to Poisson geometry and deformation quantization.
 In \S 3 we develop new methods for computing
cohomology of directed graph complexes with wheels, and prove
several theorems on cohomology of wheeled completions of minimal
resolutions of dioperads. In \S 4 we apply these methods and results to compute cohomology
of several concrete graph complexes.
In \S 5 we  prove the main theorem formulated above.
In the appendix  we use ideas of cyclic homology to construct a cyclic bicomplex
computing cohomology of wheeled completions of dg
operads.

\bip

\bip
%%%%%%%%%%%%%%%%%%%%%%%%%%%%%%%%%%%%%%%%%%%%%%%%%%%%%%%%%
%%%%%%%%%%%%%%%%%%%%%%%%%%%%%%%%%%%%%%%%%%%%%%%%%%%%%%%%%
\begin{center}
\bf \S 2. Dg props versus sheaves of dg Lie algebras
\end{center}

\bip

%%%%%%%%%%%%%%%%%%%%%%%%%%%%%%%%%%%%%%%%%%%%%%%%%%%%%%%%%

\no{\bf 2.1. Props.} An $\bS$-{\em bimodule}, $E$, is, by definition, a collection
of graded vector spaces, $\{E(m,n)\}_{m,n\geq 0}$, equipped with a left
action of the group $\bS_m$  and with a right action of $\bS_n$ which commute
with each other. For any graded vector space $M$ the collection,
${\sf End}\langle M\rangle
=\{{\sf End}\langle M\rangle(m,n):= \Hom(M^{\ot n}, M^{\ot m})\}_{m,n\geq 0}$,
is naturally an $\bS$-bimodule. A {\em morphism}\,
of $\bS$-bimodules, $\phi: E_1\rar E_2$, is a collection of equivariant
linear maps, $\{\phi(m,n):E_1(m,n)\rar E_2(m,n)\}_{m,n\geq 0}$.
A morphism $\phi: E\rar {\sf End}\langle M\rangle$ is called a
{\em representation}\, of an $\bS$-bimodule $E$ in a graded vector space $M$.

\sip

\no
There are two natural
 associative binary operations on the $\bS$-bimodule $\EndM$,
$$
{\bigotimes}: \EndM(m_1,n_1)\ot \EndM(m_2,n_2) \lon \EndM(m_1+m_2, n_1+n_2),
$$
$$
\circ: \EndM(p,m)\ot \EndM(m,n) \lon \EndM(p,n),
$$
and a distinguished element, the identity map ${\bf 1}\in \EndM(1,1)$.

\sip

\no
Axioms of prop (``{\em pro}ducts and {\em p}ermutations") are modelled  on the properties
of $(\bigotimes,\circ, {\bf 1})$ in $\EndM$ (see \cite{Mc}).

\bip
%%%%%%%%%%%%%%%%%%%%%%%%%%%%%%%%%%%%%%%%
\no{\bf 2.1.1. Definition.}  A {\em prop}, $E$, is an
$\bS$-bimodule, $E=\{E(m,n)\}_{m,n\geq 0}$, equipped
with the following data,
\Bi
\item a linear map called {\em horizontal composition},
$$
\Ba{rccc}
{\bigotimes}:& E(m_1,n_1)\ot E(m_2,n_2) & \lon & E(m_1+m_2, n_1+n_2) \\
              &{\mathfrak e}_1\ot  {\mathfrak e}_2 & \lon &{\mathfrak e}_1\bigotimes
                {\mathfrak e}_2
 \Ea
$$
such that
$({\mathfrak e}_1\bigotimes{\mathfrak e}_2)\bigotimes {\mathfrak e}_3=
{\mathfrak e}_1\bigotimes({\mathfrak e}_2\bigotimes {\mathfrak e}_3)$ and
 ${\mathfrak e}_1\bigotimes{\mathfrak e}_2=$ $(-1)^{|{\mathfrak e}_1||{\mathfrak e}_1|}
 \sigma_{m_1,m_2}(
 {\mathfrak e}_2\bigotimes{\mathfrak e}_1)\sigma_{n_2,n_1}$ where
 $\sigma_{m_1,m_2}$ is the following permutation in $\bS_{m_1+m_2}$,

 $$
 \left(\Ba{ccccccccccc}
 1 &,&\ldots&,&m_2&,& m_2+1&,&\ldots&,& m_2+m_1 \\
 1+m_1&,&\ldots&,& m_2+m_1&,&1&,&\ldots&,& m_1
 \Ea\right);
 $$
\item a linear map called {\em vertical composition},
$$
\Ba{rccc}
{\circ}:& E(p,m)\ot E(m,n) & \lon & E(p,n) \\
              &{\mathfrak e}_1\ot  {\mathfrak e}_2 & \lon &{\mathfrak e}_1\circ
                {\mathfrak e}_2
 \Ea
$$
such that $({\mathfrak e}_1\circ{\mathfrak e}_2)\circ {\mathfrak e}_3=
{\mathfrak e}_1\circ({\mathfrak e}_2\circ {\mathfrak e}_3)$ whenever both sides
 are defined;
\item an algebra morphism, $i_n: k[\bS_n] \rar (E(n,n),\circ)$, such
that (i) for any $\sigma_1\in \bS_{n_1}$, $\sigma_2\in \bS_{n_2}$
one has
$i_{n_1+n_2}(\sigma_1\times \sigma_2)=i_{n_1}(\sigma_1)\bigotimes
i_{n_2}(\sigma_2)$, and (ii) for any ${\mathfrak e}\in E(m,n)$ one has
${\bf 1}^{\otimes m}\circ{\mathfrak e}={\mathfrak e}\circ {\bf 1}^{\otimes n}
={\mathfrak e}$ where ${\bf 1}:= i_1(\Id)$.
\Ei

\no
A {\em morphism}\, of props, $\phi: E_1\rar E_2$, is a morphism of the associated
$\bS$-bimodules which respects, in the obvious sense, all the prop data.

\bip

\no
A {\em differential}\, in a prop $E$ is a collection of degree 1 linear
maps, $\{\delta: E(m,n)\rar E(m,n)\}_{m,n\geq 0}$, such that $\delta^2=0$ and
\Beqrn
\delta({\mathfrak e}_1\bigotimes{\mathfrak e}_2)&=& (\delta
{\mathfrak e}_1)\bigotimes{\mathfrak e}_2
+ (-1)^{|{\mathfrak e}_1|}{\mathfrak e}_1\bigotimes \delta{\mathfrak e}_2,\\
\delta({\mathfrak e}_3\circ{\mathfrak e}_4)&=& (\delta{\mathfrak e}_3)\circ{\mathfrak e}_4
+ (-1)^{|{\mathfrak e}_3|}{\mathfrak e}_3\circ \delta{\mathfrak e}_4,
\Eeqrn
for any ${\mathfrak e}_1,{\mathfrak e}_2\in E$ and any
${\mathfrak e}_3,{\mathfrak e}_4\in E$ such that ${\mathfrak e}_3\circ{\mathfrak e}_4$
makes sense. Note that $d{\bf 1}=0$.

\bip

\no
For any dg vector space $(M, d)$ the associated prop $ \EndM$
has a canonically induced differential which we always denote by the same
symbol $d$.

\bip

\no
A {\em representation} of a dg prop $(E, \delta)$ in a dg vector space
$(M,d)$ is, by
definition, a morphism of props, $\phi: E \rar \EndM$, which commutes
with differentials, $\phi \circ\delta = d\circ \phi$. (Here and
elsewhere $\circ$ stands for the composition of maps; it will always be clear from
the context whether $\circ$ stands for the composition of maps or for
the vertical composition in props.)

\bip

%%%%%%%%%%%%%%%%%%%%%%%%%%%%%%%%
\no
{\bf 2.1.2. Remark.}
If $\psi: (E_1,\delta) \rar (E_2,\delta)$ is a morphism of dg props,
and $\phi: (E_2,\delta) \rar (\EndM, d)$ is a representation of $E_2$,
then the composition, $\phi\circ \psi$, is a representation of $E_1$.
Thus representations  can be ``pulled back".

\bip
%%%%%%%%%%%%%%%%%%%%%%%%%%%%%%%%%%%%%%%
\no{\bf 2.1.3. Free props.}
Let ${\fG}^\uparrow(m,n)$,  $m,n\geq 0$,  be the space of
{\em directed  $(m,n)$-graphs}, $G$, that is, connected $1$-dimensional
$CW$ complexes satisfying the following conditions:
\Bi
\item[(i)] each edge (that is, 1-dimensional cell) is equipped with a
direction;
\item[(ii)] if we split the set of all vertices (that is,
0-dimensional cells) which have exactly one adjacent edge into a disjoint
union, $V_{in}\sqcup V_{out}$,
\Bi
\item[] with $V_{in}$ being the subset of vertices with the adjacent edge directed
from the vertex,
\item[] and $V_{out}$ the subset of vertices with the adjacent  edge directed
towards the vertex,
\Ei
then $|V_{in}|\geq n$ and $|V_{out}|\geq m$;
\item[(iii)] precisely $n$ of vertices from $V_{in}$ are labelled by
$\{1,\ldots,n\}$ and are called %(together with the attached edges)
{\em inputs};
\item[(iv)] precisely $m$ of vertices from $V_{out}$ are labelled by
$\{1,\ldots,m\}$ and are called {\em outputs};
\item[(v)] there are no oriented {\em wheels}, i.e.\
directed paths which begin and end at the same vertex; in particular, there
are no {\em loops} (oriented wheels consisting of one internal edge).
 Put another way,
directed edges generate a continuous flow on the graph which we always
assume in our pictures to go from bottom to the top.
\Ei

\bip

\no
Note that $G\in \fG^\uparrow(m,n)$ may not be connected. Vertices in the complement,
$$
v(G):=\overline{inputs\sqcup outputs},
$$
 are called {\em internal vertices}. For each internal vertex $v$ we denote
 by $In(v)$ (resp., by $Out(v)$)  the set of those adjacent half-edges whose
 orientation is directed towards (resp., from) the vertex.
 Input (resp., output) vertices together with
 adjacent edges are called {\em input} (resp., {\em output}) {\em legs}.
 The graph with one internal vertex, $n$ input legs and $m$ output legs
 is called the $(m,n)$-{\em corolla}.

\sip

 \no We set $\fG^\uparrow:=\sqcup_{m,n} \fG^\uparrow(m,n)$.

\sip

\no The {\em free}\, prop, $\PROP \langle E\rangle$, generated by an $\bS$-module,
 $E=\{E(m,n)\}_{m,n\geq 0}$, is defined by (see, e.g., \cite{MV,V})
$$
\PROP \langle E\rangle (m,n):= \bigoplus_{G\in \fG^\uparrow(m,n)} \left(
\bigotimes_{v\in v(G)} E(Out(v), In(v))\right)_{Aut G}
$$
where
\Bi
\item $
 E(Out(v), In(v)):= {\rm Bij}([m],Out(v))\times_{\bS_m} E(m,n)
 \times_{\bS_n} {\rm Bij}(In(v), [n])$
  with $\rm Bij$  standing for the set of bijections,
 \item $Aut(G)$ stands for
the automorphism group of the graph $G$.
\Ei

\noindent An element of the summand above, $G \langle E\rangle :=\left(
\bigotimes_{v\in v(G)} E(Out(v), In(v))\right)_{Aut G}$, is often called
a {\em graph $G$ with internal vertices decorated by elements of $E$},
or just a {\em decorated graph}.

\bip

\no
A differential, $\delta$, in a free prop
  $\PROP\langle E\rangle$ is completely determined
by its values,
$$
\delta:  E(Out(v), In(v)) \lon
\PROP\langle E \rangle(|Out(v)|,|In(v)|),
$$
on decorated corollas (whose unique internal vertex is
denoted by $v$).

\bip

\no
Prop structure on an $\bS$-bimodule $E=\{E(m,n)\}_{m,n\geq 0}$ provides us, for any graph
$G\in \fG^\uparrow(m,n)$, with a well-defined {\em evaluation}\, morphism of
$\bS$-bimodules,
$$
{\sf ev}: G \langle E\rangle \lon E(m,n).
$$
In particular, if a decorated graph $C\in \PROP\langle E\rangle$ is built from
two corollas, $C_1\in \fG(m_1,n_1)$ and  $C_2\in \fG(m_2,n_2)$
by gluing $j$th output leg of $C_2$ with $i$th input leg
of $C_1$, and if the vertices of these corollas are decorated, respectively,
by elements  $a\in E(m_1,n_1)$
and $b\in E(m_2,n_2)$, then we reserve a special notation,
$$
a\hspace{1mm}_i\hspace{-0.2mm}\circ_j b:= {\mathsf ev}(C)\in
E(m_1+m_2-1,n_1+n_2-1),
$$
for the resulting evaluation map.

\bip

%%%%%%%%%%%%%%%%%%%%%%%%%%%%%%%%%%%%%%%%%%
\no{\bf 2.1.4. Completions}. Any free prop $\PROP\langle E\rangle$ is
naturally a direct sum,
 $\PROP\langle E\rangle = \oplus_{n\geq 0}\PROP_n\langle E\rangle$, of
 subspaces spanned by decorated graphs with $n$ vertices.
Each summand $\PROP_n\langle E\rangle$ has a natural filtration by
the genus, $g$, of the underlying graphs (which is, by definition,
equal to the first Betti number of the associated CW complex). Hence
each $\PROP_n\langle E\rangle$ can be completed with respect to this
filtration. Similarly, there is a filtration by the number of
vertices. We shall  always work in this paper with completed with
respect to these filtrations free props and hence use the same
notation,  $\PROP\langle E\rangle$, and the same name, {\em free
prop}, for the completed version. Note that {\em not}\ every pair of morphisms
of props, $f: P\rar Q$ and $g:Q\rar R$, can be composed into $g\circ f: P\rar R$.
In concrete examples one must be careful to check that no divergences occur.

\bip

%%%%%%%%%%%%%%%%%%%%%%%%%%%%%%%%%%%%%
\no{\bf 2.2. Dioperads and $\frac{1}{2}$props}. A {\em dioperad}\, is an
$\bS$-bimodule, $E=\{E(m,n)\}_{m,n\geq 1\atop m+n\geq 3}$, equipped with a set
of compositions,
$$
\left\{
\hspace{1mm}_i\hspace{-0.2mm}\circ_j: E(m_1,n_1)\ot E(m_2,n_2)\lon E(m_1+m_2-1, n_1+n_2-1)
\right\}_{1\leq i \leq n_1 \atop 1\leq j\leq m_2},
$$
which satisfy the axioms imitating the  properties of the compositions
$\hspace{1mm}_i\hspace{-0.2mm}\circ_j$ in a generic prop. We refer to \cite{G}, where this
notion was introduced,
for a detailed list
of these axioms. The {\em free}\,  dioperad generated by an $\bS$-bimodule $E$ is given by,
$$
{\mathsf D}\langle E \rangle(m,n) := \bigoplus_{G\in {\mathfrak T}(m,n)} G\langle E \rangle
$$
where ${\mathfrak T}(m,n)$ is a subset of $\fG(m,n)$ consisting of connected trees (i.e.,
connected graphs of genus 0).

\bip

\no
Another and less obvious reduction of the notion of prop
was introduced by Kontsevich in \cite{Ko2} and studied in detail in \cite{MV}:
a {\em $\frac{1}{2}$prop}\, is an $\bS$-bimodule,
$E=\{E(m,n)\}_{m,n\geq 1\atop m+n\geq 3}$, equipped with two sets of compositions,
$$
\left\{
\hspace{1mm}_1\hspace{-0.2mm}\circ_j: E(m_1,1)\ot E(m_2,n_2)\lon E(m_1+m_2-1,n_2)
\right\}_{ 1\leq j\leq m_2}
$$
and
$$
\left\{
\hspace{1mm}_i\hspace{-0.2mm}\circ_1: E(m_1,n_1)\ot E(1,n_2)\lon E(m_1+m_2-1,n_2)
\right\}_{1\leq i \leq n_1}
$$
satisfying the axioms which imitate the  properties of the compositions
$\hspace{1mm}_1\hspace{-0.2mm}\circ_j$ and $\hspace{1mm}_i\hspace{-0.2mm}
\circ_1$ in a generic dioperad.  The {\em free}\, $\frac{1}{2}$prop generated
by an $\bS$-bimodule $E$ is given by,
$$
\frac{1}{2}{\mathsf P}\langle E \rangle(m,n) :=
\bigoplus_{G\in \frac{1}{2}{\mathfrak T}(m,n)} G\langle E \rangle,
$$
where $\frac{1}{2}{\mathfrak T}(m,n)$ is a subset of ${\mathfrak T}(m,n)$
consisting of those directed trees which, for each pair of internal vertices,
$(v_1,v_2)$, connected by an  edge directed from $v_1$ to $v_2$  have either
$|Out(v_1)|=1$ or/and $|In(v_2)|=1$. Such trees have
at most one vertex $v$ with $|Out(v)|\geq 2$ and
$|In(v)|\geq 2$.
\bip

\no
Axioms of dioperad (resp., $\frac{1}{2}$prop) structure on an $\bS$-bimodule $E$ ensure
that there is a well-defined evaluation map,
$$
{\sf ev}: G \langle E\rangle \lon E(m,n),
$$
for each $G\in {\mathfrak T}(m,n)$ (resp., $G\in \frac{1}{2}{\mathfrak T}(m,n)$).

\bip

%%%%%%%%%%%%%%%%%%%%%%%%%%%%%%%%%%%%%%%%%%%%%%
\no{\bf 2.2.1. Free resolutions}. A {\em free resolution}\, of a dg prop
$P$ is, by definition, a dg free prop, $(\PROP\langle E \rangle, \delta)$,
%(resp., $({\mathsf D}\langle E \rangle, \delta$, $\frac{1}{2}{\mathsf P}\langle E$)
generated by some $\bS$-bimodule $E$ together with a morphism of dg props,
$\al: (\PROP\langle E \rangle, \delta) \rar P$, which is a homology isomorphism.

\sip

\no If the differential $\delta$ in
$\PROP\langle E \rangle$ is decomposable (with respect to prop's vertical and
/or horizontal compositions), then $\al: (\PROP\langle E \rangle, \delta) \rar P$
is called a {\em minimal model}\, of $P$.

\bip

\no
Similarly one defines free resolutions and minimal models,
$({\mathsf D}\langle E \rangle, \delta) \rar P$ and
$(\frac{1}{2}\PROP\langle E \rangle, \delta) \rar P$, of
dioperads and  $\frac{1}{2}$props .

\bip
%%%%%%%%%%%%%%%%%%%%%%%%%%%%%%%%%%%%%%%%%%%%%%
\no{\bf 2.2.2. Forgetful functors and their adjoints}. There is an
obvious chain of forgetful functors, $\Prop \lon \Diop \lon
\frac{1}{2} \Prop$. Let
$$
\Omega_{\frac{1}{2}\sP\rar \sD}: \frac{1}{2}\Prop \lon \Diop, \ \ \
\ \Omega_{\sD\rar \sP}: \Diop \lon \Prop, \ \ \ \ \Omega_{
\frac{1}{2}\sP\rar \sP}: \frac{1}{2} \Prop \lon \Prop,
$$
be the associated adjoints. The main motivation behind introducing
the notion of  $\frac{1}{2}$prop is a very useful fact that the
functor $\Omega_{\frac{1}{2}\sP\rar \sP}$ is exact \cite{Ko2,MV},
i.e., it commutes with the cohomology functor. Which in turn is due
to the fact that, for any
 $\frac{1}{2}$prop $P$, there exists a kind of PBW lemma which represents
 $\Omega_{ \frac{1}{2}\sP\rar \sP}\langle P\rangle$ as a vector
 space {\em freely}\, generated by a
 family decorated graphs,
 $$
\Omega_{\frac{1}{2}\sP\rar \sP}\langle P \rangle(m,n):=
\bigoplus_{G\in \overline{\fG}(m,n)} G\langle E \rangle,
$$
where $\overline{\fG}(m,n)$ is a subset of $\fG(m,n)$ consisting of
so called {\em reduced} \, graphs, $G$, which satisfy the following
defining property \cite{MV}: for each pair of internal vertices,
$(v_1,v_2)$, of $G$ which are connected by a single edge directed
from $v_1$ to $v_2$ one has $|Out(v_1)|\geq 2$ and $|In(v_2)|\geq
2$. The prop structure on $\Omega_{\frac{1}{2}\sP\rar \sP}\langle P
\rangle$ is given by \Bi
\item[(i)] horizontal compositions $:=$ disjoint unions of decorated graphs,
\item[(ii)] vertical compositions $:=$ graftings followed by  $\frac{1}{2}$prop
compositions of all  those pairs of vertices $(v_1,v_2)$ which are connected by a single
edge directed from $v_1$ to $v_2$ and
have either $|Out(v_1)|=1$ or/and
$|In(v_2)|=1$
(if there are any).
\Ei

\bip

%%%%%%%%%%%%%%%%%%%%%%%%%%%%%%%%%%%%%%%%%%%
%%%%%%%%%%%%%%%%%%%%
\no{\bf 2.3. Graph complexes with wheels}.
 Let $\fG^\circlearrowright(m,n)$ be the set of all directed $(m,n)$-graphs which satisfy
 conditions 2.1.3(i)-(iv), and set  $\fG^\circlearrowright:=\sqcup_{m,n} \fG(m,n)$. A vertex
 (resp., edge or half-edge) of a graph $G\in \fG^\circlearrowright$ which belongs to an oriented wheel
is called a {\em cyclic}\, vertex (resp., edge or half-edge).
%Half-edges attached to cyclic vertices (including those that
%do not form a part of a cyclic edge) are called
%{\em cyclic half-edges}.
\bip

\no
Note that for each internal vertex of $G\in \fG^\circlearrowright(m,n)$
there is still
 a well defined separation of adjacent half-edges into input and output
ones, as well as a well defined separation of legs into input and
output ones.

\bip

\no
For any $\bS$-bimodule $E=\{E(m,n)\}_{m,n\geq 0}$, we define an $\bS$-bimodule,
$$
\PROP^\circlearrowright\langle E\rangle(m,n):= \bigoplus_{G\in \fG(m,n)} \left(
\bigotimes_{v\in v(G)} E(Out(v), In(v))\right)_{Aut G},
$$
and notice that $\PROP^\circlearrowright\langle E\rangle$
 has a natural prop structure
with respect to disjoint union and grafting of graphs. Clearly, this prop
contains the free prop
$\PROP\langle E\rangle$ as a natural sub-prop.

\bip

\no
A {\em derivation}\, in $\PROP^\circlearrowright\langle E\rangle$ is,
by definition,
 a collection of linear maps, $\delta: \PROP^\circlearrowright
\langle E\rangle(m,n)
 \rar \PROP^\circlearrowright\langle E\rangle(m,n)$ such that,
 for any $G\in \fG$ and any element of
 $\Graph\langle E\rangle(m,n)$ of the form,
 $$
 {\mathfrak e}=\underset{\rm orderings\ of\ v(G)}{\sf coequalizer}
\left( {\mathfrak e}_1\ot
  {\mathfrak e}_2\ot\ldots \ot {\mathfrak e}_{|v(G)|}\right), \ \ \ {\mathfrak e}_k\in
 E(Out(v_k), In(v_k))\ \mbox{for}\ 1\leq k\leq |v(G)|,
 $$
one has
\[
 \delta{\mathfrak e}= \underset{\rm orderings\ of\ v(G)}{\sf coequalizer}
\left(\sum_{k=1}^{|v(G)|}
 (-1)^{| {\mathfrak e}_1|+\ldots+ |{\mathfrak e}_{k-1}|}
  {\mathfrak e}_1\ot \ldots \ot \delta{\mathfrak e}_k\ot
  \ldots \ot {\mathfrak e}_{|v(G)|}\right).
\]
 Put another way, a graph derivation is completely determined by its
 values on decorated corollas,
$$
\begin{xy}
 <0mm,0mm>*{\mbox{$\xy *=<4mm,3mm>
 \txt{$a$}*\frm{-}\endxy$}};<0mm,0mm>*{}**@{},
%
 %<0mm,0mm>*{\bullet};<0mm,0mm>*{}**@{},
 <-2mm,1.5mm>*{};<-8mm,5mm>*{}**@{-},
 <-1mm,1.5mm>*{};<-4.5mm,5mm>*{}**@{-},
 <1mm,1.5mm>*{};<1mm,5mm>*{\ldots}**@{},
 <1mm,1.5mm>*{};<4.5mm,5mm>*{}**@{-},
 <2mm,1.5mm>*{};<8mm,5mm>*{}**@{-},
   <0mm,0mm>*{};<-8.7mm,6.3mm>*{_1}**@{},
   <0mm,0mm>*{};<-4.6mm,6.3mm>*{_2}**@{},
   <0mm,0mm>*{};<9.9mm,6.3mm>*{_m}**@{},
 <-2mm,-1.55mm>*{};<-8mm,-5mm>*{}**@{-},
 <-1mm,-1.55mm>*{};<-4.5mm,-5mm>*{}**@{-},
 <0mm,0mm>*{};<1mm,-5mm>*{\ldots}**@{},
 <1mm,-1.55mm>*{};<4.5mm,-5mm>*{}**@{-},
 <2mm,-1.55mm>*{};<8mm,-5mm>*{}**@{-},
   <0mm,0mm>*{};<-8.7mm,-6.9mm>*{_1}**@{},
   <0mm,0mm>*{};<-4.5mm,-6.9mm>*{_2}**@{},
   <0mm,0mm>*{};<9.9mm,-6.9mm>*{_n}**@{},
 \end{xy} \hspace{10mm} a\in E(m,n),
$$
that is, by linear maps,
 $$
\delta: E(m,n) \lon \Graph\langle E\rangle(m,n).
 $$

\bip

\no
A {\em  differential}\, in $\Graph\langle E\rangle$ is,
by definition, a degree 1  derivation $\delta$ satisfying the condition
$\delta^2=0$.

\bip

%%%%%%%%%%%%%%%%%%%%%%%%%%%%%%%%%%
\no{\bf 2.3.1. Remark.}
If  $(\PROP\langle E\rangle, \delta)$ is a dg free prop generated by
an $\bS$-bimodule $E$, then $\delta$ extends naturally to a differential on
 $\Graph\langle E\rangle$ which we denote by the same symbol $\delta$. It is
 worth pointing out that such an induced differential may {\em not}\, preserve the number
of oriented wheels. For example, if $\delta$
 applied to an element $a\in E(m,n)$ (which we identify with the $a$-decorated
  $(m,n)$-corolla) contains a summand of the form,
$$
\delta\left(
\begin{xy}
 <0mm,0mm>*{\mbox{$\xy *=<4mm,3mm>
 \txt{$a$}*\frm{-}\endxy$}};<0mm,0mm>*{}**@{},
%
 %<0mm,0mm>*{\bullet};<0mm,0mm>*{}**@{},
 <-2mm,1.5mm>*{};<-8mm,4mm>*{}**@{-},
 <-1mm,1.5mm>*{};<-3.5mm,4mm>*{}**@{-},
 <1mm,1.5mm>*{};<1mm,4mm>*{\ldots}**@{},
 <1mm,1.5mm>*{};<4.5mm,4mm>*{}**@{-},
 <2mm,1.5mm>*{};<8mm,4mm>*{}**@{-},
   <0mm,0mm>*{};<-8.7mm,5.3mm>*{_{i_1}}**@{},
   <0mm,0mm>*{};<-4.6mm,5.3mm>*{_{i_2}}**@{},
   %<0mm,0mm>*{};<9.9mm,6.3mm>*{_{i_m}}**@{},
 <-2mm,-1.55mm>*{};<-8mm,-4mm>*{}**@{-},
 <-1mm,-1.55mm>*{};<-3.5mm,-4mm>*{}**@{-},
 <0mm,0mm>*{};<1mm,-4mm>*{\ldots}**@{},
 <1mm,-1.55mm>*{};<4.5mm,-4mm>*{}**@{-},
 <2mm,-1.55mm>*{};<8mm,-4mm>*{}**@{-},
   <0mm,0mm>*{};<-8.7mm,-5.9mm>*{_{j_1}}**@{},
   <0mm,0mm>*{};<-4.5mm,-5.9mm>*{_{j_2}}**@{},
   %<0mm,0mm>*{};<9.9mm,-6.9mm>*{_{j_n}}**@{},
 \end{xy}
 \right)= \ldots +
\begin{xy}
 <0mm,0mm>*{\mbox{$\xy *=<4mm,3mm>
 \txt{$b$}*\frm{-}\endxy$}};<0mm,0mm>*{}**@{},
%
 %<0mm,0mm>*{\bullet};<0mm,0mm>*{}**@{},
 <-2mm,1.5mm>*{};<-8mm,4mm>*{}**@{-},
 <-1mm,1.5mm>*{};<-4.5mm,4mm>*{}**@{-},
 <1mm,1.5mm>*{};<-1mm,4mm>*{...}**@{},
 <1mm,1.5mm>*{};<4mm,4mm>*{...}**@{},
 <1mm,1.5mm>*{};<2.5mm,10.9mm>*{}**@{-},
 <2mm,1.5mm>*{};<8mm,4mm>*{}**@{-},
   <0mm,0mm>*{};<-8.7mm,5.3mm>*{_{i_1}}**@{},
   <0mm,0mm>*{};<-4.6mm,5.3mm>*{_{i_2}}**@{},
   %<0mm,0mm>*{};<9.9mm,6.3mm>*{_{i_m}}**@{},
 <-2mm,-1.55mm>*{};<-8mm,-4mm>*{}**@{-},
 <-1mm,-1.55mm>*{};<-4.5mm,-4mm>*{}**@{-},
 <0mm,0mm>*{};<1mm,-4mm>*{...}**@{},
 <1mm,-1.55mm>*{};<4.5mm,-4mm>*{}**@{-},
 <2mm,-1.55mm>*{};<8mm,-4mm>*{}**@{-},
 <0mm,0mm>*{};<-8.7mm,-5.9mm>*{_{j_2}}**@{},
   %<0mm,0mm>*{};<-8.7mm,-6.9mm>*{_1}**@{},
   %<0mm,0mm>*{};<-4.5mm,-6.9mm>*{_2}**@{},
   %<0mm,0mm>*{};<9.9mm,-6.9mm>*{_n}**@{},
%
 <2.5mm,12.5mm>*{\mbox{$\xy *=<4mm,3mm>
 \txt{$c$}*\frm{-}\endxy$}};<0mm,0mm>*{}**@{},
 <0.5mm,11mm>*{};<-3mm,9mm>*{}**@{-},
 <4.5mm,11mm>*{};<6mm,9mm>*{}**@{-},
 <4.5mm,14mm>*{};<6mm,16mm>*{}**@{-},
 <0.5mm,14mm>*{};<-3mm,16mm>*{}**@{-},
 <0mm,0mm>*{};<-3.1mm,8.2mm>*{_{j_1}}**@{},
 <1mm,1.5mm>*{};<0.1mm,9mm>*{...}**@{},
 <2.5mm,14mm>*{};<2.5mm,15.8mm>*{...}**@{},
 <1mm,1.5mm>*{};<4mm,9mm>*{...}**@{},
 \end{xy}
 \ \ + \ \ldots
$$
then the value of $\delta$ on the graph obtained from this corolla by gluing
output $i_1$ with input $j_1$ into a loop,
$$
\delta\left(
\xy
 <0mm,0mm>*{\mbox{$\xy *=<4mm,3mm>
 \txt{$a$}*\frm{-}\endxy$}};<0mm,0mm>*{}**@{},
%
 %<0mm,0mm>*{\bullet};<0mm,0mm>*{}**@{},
 %<-2mm,1.5mm>*{};<-8mm,4mm>*{}**@{-},
 <-1mm,1.5mm>*{};<-2.5mm,4mm>*{}**@{-},
 <1mm,1.5mm>*{};<0.2mm,4mm>*{...}**@{},
 <1mm,1.5mm>*{};<3.5mm,4mm>*{}**@{-},
 <2mm,1.5mm>*{};<7mm,4mm>*{}**@{-},
   %<0mm,0mm>*{};<-8.7mm,5.3mm>*{_{i_1}}**@{},
   <0mm,0mm>*{};<-2.6mm,5.3mm>*{_{i_2}}**@{},
   %<0mm,0mm>*{};<9.9mm,6.3mm>*{_{i_m}}**@{},
 %<-2mm,-1.55mm>*{};<-8mm,-4mm>*{}**@{-},
 <-1mm,-1.55mm>*{};<-2.5mm,-4mm>*{}**@{-},
 <0mm,0mm>*{};<0.2mm,-4mm>*{...}**@{},
 <1mm,-1.55mm>*{};<3.5mm,-4mm>*{}**@{-},
 <2mm,-1.55mm>*{};<7mm,-4mm>*{}**@{-},
   %<0mm,0mm>*{};<-8.7mm,-5.9mm>*{_{j_1}}**@{},
   <0mm,0mm>*{};<-2.4mm,-5.3mm>*{_{j_2}}**@{},
   %<0mm,0mm>*{};<9.9mm,-6.9mm>*{_{j_n}}**@{},
   %
   (-2,1.5)*{}
   \ar@{->}@(ul,dl) (-2,-1.5)*{}
 \endxy
 \right)= \ldots + \ \
\begin{xy}
 <0mm,0mm>*{\mbox{$\xy *=<4mm,3mm>
 \txt{$b$}*\frm{-}\endxy$}};<0mm,0mm>*{}**@{},
%
 %<0mm,0mm>*{\bullet};<0mm,0mm>*{}**@{},
 <-2mm,1.5mm>*{};<-10mm,6mm>*{}**@{-},
 <-1mm,1.5mm>*{};<-4.5mm,4mm>*{}**@{-},
 <1mm,1.5mm>*{};<-1mm,4mm>*{...}**@{},
 <1mm,1.5mm>*{};<4mm,4mm>*{...}**@{},
 <1mm,1.5mm>*{};<2.5mm,10.9mm>*{}**@{-},
 <2mm,1.5mm>*{};<8mm,4mm>*{}**@{-},
   %<0mm,0mm>*{};<-8.7mm,5.3mm>*{_{i_1}}**@{},
   <0mm,0mm>*{};<-4.6mm,5.3mm>*{_{i_2}}**@{},
   %<0mm,0mm>*{};<9.9mm,6.3mm>*{_{i_m}}**@{},
 <-2mm,-1.55mm>*{};<-8mm,-4mm>*{}**@{-},
 <-1mm,-1.55mm>*{};<-4.5mm,-4mm>*{}**@{-},
 <0mm,0mm>*{};<1mm,-4mm>*{...}**@{},
 <1mm,-1.55mm>*{};<4.5mm,-4mm>*{}**@{-},
 <2mm,-1.55mm>*{};<8mm,-4mm>*{}**@{-},
 <0mm,0mm>*{};<-8.7mm,-5.9mm>*{_{j_2}}**@{},
   %<0mm,0mm>*{};<-8.7mm,-6.9mm>*{_1}**@{},
   %<0mm,0mm>*{};<-4.5mm,-6.9mm>*{_2}**@{},
   %<0mm,0mm>*{};<9.9mm,-6.9mm>*{_n}**@{},
%
 <2.5mm,12.5mm>*{\mbox{$\xy *=<4mm,3mm>
 \txt{$c$}*\frm{-}\endxy$}};<0mm,0mm>*{}**@{},
 <0.5mm,11mm>*{};<-10mm,6mm>*{}**@{-},
 <4.5mm,11mm>*{};<6mm,9mm>*{}**@{-},
 <4.5mm,14mm>*{};<6mm,16mm>*{}**@{-},
 <0.5mm,14mm>*{};<-3mm,16mm>*{}**@{-},
 %<0mm,0mm>*{};<-3.1mm,8.2mm>*{_{j_1}}**@{},
 <1mm,1.5mm>*{};<0.1mm,9mm>*{...}**@{},
 <2.5mm,14mm>*{};<2.5mm,15.8mm>*{...}**@{},
 <1mm,1.5mm>*{};<4mm,9mm>*{...}**@{},
 \end{xy}
 \ \ + \ \ldots \ .
$$
contains a term with no oriented wheels at all. Thus propic differential  can, in general, {\em decrease}\, the number of
wheels. Notice in this connection that if $\delta$ is
induced on $\Graph\langle E\rangle$ from the minimal model of a
$\frac{1}{2}$prop, then such summands are impossible and hence
the differential preserves the number of wheels.

\bip

\no
Vector spaces $\PROP\langle E\rangle$ and
%$\PROP^+\hspace{-0.5mm}\langle E\rangle$
$\Graph\langle E\rangle$ have a  natural positive gradation,
$$
\PROP\langle E\rangle
=\bigoplus_{k\geq 1} \PROP_k\langle E\rangle, \ \ \
%
%\PROP^+\hspace{-0.5mm}\langle E\rangle
%=\bigoplus_{k\geq 1} \PROP_k^+\hspace{-0.5mm}\langle E\rangle, \ \
%
\Graph\langle E\rangle
=\bigoplus_{k\geq 1} \Graph\langle E\rangle, \ \
$$
by the number, $k$, of internal vertices of underlying graphs. In particular,
$ \PROP_1\langle E\rangle(m,n)$
 is spanned by decorated
$(m,n)$-corollas and can be identified with $E(m,n)$.

\bip

%%%%%%%%%%%%%%%%%%%%%%%%%%%%%%%%%%%%%%%%
\no{\bf 2.3.2. Representations of $\Graph \langle E \rangle$}.
Any representation, $\phi: E \rar \End\langle M \rangle$, of an $\bS$-bimodule $E$
in a finite dimensional vector space $M$ can be naturally
extended to  representations of props,
$\PROP\langle E\rangle \rar \End\langle M \rangle$ and
$\Graph\langle E\rangle \rar \End\langle M \rangle$. In the latter case
 decorated graphs with oriented wheels are mapped into appropriate traces.

\bip
%%%%%%%%%%%%%%%%%%%%%%%%%%%%%%%%%%%%%%%%%%%%%%%%
\no{\bf 2.3.3. Remark.} Prop structure on an $\bS$-bimodule  $E=\{E(m,n)\}$ can be defined as a family of evaluation
linear maps,
$$
{\mu}_G: G \langle E\rangle \lon E(m,n), \ \ \forall\ G\in \fG^\uparrow,
$$
satisfying certain associativity axiom (cf.\ \S 2.1.3). Analogously, one can define
a {\em wheeled prop struture}\, on $E$ as a family of linear maps,
$$
{\mu}_G: G \langle E\rangle \lon E(m,n), \ \ \forall\ G\in \fGG,
$$
 such that
\Bi
\item[(i)] $\mu_{(m,n)-corolla}=\Id$,
\item[(ii)] $\mu_G=\mu_{G/H}\circ\mu_H$ for
 every subgraph $H\in \fGG$ of $G$, where $G/H$ is obtained from $G$ by collapsing to the single vertex
 every connected component of $H$, and   $\mu_H:  G \langle E\rangle \rar  G/H \langle E\rangle$
 is the evaluation map on the subgraph $H$ and identity on its complement.
\Ei

\no{\bf Claim.} {\em For every finite-dimensional vector space $M$ the associated endomorphism prop
$\End\langle M\rangle$ has a natural structure of wheeled prop.}

\bip

\no
The notion of representation of  $\Graph \langle E \rangle$ in a finite dimensional vector space $M$
introduced above is just a morphism of {\em wheeled}\, props, $\Graph \langle E \rangle\rar \End\langle M\rangle$.
We shall discuss these issues in detail elsewhere as
in the present paper we are most interested in computing cohomology of dg {\em free}\, wheeled props
$(\Graph \langle E \rangle, \delta)$, where the composition maps $\mu_G$ are tautological.

\bip
%%%%%%%%%%%%%%%%%%%%%%%%%%%%%%%%%%%%%%%%%%%%%%%%
\no{\bf 2.4. Formal graded manifolds.}
For a finite-dimensional vector space $M$ we denote by $\cM$
the associated
formal graded manifold. The distinguished point of the latter is always
denoted by $*$. The structure sheaf, $\f_\cM$, is (non-canonically) isomorphic
to the completed graded symmetric tensor algebra, $\hat{\odot} M^*$.
A choice of a particular isomorphism, $\phi: \f_\cM\rar \hat{\odot} M^*$,
is called a choice of a local coordinate system on $\cM$. If
$\{e_\al\}_{\al\in I}$ is a basis in $M$ and $\{t^\al\}_{\al\in I}$
the associated dual basis in $M^*$, then
one may identify  $\f_\cM$ with the graded commutative
formal power series ring $\R[[t^\al]]$.

\bip

\no
Free modules over the ring $\f_\cM$ are called locally free sheaves
(=vector bundles) on $\cM$.
The $\f_\cM$-module, $\cT_\cM$, of derivations of $\f_\cM$ is called
the tangent sheaf on $\cM$. Its dual, $\Omega_\cM$, is called the
cotangent sheaf. One can form their (graded skewsymmetric) tensor products
such as the sheaf of polyvector fields, $\wedge^\bullet \cT_\cM$, and
the sheaf of differential forms,
$\Omega^\bullet_\cM=\wedge^\bullet \Omega_\cM$. The first sheaf is
naturally
a sheaf of Lie algebras on $\cM$ with respect to the Schouten bracket.

\bip

\no
One can also define a sheaf of polydifferential operators, $\caD_\cM\subset
\oplus_{i\geq 0}\Hom_\R(\f_\cM^{\ot i}, \f_\cM)$. The latter is naturally
a sheaf of dg Lie algebras on $\cM$ with respect to the Hochschild
differential, $d_H$, and brackets, $[\ ,\ ]_H$.

\bip
%%%%%%%%%%%%%%%%%%%%%%%%%%%%%%%%%%%%%%%%%%%%%%%%%%%%%%
\no{\bf 2.5. Geometry $\Rightarrow$ graph complexes.} We shall sketch here a
simple trick
which associates a dg free prop, $\Graph\langle E_{\cG}\rangle$,
 to a
sheaf of dg Lie algebras, $\cG_\cM$, over a smooth
graded formal manifold $\cM$.
%We often abbreviate $\Graph\langle E_{\cG}\rangle$ to $\Graph\langle\cG\rangle$.

\bip

\no
We assume that
\Bi
 \item[(i)] $\cG_\cM$ is built from direct sums and  tensor products of
 (any order)
 jets of the sheaves $\cT_\cM^{\ot \bullet}\ot \Omega_\cM^{\ot \bullet}$
  and their duals (thus $\cG_\cM$ can be defined for
{\em any}\,
 formal smooth manifold $\cM$, i.e., its definition does not depend
on the dimension of $\cM$),
 \item[(ii)] the differential and the Lie bracket in $\cG_M$ can be
 represented, in a local coordinate system, by polydifferential operators
 and natural contractions between the duals.
 \Ei
The motivating examples are  $\wedge^\bullet \cT_\cM$, $\caD_\cM$ and the
sheaf of Nijenhuis dg Lie algebras on $\cM$ (see \cite{Me2}).

\bip

\no
By assumption (i), a choice of a local coordinate system on $\cM$,
identifies $\cG_\cM$ with a subspace in
$$
\bigoplus_{q,m\geq 0}\f_\cM\ot \Hom(M^{\ot q}, M^{\ot m})
=
\prod_{p,q,m\geq 0}\Hom(\odot^p M\ot M^{\ot q}, M^{\ot m})
\subset
\prod_{m,n\geq 0} \Hom(M^{\ot n}, M^{\ot m}).
$$

\no
Let $\Gamma$ be a degree 1 element in $\cG_\cM$.
Denote by $\Gamma^m_{p,q}$ the bit of $\Gamma$ which lies in
$\Hom(\odot^p M\ot M^{\ot q}, M^{\ot m})$ and set
$\Gamma^m_{n} := \oplus_{p+q=n}\Gamma^m_{p,q}\in  \Hom(M^{\ot n}, M^{\ot m})$.

\bip

\no
There exists a uniquely defined finite-dimensional $\bS$-bimodule,
$E_\cG=\{E_\cG(m,n)\}_{m,n\geq 0}$,
whose representations in the vector space $M$ are in
one-to-one correspondence with Taylor components,
$\Gamma^m_{n}\in  \Hom(M^{\ot n}, M^{\ot m})$, of a degree 1 element
$\Gamma$
in $\cG_\cM$.
Set $\Graph\langle\cG\rangle:= \Graph\langle E_\cG\rangle$ (see Sect.\ 2.3).

\bip

\no
Next we employ the dg Lie algebra structure in $\cG_\cM$ to introduce a
differential, $\delta$,
in  $\Graph\langle\cG\rangle$. The latter is completely determined
by its restriction to the subspace of $\Graph_1\langle\cG\rangle$ spanned by decorated
corollas (without attached loops).

\bip

\no
 First we replace the Taylor
coefficients, $\Gamma^m_{n}$,  of the section $\Gamma$
by the decorated $(m,n)$-corollas
\Bi
\item
with the unique internal vertex decorated by a basis element,
$\{{\mathfrak e}_r\}_{r\in J}$,
  of $E_\cG(m,n)$,
\item with input legs labeled by basis elements, $\{e_\al\}$, of the vector
space $M$  and output legs labeled by the elements of the dual basis,
$\{t^\be\}$.
\Ei

\bip

\no
Next we consider a formal linear combination,
$$
\overline{\Gamma}_n^m=\sum_r\sum_{\al_1,...,\al_n \atop
\be_1,...,\be_m}
\begin{xy}
 <0mm,0mm>*{\mbox{$\xy *=<4mm,3mm>
 \txt{${\mathfrak e}_r$}*\frm{-}\endxy$}};<0mm,0mm>*{}**@{},
%
 %<0mm,0mm>*{\bullet};<0mm,0mm>*{}**@{},
 <-2mm,1.5mm>*{};<-8mm,5mm>*{}**@{-},
 <-1mm,1.5mm>*{};<-4.5mm,5mm>*{}**@{-},
 <1mm,1.5mm>*{};<1mm,5mm>*{\ldots}**@{},
 <1mm,1.5mm>*{};<4.5mm,5mm>*{}**@{-},
 <2mm,1.5mm>*{};<8mm,5mm>*{}**@{-},
   <0mm,0mm>*{};<-8.7mm,6.3mm>*{^{t^{\be_1}}}**@{},
   <0mm,0mm>*{};<-4.6mm,6.3mm>*{^{t^{\be_2}}}**@{},
   %<0mm,0mm>*{};<4.5mm,5.5mm>*{^{m\hspace{-0.5mm}-\hspace{-0.5mm}1}}**@{},
   <0mm,0mm>*{};<9.2mm,6.3mm>*{^{t^{\be_m}}}**@{},
 <-2mm,-1.55mm>*{};<-8mm,-5mm>*{}**@{-},
 <-1mm,-1.55mm>*{};<-4.5mm,-5mm>*{}**@{-},
 <0mm,0mm>*{};<1mm,-5mm>*{\ldots}**@{},
 <1mm,-1.55mm>*{};<4.5mm,-5mm>*{}**@{-},
 <2mm,-1.55mm>*{};<8mm,-5mm>*{}**@{-},
   <0mm,0mm>*{};<-8.7mm,-6.9mm>*{_{e_{\al_1}}}**@{},
   <0mm,0mm>*{};<-4.5mm,-6.9mm>*{_{e_{\al_2}}}**@{},
   %<0mm,0mm>*{};<4.5mm,-6.9mm>*{^{n\hspace{-0.5mm}-\hspace{-0.5mm}1}}**@{},
   <0mm,0mm>*{};<9.2mm,-6.9mm>*{_{e_{\al_n}}}**@{},
 \end{xy}
 \ \ t^{\al_1}\ot\ldots\ot t^{\al_n}\ot e_{\be_1}\ot\ldots\ot e_{\be_m}
 \ \in \Graph_1\langle\cG\rangle\ot
\Hom(M^{\ot n}, M^{\ot m}).
$$
This expression is essentially a component of the Taylor decomposition of $\Gamma$,
$$
\Gamma^m_n = \sum_{\al_1\ldots \al_n
\atop \be_1\ldots \be_m} \Gamma^{\be_1\ldots\be_m}_{\al_1\ldots\al_n}
 t^{\al_1}\ot\ldots\ot t^{\al_n}\ot
e_{\be_1}\ot\ldots\ot  e_{\be_m},
$$
in which the numerical coefficient
$ \Gamma^{\be_1\ldots\be_m}_{\al_1\ldots\al_n}$ is substituted by a
 decorated labeled graph. More precisely, the interrelation
between $\overline{\Gamma}=\oplus_{m,n\geq 0} \overline{\Gamma}_n^m$ and
${\Gamma}=\oplus_{m,n\geq 0}{\Gamma}_n^m\in \cG_\cM$ can be described
as follows:
 a choice of any
particular representation of the $\bS$-bimodule $E_\cG$,
$$
\phi: \left\{E_\cG(m,n) \rar \Hom(M^{\ot n}, M^{\ot m})\right\}_{m,n\geq 0},
$$
defines an element $\Gamma=\phi(\overline{\Gamma})\in \cG_\cM$ which is
obtained from
$\overline{\Gamma}$ by replacing
each graph,
$$
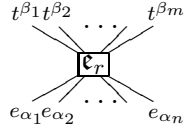
\begin{xy}
 <0mm,0mm>*{\mbox{$\xy *=<4mm,3mm>
 \txt{${\mathfrak e}_r$}*\frm{-}\endxy$}};<0mm,0mm>*{}**@{},
%
 %<0mm,0mm>*{\bullet};<0mm,0mm>*{}**@{},
 <-2mm,1.5mm>*{};<-8mm,5mm>*{}**@{-},
 <-1mm,1.5mm>*{};<-4.5mm,5mm>*{}**@{-},
 <1mm,1.5mm>*{};<1mm,5mm>*{\ldots}**@{},
 <1mm,1.5mm>*{};<4.5mm,5mm>*{}**@{-},
 <2mm,1.5mm>*{};<8mm,5mm>*{}**@{-},
   <0mm,0mm>*{};<-8.7mm,6.3mm>*{^{t^{\be_1}}}**@{},
   <0mm,0mm>*{};<-4.6mm,6.3mm>*{^{t^{\be_2}}}**@{},
   %<0mm,0mm>*{};<4.5mm,5.5mm>*{^{m\hspace{-0.5mm}-\hspace{-0.5mm}1}}**@{},
   <0mm,0mm>*{};<9.9mm,6.3mm>*{^{t^{\be_m}}}**@{},
 <-2mm,-1.55mm>*{};<-8mm,-5mm>*{}**@{-},
 <-1mm,-1.55mm>*{};<-4.5mm,-5mm>*{}**@{-},
 <0mm,0mm>*{};<1mm,-5mm>*{\ldots}**@{},
 <1mm,-1.55mm>*{};<4.5mm,-5mm>*{}**@{-},
 <2mm,-1.55mm>*{};<8mm,-5mm>*{}**@{-},
   <0mm,0mm>*{};<-8.7mm,-6.9mm>*{_{e_{\al_1}}}**@{},
   <0mm,0mm>*{};<-4.5mm,-6.9mm>*{_{e_{\al_2}}}**@{},
   %<0mm,0mm>*{};<4.5mm,-6.9mm>*{^{n\hspace{-0.5mm}-\hspace{-0.5mm}1}}**@{},
   <0mm,0mm>*{};<9.9mm,-6.9mm>*{_{e_{\al_n}}}**@{},
 \end{xy}
$$
by the value,
$\stackrel{r}{\Gamma}^{\be_1\ldots\be_m}_{\al_1\ldots\al_n}\in \R$,
of $\phi({\mathfrak e}_r\}\in \Hom(M^{\ot n}, M^{\ot m})$
on the basis vector
$e_{\al_1}\ot\ldots\ot  e_{\al_n}\ot t^{\be_1}\ot\ldots\ot t^{\be_m}$
(so that  $\Gamma^{\be_1\ldots\be_m}_{\al_1\ldots\al_n}=
\sum_r \stackrel{r}{\Gamma}^{\be_1\ldots\be_m}_{\al_1\ldots\al_n}$) .

\bip

\no
 In a similar way one can define an element,
$$
\overline{[\cdots [[d\Gamma,\Gamma], \Gamma]\cdots]}\in
\Graph_n\langle\cG\rangle\ot
\Hom(M^{\ot \bullet}, M^{\ot \bullet})
$$
for any Lie word,
$$
[\cdots[[d\Gamma,\Gamma], \Gamma]\ldots],
$$
 built from $\Gamma$, $d\Gamma$ and $n-1$ Lie brackets.
 In particular, there are uniquely defined elements,
 $$
 \overline{d\Gamma}\in   \Graph_1\langle\cG\rangle\ot
\Hom(M^{\ot \bullet}, M^{\ot \bullet})       , \ \ \
 \overline{\frac{1}{2}[{\Gamma},{\Gamma}]}\in
 \Graph_2\langle\cG\rangle\ot
\Hom(M^{\ot \bullet}, M^{\ot \bullet}),
$$
whose values, $\phi(\overline{d{\Gamma}})$ and
$\phi(\overline{\frac{1}{2}[{\Gamma},{\Gamma}]})$, for any
particular choice of representation $\phi$
of the $\bS$-bimodule $E_\cG$, coincide respectively with $d\Gamma$
and $ \frac{1}{2}[{\Gamma},{\Gamma}]$.

 \bip

\no
Finally one defines a differential $\delta$ in the graded space
$\Graph\langle\cG\rangle$
by setting
$$
\hspace{57mm}
\delta \overline{\Gamma}= \overline{d{\Gamma}} + \overline{
\frac{1}{2}[{\Gamma},{\Gamma}]},
\hspace{60mm}  \ (\star\star)
$$
i.e.\ by equating the graph coefficients of both sides.
That $\delta^2=0$ is clear from the following
calculation,
\Beqrn
\delta^2 \overline{\Gamma} &=& \delta\left(
\overline{d{\Gamma}} + \overline{
\frac{1}{2}[{\Gamma},{\Gamma}]}\right)\\
&=& \delta\overline{d\Gamma} +
\overline{[d\Gamma + \frac{1}{2}[\Gamma,\Gamma], \Gamma]}\\
&=&-\overline{d (d\Gamma + \frac{1}{2}[\Gamma,\Gamma])} +
\overline{[d\Gamma , \Gamma]} + \frac{1}{2}
\overline{[[\Gamma,\Gamma], \Gamma]}\\
&=& -\overline{[d\Gamma , \Gamma]} + \overline{[d\Gamma , \Gamma]}\\
&=& 0,
\Eeqrn
where we used both the
axioms of dg Lie algebra in $\cG_\cM$ and the axioms of the
differential in $\Graph\langle\cG\rangle$. This completes the construction
of  $(\Graph\langle\cG \rangle, \delta)$\footnote{As a first approximation to the
propic translation of {\em non}-flat geometries (Yang-Mills, Riemann, etc.)
 one might consider the following
version of the ``trick'': in addition to generic element
$\Gamma\in \cG_\cM$ of degree 1 take into consideration
(probably, {\em non}\, generic) element
of degree 2, $F\in \cG_\cM$, extend appropriately the $\bS$-bimodule $E_\cG$
to accommodate the associated ``curvature'' $F$-corollas, and then (attempt to)
 define the
differential $\delta$ in $\Graph\langle E_\cG \rangle$ by equating graph
coefficients in the expressions,
$\delta \overline{\Gamma}= \overline{F}+ \overline{d{\Gamma}} + \overline{
\frac{1}{2}[{\Gamma},{\Gamma}]}$ and
$\delta \overline{F}= \overline{d{F}} + \overline{
[{\Gamma},{F}]}$.}.

\bip

%%%%%%%%%%%%%%%%%%%%%%%%%%%%%%
\no{\bf 2.5.1. Remarks.}
(i) If the differential and Lie brackets in $\cG_\cM$ contain no traces,
then the expression  $\overline{d\Gamma}+\overline{
\frac{1}{2}[{\Gamma},{\Gamma}]}$ does not contain graphs with oriented wheels.
Hence formula  $(\star\star)$ can be used to introduce a differential
in the free prop, $\PROP\langle\cG\rangle$,
generated by the $\bS$-bimodule $E_\cG$.

\sip

(ii)  If the differential and Lie brackets in  $\cG_\cM$ contain no traces and are given
by first order differential operators, then the expression
 $\overline{d\Gamma}+\overline{
\frac{1}{2}[{\Gamma},{\Gamma}]}$ is a tree. Therefore formula  $(\star\star)$ can be used to
introduce a differential
in the free dioperad, ${\mathsf D} \langle\cG\rangle$.

\bip

%%%%%%%%%%%%%%%%%%%%%%%%%%%%%%
\no{\bf 2.5.2. Remark.} The above trick  works also for sheaves,
$$
\left(\cG_\cM, \mu_n: \wedge^n \cG_\cM\rar \cG_\cM[2-n], n=1,2,\ldots\right),
$$
 of $L_\infty$
algebras over $\cM$. The differential in $\PROP^\circlearrowright\langle\cG\rangle$
(or in $\PROP\langle\cG\rangle$,
if appropriate) is defined by,
$$
\delta\overline{\Gamma} = \sum_{n=1}^\infty \frac{1}{n!}
\overline{\mu_n(\Gamma,\ldots,\Gamma)}.
$$
The term $\overline{\mu_n(\Gamma,\ldots,\Gamma)}$ corresponds to
decorated graphs with $n$ internal vertices.

\bip

%%%%%%%%%%%%%%%%%%%%%%%%%%%%%%
\no{\bf 2.5.3. Remark.} Any sheaf of dg Lie subalgebras, $\cG'_\cM\subset
\cG_\cM$, defines a dg  prop,
 $(\Graph\langle\cG'\rangle,
\delta)$,  which is a quotient of $(\Graph\langle\cG\rangle,
\delta)$ by the ideal generated by decorated graphs lying in the complement,
$\Graph\langle\cG\rangle\setminus
\Graph\langle\cG'\rangle$.
Similar observation holds true for
$\PROP\langle\cG\rangle$ and $\PROP\langle\cG'\rangle$
(if they are defined).

\bip

%%%%%%%%%%%%%%%%%%%%%%%%%%%%%%%%%%%%%
\no{\bf 2.6. Prop profile of Poisson structures.} Let us consider the sheaf of
polyvector fields, $\wedge^\bullet\cT_\cM:=\sum_{i\geq 0}
\wedge^i \cT_\cM [1-i]$, equipped with the Schouten Lie bracket, $[\ ,\ ]_S$,
and vanishing
differential. A degree one section, $\Gamma$, of
$\wedge^\bullet\cT_\cM$ decomposes
into a direct sum, $\oplus_{i\geq 0} \Gamma_i$, with $\Gamma_i\in
\wedge^{i}\cT_\cM$ having degree $2-i$ with respect to the grading of the underlying
manifold. In a local coordinate system $\Gamma$
can be represented as a Taylor series,
$$
\Gamma = \sum_{m,n\geq 0}\sum_{\al_1\ldots \al_n
\atop \be_1\ldots \be_m} \Gamma^{\be_1\ldots\be_m}_{\al_1\ldots\al_n}
(e_{\be_1}\wedge \ldots\wedge  e_{\be_m})\ot
(t^{\al_1}\odot\ldots\odot t^{\al_n}).
$$
As $\Gamma^{\be_1\ldots\be_m}_{\al_1\ldots\al_n}=
\Gamma^{[\be_1\ldots\be_m]}_{(\al_1\ldots\al_n)}$ has degree $2-m$,
we conclude that
the associated $\bS$-bimodule $E_{\wedge^\bullet \cT}$ is given by
$$
E_{\wedge^\bullet \cT}(m,n)={\bf sgn_m}\ot {\bf 1_n}[m-2], \ \ m,n\geq 0,
$$
where $\bf sgn_m$ stands for the one dimensional sign representation of $\Sigma_m$
and ${\bf 1_n}$ stands for the trivial one-dimensional representation
of $\Sigma_n$.  Then a generator of $\PROP\langle
\wedge^\bullet \cT\rangle$
can be represented
by the  directed
planar $(m,n)$-corolla,
$$
 
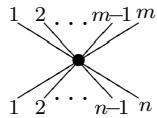
\begin{xy}
 <0mm,0mm>*{\bullet};<0mm,0mm>*{}**@{},
 <0mm,0mm>*{};<-8mm,5mm>*{}**@{-},
 <0mm,0mm>*{};<-4.5mm,5mm>*{}**@{-},
 <0mm,0mm>*{};<-1mm,5mm>*{\ldots}**@{},
 <0mm,0mm>*{};<4.5mm,5mm>*{}**@{-},
 <0mm,0mm>*{};<8mm,5mm>*{}**@{-},
   <0mm,0mm>*{};<-8.5mm,5.5mm>*{^1}**@{},
   <0mm,0mm>*{};<-5mm,5.5mm>*{^2}**@{},
   <0mm,0mm>*{};<4.5mm,5.5mm>*{^{m\hspace{-0.5mm}-\hspace{-0.5mm}1}}**@{},
   <0mm,0mm>*{};<9.0mm,5.5mm>*{^m}**@{},
 <0mm,0mm>*{};<-8mm,-5mm>*{}**@{-},
 <0mm,0mm>*{};<-4.5mm,-5mm>*{}**@{-},
 <0mm,0mm>*{};<-1mm,-5mm>*{\ldots}**@{},
 <0mm,0mm>*{};<4.5mm,-5mm>*{}**@{-},
 <0mm,0mm>*{};<8mm,-5mm>*{}**@{-},
   <0mm,0mm>*{};<-8.5mm,-6.9mm>*{^1}**@{},
   <0mm,0mm>*{};<-5mm,-6.9mm>*{^2}**@{},
   <0mm,0mm>*{};<4.5mm,-6.9mm>*{^{n\hspace{-0.5mm}-\hspace{-0.5mm}1}}**@{},
   <0mm,0mm>*{};<9.0mm,-6.9mm>*{^n}**@{},
 \end{xy}
$$
with skew-symmetric outgoing legs and symmetric ingoing legs. The formula
$(\star\star)$ in Sect.\ 2.5 gives the following explicit expression
for the induced differential, $\delta$, in
$\PROP\langle\wedge^\bullet \cT\rangle$,
$$
\delta \begin{xy}
 <0mm,0mm>*{\bullet};<0mm,0mm>*{}**@{},
 <0mm,0mm>*{};<-8mm,5mm>*{}**@{-},
 <0mm,0mm>*{};<-4.5mm,5mm>*{}**@{-},
 <0mm,0mm>*{};<-1mm,5mm>*{\ldots}**@{},
 <0mm,0mm>*{};<4.5mm,5mm>*{}**@{-},
 <0mm,0mm>*{};<8mm,5mm>*{}**@{-},
   <0mm,0mm>*{};<-8.5mm,5.5mm>*{^1}**@{},
   <0mm,0mm>*{};<-5mm,5.5mm>*{^2}**@{},
   <0mm,0mm>*{};<4.5mm,5.5mm>*{^{m\hspace{-0.5mm}-\hspace{-0.5mm}1}}**@{},
   <0mm,0mm>*{};<9.0mm,5.5mm>*{^m}**@{},
 <0mm,0mm>*{};<-8mm,-5mm>*{}**@{-},
 <0mm,0mm>*{};<-4.5mm,-5mm>*{}**@{-},
 <0mm,0mm>*{};<-1mm,-5mm>*{\ldots}**@{},
 <0mm,0mm>*{};<4.5mm,-5mm>*{}**@{-},
 <0mm,0mm>*{};<8mm,-5mm>*{}**@{-},
   <0mm,0mm>*{};<-8.5mm,-6.9mm>*{^1}**@{},
   <0mm,0mm>*{};<-5mm,-6.9mm>*{^2}**@{},
   <0mm,0mm>*{};<4.5mm,-6.9mm>*{^{n\hspace{-0.5mm}-\hspace{-0.5mm}1}}**@{},
   <0mm,0mm>*{};<9.0mm,-6.9mm>*{^n}**@{},
 \end{xy}
\ \ = \ \
 \sum_{I_1\sqcup I_2=(1,\ldots,m)\atop {J_1\sqcup J_2=(1,\ldots,n)\atop
 {|I_1|\geq 0, |I_2|\geq 1 \atop
 |J_1|\geq 1, |J_2|\geq 0}}
}\hspace{0mm}
(-1)^{\sigma(I_1\sqcup I_2) + |I_1||I_2|}
 \begin{xy}
 <0mm,0mm>*{\bullet};<0mm,0mm>*{}**@{},
 <0mm,0mm>*{};<-8mm,5mm>*{}**@{-},
 <0mm,0mm>*{};<-4.5mm,5mm>*{}**@{-},
 <0mm,0mm>*{};<0mm,5mm>*{\ldots}**@{},
 <0mm,0mm>*{};<4.5mm,5mm>*{}**@{-},
 <0mm,0mm>*{};<13mm,5mm>*{}**@{-},
     <0mm,0mm>*{};<-2mm,7mm>*{\overbrace{\ \ \ \ \ \ \ \ \ \ \ \ }}**@{},
     <0mm,0mm>*{};<-2mm,9mm>*{^{I_1}}**@{},
 <0mm,0mm>*{};<-8mm,-5mm>*{}**@{-},
 <0mm,0mm>*{};<-4.5mm,-5mm>*{}**@{-},
 <0mm,0mm>*{};<-1mm,-5mm>*{\ldots}**@{},
 <0mm,0mm>*{};<4.5mm,-5mm>*{}**@{-},
 <0mm,0mm>*{};<8mm,-5mm>*{}**@{-},
      <0mm,0mm>*{};<0mm,-7mm>*{\underbrace{\ \ \ \ \ \ \ \ \ \ \ \ \ \ \
      }}**@{},
      <0mm,0mm>*{};<0mm,-10.6mm>*{_{J_1}}**@{},
 <13mm,5mm>*{};<13mm,5mm>*{\bullet}**@{},
 <13mm,5mm>*{};<5mm,10mm>*{}**@{-},
 <13mm,5mm>*{};<8.5mm,10mm>*{}**@{-},
 <13mm,5mm>*{};<13mm,10mm>*{\ldots}**@{},
 <13mm,5mm>*{};<16.5mm,10mm>*{}**@{-},
 <13mm,5mm>*{};<20mm,10mm>*{}**@{-},
      <13mm,5mm>*{};<13mm,12mm>*{\overbrace{\ \ \ \ \ \ \ \ \ \ \ \ \ \ }}**@{},
      <13mm,5mm>*{};<13mm,14mm>*{^{I_2}}**@{},
 <13mm,5mm>*{};<8mm,0mm>*{}**@{-},
 <13mm,5mm>*{};<12mm,0mm>*{\ldots}**@{},
 <13mm,5mm>*{};<16.5mm,0mm>*{}**@{-},
 <13mm,5mm>*{};<20mm,0mm>*{}**@{-},
     <13mm,5mm>*{};<14.3mm,-2mm>*{\underbrace{\ \ \ \ \ \ \ \ \ \ \ }}**@{},
     <13mm,5mm>*{};<14.3mm,-4.5mm>*{_{J_2}}**@{},
 \end{xy}
$$
where $\sigma(I_1\sqcup I_2)$ is the sign of the shuffle
$I_1\sqcup I_2=(1,\ldots, m)$.

\bip

%%%%%%%%%%%%%%%%%%%%%%%%%
\no{\bf 2.6.1. Proposition.} {\em There is a one-to-one correspondence between
representations,
$$
\phi: (\PROP\langle\wedge^\bullet \cT\rangle,\delta) \lon
({\sf End}\langle M \rangle,d),
$$
of $(\PROP\langle \wedge^\bullet \cT\rangle,\delta)$
 in a dg vector space
$(M,d)$
and Maurer-Cartan elements, $\gamma$, in $\wedge^\bullet \cT_\cM$, that is, degree
one elements satisfying the equation, $[\ga,\ga]_S=0$.}

\bip

\Proof Let $\phi$ be a representation. Images of the above $(m,n)$-corollas
under $\phi$ provide us with a collection of linear maps,
$\Gamma_n^m: \odot^n M \rar \wedge^m M[2-m]$ which we assemble, as in
Sect.\ 2.5, into a section, $\Gamma=\sum_{m,n}\Gamma_n^m$, of
$\wedge^\bullet \cT_\cM$.

\bip

\no
The differential $d$ in $M$ can be interpreted as a linear (in the
coordinates $\{t^\al\}$) degree one section of  $\cT_\cM$ which we denote
by the same symbol.

\bip

\no
Finally, the commutativity of $\phi$ with the differentials implies
$$
[-d + \Gamma, -d + \Gamma]_S=0.
$$
Thus setting $\ga=-d + \Gamma$ one gets a Maurer-Cartan element
in $\wedge^\bullet \cT_\cM$.

\bip

\no
Reversely, if $\ga$ is a Maurer-Cartan element in $\wedge^\bullet \cT_\cM$,
then decomposing the sum $d+\ga$ into a collection of its
Taylor series components
 as in Sect.\ 2.5, one gets a representation $\phi$. \hfill $\Box$

 \bip

\no
 Let $\wedge^\bullet_0\cT_\cM =\sum_{i\geq 1}
\wedge^i_0 \cT_\cM [1-i]$ be a sheaf of Lie subalgebras of
$\wedge^\bullet\cT_\cM$ consisting of those elements which vanish
at the distinguished point $*\in \cM$, have no
$\wedge^0\cT_\cM[2]$-component, and whose $\wedge^1\cT_\cM[1]$-component
is at least quadratic in the coordinates  $\{t^\al\}$.  The associated dg free
 prop, $\PROP\langle\wedge^\bullet_0\cT\rangle$, is generated by
 $(m,n)$-corollas with $m,n\geq 1$, $m+n\geq 1$, and
 has a surprisingly small cohomology, a fact
 which is of key importance for our proof of the deformation
 quantization theorem.

 \bip

%%%%%%%%%%%%%%%%%%%%%%%%%
\no{\bf 2.6.2. Theorem.} {\em The cohomology of
$(\PROP\langle\wedge^\bullet_0\cT\rangle,\delta)$
is equal to a quadratic prop, $\Liebi$, which is  a quotient,
$$
 {\Liebi} = \frac{ \PROP\langle A \rangle}{{\sf Ideal} <R>}\ ,
$$
of the free prop generated by the following $\bS$-bimodule $A$,
\Bi
\item all $A(m,n)$ vanish except $A(2,1)$ and $A(1,2)$,
\item
$
A(2,1):= {\bf sgn_2}\ot {\bf 1_1} = span\left(
\begin{xy}
 <0mm,-0.55mm>*{};<0mm,-2.5mm>*{}**@{-},
 <0.5mm,0.5mm>*{};<2.2mm,2.2mm>*{}**@{-},
 <-0.48mm,0.48mm>*{};<-2.2mm,2.2mm>*{}**@{-},
 <0mm,0mm>*{\circ};<0mm,0mm>*{}**@{},
   <0mm,-0.55mm>*{};<0mm,-3.8mm>*{_1}**@{},
   <0.5mm,0.5mm>*{};<2.7mm,2.8mm>*{^2}**@{},
   <-0.48mm,0.48mm>*{};<-2.7mm,2.8mm>*{^1}**@{},
 \end{xy}
 = -\,
 \begin{xy}
 <0mm,-0.55mm>*{};<0mm,-2.5mm>*{}**@{-},
 <0.5mm,0.5mm>*{};<2.2mm,2.2mm>*{}**@{-},
 <-0.48mm,0.48mm>*{};<-2.2mm,2.2mm>*{}**@{-},
 <0mm,0mm>*{\circ};<0mm,0mm>*{}**@{},
   <0mm,-0.55mm>*{};<0mm,-3.8mm>*{_1}**@{},
   <0.5mm,0.5mm>*{};<2.7mm,2.8mm>*{^1}**@{},
   <-0.48mm,0.48mm>*{};<-2.7mm,2.8mm>*{^2}**@{},
 \end{xy}
 \right)
 $
\item
$A(1,2):= {\bf 1_1}\ot {\bf 1_2[-1]}=span\left(
\begin{xy}
 <0mm,0.66mm>*{};<0mm,3mm>*{}**@{-},
 <0.39mm,-0.39mm>*{};<2.2mm,-2.2mm>*{}**@{-},
 <-0.35mm,-0.35mm>*{};<-2.2mm,-2.2mm>*{}**@{-},
 <0mm,0mm>*{\bullet};<0mm,0mm>*{}**@{},
   <0mm,0.66mm>*{};<0mm,3.4mm>*{^1}**@{},
   <0.39mm,-0.39mm>*{};<2.9mm,-4mm>*{^2}**@{},
   <-0.35mm,-0.35mm>*{};<-2.8mm,-4mm>*{^1}**@{},
\end{xy}
 =
 \begin{xy}
 <0mm,0.66mm>*{};<0mm,3mm>*{}**@{-},
 <0.39mm,-0.39mm>*{};<2.2mm,-2.2mm>*{}**@{-},
 <-0.35mm,-0.35mm>*{};<-2.2mm,-2.2mm>*{}**@{-},
 <0mm,0mm>*{\bullet};<0mm,0mm>*{}**@{},
   <0mm,0.66mm>*{};<0mm,3.4mm>*{^1}**@{},
   <0.39mm,-0.39mm>*{};<2.9mm,-4mm>*{^1}**@{},
   <-0.35mm,-0.35mm>*{};<-2.8mm,-4mm>*{^2}**@{},
\end{xy}\right)
$
\Ei
modulo the ideal  generated by the following relations, $R$,
%%%%%%%%%%%%%%%%%%%% coLie %%%%%%%%%%%%%%%%%%
\Beqrn
R_1: &&
\begin{xy}
 <0mm,0mm>*{\circ};<0mm,0mm>*{}**@{},
 <0mm,-0.49mm>*{};<0mm,-3.0mm>*{}**@{-},
 <0.49mm,0.49mm>*{};<1.9mm,1.9mm>*{}**@{-},
 <-0.5mm,0.5mm>*{};<-1.9mm,1.9mm>*{}**@{-},
 <-2.3mm,2.3mm>*{\circ};<-2.3mm,2.3mm>*{}**@{},
 <-1.8mm,2.8mm>*{};<0mm,4.9mm>*{}**@{-},
 <-2.8mm,2.9mm>*{};<-4.6mm,4.9mm>*{}**@{-},
   <0.49mm,0.49mm>*{};<2.7mm,2.3mm>*{^3}**@{},
   <-1.8mm,2.8mm>*{};<0.4mm,5.3mm>*{^2}**@{},
   <-2.8mm,2.9mm>*{};<-5.1mm,5.3mm>*{^1}**@{},
 \end{xy}
\ + \
\begin{xy}
 <0mm,0mm>*{\circ};<0mm,0mm>*{}**@{},
 <0mm,-0.49mm>*{};<0mm,-3.0mm>*{}**@{-},
 <0.49mm,0.49mm>*{};<1.9mm,1.9mm>*{}**@{-},
 <-0.5mm,0.5mm>*{};<-1.9mm,1.9mm>*{}**@{-},
 <-2.3mm,2.3mm>*{\circ};<-2.3mm,2.3mm>*{}**@{},
 <-1.8mm,2.8mm>*{};<0mm,4.9mm>*{}**@{-},
 <-2.8mm,2.9mm>*{};<-4.6mm,4.9mm>*{}**@{-},
   <0.49mm,0.49mm>*{};<2.7mm,2.3mm>*{^2}**@{},
   <-1.8mm,2.8mm>*{};<0.4mm,5.3mm>*{^1}**@{},
   <-2.8mm,2.9mm>*{};<-5.1mm,5.3mm>*{^3}**@{},
 \end{xy}
\ + \
\begin{xy}
 <0mm,0mm>*{\circ};<0mm,0mm>*{}**@{},
 <0mm,-0.49mm>*{};<0mm,-3.0mm>*{}**@{-},
 <0.49mm,0.49mm>*{};<1.9mm,1.9mm>*{}**@{-},
 <-0.5mm,0.5mm>*{};<-1.9mm,1.9mm>*{}**@{-},
 <-2.3mm,2.3mm>*{\circ};<-2.3mm,2.3mm>*{}**@{},
 <-1.8mm,2.8mm>*{};<0mm,4.9mm>*{}**@{-},
 <-2.8mm,2.9mm>*{};<-4.6mm,4.9mm>*{}**@{-},
   <0.49mm,0.49mm>*{};<2.7mm,2.3mm>*{^1}**@{},
   <-1.8mm,2.8mm>*{};<0.4mm,5.3mm>*{^3}**@{},
   <-2.8mm,2.9mm>*{};<-5.1mm,5.3mm>*{^2}**@{},
 \end{xy}
\ \ \ \in \PROP\langle A \rangle(3,1)\\
&&\\
R_2:&&
%%%%%%%%%%%%%% Lie %%%%%%%%%%%%%%%%%%%%%%%%
 \begin{xy}
 <0mm,0mm>*{\bullet};<0mm,0mm>*{}**@{},
 <0mm,0.69mm>*{};<0mm,3.0mm>*{}**@{-},
 <0.39mm,-0.39mm>*{};<2.4mm,-2.4mm>*{}**@{-},
 <-0.35mm,-0.35mm>*{};<-1.9mm,-1.9mm>*{}**@{-},
 <-2.4mm,-2.4mm>*{\bullet};<-2.4mm,-2.4mm>*{}**@{},
 <-2.0mm,-2.8mm>*{};<0mm,-4.9mm>*{}**@{-},
 <-2.8mm,-2.9mm>*{};<-4.7mm,-4.9mm>*{}**@{-},
    <0.39mm,-0.39mm>*{};<3.3mm,-4.0mm>*{^3}**@{},
    <-2.0mm,-2.8mm>*{};<0.5mm,-6.7mm>*{^2}**@{},
    <-2.8mm,-2.9mm>*{};<-5.2mm,-6.7mm>*{^1}**@{},
 \end{xy}
\ + \
 \begin{xy}
 <0mm,0mm>*{\bullet};<0mm,0mm>*{}**@{},
 <0mm,0.69mm>*{};<0mm,3.0mm>*{}**@{-},
 <0.39mm,-0.39mm>*{};<2.4mm,-2.4mm>*{}**@{-},
 <-0.35mm,-0.35mm>*{};<-1.9mm,-1.9mm>*{}**@{-},
 <-2.4mm,-2.4mm>*{\bullet};<-2.4mm,-2.4mm>*{}**@{},
 <-2.0mm,-2.8mm>*{};<0mm,-4.9mm>*{}**@{-},
 <-2.8mm,-2.9mm>*{};<-4.7mm,-4.9mm>*{}**@{-},
    <0.39mm,-0.39mm>*{};<3.3mm,-4.0mm>*{^2}**@{},
    <-2.0mm,-2.8mm>*{};<0.5mm,-6.7mm>*{^1}**@{},
    <-2.8mm,-2.9mm>*{};<-5.2mm,-6.7mm>*{^3}**@{},
 \end{xy}
\ + \
 \begin{xy}
 <0mm,0mm>*{\bullet};<0mm,0mm>*{}**@{},
 <0mm,0.69mm>*{};<0mm,3.0mm>*{}**@{-},
 <0.39mm,-0.39mm>*{};<2.4mm,-2.4mm>*{}**@{-},
 <-0.35mm,-0.35mm>*{};<-1.9mm,-1.9mm>*{}**@{-},
 <-2.4mm,-2.4mm>*{\bullet};<-2.4mm,-2.4mm>*{}**@{},
 <-2.0mm,-2.8mm>*{};<0mm,-4.9mm>*{}**@{-},
 <-2.8mm,-2.9mm>*{};<-4.7mm,-4.9mm>*{}**@{-},
    <0.39mm,-0.39mm>*{};<3.3mm,-4.0mm>*{^1}**@{},
    <-2.0mm,-2.8mm>*{};<0.5mm,-6.7mm>*{^3}**@{},
    <-2.8mm,-2.9mm>*{};<-5.2mm,-6.7mm>*{^2}**@{},
 \end{xy}
\  \ \ \in   \PROP\langle A \rangle(1,3)  \\
&&\\
%%%%%%%%%%%%%%%%%%%%% %% Lie[1]Bi %%%%%%%%%%%%%%%
R_3:&&
 \begin{xy}
 <0mm,2.47mm>*{};<0mm,-0.5mm>*{}**@{-},
 <0.5mm,3.5mm>*{};<2.2mm,5.2mm>*{}**@{-},
 <-0.48mm,3.48mm>*{};<-2.2mm,5.2mm>*{}**@{-},
 <0mm,3mm>*{\circ};<0mm,3mm>*{}**@{},
  <0mm,-0.8mm>*{\bullet};<0mm,-0.8mm>*{}**@{},
<0mm,-0.8mm>*{};<-2.2mm,-3.5mm>*{}**@{-},
 <0mm,-0.8mm>*{};<2.2mm,-3.5mm>*{}**@{-},
     <0.5mm,3.5mm>*{};<2.8mm,5.7mm>*{^2}**@{},
     <-0.48mm,3.48mm>*{};<-2.8mm,5.7mm>*{^1}**@{},
   <0mm,-0.8mm>*{};<-2.7mm,-5.2mm>*{^1}**@{},
   <0mm,-0.8mm>*{};<2.7mm,-5.2mm>*{^2}**@{},
\end{xy}
\  - \
\begin{xy}
 <0mm,-1.3mm>*{};<0mm,-3.5mm>*{}**@{-},
 <0.38mm,-0.2mm>*{};<2.2mm,2.2mm>*{}**@{-},
 <-0.38mm,-0.2mm>*{};<-2.2mm,2.2mm>*{}**@{-},
<0mm,-0.8mm>*{\circ};<0mm,0.8mm>*{}**@{},
% <-2.25mm,2.2mm>*{};<-2.2mm,5.2mm>*{}**@{-},
 <2.4mm,2.4mm>*{\bullet};<2.4mm,2.4mm>*{}**@{},
 <2.5mm,2.3mm>*{};<4.4mm,-0.8mm>*{}**@{-},
% <4.4mm,-0.8mm>*{};<4.4mm,-3.5mm>*{}**@{-},
 <2.4mm,2.5mm>*{};<2.4mm,5.2mm>*{}**@{-},
     <0mm,-1.3mm>*{};<0mm,-5.3mm>*{^1}**@{},
     <2.5mm,2.3mm>*{};<5.1mm,-2.6mm>*{^2}**@{},
    <2.4mm,2.5mm>*{};<2.4mm,5.7mm>*{^2}**@{},
    <-0.38mm,-0.2mm>*{};<-2.8mm,2.5mm>*{^1}**@{},
    \end{xy}
\  + \
\begin{xy}
 <0mm,-1.3mm>*{};<0mm,-3.5mm>*{}**@{-},
 <0.38mm,-0.2mm>*{};<2.2mm,2.2mm>*{}**@{-},
 <-0.38mm,-0.2mm>*{};<-2.2mm,2.2mm>*{}**@{-},
<0mm,-0.8mm>*{\circ};<0mm,0.8mm>*{}**@{},
% <-2.25mm,2.2mm>*{};<-2.2mm,5.2mm>*{}**@{-},
 <2.4mm,2.4mm>*{\bullet};<2.4mm,2.4mm>*{}**@{},
 <2.5mm,2.3mm>*{};<4.4mm,-0.8mm>*{}**@{-},
% <4.4mm,-0.8mm>*{};<4.4mm,-3.5mm>*{}**@{-},
 <2.4mm,2.5mm>*{};<2.4mm,5.2mm>*{}**@{-},
     <0mm,-1.3mm>*{};<0mm,-5.3mm>*{^1}**@{},
     <2.5mm,2.3mm>*{};<5.1mm,-2.6mm>*{^2}**@{},
    <2.4mm,2.5mm>*{};<2.4mm,5.7mm>*{^1}**@{},
    <-0.38mm,-0.2mm>*{};<-2.8mm,2.5mm>*{^2}**@{},
    \end{xy}
\  - \
\begin{xy}
 <0mm,-1.3mm>*{};<0mm,-3.5mm>*{}**@{-},
 <0.38mm,-0.2mm>*{};<2.2mm,2.2mm>*{}**@{-},
 <-0.38mm,-0.2mm>*{};<-2.2mm,2.2mm>*{}**@{-},
<0mm,-0.8mm>*{\circ};<0mm,0.8mm>*{}**@{},
% <-2.25mm,2.2mm>*{};<-2.2mm,5.2mm>*{}**@{-},
 <2.4mm,2.4mm>*{\bullet};<2.4mm,2.4mm>*{}**@{},
 <2.5mm,2.3mm>*{};<4.4mm,-0.8mm>*{}**@{-},
% <4.4mm,-0.8mm>*{};<4.4mm,-3.5mm>*{}**@{-},
 <2.4mm,2.5mm>*{};<2.4mm,5.2mm>*{}**@{-},
     <0mm,-1.3mm>*{};<0mm,-5.3mm>*{^2}**@{},
     <2.5mm,2.3mm>*{};<5.1mm,-2.6mm>*{^1}**@{},
    <2.4mm,2.5mm>*{};<2.4mm,5.7mm>*{^2}**@{},
    <-0.38mm,-0.2mm>*{};<-2.8mm,2.5mm>*{^1}**@{},
    \end{xy}
\  + \
\begin{xy}
 <0mm,-1.3mm>*{};<0mm,-3.5mm>*{}**@{-},
 <0.38mm,-0.2mm>*{};<2.2mm,2.2mm>*{}**@{-},
 <-0.38mm,-0.2mm>*{};<-2.2mm,2.2mm>*{}**@{-},
<0mm,-0.8mm>*{\circ};<0mm,0.8mm>*{}**@{},
% <-2.25mm,2.2mm>*{};<-2.2mm,5.2mm>*{}**@{-},
 <2.4mm,2.4mm>*{\bullet};<2.4mm,2.4mm>*{}**@{},
 <2.5mm,2.3mm>*{};<4.4mm,-0.8mm>*{}**@{-},
% <4.4mm,-0.8mm>*{};<4.4mm,-3.5mm>*{}**@{-},
 <2.4mm,2.5mm>*{};<2.4mm,5.2mm>*{}**@{-},
     <0mm,-1.3mm>*{};<0mm,-5.3mm>*{^2}**@{},
     <2.5mm,2.3mm>*{};<5.1mm,-2.6mm>*{^1}**@{},
    <2.4mm,2.5mm>*{};<2.4mm,5.7mm>*{^1}**@{},
    <-0.38mm,-0.2mm>*{};<-2.8mm,2.5mm>*{^2}**@{},
    \end{xy}
\ \ \ \in\PROP\langle A \rangle(2,2). \Eeqrn }

\sip

\Proof The cohomology of
$(\PROP\langle\wedge_0^\bullet\cT\rangle,\delta)$ can not be
computed directly. At the dioperadic level the theorem was
established in \cite{Me1}. That this result extends to the level of
props  can be easily shown using either ideas of perturbations of
1/2props and  path filtrations developed in \cite{Ko2,MV} or the
idea of Koszul duality for props developed in \cite{V}. One can
argue, for example,
 as follows: for any $f\in \PROP\langle\wedge_0^\bullet\cT\rangle$, define the natural number,
$$
|f| :=\Ba{l}\mbox{number of  directed paths in the graph}\ f\ \\
       \mbox{which connect
     input legs with output ones.}
     \Ea
$$
and notice that the differential $\delta$ preserves the filtration,
$$
F_p \PROP\langle\wedge_0^\bullet\cT\rangle:= \left\{\mbox{span}\  f\in
  \PROP\langle\wedge_0^\bullet\cT\rangle  : \ |f|\leq p
 \right\}.
$$
The associated spectral sequence, $\{E_r\PROP\langle\wedge^\bullet_0\cT\rangle,\delta_r\}_{r\geq 0}$, is
exhaustive and bounded below
so that it  converges to
the cohomology of
$\left( \PROP\langle\wedge_0^\bullet\cT\rangle, \delta\right)$.

\bip

\no
By Koszulness of the operad ${\Lie}$ and exactness of the
functor $\Omega_{\frac{1}{2}\sP\rar \sP}$, the zeroth term of the
spectral sequence,
$(E_0\PROP\langle\wedge^\bullet_0\cT\rangle,\delta_0)$, is precisely
the minimal resolution of a quadratic prop, $\Liebi'$, generated by
the $\bS$-bimodule $A$ modulo the ideal generated by relations
$R_1$, $R_2$ and the following one,
$$
R_3':\ \
 \begin{xy}
 <0mm,2.47mm>*{};<0mm,-0.5mm>*{}**@{-},
 <0.5mm,3.5mm>*{};<2.2mm,5.2mm>*{}**@{-},
 <-0.48mm,3.48mm>*{};<-2.2mm,5.2mm>*{}**@{-},
 <0mm,3mm>*{\circ};<0mm,3mm>*{}**@{},
  <0mm,-0.8mm>*{\bullet};<0mm,-0.8mm>*{}**@{},
<0mm,-0.8mm>*{};<-2.2mm,-3.5mm>*{}**@{-},
 <0mm,-0.8mm>*{};<2.2mm,-3.5mm>*{}**@{-},
     <0.5mm,3.5mm>*{};<2.8mm,5.7mm>*{^2}**@{},
     <-0.48mm,3.48mm>*{};<-2.8mm,5.7mm>*{^1}**@{},
   <0mm,-0.8mm>*{};<-2.7mm,-5.2mm>*{^1}**@{},
   <0mm,-0.8mm>*{};<2.7mm,-5.2mm>*{^2}**@{},
\end{xy} =0.
$$
As the differential $\delta$ vanishes on the generators of $A$,
 this spectral sequence degenerates at the first term, $(E_1\PROP\langle\wedge^\bullet_0\cT\rangle, d_1=0)$, implying the isomorphism
$$
\bigoplus_{p\geq 1} \frac{F_{p+1}H(\PROP\langle\wedge^\bullet_0\cT\rangle,\delta)}
{F_{p}H(\PROP\langle\wedge^\bullet_0\cT\rangle,\delta)} = \Liebi'.
$$
There is a natural surjective morphism of dg props,  $p:(\PROP\langle\wedge^\bullet_0\cT\rangle,\delta) \rar \Liebi$.
Define the dg prop $(X,\delta)$ via an exact sequence,
$$
0\lon (X,\delta) \stackrel{i}{\lon} (\PROP\langle\wedge^\bullet_0\cT\rangle,\delta) \stackrel{p}{\lon}
(\Liebi, 0)\lon 0.
$$
The filtration on $(\PROP\langle\wedge^\bullet_0\cT\rangle,\delta)$ induces filtrations on sub- and quotient
complexes,
$$
0\lon (F_pX,\delta) \stackrel{i}{\lon} (F_p\PROP\langle\wedge^\bullet_0\cT\rangle,\delta) \stackrel{p}{\lon}
(F_p\Liebi, 0)\lon 0,
$$
and hence an exact sequence of 0th terms of the associated spectral sequences,
$$
0\lon (E_0X,\delta_0) \stackrel{i_0}{\lon} (E_0\PROP\langle\wedge^\bullet_0\cT\rangle,\delta_0) \stackrel{p_0}{\lon}
(\bigoplus_{p\geq 1} \frac{F_{p+1}\Liebi}{F_p\Liebi}, 0)\lon 0.
$$
By the above observation,
$$
E_1\PROP\langle\wedge^\bullet_0\cT\rangle=\bigoplus_{p\geq 1} \frac{F_{p+1}H(\PROP\langle\wedge^\bullet_0\cT\rangle,\delta)}
{F_{p}H(\PROP\langle\wedge^\bullet_0\cT\rangle,\delta)} = \Liebi'.
$$
On the other hand, it is not hard to check that
$$
\bigoplus_{p\geq 1} \frac{F_{p+1}\Liebi}{F_p\Liebi}= \Liebi'.
$$
Thus the map $p_0$ is a quasi-isomorphism implying vanishing of $E_1X$ and hence acyclicity of $(X,\delta)$.
Thus the projection map $p$ is a quasi-isomorphism.
 \hfill
$\Box$

\bip

%%%%%%%%%%%%%%%%%%%%%%%%%
\no{\bf 2.6.3. Corollary.} The dg prop
 $(\PROP\langle\wedge^\bullet_0\cT\rangle, \delta)$
is a minimal model
 of the prop ${\Liebi}$: the natural morphism of dg props,
$$
p: (\PROP\langle\wedge^\bullet_0\cT\rangle, \delta)\lon ({\Liebi},
\mbox{vanishing differential}).
$$
 which sends to zero
all generators of $ \PROP\langle\wedge^\bullet_0\cT\rangle$ except those
in $A(2,1)$ and $A(1,2)$, is a quasi-isomorphism. Hence we can and shall re-denote
 $ \PROP\langle\wedge^\bullet_0\cT\rangle$ as ${\Liebi_\infty}$.

\bip

\no{\bf 2.6.4. Definition.} A representation of the dg prop $\Liebi_\infty$ in a dg space $V$ is called
a {\em Poisson structure}\, on the formal graded manifold $V$.

\bip

\no
Thus Poisson structure on $V$ is the same as a Maurer-Cartan element, $\ga\in \wedge^\bullet\cT_V$,
in the Lie algebra of formal polyvector fields on $V$ satisfying the conditions. $[\ga,\ga]=0$,
and $\ga|_{x=0}=0$.

\bip

\no{\bf 2.6.5. Remark.} The condition $\ga|_{x=0}=0$ above is no serious restriction: given an arbitrary Poisson
structure $\pi$ on $V$ (not necessary  vanishing at $0\in V$), then, for any parameter $\hbar$ viewed as a coordinate
on $\K$, the product
$\ga:=\hbar\pi$ is a Poisson structure on $V=V\times \K$ vanishing at zero $0\in V$ and hence
is a representation of the prop $\Liebi_\infty$.

\bip

%%%%%%%%%%%%%%%%%%%%%%%%%%%%%%%%%%%%%
\no{\bf 2.7. Prop profile of ``star products".} Let us consider a sheaf
of dg Lie algebras,
$$
\caD_\cM\subset \bigoplus_{k\geq 0} \Hom(\f_\cM^{\ot k},\f_\cM)[1-k],
$$
consisting of polydifferential operators on $\f_\cM$ which, for $k\geq 1$,
vanish on every
element $f_1\ot\ldots \ot f_k\in \f_\cM^{\ot k}$ with at least one function,
 $f_i$,
$i=1,\ldots,k$, constant. A degree one
section, $\Gamma\in \caD_\cM$, decomposes into a sum, $\sum_{k\geq 0}
\Gamma_k$,  with $\Gamma_k\in \Hom_{2-k}(\f_\cM^{\ot k},\f_\cM)$. In a local
coordinate system, $(t^\al, \p/\p t^\al\simeq e_\al)$ on $\cM$, $\Gamma$
can be represented as a Taylor series,
$$
\Gamma= \sum_{k\geq 0}\sum_{I_1,\ldots, I_k, J}
\Gamma^{I_1,\ldots, I_k}_J e_{I_1}\ot \ldots e_{I_k} \ot t^J,
$$
where for each fixed $k$ and $|J|$ only a finite number of coefficients
$\Gamma^{I_1,\ldots, I_k}_J$ is non-zero.
The summation runs over  multi-indices, $I:=\al_1\al_2\ldots
\al_{|I|}$, and $e_{I}:=e_{\al_1}\odot \ldots \odot e_{\al_{|I|}}$,
$t^{I}:= t^{\al_1}\odot \ldots t^{\al_{|I|}}$. Hence the associated
$\bS$-bimodule $E_\caD$ is given by
$$
E_\caD(m,n)= {\bf E}(m)\ot {\bf 1_n}, \ \ m, n\geq 0,
$$
where
\Beqrn
{\bf E}(0)&:=& \R[-2], \\
{\bf E}(m\geq 1)&:=&\bigoplus_{k\geq 1}\bigoplus_{
[m]=I_1\sqcup \ldots \sqcup I_k \atop |I_1|,\ldots,|I_k|\geq 1}
{\rm Ind}^{\bS_m}_{\bS_{|I_1|}\times
\ldots \times \bS_{|I_k|}} {\bf 1}_{|I_1|}\ot \ldots \ot {\bf 1}_{|I_k|}
[k-2]
\Eeqrn
Let $\DefQ$ stand for the free prop, $\PROP\langle E_\caD\rangle$, generated by the
above bimodule. The generators of $\DefQ$ can be identified with
directed planar
corollas of the form,
$$
\xy
%\begin{xy}
 <0mm,0mm>*{\mbox{$\xy *=<20mm,3mm>\txt{}*\frm{-}\endxy$}};<0mm,0mm>*{}**@{},
  <-10mm,1.5mm>*{};<-12mm,7mm>*{}**@{-},
  <-10mm,1.5mm>*{};<-11mm,7mm>*{}**@{-},
  <-10mm,1.5mm>*{};<-9.5mm,6mm>*{}**@{-},
  <-10mm,1.5mm>*{};<-8mm,7mm>*{}**@{-},
 <-10mm,1.5mm>*{};<-9.5mm,6.6mm>*{.\hspace{-0.4mm}.\hspace{-0.4mm}.}**@{},
 <0mm,0mm>*{};<-6.5mm,3.6mm>*{.\hspace{-0.1mm}.\hspace{-0.1mm}.}**@{},
  <-3mm,1.5mm>*{};<-5mm,7mm>*{}**@{-},
  <-3mm,1.5mm>*{};<-4mm,7mm>*{}**@{-},
  <-3mm,1.5mm>*{};<-2.5mm,6mm>*{}**@{-},
  <-3mm,1.5mm>*{};<-1mm,7mm>*{}**@{-},
 <-3mm,1.5mm>*{};<-2.5mm,6.6mm>*{.\hspace{-0.4mm}.\hspace{-0.4mm}.}**@{},
  <2mm,1.5mm>*{};<0mm,7mm>*{}**@{-},
  <2mm,1.5mm>*{};<1mm,7mm>*{}**@{-},
  <2mm,1.5mm>*{};<2.5mm,6mm>*{}**@{-},
  <2mm,1.5mm>*{};<4mm,7mm>*{}**@{-},
 <2mm,1.5mm>*{};<2.5mm,6.6mm>*{.\hspace{-0.4mm}.\hspace{-0.4mm}.}**@{},
 <0mm,0mm>*{};<6mm,3.6mm>*{.\hspace{-0.1mm}.\hspace{-0.1mm}.}**@{},
<10mm,1.5mm>*{};<8mm,7mm>*{}**@{-},
  <10mm,1.5mm>*{};<9mm,7mm>*{}**@{-},
  <10mm,1.5mm>*{};<10.5mm,6mm>*{}**@{-},
  <10mm,1.5mm>*{};<12mm,7mm>*{}**@{-},
 <10mm,1.5mm>*{};<10.5mm,6.6mm>*{.\hspace{-0.4mm}.\hspace{-0.4mm}.}**@{},
 <-10mm,-1.5mm>*{};<-12mm,-6mm>*{}**@{-},
 <-7mm,-1.5mm>*{};<-8mm,-6mm>*{}**@{-},
 <-4mm,-1.5mm>*{};<-4.5mm,-6mm>*{}**@{-},
 <0mm,0mm>*{};<0mm,-4.6mm>*{.\hspace{0.1mm}.\hspace{0.1mm}.}**@{},
<10mm,-1.5mm>*{};<12mm,-6mm>*{}**@{-},
 <7mm,-1.5mm>*{};<8mm,-6mm>*{}**@{-},
  <4mm,-1.5mm>*{};<4.5mm,-6mm>*{}**@{-},
  %%%%
<0mm,0mm>*{};<-9.5mm,8.2mm>*{^{I_{ 1}}}**@{},
<0mm,0mm>*{};<-3mm,8.2mm>*{^{I_{ i}}}**@{},
<0mm,0mm>*{};<2mm,8.2mm>*{^{I_{ i+1}}}**@{},
<0mm,0mm>*{};<10mm,8.2mm>*{^{I_{ k}}}**@{},
 %%%
<0mm,0mm>*{};<-12mm,-7.4mm>*{_1}**@{},
<0mm,0mm>*{};<-8mm,-7.4mm>*{_2}**@{},
<0mm,0mm>*{};<-4mm,-7.4mm>*{_3}**@{},
<0mm,0mm>*{};<6mm,-7.4mm>*{\ldots}**@{},
%<0mm,0mm>*{};<8.5mm,-7.4mm>*{_{n-1}}**@{},
<0mm,0mm>*{};<12.5mm,-7.4mm>*{_n}**@{},
\endxy
\ \ \ \simeq \ \ \ \Gamma^{I_1,\ldots, I_k}_J,
$$
where
\Bi
\item the input legs are labeled by the set
$[n]:=\{1,2,\ldots,n\}$ and  are symmetric (so that it does not matter
how labels from $[n]$ are distributed over them),
\item the output legs (if there are any) are labeled by the set $[m]$
partitioned into $k$ disjoint non-empty subsets,
$$
[m]=I_1\sqcup \ldots\sqcup I_i\sqcup I_{i+1}\sqcup\ldots \sqcup I_k,
$$
and legs in each $I_i$-bunch are symmetric (so that it does not matter
how labels from the set $I_i$ are distributed over legs in $I_i$th bunch).
\Ei
 The
$\Z$-grading in $\DefQ$ is defined by associating degree
 $2-k$ to such a corolla. The formula $(\star\star)$ in Sect.\ 2.5 provides
 us with the following explicit expression for the differential, $\delta$,
 in $\DefQ$,
\Beqrn
\delta\left(\xy
%\begin{xy}
 <0mm,0mm>*{\mbox{$\xy *=<20mm,3mm>\txt{}*\frm{-}\endxy$}};<0mm,0mm>*{}**@{},
  <-10mm,1.5mm>*{};<-12mm,7mm>*{}**@{-},
  <-10mm,1.5mm>*{};<-11mm,7mm>*{}**@{-},
  <-10mm,1.5mm>*{};<-9.5mm,6mm>*{}**@{-},
  <-10mm,1.5mm>*{};<-8mm,7mm>*{}**@{-},
 <-10mm,1.5mm>*{};<-9.5mm,6.6mm>*{.\hspace{-0.4mm}.\hspace{-0.4mm}.}**@{},
 <0mm,0mm>*{};<-6.4mm,3.6mm>*{.\hspace{-0.1mm}.\hspace{-0.1mm}.}**@{},
  <-3mm,1.5mm>*{};<-5mm,7mm>*{}**@{-},
  <-3mm,1.5mm>*{};<-4mm,7mm>*{}**@{-},
  <-3mm,1.5mm>*{};<-2.5mm,6mm>*{}**@{-},
  <-3mm,1.5mm>*{};<-1mm,7mm>*{}**@{-},
 <-3mm,1.5mm>*{};<-2.5mm,6.6mm>*{.\hspace{-0.4mm}.\hspace{-0.4mm}.}**@{},
  <2mm,1.5mm>*{};<0mm,7mm>*{}**@{-},
  <2mm,1.5mm>*{};<1mm,7mm>*{}**@{-},
  <2mm,1.5mm>*{};<2.5mm,6mm>*{}**@{-},
  <2mm,1.5mm>*{};<4mm,7mm>*{}**@{-},
 <2mm,1.5mm>*{};<2.5mm,6.6mm>*{.\hspace{-0.4mm}.\hspace{-0.4mm}.}**@{},
 <0mm,0mm>*{};<6mm,3.6mm>*{.\hspace{-0.1mm}.\hspace{-0.1mm}.}**@{},
<10mm,1.5mm>*{};<8mm,7mm>*{}**@{-},
  <10mm,1.5mm>*{};<9mm,7mm>*{}**@{-},
  <10mm,1.5mm>*{};<10.5mm,6mm>*{}**@{-},
  <10mm,1.5mm>*{};<12mm,7mm>*{}**@{-},
 <10mm,1.5mm>*{};<10.5mm,6.6mm>*{.\hspace{-0.4mm}.\hspace{-0.4mm}.}**@{},
 <-10mm,-1.5mm>*{};<-12mm,-6mm>*{}**@{-},
 <-7mm,-1.5mm>*{};<-8mm,-6mm>*{}**@{-},
 <-4mm,-1.5mm>*{};<-4.5mm,-6mm>*{}**@{-},
 <0mm,0mm>*{};<0mm,-4.6mm>*{.\hspace{0.1mm}.\hspace{0.1mm}.}**@{},
<10mm,-1.5mm>*{};<12mm,-6mm>*{}**@{-},
 <7mm,-1.5mm>*{};<8mm,-6mm>*{}**@{-},
  <4mm,-1.5mm>*{};<4.5mm,-6mm>*{}**@{-},
  %%%%
<0mm,0mm>*{};<-9.5mm,8.2mm>*{^{I_{ 1}}}**@{},
<0mm,0mm>*{};<-3mm,8.2mm>*{^{I_{ i}}}**@{},
<0mm,0mm>*{};<2mm,8.2mm>*{^{I_{ i+1}}}**@{},
<0mm,0mm>*{};<10mm,8.2mm>*{^{I_{ k}}}**@{},
 %%%
<0mm,0mm>*{};<-12mm,-7.4mm>*{_1}**@{},
<0mm,0mm>*{};<-8mm,-7.4mm>*{_2}**@{},
<0mm,0mm>*{};<-4mm,-7.4mm>*{_3}**@{},
<0mm,0mm>*{};<6mm,-7.4mm>*{\ldots}**@{},
%<0mm,0mm>*{};<8.5mm,-7.4mm>*{_{n-1}}**@{},
<0mm,0mm>*{};<12.5mm,-7.4mm>*{_n}**@{},
\endxy
\right) &=& \sum_{i=1}^k(-1)^{i+1}
%%%%%%%%%%%%%%%%%%%%%%%%%%%%%%%%%%%%%%%%%%%%%%%%%%%%
\xy
%\begin{xy}
 <0mm,0mm>*{\mbox{$\xy *=<20mm,3mm>\txt{}*\frm{-}\endxy$}};<0mm,0mm>*{}**@{},
  <-10mm,1.5mm>*{};<-12mm,7mm>*{}**@{-},
  <-10mm,1.5mm>*{};<-11mm,7mm>*{}**@{-},
  <-10mm,1.5mm>*{};<-9.5mm,6mm>*{}**@{-},
  <-10mm,1.5mm>*{};<-8mm,7mm>*{}**@{-},
 <-10mm,1.5mm>*{};<-9.5mm,6.6mm>*{.\hspace{-0.4mm}.\hspace{-0.4mm}.}**@{},
 <0mm,0mm>*{};<-5.5mm,3.6mm>*{.\hspace{-0.1mm}.\hspace{-0.1mm}.}**@{},
%
  %<0mm,1.5mm>*{};<-5mm,7mm>*{}**@{-},
  <0mm,1.5mm>*{};<-4mm,7mm>*{}**@{-},
  <0mm,1.5mm>*{};<-2.0mm,6mm>*{}**@{-},
  <0mm,1.5mm>*{};<-1mm,7mm>*{}**@{-},
 <0mm,1.5mm>*{};<-2.3mm,6.6mm>*{.\hspace{-0.4mm}.\hspace{-0.4mm}.}**@{},
%
  %<0mm,1.5mm>*{};<0mm,7mm>*{}**@{-},
  <0mm,1.5mm>*{};<1mm,7mm>*{}**@{-},
  <0mm,1.5mm>*{};<2.0mm,6mm>*{}**@{-},
  <0mm,1.5mm>*{};<4mm,7mm>*{}**@{-},
 <0mm,1.5mm>*{};<2.3mm,6.6mm>*{.\hspace{-0.4mm}.\hspace{-0.4mm}.}**@{},
 <0mm,0mm>*{};<6mm,3.6mm>*{.\hspace{-0.1mm}.\hspace{-0.1mm}.}**@{},
<10mm,1.5mm>*{};<8mm,7mm>*{}**@{-},
  <10mm,1.5mm>*{};<9mm,7mm>*{}**@{-},
  <10mm,1.5mm>*{};<10.5mm,6mm>*{}**@{-},
  <10mm,1.5mm>*{};<12mm,7mm>*{}**@{-},
 <10mm,1.5mm>*{};<10.5mm,6.6mm>*{.\hspace{-0.4mm}.\hspace{-0.4mm}.}**@{},
 <-10mm,-1.5mm>*{};<-12mm,-6mm>*{}**@{-},
 <-7mm,-1.5mm>*{};<-8mm,-6mm>*{}**@{-},
 <-4mm,-1.5mm>*{};<-4.5mm,-6mm>*{}**@{-},
 <0mm,0mm>*{};<0mm,-4.6mm>*{.\hspace{0.1mm}.\hspace{0.1mm}.}**@{},
<10mm,-1.5mm>*{};<12mm,-6mm>*{}**@{-},
 <7mm,-1.5mm>*{};<8mm,-6mm>*{}**@{-},
  <4mm,-1.5mm>*{};<4.5mm,-6mm>*{}**@{-},
  %%%%
<0mm,0mm>*{};<-9.5mm,8.2mm>*{^{I_{ 1}}}**@{},
<0mm,0mm>*{};<-2.5mm,8.2mm>*{{^{I_{ i}\sqcup}}}**@{},
<0mm,0mm>*{};<2.7mm,8.2mm>*{^{I_{ i+1}}}**@{},
<0mm,0mm>*{};<10mm,8.2mm>*{^{I_{ k}}}**@{},
 %%%
<0mm,0mm>*{};<-12mm,-7.4mm>*{_1}**@{},
<0mm,0mm>*{};<-8mm,-7.4mm>*{_2}**@{},
<0mm,0mm>*{};<-4mm,-7.4mm>*{_3}**@{},
<0mm,0mm>*{};<6mm,-7.4mm>*{\ldots}**@{},
%<0mm,0mm>*{};<8.5mm,-7.4mm>*{_{n-1}}**@{},
<0mm,0mm>*{};<12.5mm,-7.4mm>*{_n}**@{},
\endxy
\\
&& +\ \sum_{p+q=k+1\atop p\geq 1,q\geq 0}\sum_{i=0}^{p-1}
\sum_{  {I_{i+1}=I_{i+1}'\sqcup I''_{i+1}\atop .......................}
\atop
I_{i+q}=I_{i+q}'\sqcup I''_{i+q}}
\sum_{[n]=J_1\sqcup J_2}\sum_{s\geq 0}(-1)^{(p+1)q + i(q-1)}\\
&&
%%%%%%%%%%%%%%%%%%%%%%%%%%%%%%%%%%%%%%%%%%%%%%%%%%%%%%%%%%%
\frac{1}{s!}\ \
\xy
%\begin{xy}
 <19mm,0mm>*{\mbox{$\xy *=<58mm,3mm>\txt{}*\frm{-}\endxy$}};<0mm,0mm>*{}**@{},
  <-10mm,1.5mm>*{};<-12mm,7mm>*{}**@{-},
  <-10mm,1.5mm>*{};<-11mm,7mm>*{}**@{-},
  <-10mm,1.5mm>*{};<-9.5mm,6mm>*{}**@{-},
  <-10mm,1.5mm>*{};<-8mm,7mm>*{}**@{-},
 <-10mm,1.5mm>*{};<-9.5mm,6.6mm>*{.\hspace{-0.4mm}.\hspace{-0.4mm}.}**@{},
 <0mm,0mm>*{};<-6.5mm,3.6mm>*{.\hspace{-0.1mm}.\hspace{-0.1mm}.}**@{},
  <-3mm,1.5mm>*{};<-5mm,7mm>*{}**@{-},
  <-3mm,1.5mm>*{};<-4mm,7mm>*{}**@{-},
  <-3mm,1.5mm>*{};<-2.5mm,6mm>*{}**@{-},
  <-3mm,1.5mm>*{};<-1mm,7mm>*{}**@{-},
 <-3mm,1.5mm>*{};<-2.5mm,6.6mm>*{.\hspace{-0.4mm}.\hspace{-0.4mm}.}**@{},
%%%%%
  <10mm,1.5mm>*{};<0mm,7mm>*{}**@{-},
  <10mm,1.5mm>*{};<4mm,7mm>*{}**@{-},
<10mm,1.5mm>*{};<7.3mm,5.9mm>*{.\hspace{-0.0mm}.\hspace{-0.0mm}.}**@{},
   <10mm,1.5mm>*{};<3.8mm,6.0mm>*{}**@{-},
 <10mm,1.5mm>*{};<2.5mm,6.6mm>*{.\hspace{-0.4mm}.\hspace{-0.4mm}.}**@{},
%
% <0mm,0mm>*{};<6mm,3.6mm>*{.\hspace{-0.1mm}.\hspace{-0.1mm}.}**@{},
%
<10mm,1.5mm>*{};<9mm,7mm>*{}**@{-},
  <10mm,1.5mm>*{};<10.5mm,6mm>*{}**@{-},
  <10mm,1.5mm>*{};<12mm,7mm>*{}**@{-},
 <10mm,1.5mm>*{};<10.5mm,6.6mm>*{.\hspace{-0.4mm}.\hspace{-0.4mm}.}**@{},
 %%%%%%%%%%%%%%%%%%%%%%%%%%%%%%%%%%%%%%%
%%%% lower legs of 1st vertex %%%%%%%
 <-10mm,-1.5mm>*{};<-12mm,-6mm>*{}**@{-},
 <-7mm,-1.5mm>*{};<-8mm,-6mm>*{}**@{-},
 <-4mm,-1.5mm>*{};<-4.5mm,-6mm>*{}**@{-},
 <-1mm,-1.5mm>*{};<-1.1mm,-6mm>*{}**@{-},
 <2mm,-1.5mm>*{};<2.0mm,-6mm>*{}**@{-},
 <0mm,0mm>*{};<20mm,-4.6mm>*{.\hspace{2mm}.\hspace{2mm}.\hspace{2mm}
 .\hspace{2mm}.\hspace{2mm}.\hspace{2mm}.\hspace{2mm}
 .\hspace{2mm}.\hspace{2mm}}**@{},
<48mm,-1.5mm>*{};<50mm,-6mm>*{}**@{-},
 <45mm,-1.5mm>*{};<46mm,-6mm>*{}**@{-},
  <42mm,-1.5mm>*{};<42.5mm,-6mm>*{}**@{-},
  <39mm,-1.5mm>*{};<39.2mm,-6mm>*{}**@{-},
   <36mm,-1.5mm>*{};<36mm,-6mm>*{}**@{-},
 <20mm,-1.5mm>*{};<20.0mm,-8mm>*{\underbrace{\hspace{66mm}}}**@{},
 <20mm,-1.5mm>*{};<20.0mm,-11mm>*{_{J_1}}**@{},
  %%%%
<0mm,0mm>*{};<-9.5mm,8.4mm>*{^{I_{ 1}}}**@{},
<0mm,0mm>*{};<-3mm,8.4mm>*{^{I_{ i}}}**@{},
<0mm,0mm>*{};<2mm,8.6mm>*{^{I_{ i+1}'}}**@{},
<0mm,0mm>*{};<10.5mm,8.6mm>*{^{I_{ i+q}'}}**@{},
%%%% edges connecting two vertices %%%%%%%%%%
<10mm,1.5mm>*{};<18mm,12mm>*{}**@{-},
<10mm,1.5mm>*{};<20.0mm,12mm>*{}**@{-},
<10mm,1.5mm>*{};<25mm,12mm>*{}**@{-},
<10mm,1.5mm>*{};<18.7mm,8.6mm>*{.\hspace{-0.4mm}.\hspace{-0.4mm}.}**@{},
<10mm,1.5mm>*{};<19.6mm,9.0mm>*{^s}**@{},
%%%%%%%%%%%%%%%%%%%%%%%%%%%%%%%%%%%%%%%%%%%%%%
%%%% second vertex %%%%%%%%%%%%%%%%%%%%%%%%%%
%%%%%%%%%%%%%%%%%%%%%%%%%%%%%%%%%%%%%%%%%%%%%
<25mm,13.75mm>*{\mbox{$\xy *=<14mm,3mm>\txt{}*\frm{-}\endxy$}};
<0mm,0mm>*{}**@{},
 <18mm,15mm>*{};<16mm,20.5mm>*{}**@{-},
 <18mm,15mm>*{};<17mm,20.5mm>*{}**@{-},
 <18mm,15mm>*{};<18.5mm,19.6mm>*{}**@{-},
 <18mm,15mm>*{};<20mm,20.5mm>*{}**@{-},
 <18mm,15mm>*{};<18.6mm,20.3mm>*{.\hspace{-0.4mm}.\hspace{-0.4mm}.}**@{},
<22mm,15mm>*{};<25.5mm,17.7mm>*{\cdots}**@{},
 <32mm,15.2mm>*{};<30mm,20.5mm>*{}**@{-},
 <32mm,15.2mm>*{};<31mm,20.5mm>*{}**@{-},
 <32mm,15.2mm>*{};<32.5mm,19.6mm>*{}**@{-},
 <32mm,15mm>*{};<34mm,20.5mm>*{}**@{-},
 <32mm,15mm>*{};<32.3mm,20.3mm>*{.\hspace{-0.4mm}.\hspace{-0.4mm}.}**@{},
%
  %%%%
<0mm,0mm>*{};<18mm,22.6mm>*{^{I_{ i+1}''}}**@{},
<0mm,0mm>*{};<32.5mm,22.6mm>*{^{I_{ i+q}''}}**@{},
  %%%%
%%%% lower legs of 2nd vertex %%%%%%
 <26mm,12mm>*{};<25mm,9mm>*{}**@{-},
 <27mm,12mm>*{};<26.8mm,9mm>*{}**@{-},
 <29mm,12mm>*{};<29.2mm,10mm>*{.\hspace{-0.1mm}.\hspace{-0.1mm}.}**@{},
 <32mm,12mm>*{};<33mm,9mm>*{}**@{-},
  <31mm,12mm>*{};<31.5mm,9mm>*{}**@{-},
 <29mm,12mm>*{};<29mm,8mm>*{\underbrace{\ \ \ \ \ \ \ \  }}**@{},
 <29mm,12mm>*{};<29mm,5.3mm>*{_{J_2}}**@{},
%%%%%%%%%%%%%%%%%% 1st vertex contd %%%%%%%%%%%
<38mm,1.5mm>*{};<36mm,7mm>*{}**@{-},
<38mm,1.5mm>*{};<37mm,7mm>*{}**@{-},
<38mm,1.5mm>*{};<38.5mm,6mm>*{}**@{-},
<38mm,1.5mm>*{};<40mm,7mm>*{}**@{-},
<38mm,1.5mm>*{};<38.5mm,6.6mm>*{.\hspace{-0.4mm}.\hspace{-0.4mm}.}**@{},
<38mm,1.5mm>*{};<43mm,4mm>*{.\hspace{-0.0mm}.\hspace{-0.0mm}.}**@{},
<48mm,1.5mm>*{};<46mm,7mm>*{}**@{-},
<48mm,1.5mm>*{};<47mm,7mm>*{}**@{-},
<48mm,1.5mm>*{};<48.5mm,6mm>*{}**@{-},
<48mm,1.5mm>*{};<50mm,7mm>*{}**@{-},
<48mm,1.5mm>*{};<48.5mm,6.6mm>*{.\hspace{-0.4mm}.\hspace{-0.4mm}.}**@{},
<0mm,0mm>*{};<40.3mm,8.6mm>*{^{I_{ i+q+1}}}**@{},
<0mm,0mm>*{};<48.5mm,8.6mm>*{^{I_{ k}}}**@{},
  %%%%
\endxy
\Eeqrn where the first sum comes from the Hochschild differential
$d_{\rm H}$ and the second sum comes from the Hochschild brackets
$[\ ,\ ]_{\rm H}$. The $s$-summation in the latter runs over the
number, $s$, of  edges connecting the two internal vertices. As $s$
can be zero, the r.h.s. above contains disconnected graphs (more
precisely, disjoint unions of two corollas).

\bip

%%%%%%%%%%%%%%%%%%%%%%%%%
\no{\bf 2.7.1. Proposition.} {\em There is a one-to-one correspondence between
representations,
$$
\phi: (\DefQ,\delta) \lon ({\sf End}\langle M\rangle,d),
$$
of $(\DefQ,\delta)$
 in a dg vector space
$(M,d)$
and Maurer-Cartan elements, $\gamma$, in $\caD_\cM$, that is, degree
one elements satisfying the equation, $d_{\rm H}\ga+
\frac{1}{2}[\ga,\ga]_{\rm H}=0$.}

\bip

\no{\bf Proof} is similar to the proof of Proposition 2.6.1.

\bip

\no{\bf 2.8. Remark.} Kontsevich's formality map \cite{Ko} can be
interpreted as a morphism of dg props,
$$
F_{\infty}: (\DefQ,\delta) \lon
(\PROP\langle\wedge^\bullet \cT \rangle^\circlearrowright, \delta).
$$
Vice versa, any morphism of the above dg props gives rise
 to a {\em universal}\, formality map in the sense of
\cite{Ko}. Note that the dg prop of polyvector fields,
$\PROP\langle\wedge^\bullet \cT \rangle$, appears above in the wheel
extended  form, $\PROP\langle\wedge^\bullet \cT
\rangle^\circlearrowright$. This is not accidental: we shall show below in \S 5 (by quantizing a pair consisting
of a linear Poisson structure and quadratic homological vector field) that there does {\em not}\, exist a morphism between
ordinary (i.e.\ unwheeled) dg props, $(\DefQ,\delta) \lon
(\PROP\langle\wedge^\bullet \cT \rangle, \delta)$, satisfying the
quasi-classical limit condition.

\bip

%%%%%%%%%%%%%%%%%%%%%%%%%%%%%%%%%%%%%%%%%%%%%%%%%%%%%%%%%%%%%%%%%%%%%%%%%%
\no{\bf 2.9. Prop profile of perturbative ``star products"}.
Let $\K[[\hbar]]$ be the formal power series in a formal parameter
$\hbar$, and let $\hbar\K[[\hbar]]$ be its (maximal) ideal spanned
by series which vanish at $\hbar=0$. Then $\caD^\hbar_V:= \caD_V\ot
\hbar\K[[\hbar]]$ is a dg Lie algebra of polydifferential operators
on $\f_V[[\hbar]]$ which vanish at $\hbar=0$. A solution, $\Ga$, of
Maurer-Cartan equations in $\caD^\hbar_V$ gives a (generalized)
$A_\infty$-structure on $\f_V[[\hbar]]$ which at $\hbar=0$ reduces
to the ordinary graded commutative multiplication in $\f_V$. Thus
Maurer-Cartan elements in this  Lie algebra describe perturbative
deformations of the ordinary product in $\f_V$. Again, the set of
all possible Maurer-Cartan elements in $\caD^\hbar_V$ can be
understood as the set of all possible representations of a certain
dg free prop, $\DefQh:=\PROP\langle\caD^\hbar \rangle$, defined as
follows:

\Bi
\item the $\bS$-bimodule of generators, $E_{\caD^\hbar}=\{E_{\caD^\hbar}(m,n)\}$, is a direct sum,
$E_{\caD^\hbar}(m,n)= \bigoplus_{a=1}^\infty E_{\caD}^a(m,n)$, $m,
n\geq 0$, of labeled by  $a\in \N^*$ copies,
$E_{\caD}^a(m,n):={E_{\caD}}(m,n)$, of the $\bS$-module
corresponding to the sheaf $\caD_V$; each copy corresponds to the ``$\hbar^a$" coefficient
 in the formal power series;
\item generators of
$\DefQh$ can be identified with planar corollas, \
$
\xy
%\begin{xy}
 <0mm,0mm>*{\mbox{$\xy *=<20mm,3mm>\txt{\em a}*\frm{-}\endxy$}};<0mm,0mm>*{}**@{},
  <-10mm,1.5mm>*{};<-12mm,7mm>*{}**@{-},
  <-10mm,1.5mm>*{};<-11mm,7mm>*{}**@{-},
  <-10mm,1.5mm>*{};<-9.5mm,6mm>*{}**@{-},
  <-10mm,1.5mm>*{};<-8mm,7mm>*{}**@{-},
 <-10mm,1.5mm>*{};<-9.5mm,6.6mm>*{.\hspace{-0.4mm}.\hspace{-0.4mm}.}**@{},
 <0mm,0mm>*{};<-6.5mm,3.6mm>*{.\hspace{-0.1mm}.\hspace{-0.1mm}.}**@{},
  <-3mm,1.5mm>*{};<-5mm,7mm>*{}**@{-},
  <-3mm,1.5mm>*{};<-4mm,7mm>*{}**@{-},
  <-3mm,1.5mm>*{};<-2.5mm,6mm>*{}**@{-},
  <-3mm,1.5mm>*{};<-1mm,7mm>*{}**@{-},
 <-3mm,1.5mm>*{};<-2.5mm,6.6mm>*{.\hspace{-0.4mm}.\hspace{-0.4mm}.}**@{},
  <2mm,1.5mm>*{};<0mm,7mm>*{}**@{-},
  <2mm,1.5mm>*{};<1mm,7mm>*{}**@{-},
  <2mm,1.5mm>*{};<2.5mm,6mm>*{}**@{-},
  <2mm,1.5mm>*{};<4mm,7mm>*{}**@{-},
 <2mm,1.5mm>*{};<2.5mm,6.6mm>*{.\hspace{-0.4mm}.\hspace{-0.4mm}.}**@{},
 <0mm,0mm>*{};<6mm,3.6mm>*{.\hspace{-0.1mm}.\hspace{-0.1mm}.}**@{},
<10mm,1.5mm>*{};<8mm,7mm>*{}**@{-},
  <10mm,1.5mm>*{};<9mm,7mm>*{}**@{-},
  <10mm,1.5mm>*{};<10.5mm,6mm>*{}**@{-},
  <10mm,1.5mm>*{};<12mm,7mm>*{}**@{-},
 <10mm,1.5mm>*{};<10.5mm,6.6mm>*{.\hspace{-0.4mm}.\hspace{-0.4mm}.}**@{},
 <-10mm,-1.5mm>*{};<-12mm,-6mm>*{}**@{-},
 <-7mm,-1.5mm>*{};<-8mm,-6mm>*{}**@{-},
 <-4mm,-1.5mm>*{};<-4.5mm,-6mm>*{}**@{-},
 <0mm,0mm>*{};<0mm,-4.6mm>*{.\hspace{0.1mm}.\hspace{0.1mm}.}**@{},
<10mm,-1.5mm>*{};<12mm,-6mm>*{}**@{-},
 <7mm,-1.5mm>*{};<8mm,-6mm>*{}**@{-},
  <4mm,-1.5mm>*{};<4.5mm,-6mm>*{}**@{-},
  %%%%
<0mm,0mm>*{};<-9.5mm,8.2mm>*{^{I_{ 1}}}**@{},
<0mm,0mm>*{};<-3mm,8.2mm>*{^{I_{ i}}}**@{},
<0mm,0mm>*{};<2mm,8.2mm>*{^{I_{ i+1}}}**@{},
<0mm,0mm>*{};<10mm,8.2mm>*{^{I_{ k}}}**@{},
 %%%
<0mm,0mm>*{};<-12mm,-7.4mm>*{_1}**@{},
<0mm,0mm>*{};<-8mm,-7.4mm>*{_2}**@{},
<0mm,0mm>*{};<-4mm,-7.4mm>*{_3}**@{},
<0mm,0mm>*{};<6mm,-7.4mm>*{\ldots}**@{},
%<0mm,0mm>*{};<8.5mm,-7.4mm>*{_{n-1}}**@{},
<0mm,0mm>*{};<12.5mm,-7.4mm>*{_n}**@{},
\endxy
$\ ,
which are exactly the same as in \S 2.7 except that now the vertex gets a numerical label $a\in \N^*$;
\item the differential $\delta$ is given on generators by
\Beqr
\delta\left(
\xy
%\begin{xy}
 <0mm,0mm>*{\mbox{$\xy *=<20mm,3mm>\txt{\em a}*\frm{-}\endxy$}};<0mm,0mm>*{}**@{},
  <-10mm,1.5mm>*{};<-12mm,7mm>*{}**@{-},
  <-10mm,1.5mm>*{};<-11mm,7mm>*{}**@{-},
  <-10mm,1.5mm>*{};<-9.5mm,6mm>*{}**@{-},
  <-10mm,1.5mm>*{};<-8mm,7mm>*{}**@{-},
 <-10mm,1.5mm>*{};<-9.5mm,6.6mm>*{.\hspace{-0.4mm}.\hspace{-0.4mm}.}**@{},
 <0mm,0mm>*{};<-6.4mm,3.6mm>*{.\hspace{-0.1mm}.\hspace{-0.1mm}.}**@{},
  <-3mm,1.5mm>*{};<-5mm,7mm>*{}**@{-},
  <-3mm,1.5mm>*{};<-4mm,7mm>*{}**@{-},
  <-3mm,1.5mm>*{};<-2.5mm,6mm>*{}**@{-},
  <-3mm,1.5mm>*{};<-1mm,7mm>*{}**@{-},
 <-3mm,1.5mm>*{};<-2.5mm,6.6mm>*{.\hspace{-0.4mm}.\hspace{-0.4mm}.}**@{},
  <2mm,1.5mm>*{};<0mm,7mm>*{}**@{-},
  <2mm,1.5mm>*{};<1mm,7mm>*{}**@{-},
  <2mm,1.5mm>*{};<2.5mm,6mm>*{}**@{-},
  <2mm,1.5mm>*{};<4mm,7mm>*{}**@{-},
 <2mm,1.5mm>*{};<2.5mm,6.6mm>*{.\hspace{-0.4mm}.\hspace{-0.4mm}.}**@{},
 <0mm,0mm>*{};<6mm,3.6mm>*{.\hspace{-0.1mm}.\hspace{-0.1mm}.}**@{},
<10mm,1.5mm>*{};<8mm,7mm>*{}**@{-},
  <10mm,1.5mm>*{};<9mm,7mm>*{}**@{-},
  <10mm,1.5mm>*{};<10.5mm,6mm>*{}**@{-},
  <10mm,1.5mm>*{};<12mm,7mm>*{}**@{-},
 <10mm,1.5mm>*{};<10.5mm,6.6mm>*{.\hspace{-0.4mm}.\hspace{-0.4mm}.}**@{},
 <-10mm,-1.5mm>*{};<-12mm,-6mm>*{}**@{-},
 <-7mm,-1.5mm>*{};<-8mm,-6mm>*{}**@{-},
 <-4mm,-1.5mm>*{};<-4.5mm,-6mm>*{}**@{-},
 <0mm,0mm>*{};<0mm,-4.6mm>*{.\hspace{0.1mm}.\hspace{0.1mm}.}**@{},
<10mm,-1.5mm>*{};<12mm,-6mm>*{}**@{-},
 <7mm,-1.5mm>*{};<8mm,-6mm>*{}**@{-},
  <4mm,-1.5mm>*{};<4.5mm,-6mm>*{}**@{-},
  %%%%
<0mm,0mm>*{};<-9.5mm,8.2mm>*{^{I_{ 1}}}**@{},
<0mm,0mm>*{};<-3mm,8.2mm>*{^{I_{ i}}}**@{},
<0mm,0mm>*{};<2mm,8.2mm>*{^{I_{ i+1}}}**@{},
<0mm,0mm>*{};<10mm,8.2mm>*{^{I_{ k}}}**@{},
 %%%
<0mm,0mm>*{};<-12mm,-7.4mm>*{_1}**@{},
<0mm,0mm>*{};<-8mm,-7.4mm>*{_2}**@{},
<0mm,0mm>*{};<-4mm,-7.4mm>*{_3}**@{},
<0mm,0mm>*{};<6mm,-7.4mm>*{\ldots}**@{},
%<0mm,0mm>*{};<8.5mm,-7.4mm>*{_{n-1}}**@{},
<0mm,0mm>*{};<12.5mm,-7.4mm>*{_n}**@{},
\endxy
\right) &=& \sum_{i=1}^k(-1)^{i+1}
%%%%%%%%%%%%%%%%%%%%%%%%%%%%%%%%%%%%%%%%%%%%%%%%%%%%
\xy
%\begin{xy}
 <0mm,0mm>*{\mbox{$\xy *=<20mm,3mm>\txt{\em a}*\frm{-}\endxy$}};<0mm,0mm>*{}**@{},
  <-10mm,1.5mm>*{};<-12mm,7mm>*{}**@{-},
  <-10mm,1.5mm>*{};<-11mm,7mm>*{}**@{-},
  <-10mm,1.5mm>*{};<-9.5mm,6mm>*{}**@{-},
  <-10mm,1.5mm>*{};<-8mm,7mm>*{}**@{-},
 <-10mm,1.5mm>*{};<-9.5mm,6.6mm>*{.\hspace{-0.4mm}.\hspace{-0.4mm}.}**@{},
 <0mm,0mm>*{};<-5.5mm,3.6mm>*{.\hspace{-0.1mm}.\hspace{-0.1mm}.}**@{},
%
  %<0mm,1.5mm>*{};<-5mm,7mm>*{}**@{-},
  <0mm,1.5mm>*{};<-4mm,7mm>*{}**@{-},
  <0mm,1.5mm>*{};<-2.0mm,6mm>*{}**@{-},
  <0mm,1.5mm>*{};<-1mm,7mm>*{}**@{-},
 <0mm,1.5mm>*{};<-2.3mm,6.6mm>*{.\hspace{-0.4mm}.\hspace{-0.4mm}.}**@{},
%
  %<0mm,1.5mm>*{};<0mm,7mm>*{}**@{-},
  <0mm,1.5mm>*{};<1mm,7mm>*{}**@{-},
  <0mm,1.5mm>*{};<2.0mm,6mm>*{}**@{-},
  <0mm,1.5mm>*{};<4mm,7mm>*{}**@{-},
 <0mm,1.5mm>*{};<2.3mm,6.6mm>*{.\hspace{-0.4mm}.\hspace{-0.4mm}.}**@{},
 <0mm,0mm>*{};<6mm,3.6mm>*{.\hspace{-0.1mm}.\hspace{-0.1mm}.}**@{},
<10mm,1.5mm>*{};<8mm,7mm>*{}**@{-},
  <10mm,1.5mm>*{};<9mm,7mm>*{}**@{-},
  <10mm,1.5mm>*{};<10.5mm,6mm>*{}**@{-},
  <10mm,1.5mm>*{};<12mm,7mm>*{}**@{-},
 <10mm,1.5mm>*{};<10.5mm,6.6mm>*{.\hspace{-0.4mm}.\hspace{-0.4mm}.}**@{},
 <-10mm,-1.5mm>*{};<-12mm,-6mm>*{}**@{-},
 <-7mm,-1.5mm>*{};<-8mm,-6mm>*{}**@{-},
 <-4mm,-1.5mm>*{};<-4.5mm,-6mm>*{}**@{-},
 <0mm,0mm>*{};<0mm,-4.6mm>*{.\hspace{0.1mm}.\hspace{0.1mm}.}**@{},
<10mm,-1.5mm>*{};<12mm,-6mm>*{}**@{-},
 <7mm,-1.5mm>*{};<8mm,-6mm>*{}**@{-},
  <4mm,-1.5mm>*{};<4.5mm,-6mm>*{}**@{-},
  %%%%
<0mm,0mm>*{};<-9.5mm,8.2mm>*{^{I_{ 1}}}**@{},
<0mm,0mm>*{};<-2.5mm,8.2mm>*{{^{I_{ i}\sqcup}}}**@{},
<0mm,0mm>*{};<2.7mm,8.2mm>*{^{I_{ i+1}}}**@{},
<0mm,0mm>*{};<10mm,8.2mm>*{^{I_{ k}}}**@{},
 %%%
<0mm,0mm>*{};<-12mm,-7.4mm>*{_1}**@{},
<0mm,0mm>*{};<-8mm,-7.4mm>*{_2}**@{},
<0mm,0mm>*{};<-4mm,-7.4mm>*{_3}**@{},
<0mm,0mm>*{};<6mm,-7.4mm>*{\ldots}**@{},
%<0mm,0mm>*{};<8.5mm,-7.4mm>*{_{n-1}}**@{},
<0mm,0mm>*{};<12.5mm,-7.4mm>*{_n}**@{},
\endxy
\nonumber
\\
&& +\ \sum_{b+c=a\atop b,c\geq 1}\sum_{p+q=k+1\atop p\geq 1,q\geq
0}\sum_{i=0}^{p-1} \sum_{  {I_{i+1}=I_{i+1}'\sqcup I''_{i+1}\atop
.......................} \atop I_{i+q}=I_{i+q}'\sqcup I''_{i+q}}
\sum_{[n]=J_1\sqcup J_2}\sum_{s\geq 0}(-1)^{(p+1)q + i(q-1)}\nonumber
\\
&&
%%%%%%%%%%%%%%%%%%%%%%%%%%%%%%%%%%%%%%%%%%%%%%%%%%%%%%%%%%%
\frac{1}{s!}\ \ \xy
%\begin{xy}
 <19mm,0mm>*{\mbox{$\xy *=<58mm,3mm>\txt{\em b}*\frm{-}\endxy$}};
 <0mm,0mm>*{}**@{},
  <-10mm,1.5mm>*{};<-12mm,7mm>*{}**@{-},
  <-10mm,1.5mm>*{};<-11mm,7mm>*{}**@{-},
  <-10mm,1.5mm>*{};<-9.5mm,6mm>*{}**@{-},
  <-10mm,1.5mm>*{};<-8mm,7mm>*{}**@{-},
 <-10mm,1.5mm>*{};<-9.5mm,6.6mm>*{.\hspace{-0.4mm}.\hspace{-0.4mm}.}**@{},
 <0mm,0mm>*{};<-6.5mm,3.6mm>*{.\hspace{-0.1mm}.\hspace{-0.1mm}.}**@{},
  <-3mm,1.5mm>*{};<-5mm,7mm>*{}**@{-},
  <-3mm,1.5mm>*{};<-4mm,7mm>*{}**@{-},
  <-3mm,1.5mm>*{};<-2.5mm,6mm>*{}**@{-},
  <-3mm,1.5mm>*{};<-1mm,7mm>*{}**@{-},
 <-3mm,1.5mm>*{};<-2.5mm,6.6mm>*{.\hspace{-0.4mm}.\hspace{-0.4mm}.}**@{},
%%%%%
  <10mm,1.5mm>*{};<0mm,7mm>*{}**@{-},
  <10mm,1.5mm>*{};<4mm,7mm>*{}**@{-},
<10mm,1.5mm>*{};<7.3mm,5.9mm>*{.\hspace{-0.0mm}.\hspace{-0.0mm}.}**@{},
   <10mm,1.5mm>*{};<3.8mm,6.0mm>*{}**@{-},
 <10mm,1.5mm>*{};<2.5mm,6.6mm>*{.\hspace{-0.4mm}.\hspace{-0.4mm}.}**@{},
%
% <0mm,0mm>*{};<6mm,3.6mm>*{.\hspace{-0.1mm}.\hspace{-0.1mm}.}**@{},
%
<10mm,1.5mm>*{};<9mm,7mm>*{}**@{-},
  <10mm,1.5mm>*{};<10.5mm,6mm>*{}**@{-},
  <10mm,1.5mm>*{};<12mm,7mm>*{}**@{-},
 <10mm,1.5mm>*{};<10.5mm,6.6mm>*{.\hspace{-0.4mm}.\hspace{-0.4mm}.}**@{},
 %%%%%%%%%%%%%%%%%%%%%%%%%%%%%%%%%%%%%%%
%%%% lower legs of 1st vertex %%%%%%%
 <-10mm,-1.5mm>*{};<-12mm,-6mm>*{}**@{-},
 <-7mm,-1.5mm>*{};<-8mm,-6mm>*{}**@{-},
 <-4mm,-1.5mm>*{};<-4.5mm,-6mm>*{}**@{-},
 <-1mm,-1.5mm>*{};<-1.1mm,-6mm>*{}**@{-},
 <2mm,-1.5mm>*{};<2.0mm,-6mm>*{}**@{-},
 <0mm,0mm>*{};<20mm,-4.6mm>*{.\hspace{2mm}.\hspace{2mm}.\hspace{2mm}
 .\hspace{2mm}.\hspace{2mm}.\hspace{2mm}.\hspace{2mm}
 .\hspace{2mm}.\hspace{2mm}}**@{},
<48mm,-1.5mm>*{};<50mm,-6mm>*{}**@{-},
 <45mm,-1.5mm>*{};<46mm,-6mm>*{}**@{-},
  <42mm,-1.5mm>*{};<42.5mm,-6mm>*{}**@{-},
  <39mm,-1.5mm>*{};<39.2mm,-6mm>*{}**@{-},
   <36mm,-1.5mm>*{};<36mm,-6mm>*{}**@{-},
 <20mm,-1.5mm>*{};<20.0mm,-8mm>*{\underbrace{\hspace{66mm}}}**@{},
 <20mm,-1.5mm>*{};<20.0mm,-11mm>*{_{J_1}}**@{},
  %%%%
<0mm,0mm>*{};<-9.5mm,8.4mm>*{^{I_{ 1}}}**@{},
<0mm,0mm>*{};<-3mm,8.4mm>*{^{I_{ i}}}**@{},
<0mm,0mm>*{};<2mm,8.6mm>*{^{I_{ i+1}'}}**@{},
<0mm,0mm>*{};<10.5mm,8.6mm>*{^{I_{ i+q}'}}**@{},
%%%% edges connecting two vertices %%%%%%%%%%
<10mm,1.5mm>*{};<18mm,12mm>*{}**@{-},
<10mm,1.5mm>*{};<20.0mm,12mm>*{}**@{-},
<10mm,1.5mm>*{};<25mm,12mm>*{}**@{-},
<10mm,1.5mm>*{};<18.7mm,8.6mm>*{.\hspace{-0.4mm}.\hspace{-0.4mm}.}**@{},
<10mm,1.5mm>*{};<19.6mm,9.0mm>*{^s}**@{},
%%%%%%%%%%%%%%%%%%%%%%%%%%%%%%%%%%%%%%%%%%%%%%
%%%% second vertex %%%%%%%%%%%%%%%%%%%%%%%%%%
%%%%%%%%%%%%%%%%%%%%%%%%%%%%%%%%%%%%%%%%%%%%%
<25mm,13.75mm>*{\mbox{$\xy *=<14mm,3mm>\txt{\em c}*\frm{-}\endxy$}};
<0mm,0mm>*{}**@{},
 <18mm,15mm>*{};<16mm,20.5mm>*{}**@{-},
 <18mm,15mm>*{};<17mm,20.5mm>*{}**@{-},
 <18mm,15mm>*{};<18.5mm,19.6mm>*{}**@{-},
 <18mm,15mm>*{};<20mm,20.5mm>*{}**@{-},
 <18mm,15mm>*{};<18.6mm,20.3mm>*{.\hspace{-0.4mm}.\hspace{-0.4mm}.}**@{},
<22mm,15mm>*{};<25.5mm,17.7mm>*{\cdots}**@{},
 <32mm,15.2mm>*{};<30mm,20.5mm>*{}**@{-},
 <32mm,15.2mm>*{};<31mm,20.5mm>*{}**@{-},
 <32mm,15.2mm>*{};<32.5mm,19.6mm>*{}**@{-},
 <32mm,15mm>*{};<34mm,20.5mm>*{}**@{-},
 <32mm,15mm>*{};<32.3mm,20.3mm>*{.\hspace{-0.4mm}.\hspace{-0.4mm}.}**@{},
%
  %%%%
<0mm,0mm>*{};<18mm,22.6mm>*{^{I_{ i+1}''}}**@{},
<0mm,0mm>*{};<32.5mm,22.6mm>*{^{I_{ i+q}''}}**@{},
  %%%%
%%%% lower legs of 2nd vertex %%%%%%
 <26mm,12mm>*{};<25mm,9mm>*{}**@{-},
 <27mm,12mm>*{};<26.8mm,9mm>*{}**@{-},
 <29mm,12mm>*{};<29.2mm,10mm>*{.\hspace{-0.1mm}.\hspace{-0.1mm}.}**@{},
 <32mm,12mm>*{};<33mm,9mm>*{}**@{-},
  <31mm,12mm>*{};<31.5mm,9mm>*{}**@{-},
 <29mm,12mm>*{};<29mm,8mm>*{\underbrace{\ \ \ \ \ \ \ \  }}**@{},
 <29mm,12mm>*{};<29mm,5.3mm>*{_{J_2}}**@{},
%%%%%%%%%%%%%%%%%% 1st vertex contd %%%%%%%%%%%
<38mm,1.5mm>*{};<36mm,7mm>*{}**@{-},
<38mm,1.5mm>*{};<37mm,7mm>*{}**@{-},
<38mm,1.5mm>*{};<38.5mm,6mm>*{}**@{-},
<38mm,1.5mm>*{};<40mm,7mm>*{}**@{-},
<38mm,1.5mm>*{};<38.5mm,6.6mm>*{.\hspace{-0.4mm}.\hspace{-0.4mm}.}**@{},
<38mm,1.5mm>*{};<43mm,4mm>*{.\hspace{-0.0mm}.\hspace{-0.0mm}.}**@{},
<48mm,1.5mm>*{};<46mm,7mm>*{}**@{-},
<48mm,1.5mm>*{};<47mm,7mm>*{}**@{-},
<48mm,1.5mm>*{};<48.5mm,6mm>*{}**@{-},
<48mm,1.5mm>*{};<50mm,7mm>*{}**@{-},
<48mm,1.5mm>*{};<48.5mm,6.6mm>*{.\hspace{-0.4mm}.\hspace{-0.4mm}.}**@{},
<0mm,0mm>*{};<40.3mm,8.6mm>*{^{I_{ i+q+1}}}**@{},
<0mm,0mm>*{};<48.5mm,8.6mm>*{^{I_{ k}}}**@{},
  %%%%
\endxy
\label{deltah}
\Eeqr
\Ei

\no
Note that $\DefQh$ is spanned by graphs
 over $\K$, {\em not}\, over $\K[[\hbar]]$. At the prop level the only
remnant of the presence of the Plank constant $\hbar$ in the input
geometry
 is in the decoration of vertices
by a natural number $a\in \N^*$. In this respect the notation
$\DefQh$ could be misleading.

\bip

%%%%%%%%%%%%%%%%%%%%%%%%%%%%%%%%%%%%%%%%%%%%%%%
\no{\bf 2.9.1. Proposition.} {\em There is a one-to-one correspondence
between representations,
$$
\phi: \left(\DefQh,\delta\right) \lon \left({\sf
End}_{\K[[\hbar]]}\langle V[[\hbar]]\rangle,d\right),
$$
of $\DefQh$
 in an $\K[[\hbar]]$-extension of a  dg vector space
$(V,d)$, and star products on $\f_V[[\hbar]]$, i.e.\ Maurer-Cartan elements,
$\Ga$, in $\sD_V[[\hbar]]$
satisfying the equation, $d_{\rm H}\Ga+ \frac{1}{2}[\Ga,\Ga]_{\rm
H}=0$ and the condition $\Ga|_{\hbar=0}=d$.}

\bip

\no{Proof} is similar to the proof of Proposition 2.6.1.

\bip

%%%%%%%%%%%%%%%%%%%%%%%%%%%%%%%%%%%%%%%%%%%%%%%
\no{\bf 2.9.2. Remark.} By constructions of
 $(\DefQ, \delta)$ and  $(\DefQh, \delta)$,
 there is a canonical morphism of dg
props,
\[
\Ba{rccc} \chi_\hbar: & \DefQ & \lon
&
\DefQh[[\hbar]] \vspace{2mm}\\
& \xy
%\begin{xy}
 <0mm,0mm>*{\mbox{$\xy *=<20mm,3mm>\txt{}*\frm{-}\endxy$}};<0mm,0mm>*{}**@{},
  <-10mm,1.5mm>*{};<-12mm,7mm>*{}**@{-},
  <-10mm,1.5mm>*{};<-11mm,7mm>*{}**@{-},
  <-10mm,1.5mm>*{};<-9.5mm,6mm>*{}**@{-},
  <-10mm,1.5mm>*{};<-8mm,7mm>*{}**@{-},
 <-10mm,1.5mm>*{};<-9.5mm,6.6mm>*{.\hspace{-0.4mm}.\hspace{-0.4mm}.}**@{},
 <0mm,0mm>*{};<-6.5mm,3.6mm>*{.\hspace{-0.1mm}.\hspace{-0.1mm}.}**@{},
  <-3mm,1.5mm>*{};<-5mm,7mm>*{}**@{-},
  <-3mm,1.5mm>*{};<-4mm,7mm>*{}**@{-},
  <-3mm,1.5mm>*{};<-2.5mm,6mm>*{}**@{-},
  <-3mm,1.5mm>*{};<-1mm,7mm>*{}**@{-},
 <-3mm,1.5mm>*{};<-2.5mm,6.6mm>*{.\hspace{-0.4mm}.\hspace{-0.4mm}.}**@{},
  <2mm,1.5mm>*{};<0mm,7mm>*{}**@{-},
  <2mm,1.5mm>*{};<1mm,7mm>*{}**@{-},
  <2mm,1.5mm>*{};<2.5mm,6mm>*{}**@{-},
  <2mm,1.5mm>*{};<4mm,7mm>*{}**@{-},
 <2mm,1.5mm>*{};<2.5mm,6.6mm>*{.\hspace{-0.4mm}.\hspace{-0.4mm}.}**@{},
 <0mm,0mm>*{};<6mm,3.6mm>*{.\hspace{-0.1mm}.\hspace{-0.1mm}.}**@{},
<10mm,1.5mm>*{};<8mm,7mm>*{}**@{-},
  <10mm,1.5mm>*{};<9mm,7mm>*{}**@{-},
  <10mm,1.5mm>*{};<10.5mm,6mm>*{}**@{-},
  <10mm,1.5mm>*{};<12mm,7mm>*{}**@{-},
 <10mm,1.5mm>*{};<10.5mm,6.6mm>*{.\hspace{-0.4mm}.\hspace{-0.4mm}.}**@{},
 <-10mm,-1.5mm>*{};<-12mm,-6mm>*{}**@{-},
 <-7mm,-1.5mm>*{};<-8mm,-6mm>*{}**@{-},
 <-4mm,-1.5mm>*{};<-4.5mm,-6mm>*{}**@{-},
 <0mm,0mm>*{};<0mm,-4.6mm>*{.\hspace{0.1mm}.\hspace{0.1mm}.}**@{},
<10mm,-1.5mm>*{};<12mm,-6mm>*{}**@{-},
 <7mm,-1.5mm>*{};<8mm,-6mm>*{}**@{-},
  <4mm,-1.5mm>*{};<4.5mm,-6mm>*{}**@{-},
  %%%%
<0mm,0mm>*{};<-9.5mm,8.2mm>*{^{I_{ 1}}}**@{},
<0mm,0mm>*{};<-3mm,8.2mm>*{^{I_{ i}}}**@{},
<0mm,0mm>*{};<2mm,8.2mm>*{^{I_{ i+1}}}**@{},
<0mm,0mm>*{};<10mm,8.2mm>*{^{I_{ k}}}**@{},
 %%%
<0mm,0mm>*{};<-12mm,-7.4mm>*{_1}**@{},
<0mm,0mm>*{};<-8mm,-7.4mm>*{_2}**@{},
<0mm,0mm>*{};<-4mm,-7.4mm>*{_3}**@{},
<0mm,0mm>*{};<6mm,-7.4mm>*{\ldots}**@{},
%<0mm,0mm>*{};<8.5mm,-7.4mm>*{_{n-1}}**@{},
<0mm,0mm>*{};<12.5mm,-7.4mm>*{_n}**@{},
\endxy
& \lon & \underset{a\geq 1}{\bigoplus} \ \ \hbar^a
\xy
%\begin{xy}
 <0mm,0mm>*{\mbox{$\xy *=<20mm,3mm>\txt{\em a}*\frm{-}\endxy$}};<0mm,0mm>*{}**@{},
  <-10mm,1.5mm>*{};<-12mm,7mm>*{}**@{-},
  <-10mm,1.5mm>*{};<-11mm,7mm>*{}**@{-},
  <-10mm,1.5mm>*{};<-9.5mm,6mm>*{}**@{-},
  <-10mm,1.5mm>*{};<-8mm,7mm>*{}**@{-},
 <-10mm,1.5mm>*{};<-9.5mm,6.6mm>*{.\hspace{-0.4mm}.\hspace{-0.4mm}.}**@{},
 <0mm,0mm>*{};<-6.5mm,3.6mm>*{.\hspace{-0.1mm}.\hspace{-0.1mm}.}**@{},
  <-3mm,1.5mm>*{};<-5mm,7mm>*{}**@{-},
  <-3mm,1.5mm>*{};<-4mm,7mm>*{}**@{-},
  <-3mm,1.5mm>*{};<-2.5mm,6mm>*{}**@{-},
  <-3mm,1.5mm>*{};<-1mm,7mm>*{}**@{-},
 <-3mm,1.5mm>*{};<-2.5mm,6.6mm>*{.\hspace{-0.4mm}.\hspace{-0.4mm}.}**@{},
  <2mm,1.5mm>*{};<0mm,7mm>*{}**@{-},
  <2mm,1.5mm>*{};<1mm,7mm>*{}**@{-},
  <2mm,1.5mm>*{};<2.5mm,6mm>*{}**@{-},
  <2mm,1.5mm>*{};<4mm,7mm>*{}**@{-},
 <2mm,1.5mm>*{};<2.5mm,6.6mm>*{.\hspace{-0.4mm}.\hspace{-0.4mm}.}**@{},
 <0mm,0mm>*{};<6mm,3.6mm>*{.\hspace{-0.1mm}.\hspace{-0.1mm}.}**@{},
<10mm,1.5mm>*{};<8mm,7mm>*{}**@{-},
  <10mm,1.5mm>*{};<9mm,7mm>*{}**@{-},
  <10mm,1.5mm>*{};<10.5mm,6mm>*{}**@{-},
  <10mm,1.5mm>*{};<12mm,7mm>*{}**@{-},
 <10mm,1.5mm>*{};<10.5mm,6.6mm>*{.\hspace{-0.4mm}.\hspace{-0.4mm}.}**@{},
 <-10mm,-1.5mm>*{};<-12mm,-6mm>*{}**@{-},
 <-7mm,-1.5mm>*{};<-8mm,-6mm>*{}**@{-},
 <-4mm,-1.5mm>*{};<-4.5mm,-6mm>*{}**@{-},
 <0mm,0mm>*{};<0mm,-4.6mm>*{.\hspace{0.1mm}.\hspace{0.1mm}.}**@{},
<10mm,-1.5mm>*{};<12mm,-6mm>*{}**@{-},
 <7mm,-1.5mm>*{};<8mm,-6mm>*{}**@{-},
  <4mm,-1.5mm>*{};<4.5mm,-6mm>*{}**@{-},
  %%%%
<0mm,0mm>*{};<-9.5mm,8.2mm>*{^{I_{ 1}}}**@{},
<0mm,0mm>*{};<-3mm,8.2mm>*{^{I_{ i}}}**@{},
<0mm,0mm>*{};<2mm,8.2mm>*{^{I_{ i+1}}}**@{},
<0mm,0mm>*{};<10mm,8.2mm>*{^{I_{ k}}}**@{},
 %%%
<0mm,0mm>*{};<-12mm,-7.4mm>*{_1}**@{},
<0mm,0mm>*{};<-8mm,-7.4mm>*{_2}**@{},
<0mm,0mm>*{};<-4mm,-7.4mm>*{_3}**@{},
<0mm,0mm>*{};<6mm,-7.4mm>*{\ldots}**@{},
%<0mm,0mm>*{};<8.5mm,-7.4mm>*{_{n-1}}**@{},
<0mm,0mm>*{};<12.5mm,-7.4mm>*{_n}**@{},
\endxy
\Ea
\]

\sip

\no
As  we assume that the prop
$\DefQ^\hbar$ is completed with respect to the
number of vertices, we can set $\hbar=1$ above and get a
well-defined morphism of dg props, $\chi_{\hbar=1}:
\DefQ \rar \DefQ^\hbar$.
However, in general a morphism of dg props $\phi:
\DefQ^\hbar\rar \sP$ can not be composed with
$\chi_{\hbar=1}$ while $\phi\circ \chi_\hbar:
\DefQ  \rar \sP[[\hbar]]$ is always well defined.

\bip

We shall next investigate how wheeled completion of  directed graph
complexes affects their cohomology.

\bip

\bip

%%%%%%%%%%%%%%%%%%%%%%%%%% 3 %%%%%%%%%%%%%%%%%%%%%%%%%%%%%%%
%%%%%%%%%%%%%%%%%%%%%%%%%%%%%%%%%%%%%%%%%%%%%%%%%%%%%%%%%
\begin{center}
\bf \S 3.  Directed graph complexes with loops and wheels
\end{center}

\bip

%%%%%%%%%%%%%%%%%%%%%%%%%%%%%%%%%%%%%%%%%%%%%%%%%%%%%%%%%
\no{\bf 3.1. $\fG^\circlearrowright$ versus $\fG^\uparrow$.}
One of the most effective methods for computing cohomology of dg free props
(that is, decorated $\fG^\uparrow$-graph complexes) is based on the idea
of interpreting the differential as a perturbation of its $\frac{1}{2}$propic
part which,  in this $\fG^\uparrow$-case, can often be singled out by the path
filtration \cite{Ko2,MV}. However, one can {\em not}\, apply this idea directly
to graphs with wheels --- it is shown below that
 a filtration which singles out the $\frac{1}{2}$propic part of the
 differential does {\em not}\, exist in general even for dioperadic
 differentials. Put another way, if one takes  a $\fG^\uparrow$-graph complex,
 $(\PROP\langle E \rangle, \delta)$, enlarges it by adding decorated graphs
 with wheels  while keeping the original differential $\delta$ unchanged, then
 one ends up in a very different situation in which the idea of
 $\frac{1}{2}$props is no longer directly applicable.

\bip

\no{\bf 3.2. Graphs with back-in-time edges.}
Here we suggest the following trick: we further enlarge
our set of graphs with wheels, $\fGG\rightsquigarrow \fGGG$, by putting a mark
on one (and only one) of the edges in each wheel, and then study the natural ``forgetful"
surjection, $\fGGG \rar \fGG$. The point is that $\fGGG$-graph complexes  again admit a filtration
which singles out the  $\frac{1}{2}$prop part of the differential and
hence their cohomology are  often easily computable.

\bip

\no
More precisely,
let $\fGGG(m,n)$
be the set of all directed $(m,n)$-graphs $G$ which satisfy
 conditions 2.1.3(i)-(iv),  and the following one:
\Bi
\item[(v)] every oriented wheel in $G$ (if any) has one and only one of its
internal edges marked (say, dashed) and called {\em back-in-time}\,
edge.
\Ei
For example,
$$
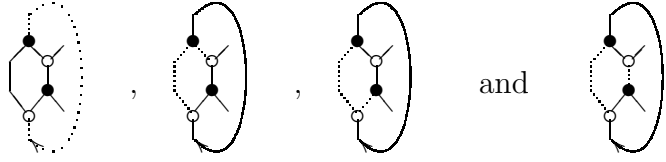
\begin{xy}
 <0mm,2.47mm>*{};<0mm,-0.5mm>*{}**@{-},
 <0.5mm,3.5mm>*{};<2.2mm,5.2mm>*{}**@{-},
 <-0.48mm,3.48mm>*{};<-2.2mm,5.2mm>*{}**@{-},
 <0mm,3mm>*{\circ};<0mm,3mm>*{}**@{},
  <0mm,-0.8mm>*{\bullet};<0mm,-0.8mm>*{}**@{},
<0mm,-0.8mm>*{};<-2.2mm,-3.5mm>*{}**@{-},
<0mm,-0.8mm>*{};<2.2mm,-3.5mm>*{}**@{-},
 <-2.5mm,5.7mm>*{\bullet};<0mm,0mm>*{}**@{},
<-2.5mm,5.7mm>*{};<-2.5mm,9.4mm>*{}**@{.},                  %1
<-2.5mm,5.7mm>*{};<-5mm,3mm>*{}**@{-},
<-5mm,3mm>*{};<-5mm,-0.8mm>*{}**@{-},
 <-2.5mm,-4.2mm>*{\circ};<0mm,3mm>*{}**@{},
 <-2.8mm,-3.6mm>*{};<-5mm,-0.8mm>*{}**@{-},
 <-2.5mm,-4.6mm>*{};<-2.5mm,-7.3mm>*{}**@{.},
%   %
 (-2.2,9.4)*{}
   \ar@{.>}@(ur,dr) (-2.2,-7.4)*{}
\end{xy}
\ \ \ , \ \ \
\begin{xy}
 <0mm,2.47mm>*{};<0mm,-0.5mm>*{}**@{-},
 <0.5mm,3.5mm>*{};<2.2mm,5.2mm>*{}**@{-},
 <-0.48mm,3.48mm>*{};<-2.2mm,5.2mm>*{}**@{.},
 <0mm,3mm>*{\circ};<0mm,3mm>*{}**@{},
  <0mm,-0.8mm>*{\bullet};<0mm,-0.8mm>*{}**@{},
<0mm,-0.8mm>*{};<-2.2mm,-3.5mm>*{}**@{-},
<0mm,-0.8mm>*{};<2.2mm,-3.5mm>*{}**@{-},                  %2
 <-2.5mm,5.7mm>*{\bullet};<0mm,0mm>*{}**@{},
<-2.5mm,5.7mm>*{};<-2.5mm,9.4mm>*{}**@{-},
<-2.5mm,5.7mm>*{};<-5mm,3mm>*{}**@{.},
<-5mm,3mm>*{};<-5mm,-0.8mm>*{}**@{.},
 <-2.5mm,-4.2mm>*{\circ};<0mm,3mm>*{}**@{},
 <-2.8mm,-3.6mm>*{};<-5mm,-0.8mm>*{}**@{.},
 <-2.5mm,-4.6mm>*{};<-2.5mm,-7.3mm>*{}**@{-},
%   %
 (-2.2,9.4)*{}
   \ar@{->}@(ur,dr) (-2.2,-7.4)*{}
\end{xy}
\ \ \ , \ \ \
\begin{xy}
 <0mm,2.47mm>*{};<0mm,-0.5mm>*{}**@{-},
 <0.5mm,3.5mm>*{};<2.2mm,5.2mm>*{}**@{-},
 <-0.48mm,3.48mm>*{};<-2.2mm,5.2mm>*{}**@{-},
 <0mm,3mm>*{\circ};<0mm,3mm>*{}**@{},
  <0mm,-0.8mm>*{\bullet};<0mm,-0.8mm>*{}**@{},
<0mm,-0.8mm>*{};<-2.2mm,-3.5mm>*{}**@{.},
<0mm,-0.8mm>*{};<2.2mm,-3.5mm>*{}**@{-},
 <-2.5mm,5.7mm>*{\bullet};<0mm,0mm>*{}**@{},
<-2.5mm,5.7mm>*{};<-2.5mm,9.4mm>*{}**@{-},
<-2.5mm,5.7mm>*{};<-5mm,3mm>*{}**@{.},
<-5mm,3mm>*{};<-5mm,-0.8mm>*{}**@{.},
 <-2.5mm,-4.2mm>*{\circ};<0mm,3mm>*{}**@{},
 <-2.8mm,-3.6mm>*{};<-5mm,-0.8mm>*{}**@{.},
 <-2.5mm,-4.6mm>*{};<-2.5mm,-7.3mm>*{}**@{-},
%   %
 (-2.2,9.4)*{}
   \ar@{->}@(ur,dr) (-2.2,-7.4)*{}
\end{xy}
\ \ \ \ \ \mbox{and}\ \ \ \ \ \
\begin{xy}
 <0mm,2.47mm>*{};<0mm,-0.5mm>*{}**@{.},
 <0.5mm,3.5mm>*{};<2.2mm,5.2mm>*{}**@{-},
 <-0.48mm,3.48mm>*{};<-2.2mm,5.2mm>*{}**@{-},
 <0mm,3mm>*{\circ};<0mm,3mm>*{}**@{},
  <0mm,-0.8mm>*{\bullet};<0mm,-0.8mm>*{}**@{},
<0mm,-0.8mm>*{};<-2.2mm,-3.5mm>*{}**@{-},
<0mm,-0.8mm>*{};<2.2mm,-3.5mm>*{}**@{-},
 <-2.5mm,5.7mm>*{\bullet};<0mm,0mm>*{}**@{},
<-2.5mm,5.7mm>*{};<-2.5mm,9.4mm>*{}**@{-},
<-2.5mm,5.7mm>*{};<-5mm,3mm>*{}**@{.},
<-5mm,3mm>*{};<-5mm,-0.8mm>*{}**@{.},
 <-2.5mm,-4.2mm>*{\circ};<0mm,3mm>*{}**@{},
 <-2.8mm,-3.6mm>*{};<-5mm,-0.8mm>*{}**@{.},
 <-2.5mm,-4.6mm>*{};<-2.5mm,-7.3mm>*{}**@{-},
%   %
 (-2.2,9.4)*{}
   \ar@{->}@(ur,dr) (-2.2,-7.4)*{}
\end{xy}
$$
are four different graphs in $\fG^+(1,1)$.

\bip

\no
Clearly, we have a natural surjection,
$$
u: \fG^+(m,n)\lon
\fG^\circlearrowright(m,n),
$$
which forgets the markings. For example, the four graphs above are mapped
under $u$ into the same graph,
$$
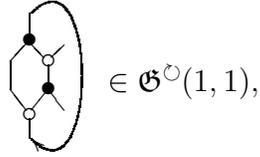
\begin{xy}
 <0mm,2.47mm>*{};<0mm,-0.5mm>*{}**@{-},
 <0.5mm,3.5mm>*{};<2.2mm,5.2mm>*{}**@{-},
 <-0.48mm,3.48mm>*{};<-2.2mm,5.2mm>*{}**@{-},
 <0mm,3mm>*{\circ};<0mm,3mm>*{}**@{},
  <0mm,-0.8mm>*{\bullet};<0mm,-0.8mm>*{}**@{},
<0mm,-0.8mm>*{};<-2.2mm,-3.5mm>*{}**@{-},
<0mm,-0.8mm>*{};<2.2mm,-3.5mm>*{}**@{-},
 <-2.5mm,5.7mm>*{\bullet};<0mm,0mm>*{}**@{},
<-2.5mm,5.7mm>*{};<-2.5mm,9.4mm>*{}**@{-},                  %1
<-2.5mm,5.7mm>*{};<-5mm,3mm>*{}**@{-},
<-5mm,3mm>*{};<-5mm,-0.8mm>*{}**@{-},
 <-2.5mm,-4.2mm>*{\circ};<0mm,3mm>*{}**@{},
 <-2.8mm,-3.6mm>*{};<-5mm,-0.8mm>*{}**@{-},
 <-2.5mm,-4.6mm>*{};<-2.5mm,-7.3mm>*{}**@{-},
%   %
 (-2.2,9.4)*{}
   \ar@{->}@(ur,dr) (-2.2,-7.4)*{}
\end{xy} \in \fG^\circlearrowright(1,1),
$$
and, in fact, span its pre-image under $u$.

\bip

%%%%%%%%%%%%%%%%%%%%%%%%%%%%%%%%%%%%%%%%%%%%%%%%%%%%%%%%%%
%%%%%%%%%%%%%%%%%%%%%%%%%%%%%%%%%%%%%%%%%%%%%%%%%%%%%%%%%
\no{\bf 3.3. Graph complexes}. Let $E=\{E(m,n)\}_{m,n\geq 1, m+n\geq 3}$ be an
$\bS$-bimodule and  let $(\PROP \langle E \rangle, \delta)$
be a dg  free prop on $E$ with the differential $\delta$ which preserves
 connectedness and genus,
that is, $\delta$ applied to any decorated $(m,n)$-corolla creates a
connected
$(m,n)$-tree. Such a differential can be called {\em dioperadic}, and from
now on we restrict ourselves to dioperadic differentials only.
This
restriction is not that dramatic: every dg free prop with a non-dioperadic but
connected\footnote
{A differential $\delta$ in $\PROP \langle E \rangle$ is called
{\em connected} if it preserves the filtration of $\PROP \langle E \rangle$
by the number of connected (in the topological sense) components.}
differential always admits a filtration which singles out its
dioperadic part \cite{MV}. Thus the technique we develop here in \S 3 can,
in principle, be applied to a wheeled extension of any dg free prop
with a connected differential.

\bip

\no
We enlarge the  $\fG^\uparrow$-graph complex
$(\PROP \langle E \rangle, \delta)$ in two ways,
$$
\PP\langle E\rangle(m,n):= \bigoplus_{G\in \fGG(m,n)} \left(
\bigotimes_{v\in v(G)} E(Out(v), In(v))\right)_{Aut G},
$$
$$
\PPP\langle E\rangle(m,n):= \bigoplus_{G\in \fGGG(m,n)} \left(
\bigotimes_{v\in v(G)} E(Out(v), In(v))\right)_{Aut G},
$$
and notice that both
$\PP\langle E\rangle:=\{\PP\langle E\rangle(m,n)\}$
and $\PPP\langle E\rangle:=\{\PP\langle E\rangle(m,n)\}$
 have a natural prop structure
with respect to disjoint unions and grafting of decorated graphs. The original differential $\delta$ in
$\PROP\langle E\rangle$ extends naturally to $\PP\langle E\rangle$
and $\PPP\langle E\rangle$ making the latter into {\em dg}\, props.
However that differential preserves, in general, neither the number of wheels in $\PP\langle E\rangle$
nor the number of marked edges in $\PPP\langle E\rangle$. Clearly, they both contain $(\PROP\langle
E\rangle, \delta)$ as a dg subprop.  There is a natural morphism of
dg props,
$$
u: (\PROP^{+} \langle E\rangle, \delta) \lon
(\PROP^\circlearrowright\langle E\rangle, \delta),
$$
which forgets the markings. Let $({\mathsf L}^+\langle E \rangle,
\delta):= \Ker u$ and denote the natural inclusion ${\mathsf
L}^+\langle E \rangle\subset \PROP^+\langle E\rangle$ by $i$.

\bip

%%%%%%%%%%%%%%%%%%%%%%%%%%%%%%%%%%%%%%%%%%%%%%%%%%%%%%%
\no{\bf 3.4. Fact.}  There is a short exact sequence of graph
complexes,
$$
0\lon (\LLL\langle E \rangle, \delta) \stackrel{i}{\lon}
 (\PPP\langle E \rangle,\delta)
\stackrel{u}{\lon}
(\PP\langle E \rangle,\delta) \lon 0,
$$
 Thus, if the natural
inclusion of complexes, $i: (\LLL\langle E \rangle, \delta) \rar
(\PPP\langle E \rangle,\delta)$, induces a monomorphism in
cohomology, $[i]: \sH(\LLL\langle E \rangle,\delta)\rar \sH(\PPP\langle
E \rangle,\delta)$, then
$$
 \sH(\PP\langle E \rangle,\delta) =\frac{ \sH(\PPP\langle E \rangle,\delta)}
 {\sH(\LLL\langle E \rangle,\delta)}.
$$
 Put another way, if $[i]$ is a monomorphism, then $\sH(\PP\langle E
\rangle,\delta)$ is obtained from $\sH(\PPP\langle E \rangle,\delta)$
simply by forgetting the markings.

\bip

%%%%%%%%%%%%%%%%%%%%%%%%%%%%%%%%%%%%%%%%
\no{\bf 3.5. Functors which adjoin wheels.}
We are interested in this paper
 in dioperads, $D$, which are either free, ${\mathsf D}\langle E \rangle$,
 on an $\bS$-bimodule $E=\{E(m,n)\}_{m,n\geq 1, m+n\geq 3}$, or are naturally
 represented as quotients of free dioperads,
 $$
D=\frac{{\mathsf D}\langle E \rangle}{<I>},
 $$
modulo the ideals  generated by some relations
 $I\subset {\mathsf D}\langle E \rangle $.
Then the free prop, $\Omega_{\mathsf D\rar P}\langle E \rangle$,
  generated by $D$ is simply the quotient of the free prop,
  $\PROP\langle E \rangle$,
$$
\Omega_{\sD\rar \sP}\langle D \rangle := \frac{\PROP\langle E
\rangle}{<I>},
$$
by the ideal generated by the same relations $I$. Now we define two
other props\footnote{The prop $\Omega_{\mathsf D\rar \PP}\langle D \rangle$  ia  a particular example of a {\em
wheeled prop}\, which will be discussed in detail elsewhere.},
$$
\Omega_{\mathsf D\rar \PP}\langle D \rangle :=
\frac{\PP\langle E \rangle}{<I>^\circlearrowright}, \ \ \ \ \ \ \
\Omega_{\mathsf D\rar \PPP}\langle D \rangle :=
\frac{\PPP\langle E \rangle}{<I>^+}, \ \ \
$$
where  $<\hspace{-0.5mm}I\hspace{-0.5mm}>^\circlearrowright$ (resp.,
$<\hspace{-0.5mm} I\hspace{-0.5mm} >^+$) is the subspace of those graphs
$G$ in $\PP\langle E \rangle$ (resp., in $\PPP\langle E \rangle$) which satisfy
the following condition: there exists a (possibly empty) set of cyclic
edges whose breaking up into two legs produces a graph lying in the ideal
$<\hspace{-0.5mm} I\hspace{-0.5mm}>$ which defines the prop
$\Omega_{\mathsf D\rar P}\langle D\rangle$.

\bip

\no
Analogously one defines functors $\Omega_{\frac{1}{2}\sP\rar \PP}$
and $\Omega_{\frac{1}{2}\sP\rar \PPP}$.

\bip

\no
From now on we abbreviate notations as follows,
$$
D^\uparrow:= \Omega_{\sD\rar \sP}\langle D \rangle, \ \ \ D^+:=
\Omega_{\sD\rar \PPP}\langle D \rangle, \ \ \ D^\circlearrowright:=
\Omega_{\sD\rar \PP}\langle D \rangle,
$$
for values of the above defined functors on dioperads, and,
respectively
$$
D_0^\uparrow:= \Omega_{\frac{1}{2}\sP\rar \sP}\langle D_0 \rangle, \ \ \ D_0^+:=
\Omega_{\frac{1}{2}\sP \rar \PP}\langle D_0\rangle,\ \ \  D_0^\circlearrowright:=
\Omega_{\frac{1}{2}\sP \rar \PP}\langle D_0 \rangle.
$$
for their values  on $\frac{1}{2}$props.

\bip

%%%%%%%%%%%%%%%%%%%%%%%%%%%%%%%%%%%%%%%%
\no{\bf 3.5.1. Facts}. {\bf (i)}
If $D$ is a dg dioperad, then both $D^\circlearrowright$ and
%$\Omega_{\mathsf D\rar P^+}\langle D\rangle$ a
% $\Omega_{\mathsf D\rar P^+}\langle D \rangle$
$D^+$ are naturally  dg props. {\bf (ii)} If $D$ is an operad, then
both $D^\circlearrowright$ and $D^+$
 may contain at most one wheel.

\bip
%%%%%%%%%%%%%%%%%%%%%%%%%
\no{\bf 3.5.2. Proposition.} {\em Any finite-dimensional representation of the dioperad $D$ lifts
to a representation of its wheeled prop extension\,
%$\Omega_{\mathsf D\rar P^+}\langle D\rangle$.
$D^\circlearrowright$}.

\bip

\Proof If $\phi:  D\rar \End\langle M\rangle$ is a representation, then
we first extend it to a representation, $\phi^\circlearrowright$, of
$\PP\langle E \rangle$ as in Sect.\ 2.3.2
and then notice that
$\phi^\circlearrowright(f)= 0$ for any $f\in <I>^\circlearrowright$.
\hfill $\Box$

\bip
%%%%%%%%%%%%%%%%%%%%%%%%%
\no{\bf 3.5.3. Definition.} Let $D$ be a Koszul dioperad with
$(D_\infty,\delta )\rar (D,0)$ being its minimal resolution. The
dioperad $D$ is called {\em stably Koszul}\, if the associated
morphism of the wheeled completions,
$$
(D_\infty^\circlearrowright, \delta) \lon (D^\circlearrowright,0)
$$
remains a quasi-isomorphism.

\bip

%%%%%%%%%%%%%%%%%%%%%%%%%
\no{\bf 3.5.4. Example.} The notion of stable Koszulness is
non-trivial. Just adding oriented wheels to a minimal resolution of
a Koszul operad while keeping the differential unchanged may
 alter the cohomology group of the resulting graph
complex as the following example shows.

\bip

\no {\bf Claim}. {\em The operad, $\sf Ass$, of associative
algebras is not stably Koszul}.

\sip

\Proof The operad $\sf Ass$ can be represented a quotient,
$$
 {\sf Ass} = \frac{ {\mathsf O\mathsf p\mathsf e\mathsf r}\langle E \rangle}{{\sf Ideal} <R>}\ ,
$$
of the free operad,  ${\mathsf O\mathsf p\mathsf e\mathsf r}\langle
E \rangle$, generated by the following $\bS$-module $E$,
$$
E(n):= \left\{ \Ba{ll} k[\bS_2]=\mbox{span}\ \left(
\begin{xy}
 <0mm,-0.55mm>*{};<0mm,-2.5mm>*{}**@{-},
 <0.5mm,0.5mm>*{};<2.2mm,2.2mm>*{}**@{-},
 <-0.48mm,0.48mm>*{};<-2.2mm,2.2mm>*{}**@{-},
 <0mm,0mm>*{\circ};<0mm,0mm>*{}**@{},
   %<0mm,-0.55mm>*{};<0mm,-3.8mm>*{_1}**@{},
   <0.5mm,0.5mm>*{};<2.7mm,2.8mm>*{^2}**@{},
   <-0.48mm,0.48mm>*{};<-2.7mm,2.8mm>*{^1}**@{},
 \end{xy}\ , \
 \begin{xy}
 <0mm,-0.55mm>*{};<0mm,-2.5mm>*{}**@{-},
 <0.5mm,0.5mm>*{};<2.2mm,2.2mm>*{}**@{-},
 <-0.48mm,0.48mm>*{};<-2.2mm,2.2mm>*{}**@{-},
 <0mm,0mm>*{\circ};<0mm,0mm>*{}**@{},
   %<0mm,-0.55mm>*{};<0mm,-3.8mm>*{_1}**@{},
   <0.5mm,0.5mm>*{};<2.7mm,2.8mm>*{^1}**@{},
   <-0.48mm,0.48mm>*{};<-2.7mm,2.8mm>*{^2}**@{},
 \end{xy}
\right) & \mbox{for}\ n=2 \\
0 & \mbox{otherwise}, \Ea \right.
$$
modulo the ideal  generated by the following relations,
%%%%%%%%%%%%%%%%%%%% coLie %%%%%%%%%%%%%%%%%%
$$
\begin{xy}
 <0mm,0mm>*{\circ};<0mm,0mm>*{}**@{},
 <0mm,-0.49mm>*{};<0mm,-3.0mm>*{}**@{-},
 <0.49mm,0.49mm>*{};<1.9mm,1.9mm>*{}**@{-},
 <-0.5mm,0.5mm>*{};<-1.9mm,1.9mm>*{}**@{-},
 <-2.3mm,2.3mm>*{\circ};<-2.3mm,2.3mm>*{}**@{},
 <-1.8mm,2.8mm>*{};<0mm,4.9mm>*{}**@{-},
 <-2.8mm,2.9mm>*{};<-4.6mm,4.9mm>*{}**@{-},
   <0.49mm,0.49mm>*{};<3.4mm,2.5mm>*{^{_{\sigma(3)}}}**@{},
   <-1.8mm,2.8mm>*{};<1.1mm,5.3mm>*{^{^{\sigma(2)}}}**@{},
   <-2.8mm,2.9mm>*{};<-5.5mm,5.3mm>*{^{^{\sigma(1)}}}**@{},
 \end{xy}
\ - \
\begin{xy}
 <0mm,0mm>*{\circ};<0mm,0mm>*{}**@{},
 <0mm,-0.49mm>*{};<0mm,-3.0mm>*{}**@{-},
 <0.49mm,0.49mm>*{};<1.9mm,1.9mm>*{}**@{-},
 <-0.5mm,0.5mm>*{};<-1.9mm,1.9mm>*{}**@{-},
 <2.3mm,2.3mm>*{\circ};<2.3mm,2.3mm>*{}**@{},
 <1.8mm,2.8mm>*{};<0mm,4.9mm>*{}**@{-},
 <2.8mm,2.9mm>*{};<4.6mm,4.9mm>*{}**@{-},
   <0.49mm,0.49mm>*{};<-3.4mm,2.5mm>*{^{_{\sigma(1)}}}**@{},
   <1.8mm,2.8mm>*{};<1.1mm,5.3mm>*{^{^{\sigma(2)}}}**@{},
   <2.8mm,2.9mm>*{};<5.5mm,5.3mm>*{^{^{\sigma(3)}}}**@{},
 \end{xy}
=0,  \ \ \ \ \forall \sigma\in \bS_3.
$$
Hence the minimal resolution, $({\sf Ass}_\infty,\delta)$ of $\sf
Ass$ contains a degree -1 corolla $
\begin{xy}
 <0mm,-0.55mm>*{};<0mm,-2.5mm>*{}**@{-},
 <0.5mm,0.5mm>*{};<2.2mm,2.2mm>*{}**@{-},
 <0mm,0.6mm>*{};<0mm,2.2mm>*{}**@{-},
 <-0.48mm,0.48mm>*{};<-2.2mm,2.2mm>*{}**@{-},
 <0mm,0mm>*{\circ};<0mm,0mm>*{}**@{},
   <0mm, 0.55mm>*{};<0mm,2.8mm>*{^2}**@{},
   <0.5mm,0.5mm>*{};<2.7mm,2.8mm>*{^3}**@{},
   <-0.48mm,0.48mm>*{};<-2.7mm,2.8mm>*{^1}**@{},
 \end{xy}
 $
 such that
$$
\delta
\begin{xy}
 <0mm,-0.55mm>*{};<0mm,-2.5mm>*{}**@{-},
 <0.5mm,0.5mm>*{};<2.2mm,2.2mm>*{}**@{-},
 <0mm,0.6mm>*{};<0mm,2.2mm>*{}**@{-},
 <-0.48mm,0.48mm>*{};<-2.2mm,2.2mm>*{}**@{-},
 <0mm,0mm>*{\circ};<0mm,0mm>*{}**@{},
   <0mm, 0.55mm>*{};<0mm,2.8mm>*{^2}**@{},
   <0.5mm,0.5mm>*{};<2.7mm,2.8mm>*{^3}**@{},
   <-0.48mm,0.48mm>*{};<-2.7mm,2.8mm>*{^1}**@{},
 \end{xy}
 =
 \begin{xy}
 <0mm,0mm>*{\circ};<0mm,0mm>*{}**@{},
 <0mm,-0.49mm>*{};<0mm,-3.0mm>*{}**@{-},
 <0.49mm,0.49mm>*{};<1.9mm,1.9mm>*{}**@{-},
 <-0.5mm,0.5mm>*{};<-1.9mm,1.9mm>*{}**@{-},
 <-2.3mm,2.3mm>*{\circ};<-2.3mm,2.3mm>*{}**@{},
 <-1.8mm,2.8mm>*{};<0mm,4.9mm>*{}**@{-},
 <-2.8mm,2.9mm>*{};<-4.6mm,4.9mm>*{}**@{-},
   <0.49mm,0.49mm>*{};<2.7mm,2.3mm>*{^{3}}**@{},
   <-1.8mm,2.8mm>*{};<0.4mm,5.3mm>*{^{2}}**@{},
   <-2.8mm,2.9mm>*{};<-5.1mm,5.3mm>*{^{1}}**@{},
 \end{xy}
 \ - \
\begin{xy}
 <0mm,0mm>*{\circ};<0mm,0mm>*{}**@{},
 <0mm,-0.49mm>*{};<0mm,-3.0mm>*{}**@{-},
 <0.49mm,0.49mm>*{};<1.9mm,1.9mm>*{}**@{-},
 <-0.5mm,0.5mm>*{};<-1.9mm,1.9mm>*{}**@{-},
 <2.3mm,2.3mm>*{\circ};<2.3mm,2.3mm>*{}**@{},
 <1.8mm,2.8mm>*{};<0mm,4.9mm>*{}**@{-},
 <2.8mm,2.9mm>*{};<4.6mm,4.9mm>*{}**@{-},
   <0.49mm,0.49mm>*{};<-3.4mm,2.5mm>*{^{1}}**@{},
   <1.8mm,2.8mm>*{};<1.1mm,5.3mm>*{^{2}}**@{},
   <2.8mm,2.9mm>*{};<5.5mm,5.3mm>*{^{3}}**@{},
 \end{xy}
$$
Therefore, in  its wheeled extension, $({\sf
Ass}_\infty^\circlearrowright,\delta)$, one has
\Beqrn \delta
\begin{xy}
 <0mm,-0.55mm>*{};<0mm,-2.5mm>*{}**@{-},
 <0.5mm,0.5mm>*{};<2.2mm,2.2mm>*{}**@{-},
 <0mm,0.6mm>*{};<0mm,3.2mm>*{}**@{-},
 <-0.48mm,0.48mm>*{};<-2.2mm,2.2mm>*{}**@{-},
 <0mm,0mm>*{\circ};<0mm,0mm>*{}**@{},
   %<0mm, 0.55mm>*{};<0mm,2.8mm>*{^1}**@{},
   <0.5mm,0.5mm>*{};<2.4mm,2.8mm>*{^2}**@{},
   <-0.48mm,0.48mm>*{};<-2.7mm,2.8mm>*{^1}**@{},
(0,3.3)*{}
   \ar@{->}@(ul,dl) (0,-2.5)*{}
 \end{xy}
& = &
\begin{xy}
 <0mm,0mm>*{\circ};<0mm,0mm>*{}**@{},
 <0mm,-0.49mm>*{};<0mm,-2.8mm>*{}**@{-},
 <0.49mm,0.49mm>*{};<1.9mm,1.9mm>*{}**@{-},
 <-0.5mm,0.5mm>*{};<-1.9mm,1.9mm>*{}**@{-},
 <-2.3mm,2.3mm>*{\circ};<-2.3mm,2.3mm>*{}**@{},
 <-1.8mm,2.8mm>*{};<0mm,4.9mm>*{}**@{-},
 <-2.8mm,2.9mm>*{};<-4.6mm,4.9mm>*{}**@{-},
   <0.49mm,0.49mm>*{};<2.7mm,2.3mm>*{^{2}}**@{},
   %<-1.8mm,2.8mm>*{};<0.4mm,5.3mm>*{^{2}}**@{},
   <-2.8mm,2.9mm>*{};<-5.1mm,5.3mm>*{^{1}}**@{},
   (0.3,5.2)*{}
   \ar@{->}@(ur,dr) (0,-2.5)*{}
 \end{xy}
 \ - \
\begin{xy}
 <0mm,0mm>*{\circ};<0mm,0mm>*{}**@{},
 <0mm,-0.49mm>*{};<0mm,-2.8mm>*{}**@{-},
 <0.49mm,0.49mm>*{};<1.9mm,1.9mm>*{}**@{-},
 <-0.5mm,0.5mm>*{};<-1.9mm,1.9mm>*{}**@{-},
 <2.3mm,2.3mm>*{\circ};<2.3mm,2.3mm>*{}**@{},
 <1.8mm,2.8mm>*{};<0mm,4.9mm>*{}**@{-},
 <2.8mm,2.9mm>*{};<4.6mm,4.9mm>*{}**@{-},
   <0.49mm,0.49mm>*{};<-3.4mm,2.5mm>*{^{1}}**@{},
   %<1.8mm,2.8mm>*{};<1.1mm,5.3mm>*{^{2}}**@{},
   <2.8mm,2.9mm>*{};<5.5mm,5.3mm>*{^{2}}**@{},
   (0,5.6)*{}
   \ar@{->}@(ul,dl) (0,-2.8)*{}
 \end{xy}
\\&=&
\begin{xy}
 <0mm,0mm>*{\circ};<0mm,0mm>*{}**@{},
 <0mm,-0.49mm>*{};<0mm,-2.8mm>*{}**@{-},
 <0.49mm,0.49mm>*{};<1.9mm,1.9mm>*{}**@{-},
 <-0.5mm,0.5mm>*{};<-1.9mm,1.9mm>*{}**@{-},
 <-2.3mm,2.3mm>*{\circ};<-2.3mm,2.3mm>*{}**@{},
 <-1.8mm,2.8mm>*{};<0mm,4.9mm>*{}**@{-},
 <-2.8mm,2.9mm>*{};<-4.6mm,4.9mm>*{}**@{-},
   <0.49mm,0.49mm>*{};<2.7mm,2.3mm>*{^{2}}**@{},
   %<-1.8mm,2.8mm>*{};<0.4mm,5.3mm>*{^{2}}**@{},
   <-2.8mm,2.9mm>*{};<-5.1mm,5.3mm>*{^{1}}**@{},
   (0.3,5.2)*{}
   \ar@{->}@(ur,dr) (0,-2.5)*{}
 \end{xy}
 \ - \
 \begin{xy}
 <0mm,0mm>*{\circ};<0mm,0mm>*{}**@{},
 <0mm,-0.49mm>*{};<0mm,-2.8mm>*{}**@{-},
 <0.49mm,0.49mm>*{};<1.9mm,1.9mm>*{}**@{-},
 <-0.5mm,0.5mm>*{};<-1.9mm,1.9mm>*{}**@{-},
 <-2.3mm,2.3mm>*{\circ};<-2.3mm,2.3mm>*{}**@{},
 <-1.8mm,2.8mm>*{};<0mm,4.9mm>*{}**@{-},
 <-2.8mm,2.9mm>*{};<-4.6mm,4.9mm>*{}**@{-},
   <0.49mm,0.49mm>*{};<2.7mm,2.3mm>*{^{2}}**@{},
   %<-1.8mm,2.8mm>*{};<0.4mm,5.3mm>*{^{2}}**@{},
   <-2.8mm,2.9mm>*{};<-5.1mm,5.3mm>*{^{1}}**@{},
   (0.3,5.2)*{}
   \ar@{->}@(ur,dr) (0,-2.5)*{}
 \end{xy}\\
 &=& 0,
\Eeqrn implying existence of a non-trivial cohomology class $
\begin{xy}
 <0mm,-0.55mm>*{};<0mm,-2.5mm>*{}**@{-},
 <0.5mm,0.5mm>*{};<2.2mm,2.2mm>*{}**@{-},
 <0mm,0.6mm>*{};<0mm,3.2mm>*{}**@{-},
 <-0.48mm,0.48mm>*{};<-2.2mm,2.2mm>*{}**@{-},
 <0mm,0mm>*{\circ};<0mm,0mm>*{}**@{},
   %<0mm, 0.55mm>*{};<0mm,2.8mm>*{^1}**@{},
   <0.5mm,0.5mm>*{};<2.4mm,2.8mm>*{^2}**@{},
   <-0.48mm,0.48mm>*{};<-2.7mm,2.8mm>*{^1}**@{},
(0,3.3)*{}
   \ar@{->}@(ul,dl) (0,-2.5)*{}
 \end{xy}
$ in  $\sH({\mathsf Ass}_\infty^\circlearrowright,\delta)$ which
does {\em not}\, belong to ${\sf Ass}^\circlearrowright$. Thus
$\mathsf Ass$ is Koszul, but {\em not}\, stably Koszul. \hfill
$\Box$

\bip

\bip

\no
It is instructive to see explicitly how the map $[i]: \sH(\LLL
\langle {\sf Ass}_\infty \rangle, \delta) \rar \sH(\PPP\langle
{\mathsf Ass}_\infty \rangle,\delta)$ fails to be a monomorphism. As
$\LLL\langle {\sf Ass}_\infty \rangle$ does not contain loops, the
element
$$
a:=\ \
\begin{xy}
 <0mm,0mm>*{\circ};<0mm,0mm>*{}**@{},
 <0mm,-0.49mm>*{};<0mm,-2.8mm>*{}**@{-},
 <0.49mm,0.49mm>*{};<1.9mm,1.9mm>*{}**@{-},
 <-0.5mm,0.5mm>*{};<-1.9mm,1.9mm>*{}**@{-},
 <-2.3mm,2.3mm>*{\circ};<-2.3mm,2.3mm>*{}**@{},
 <-1.8mm,2.8mm>*{};<0mm,4.9mm>*{}**@{-},
 <-2.8mm,2.9mm>*{};<-4.6mm,4.9mm>*{}**@{-},
   <0.49mm,0.49mm>*{};<2.7mm,2.3mm>*{^{2}}**@{},
   %<-1.8mm,2.8mm>*{};<0.4mm,5.3mm>*{^{2}}**@{},
   <-2.8mm,2.9mm>*{};<-5.1mm,5.3mm>*{^{1}}**@{},
   (0.3,5.2)*{}
   \ar@{.>}@(ur,dr) (0,-2.5)*{}
 \end{xy}
 \ - \
 \begin{xy}
 <0mm,0mm>*{\circ};<0mm,0mm>*{}**@{},
 <0mm,-0.49mm>*{};<0mm,-2.8mm>*{}**@{-},
 <0.49mm,0.49mm>*{};<1.9mm,1.9mm>*{}**@{-},
 <-0.5mm,0.5mm>*{};<-1.9mm,1.9mm>*{}**@{.},
 <-2.3mm,2.3mm>*{\circ};<-2.3mm,2.3mm>*{}**@{},
 <-1.8mm,2.8mm>*{};<0mm,4.9mm>*{}**@{-},
 <-2.8mm,2.9mm>*{};<-4.6mm,4.9mm>*{}**@{-},
   <0.49mm,0.49mm>*{};<2.7mm,2.3mm>*{^{2}}**@{},
   %<-1.8mm,2.8mm>*{};<0.4mm,5.3mm>*{^{2}}**@{},
   <-2.8mm,2.9mm>*{};<-5.1mm,5.3mm>*{^{1}}**@{},
   (0.3,5.2)*{}
   \ar@{->}@(ur,dr) (0,-2.5)*{}
 \end{xy}
$$
defines a non-trivial cohomology class, $[a]$, in $\sH(\LLL\langle
{\sf Ass}_\infty \rangle, \delta)$, whose image, $[i]([a])$, in $
\sH(\PPP\langle {\sf Ass}_\infty \rangle,\delta)$ vanishes.

\bip
%One of our main tasks in the rest of this section is to identify a class
%of stably Koszul dioperads which contains many important examples
%including the one in which we are most interested in this paper,
%$\Liebi$.

\bip
%%%%%%%%%%%%%%%%%%%%%%%%%%%%%%%%%%%%%%%%
\no{\bf 3.6. Koszul substitution laws}. Let $P=\{P(n)\}_{n\geq 1}$ and
 $Q=\{Q(n)\}_{n\geq 1}$ be two quadratic Koszul operads generated,
$$
P:= \frac{\PROP\langle E_P(2)\rangle}{<I_P>}, \ \ \ \ \ \
Q:= \frac{\PROP\langle E_Q(2)\rangle}{<I_Q>},
$$
  by
 $\bS_2$-modules $E_P(2)$, and,
 respectively, $E_Q(2)$.

 \bip

One can canonically associate \cite{MV} to such a pair the $\frac{1}{2}$prop,
 $P\diamond Q^\dag$, with
$$
P\diamond Q^\dag(m,n)=\left\{
\Ba{ll}
P(n) & {\mathrm for} \ m=1, n\geq 2,\\
Q(m) & {\mathrm for} \ n=1, m\geq 2,\\
0   & {\mathrm otherwise},
\Ea
\right.
$$
and the $\frac{1}{2}$prop compositions,
$$
\left\{
\hspace{1mm}_1\hspace{-0.2mm}\circ_j:
P\diamond Q^\dag(m_1,1)\ot P\diamond Q^\dag(m_2,n_2)\lon
P\diamond Q^\dag(m_1+m_2-1,n_2)
\right\}_{ 1\leq j\leq m_2}
$$
being zero for $n_2\geq 2$ and  coinciding with the operadic
composition in $Q$ for $n_2=1$, and
$$
\left\{
\hspace{1mm}_i\hspace{-0.2mm}\circ_1: P\diamond Q^\dag(m_1,n_1)\ot P\diamond Q^\dag(1,n_2)
\lon P\diamond Q^\dag(m_1+m_2-1,n_2)
\right\}_{1\leq i \leq n_1}
$$
being zero for $m_1\geq 2$ and otherwise coinciding with the operadic
composition in $P$ for $m_1=1$.

\bip

\no
Let $D_0=\Omega_{\frac{1}{2}\sP\rar \sD}\langle P\diamond
Q^\dag\rangle$ be the associated free dioperad, $D_0^!$ its Koszul
dual dioperad, and $(D_{0\infty}:={\mathbf D}D_0^!, \delta_0)$ the
associated cobar construction \cite{G}. As $D_0$ is Koszul
\cite{G,MV}, the latter
 provides us with the dioperadic minimal model of $D_0$. By
exactness of $\Omega_{\frac{1}{2}\sP\rar \sP}$, the dg free prop,
$(D_{0\infty}^\uparrow:=\Omega_{\sD\rar \sP}\langle D_{0\infty}\rangle, \delta_0)$, is the
minimal  model of the prop  $D_0^\uparrow:=\Omega_{\sD\rar \sP}\langle D_0\rangle\simeq
\Omega_{\frac{1}{2}\sP\rar \sP}\langle P\diamond Q^\dag\rangle$.

\bip

\no{\bf 3.6.1. Remark.} The prop  $D_0^\uparrow$ can be equivalently defined as the quotient,
$$
\frac{ P* Q^\dag}{I_0}
$$
where $P*Q^\dag$ is the free product of props associated to operads
$P$ and\footnote{the symbol $^\dag$ stands for the functor on props, $P=\{P(m,n)\}\rar P^\dag=\{P^\dag(m,n)\}$ which
reverses ``time flow", i.e.\ $P^\dag(m,n):= P(n,m)$.}  $Q^\dag$, and the ideal $I_0$ is generated by graphs of the
form,
$$
I_0={\mathsf span}\left\langle
\begin{xy}
 <0mm,2.47mm>*{};<0mm,-0.5mm>*{}**@{-},
 <0.5mm,3.5mm>*{};<2.2mm,5.2mm>*{}**@{-},
 <-0.48mm,3.48mm>*{};<-2.2mm,5.2mm>*{}**@{-},
 <0mm,3mm>*{\circ};<0mm,3mm>*{}**@{},
  <0mm,-0.8mm>*{\bullet};<0mm,-0.8mm>*{}**@{},
<0mm,-0.8mm>*{};<-2.2mm,-3.5mm>*{}**@{-},
 <0mm,-0.8mm>*{};<2.2mm,-3.5mm>*{}**@{-},
    % <0.5mm,3.5mm>*{};<2.8mm,5.7mm>*{^2}**@{},
    % <-0.48mm,3.48mm>*{};<-2.8mm,5.7mm>*{^1}**@{},
   %<0mm,-0.8mm>*{};<-2.7mm,-5.2mm>*{^1}**@{},
   %<0mm,-0.8mm>*{};<2.7mm,-5.2mm>*{^2}**@{},
\end{xy}\right\rangle \simeq  D_0(2,1)\ot D_0(1,2)=E_Q(2)\ot E_P(2)
$$
with white vertex decorated by elements of $E_Q(2)$ and black vertex decorated
by elements of $E_P(2)$.

\bip

\no
Let us consider
a morphism of $\bS_2$-bimodules,
$$
\Ba{rccc}
\lambda: & D_0(2,1)\ot D_0(1,2) & \lon & D_0(2,2) \\
&{\mathsf span}\left\langle
\begin{xy}
 <0mm,2.47mm>*{};<0mm,-0.5mm>*{}**@{-},
 <0.5mm,3.5mm>*{};<2.2mm,5.2mm>*{}**@{-},
 <-0.48mm,3.48mm>*{};<-2.2mm,5.2mm>*{}**@{-},
 <0mm,3mm>*{\circ};<0mm,3mm>*{}**@{},
  <0mm,-0.8mm>*{\bullet};<0mm,-0.8mm>*{}**@{},
<0mm,-0.8mm>*{};<-2.2mm,-3.5mm>*{}**@{-},
 <0mm,-0.8mm>*{};<2.2mm,-3.5mm>*{}**@{-},
    % <0.5mm,3.5mm>*{};<2.8mm,5.7mm>*{^2}**@{},
    % <-0.48mm,3.48mm>*{};<-2.8mm,5.7mm>*{^1}**@{},
   %<0mm,-0.8mm>*{};<-2.7mm,-5.2mm>*{^1}**@{},
   %<0mm,-0.8mm>*{};<2.7mm,-5.2mm>*{^2}**@{},
\end{xy} \right\rangle & \lon &
{\mathsf span}\left\langle
\begin{xy}
 <0mm,-1.3mm>*{};<0mm,-3.5mm>*{}**@{-},
 <0.38mm,-0.2mm>*{};<2.2mm,2.2mm>*{}**@{-},
 <-0.38mm,-0.2mm>*{};<-2.2mm,2.2mm>*{}**@{-},
<0mm,-0.8mm>*{\circ};<0mm,0.8mm>*{}**@{},
 <2.4mm,2.4mm>*{\bullet};<2.4mm,2.4mm>*{}**@{},
 <2.5mm,2.3mm>*{};<4.4mm,-0.8mm>*{}**@{-},
 <2.4mm,2.5mm>*{};<2.4mm,5.2mm>*{}**@{-},
    \end{xy}
\  , \
\begin{xy}
 <0mm,-1.3mm>*{};<0mm,-3.5mm>*{}**@{-},
 <0.38mm,-0.2mm>*{};<2.2mm,2.2mm>*{}**@{-},
 <-0.38mm,-0.2mm>*{};<-2.2mm,2.2mm>*{}**@{-},
<0mm,-0.8mm>*{\circ};<0mm,0.8mm>*{}**@{},
 <-2.4mm,2.4mm>*{\bullet};<2.4mm,2.4mm>*{}**@{},
 <-2.5mm,2.3mm>*{};<-4.4mm,-0.8mm>*{}**@{-},
 <-2.4mm,2.5mm>*{};<-2.4mm,5.2mm>*{}**@{-},
    \end{xy}
\right\rangle
\Ea
$$
and define \cite{MV} the dioperad, $D_\lambda$, as the quotient of the free
dioperad generated by the two spaces of binary operations,
$D_0(2,1)=E_Q(2)$ and $D_0(1,2)=E_P(2)$, modulo the
ideal generated  by
relations in $P$, relations in $Q$ as well as
the followings ones,
$$
I_\lambda={\mathsf span}\{ f - \lambda f : \forall f\in  D_0(2,1)\ot D_0(1,2)\}.
$$
Note that in notations of \S~3.6.1 the associated prop,
  $D_\lambda^\uparrow:=\Omega_{\mathsf D\rar P}\langle D_\lambda\rangle$,
is just the quotient,
${ P* Q^\dag}/{I_\lambda}$.

\bip

\no
The substitution law $\lambda$ is called {\em Koszul}, if $D_\lambda$ is
isomorphic
to $D_0$ as an $\bS$-bimodule. Which implies that $D_\lambda$ is Koszul
\cite{G}.
Koszul duality technique provides the space
${\mathbf D}D_0^!\simeq {\mathbf D}D_\lambda^!$ with a
perturbed differential $\delta_\lambda$ such that
 $({\mathbf D}D_0^!, \delta_\lambda)$ is the minimal model,
 $(D_{\lambda\infty}, \delta_\lambda)$,  of the
 dioperad $D_\lambda$.

\bip

 %%%%%%%%%%%%%%%%%%%%%%%%%%%%%%%%%%%%%%%%%%%%
\no{\bf 3.7. Theorem}  \cite{MV,V} .
{\em  The dg free prop
$D_{\lambda\infty}^\uparrow:=\Omega_{\mathsf D\rar P}\langle D_{\lambda\infty}\rangle$
is the minimal model of the prop $D_\lambda^\uparrow$, i.e.\ the natural
morphism
$$
(D_{\lambda\infty}^\uparrow, \delta_\lambda)
\lon (D_\lambda^\uparrow,0),
$$
which sends to zero all vertices of
$D_{\lambda\infty}^\uparrow$ except
binary ones decorated by elements of $E_P(2)$ and $E_Q(2)$, is a
quasi-isomorphism.}
\bip

\Proof The main point  is that
$$
F_p:= \left\{\mbox{span}\langle f\in  D_{\lambda\infty}^\uparrow \rangle
   \ : \ \Ba{l}\mbox{number of  directed paths in the graph}\ f\ \\
       \mbox{which connect
     input legs with output ones}
     \Ea
   \leq p
 \right\}.
$$
defines a  filtration of the complex $D_{\lambda\infty}^\uparrow$.
The associated spectral sequence, $\{E_r,d_r\}_{r\geq 0}$, is exhaustive
and bounded below so that it  converges to
the cohomology of
$(D_{\lambda\infty}^\uparrow, \delta_\lambda)$.

\bip

\no
The zeroth term of this spectral sequence is isomorphic to
$(D_{0\infty}^\uparrow,\delta_0)$ and hence, by Koszulness of the
dioperad $D_0$ and exactness of the functor $\Omega_{
\frac{1}{2}\sP\rar \sP}$, has the cohomology, $E_1$, isomorphic to
$D_0^\uparrow$.  Which, by
Koszulness of $D_\lambda$, is isomorphic to $D_\lambda^\uparrow$ as an $\bS$-bimodule. Hence
$\{d_r=0\}_{r\geq 1}$ and the result follows along the same lines as
in the second part of the proof of Theorem 2.6.2. \hfill $\Box$

\bip

 %%%%%%%%%%%%%%%%%%%%%%%%%%%%%%%%%%%%%%%%%%%%

%%%%%%%%%%%%%%%%%%%%%%%%%%%%%%%%%%%%%%%%
\no{\bf 3.8. Cohomology of graph complexes with marked wheels}.
In this section we analyze the functor
 $\Omega_{\frac{1}{2}\sP\rar \PPP}$. The following statement
 is one of the motivations for its introduction (it does {\em not}\, hold true
 for the ``unmarked" version $\Omega_{\frac{1}{2}\sP\rar \PP}$).

\bip

%%%%%%%%%%%%%%%%%%%%%%%%%%%%%%%%%%%%%%%%
\no{\bf 3.8.1. Theorem}. {\em The functor
 $\Omega_{\frac{1}{2}\sP\rar \PPP}$
is exact}.

\bip

\Proof Let $T$ be an arbitrary dg $\frac{1}{2}$prop. The main point
is that we can use $\frac{1}{2}$prop compositions and presence of
marks on cyclic edges to represent
 $\Omega_{\frac{1}{2}\sP\rar \PPP}\langle T \rangle$ as
a vector space {\em freely}\, generated by a family of decorated graphs,
$$
\Omega_{\frac{1}{2}\sP\rar \PPP}\langle T \rangle(m,n)=
\bigoplus_{G\in \overline{\fGGG}(m,n)} G\langle P \rangle,
$$
where $ \overline{\fGGG}(m,n)$
is a subset of $\fGGG(m,n)$ consisting of so called
{\em reduced}\, graphs, $G$, which
satisfy the following defining
property: for each pair of internal vertices,
$(v_1,v_2)$, of $G$ which are connected by an unmarked edge directed from
$v_1$ to $v_2$ one has $|Out(v_1)|\geq 2$ and $|In(v_2)|\geq 2$.
Put another way, given an arbitrary $T$-decorated graph with
wheels, one can perform $\frac{1}{2}$prop compositions (``contractions") along
all unmarked internal edges $(v_1,v_2)$ which do not satisfy the above
conditions. The result is a reduced decorated graph (with wheels).
Which is uniquely defined by the original one. Notice that marks are vital
for this contraction procedure, e.g.
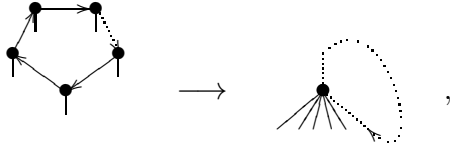
$$
\begin{xy}
<0mm,0mm>*{\bullet};<0mm,0mm>*{}**@{},
<-7mm,5mm>*{\bullet};<0mm,0mm>*{}**@{},
<-4mm,11mm>*{\bullet};<0mm,0mm>*{}**@{},
<7mm,5mm>*{\bullet};<0mm,0mm>*{}**@{},
<4mm,11mm>*{\bullet};<0mm,0mm>*{}**@{},
 <0mm,0mm>*{};<0mm,-3mm>*{}**@{-},
 <-7mm,5mm>*{};<-7mm,2mm>*{}**@{-},
  <-4mm,11mm>*{};<-4mm,8mm>*{}**@{-},
   <7mm,5mm>*{};<7mm,2mm>*{}**@{-},
    <4mm,11mm>*{};<4mm,8mm>*{}**@{-},
{(0,0)*{}\ar@{->} (-7,5)*{}}
\hspace{-7mm}\mbox{\xy (-7,5)*{}\ar@{->} (-4,11)*{}\endxy}
\hspace{-2mm}\mbox{\xy (-4,11)*{}\ar@{->} (4,11)*{}\endxy}
\hspace{-0.5mm}\mbox{\xy (4,11)*{}\ar@{.>} (7,5)*{}\endxy}
\hspace{-8.3mm}\mbox{\xy (7,5)*{}\ar@{->} (0,0)*{}\endxy}
 \end{xy}
\ \ \ \  \lon \ \ \ \
\begin{xy}
 <0mm,0mm>*{\bullet};<0mm,0mm>*{}**@{},
 <0mm,0mm>*{};<0mm,5mm>*{}**@{.},
   %<0mm,0mm>*{};<0mm,5.5mm>*{^{y}}**@{},
   %
<0mm,0mm>*{\bullet};<0mm,0mm>*{}**@{},
 <0mm,0mm>*{};<-6mm,-5mm>*{}**@{-},
 <0mm,0mm>*{};<-3.1mm,-5mm>*{}**@{-},
 <0mm,0mm>*{};<-1.2mm,-5mm>*{}**@{-},
  <0mm,0mm>*{};<1.2mm,-5mm>*{}**@{-},
 %<0mm,0mm>*{};<1.0mm,-5mm>*{}**@{-},
 <0mm,0mm>*{};<3.1mm,-5mm>*{}**@{-},
 <0mm,0mm>*{};<6mm,-5mm>*{}**@{.},
(0,5)*{}
   \ar@{.>}@(ur,dr) (6,-5)*{}
 \end{xy} \ \ \ ,
$$
to be well-defined.

\bip

\no
Then we have \Beqrn H^*\left(\Omega_{\frac{1}{2}\sP\rar \PPP}\langle
T \rangle (m,n)\right)&=& H^*\left( \bigoplus_{G\in \bar{\mathfrak
G}^+(m,n)} \left(
\bigotimes_{v\in v(G)} T(Out(v), In(v))\right)_{Aut G}\right)\\
&=& \bigoplus_{G\in \bar{\mathfrak G}^+(m,n)}H^* \left(
\bigotimes_{v\in v(G)} T(Out(v), In(v))\right)_{Aut G} \
\mbox{by Maschke's theorem}\\
&=& \bigoplus_{G\in \bar{\mathfrak G}^+(m,n)}\left(
\bigotimes_{v\in v(G)}H^*(T)(Out(v), In(v))\right)_{Aut G} \
\mbox{by K\"unneth formula}\\
&=& \Omega_{\frac{1}{2}\sP\rar \PPP}\langle H^*(T)\rangle(m,n).
\Eeqrn In the second line we used the fact that the group $AutG$ is
finite. \hfill $\Box$

\bip

\no
Another motivation for introducing graph complexes with {\em marked}\,
wheels is that they admit a filtration which singles out the
$\frac{1}{2}$propic part of the differential. A fact which we heavily use
in the proof of the following

\bip
 %%%%%%%%%%%%%%%%%%%%%%%%%%%%%%%%%%%%%%%%%%%%
\no{\bf 3.8.2. Theorem}. {\em Let $D_\lambda$ be a Koszul dioperad
with Koszul substitution law and let $(D_{\lambda\infty},\, \delta)$ be
its minimal resolution.
The natural morphism of graph complexes,
$$
(D_{\lambda\infty}^+,\, \delta_\lambda) \lon
(D_\lambda^+, 0)
%(\Omega_{\mathsf D\rar P_+^+}\langle D_\infty^\lambda\rangle, \delta_\lambda)
%\lon \Omega_{\mathsf D\rar P_+^+}\langle D_\lambda\rangle,
$$
is a quasi-isomorphism.}

\bip

\Proof Consider first a filtration of the complex
$(D_{\lambda\infty}^+,\, \delta_\lambda)$
 by the number of marked edges, and let $(D_{\lambda\infty}^+,\, b)$
denote 0th term of the associated spectral sequence (which, as we shall
show below, degenerates at the 1st term).

\bip

\no
To any decorated graph $f \in D_{\lambda\infty}^+$ one can
associate a graph without wheels, $\overline{f} \in D_{\lambda\infty}^\uparrow$, by breaking every marked
cyclic edge into two legs (one of which is input and another one is
output). Let $|\overline{f}|$ be the number of  directed paths in
the graph $\overline{f}$ which connect input legs with output ones.
Then
$$
F_p:= \left\{f\in D_{\lambda\infty}^+
   \ : \ |\overline{f}|
   \leq p
 \right\}.
$$
defines a  filtration of the complex
$(D_{\lambda\infty}^+,\, b)$.
%The associated spectral sequence, $\{E_r,d_r\}_{r\geq 0}$, is exhaustive
%and bounded below so that it  converges to
%the cohomology of
%$(\Omega_{\mathsf D\rar P}\langle D_\infty^\lambda\rangle, \delta)$.

\bip

\no
The zeroth term of the spectral sequence, $\{E_r,d_r\}_{r\geq 0}$,
associated to this filtration is isomorphic to
$(\langle D_{0\infty}^+,
\delta_0)$ and hence, by Theorem 3.8.1, has the cohomology, $E_1$,
equal to $D_0^+$. Which, by Koszulness of $D_\lambda$, is isomorphic as a
vector space to $D_{\lambda}^+$. Hence the  differentials of all higher
terms of both our spectral sequences vanish, and
 the result follows.
\hfill $\Box$

\bip

 %%%%%%%%%%%%%%%%%%%%%%%%%%%%%%%%%%%%%%%%%%%%
\no{\bf 3.8.3. Remark}. In the proof of Theorem 3.8.2 the
$\frac{1}{2}$propic part, $\delta_0$, of the differential
$\delta_\lambda$ was singled out in two steps: first we introduced a
filtration by the number of marked edges, and then a filtration by
the number of paths, $|\bar{f}|$, in the unwheeled graphs
$\overline{f}$. As the following lemma shows, one can do it in one
step. Let $w(f)$ stand for the number of marked edges in a decorated
graph $f\in D_{\lambda\infty}^+$.

\bip

 %%%%%%%%%%%%%%%%%%%%%%%%%%%%%%%%%%%%%%%%%%%%
\no{\bf 3.8.4. Lemma} {\em The sequence of vector spaces spaces, $p\in \N$,
$$
{\mathsf F}_p:= \left\{\mbox{span}\langle f\in D_{\lambda\infty}^+\rangle
%\Omega_{\mathsf D\rar P_+^+}\langle D_\infty^\lambda\rangle
   \ : \ ||f||:= 3^{w(f)}|\overline{f}|
   \leq p
 \right\},
$$
defines a  filtration of the complex
$(D_{\lambda\infty}^+, \delta_\lambda)$
%$(\Omega_{\mathsf D\rar P^+_+}\langle D_\infty^\lambda\rangle, \delta)$
whose spectral sequence has  0-th term isomorphic to
$(D_{0\infty}^+, \delta_0)$.}
%$(\Omega_{\mathsf \frac{1}{2}P\rar P^+_+}\langle D^0_\infty\rangle$

\bip

\Proof It is enough to show that for any graph $f$ in $D_{\lambda\infty}^+$
%$\Omega_{\mathsf D\rar P_+^+}\langle D_\infty^\lambda\rangle$
with $w(f)\neq 0$
one has, $||\delta_\lambda f||\leq ||f||$.

\bip

\no
We can, in general, split $\delta_\lambda f$ into two groups of summands,
$$
\delta_\lambda f =\sum_{a\in I_1} g_a + \sum_{b\in I_2} g_b
$$
where $w(g_a)=w(f)$, $\forall a\in I_1$, and  $w(g_b)=w(f)-p_b$ for some $p_b\geq 1$
and all  $b\in I_2$.

\bip

\no
For any $a\in I_1$,
$$
||g_a||=3^{w(f)}|\overline{g_a}| \leq 3^{w(f)}|\overline{f}| =||f||.
$$
So it remains to check the inequality $||g_b||\leq ||f||$, $\forall b\in I_2$.

\bip

\no
There is an associated splitting  of $\delta_\lambda \overline{f}$ into two groups of
summands,
$$
\delta_\lambda \overline{f} =\sum_{a\in I_1} h_a + \sum_{b\in I_2} h_b
$$
where $\{h_b\}_{b\in I_2}$ is the set of all those summands
which contain two-vertex subgraphs of the form,
$$
 
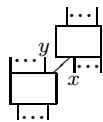
\xy <0mm,0mm>*{\mbox{$\xy *=<6mm,4mm>
\txt{}*\frm{-}\endxy$}};<0mm,0mm>*{}**@{},
   <1.5mm,2mm>*{};<1.5mm,4mm>*{}**@{-},
    <-3mm,2mm>*{};<-3mm,4mm>*{}**@{-},
  <2.7mm,2mm>*{};<5mm,4mm>*{}**@{-},
  <-2mm,-2.2mm>*{};<-2mm,-4.2mm>*{}**@{-},
  <0mm,-2.2mm>*{};<0mm,-3.9mm>*{...}**@{},
  <2mm,-2.2mm>*{};<2mm,-4.2mm>*{}**@{-},
<6mm,6.2mm>*{\mbox{$\xy *=<6mm,4mm>
\txt{}*\frm{-}\endxy$}};<0mm,6mm>*{}**@{},
<1.5mm,2mm>*{};<1.5mm,4.6mm>*{^y}**@{},
<-1mm,2mm>*{};<-1mm,3.6mm>*{...}**@{},
 <5.5mm,2mm>*{};<5.5mm,4mm>*{}**@{-},
 <5.5mm,2mm>*{};<5.5mm,1mm>*{_x}**@{},
 <9mm,2mm>*{};<9mm,4mm>*{}**@{-},
 <9mm,2mm>*{};<7.3mm,2.8mm>*{...}**@{},
<8.5mm,8.1mm>*{};<8.5mm,10.6mm>*{}**@{-},
<4.5mm,8.1mm>*{};<4.5mm,10.6mm>*{}**@{-},
<8.5mm,8.1mm>*{};<6.3mm,9.6mm>*{...}**@{},
  \endxy
$$
having half-edges of the type $x$ and $y$ corresponding to broken wheeled paths in $f$.
Every graph $g_b$ is obtained from the corresponding $h_b$ by gluing
some number of path connected to $y$ output legs  with the same number of path connected to $x$  input legs into new internal
non-cyclic edges.
This gluing operation creates $p_b$ new paths connecting some internal vertices in $\overline{h_b}$, and hence
may increase the total number of paths in $\overline{h_b}$ but no more than by the factor of $p_b+1$,
i.e.\ $|\overline{g_b}|\leq (p_b+1)|h_b|$, $\forall b\in I_2$.

\bip

\no
Finally, we have
$$
||g_b||  = 3^{w(f)-p_b}|\overline{g_b}| \leq  3^{w(f)-p_b}(p_b+1)|{h_b}|
< 3^{w(f)}|\overline{f}|= ||f||, \ \ \ \forall b\in I_2.
$$
The part of the differential $\delta_\lambda$ which preserves the filtration
 must in fact  preserve both
the number of marked edges, $w(f)$, and the number of paths, $|\overline{f}|$,
for any decorated graph $f$. Hence this is precisely $\delta_0$.
\hfill $\Box$

\bip

%%%%%%%%%%%%%%%%%%%%%%%%%%%%%%%%%%%%%%
\no\mbox{\bf 3.9.~Graph~complexes~with~unmarked wheels built on
$\frac{1}{2}$props.}

\no
Let $ \left (T= {{\frac{1}{2}\sP}\langle E
\rangle}/{<I>},\,\delta\right) $
 be a
dg $\frac{1}{2}$prop.
In \S 3.5 we defined its wheeled extension,
%$\left(T^\circlearrowright:=\Omega_{\mathsf \frac{1}{2}P\rar \PP}
%\langle H\rangle, \delta\right)$
$$
\left(T^\circlearrowright
%\Omega_{\mathsf \frac{1}{2}P\rar {P}^+}\langle H\rangle
:=
\frac{{\PP}\langle E \rangle}{<I>^\circlearrowright}, \delta\right).
$$
Now we specify a dg subprop,  $\Omega_{\mathsf no-oper}
\langle T\rangle\subset T^\circlearrowright$, whose cohomology is easy
to compute.

\bip

%%%%%%%%%%%%%%%%%%%%%%%%%%%%%%%%%%%%%%%%%
\no{\bf 3.9.1. Definition.}
Let $E=\{E(m,n)\}_{m,n\geq 1, m+n\geq 3}$ be an $\bS$-bimodule,
and $\PP\langle E \rangle$ the associated prop of decorated graphs
with  wheels. We say that a wheel $W$ in a graph
$G\in \PP\langle E \rangle$ is {\em operadic}\, if all its cyclic
vertices $v\in W$ are
%\Bi
%\item[] either of type $(1,n)$, $n\geq 2$, only
%\item[] or of type $(m,1)$, $m\geq 2$ only
%\Ei
%are
decorated either by elements of $\{E(1,n_v)\}_{n_v\geq 2}$ only, or
by elements $E(n_v,1)_{n_v\geq 2}$ only.
Vertices of operadic wheels are called {\em operadic cyclic vertices}.
Notice that operadic wheels can be of geometric genus 1 only.

\bip

\no
Let $\PP_{\mathsf no-oper}\langle E\rangle$
be the subspace of
$\PP\langle E\rangle$ consisting of graphs with
no operadic wheels, and let
$$
\Omega_{\mathsf no-oper}\langle T\rangle
=
\frac{\PP_{\mathsf no-oper}\langle E \rangle}{<I>^\circlearrowright},
$$
be the associated dg sub-prop of $(T^\circlearrowright,\, \delta)$.
%$(\Omega_{\mathsf \frac{1}{2}P\rar \PP}\langle H\rangle, \delta)$

\bip

\no
 Clearly, $\Omega_{\mathsf no-oper}$ is
a functor from the category of dg $\frac{1}{2}$props to the category
of dg props. It is worth pointing out that this functor can {\em not}\,
be extended to dg dioperads as differential can, in general, create operadic
wheels from non-operadic ones.
\bip

 %%%%%%%%%%%%%%%%%%%%%%%%%%%%%%%%%%%%%%%%%%%%
\no{\bf 3.9.2. Theorem}. {\em The functor \, $\Omega_{\mathsf no-oper}$
is exact.}

\bip

\Proof Let $(T,\delta)$ be an arbitrary dg $\frac{1}{2}$prop.
Every wheel in $\Omega_{\mathsf no-oper}\langle T\rangle$
contains at least one cyclic edge along which $\frac{1}{2}$prop
composition in $T$ is not possible. This fact allows one to
non-ambiguously perform such compositions along all those cyclic and
non-cyclic edges at which such a composition makes sense, and
hence represent  $\Omega_{\mathsf no-oper}\langle T\rangle$ as
a vector space {\em freely}\, generated by a family of decorated graphs,
$$
 \Omega_{\mathsf no-oper}\langle T\rangle(m,n)=
\bigoplus_{G\in \overline{\fGG}(m,n)} G\langle T \rangle
$$
where $ \overline{\fGG}(m,n)$
is a subset of $\fGG(m,n)$ consisting of {\em reduced}\, graphs, $G$, which
satisfy the following defining
properties: (i) for each pair of internal vertices,
$(v_1,v_2)$, of $G$ which are connected by an edge directed from
$v_1$ to $v_2$ one has $|Out(v_1)|\geq 2$ and $|In(v_2)|\geq 2$;
(ii) there are no operadic wheels in $G$.
The rest of the proof is exactly the same as in \S 3.8.1.
\hfill $\Box$

\bip

\no
Let $P$ and $Q$ be Koszul operads and let $D_0$ be the associated Koszul
dioperad (defined in \S 3.6) whose minimal resolution
is denoted by $(D_{0\infty},\delta_0)$.

\bip

%%%%%%%%%%%%%%%%%%%%%%%%%%%%%%%%%%%%%%%%%%%%
\no{\bf 3.9.3. Corollary}.
 $H(\Omega_{\mathrm no-oper}\langle D_{0\infty}\rangle, \delta_0)=
\Omega_{\frac{1}{2}\sP\rar {P}}\langle D_0\rangle$.

\bip

\Proof By Theorem~3.9.2,
$$
H(\Omega_{\mathrm no-oper}\langle D_{0\infty}\rangle, \delta_0)
= \Omega_{\mathrm no-oper}\langle H(D_{0\infty},\delta_0)\rangle
= \Omega_{\mathrm no-oper}\langle D_0\rangle.
$$
But the latter space can not have graphs with wheels as any such a wheel
would contain at least one ``non-reduced" internal cyclic edge
corresponding to composition,
$$
\circ_{1,1}: D_0(m,1)\ot D_0(1,n) \lon D_0(m,n),
$$
which is zero by the definition of $D_0$ (see \S 3.6).
\hfill $\Box$

\bip

%%%%%%%%%%%%%%%%%%%%%%%%%%%%%%%%%%%%%%%%%%%%
\no{\bf 3.10. Theorem}. {\em For any Koszul operads $P$ and $Q$,

(i)
 the natural morphism of graph complexes,
$$
(D_{0\infty}^\circlearrowright,\, \delta_0) \lon (D_0^\circlearrowright,0)
%(\Omega_{\mathsf D\rar P^+}\langle D_\infty^0\rangle, \delta_0)
%\lon \Omega_{\mathsf D\rar P^+}\langle D_0\rangle,
$$
is a quasi-isomorphism if and only if the operads $P$ and $Q$ are
stably Koszul;

\bip

(ii) there is, in general, an isomorphism of $\bS$-bimodules,
$$
H(D_{0\infty}^\circlearrowright,\, \delta_0)=
\frac{H(P^\circlearrowright_\infty)*H(Q^\circlearrowright_\infty)^\dag}{I_0}
$$
where  $H(P^\circlearrowright_\infty)$ and $H(Q^\circlearrowright_\infty)$ are
cohomologies of the wheeled completions of the minimal resolutions
of the operads $P$ and $Q$, $*$ stands for the free product of PROPs, and
the ideal $I_0$ is defined in \S~3.6.1.  }

\bip

\Proof (i) The necessity of the condition is obvious. Let us prove
its sufficiency.

\bip

\no
Let $P$ and $Q$ be stably Koszul operads so that the natural morphisms,
$$
(P_\infty^\circlearrowright, \delta_P) \rar P^\circlearrowright
\ \ \ \ \mbox{and} \ \ \ \
(Q_\infty^\circlearrowright, \delta_Q) \rar Q^\circlearrowright,
$$
are quasi-isomorphisms, where $(P_\infty,\delta_P)$ and
$(Q_\infty,\delta_Q)$
are minimal resolutions of $P$ and $Q$ respectively.

\bip

\no
 Consider a filtration of the complex
 $(D_{0\infty}^\circlearrowright,\, \delta_0)$,
 %$(\Omega_{\mathsf D\rar P^+}\langle D_\infty^0\rangle, \delta_0)$,
 $$
 F_p:= \{ \mbox{span} \langle f\in D_{0\infty}^\circlearrowright \rangle:\
 %(\Omega_{\mathsf D\rar P^+}\langle D_\infty^0\rangle:\
 |f|_2 - |f|_1\ \leq p\},
 $$
where
\Bi
\item[-]$|f|_1$ is the number of cyclic vertices in $f$ which belong
to operadic wheels;
\item[-] $|f|_2$ is the number of non-cyclic half-edges
attached to cyclic vertices in $f$ which belong
to operadic wheels.
\Ei
Note that $|f|_2 - |f|_1\geq 0$.
Let $\{E_r,d_r\}_{r\geq 0}$ be the associated spectral sequence.
The differential
$d_0$ in $E_0$ is given by its values on the vertices as follows:
\Bi
\item[(a)]  on every non-cyclic vertex and on every cyclic vertex
which does {\em not}\, belong to an operadic wheel one has $d_0=\delta_0$;
\item[(b)] on every cyclic vertex which  belongs to
an operadic wheel one has  $d_0=0$.
\Ei
Hence modulo the action of finite groups (which we can ignore by
Maschke theorem)
the complex $(E_0, d_0)$ is isomorphic to the complex
$(\Omega_{\mathrm no-oper}\langle D_{0\infty}\rangle, \delta_0)$,
tensored with a trivial complex (i.e.\ one with vanishing differential).
By Corollary~3.9.3 and K\"{u}nneth formula we obtain,
$$
E_1 = H(E_0,d_0)= W_1/h(W_2)
$$
where
\Bi
\item[-] $W_1$ is the subspace of $\PP\langle E_P\oplus E_Q^\dag \rangle$
consisting of graphs whose wheels (if any) are operadic; here
the $\bS$-bimodule $E_P\oplus E_Q^\dag$ is given by
$$
 (E_P\oplus E_Q^\dag)(m,n)=\left\{
\Ba{ll}
E_P(2), \mbox{the space of generators of}\ P, & \mbox{if}\ m=1,n=2\\
E_Q(2), \mbox{the space of generators of}\ Q, & \mbox{if}\ m=2,n=1\\
 0, & \mbox{otherwise;}\\
\Ea
 \right.
$$
\item[-] $W_2$ is the subspace  of  $\PP\langle E_P\oplus E_Q^\dag \oplus
{I}_P \oplus {I}_Q^\dag\rangle$
consisting of graphs, $G$, whose wheels (if any) are operadic and satisfy
the following condition: the elements of ${I}_P$ and $I^\dag_Q$ are used to
decorate at least one non-cyclic vertex in $G$. Here $I_P$ and
$I_Q^\dag$ are $\bS$-bimodules of relations of the quadratic operads $P$ and
$Q^\dag$ respectively.
\item[-] the map $h: W_2\rar W_1$ is defined to be
the identity on vertices
decorated by elements of $E_P\oplus E_Q^\dag$, and the tautological
(in the obvious sense) morphism on vertices decorated by elements
of  $I_P$ and $I_Q^\dag$.
\Ei

\sip

\no
To understand all the remaining terms $\{E_r,d_r\}_{r\geq 1}$ of the
spectral sequence we step aside and contemplate for a moment
 a purely operadic graph complex with wheels, say,
$(P_\infty^\circlearrowright, \delta_P)$.
%and $(\Omega_{\mathsf D\rar {P}^+}\langle Q_\infty\rangle, \delta_Q)$.

\bip

\no
The  complex  $(P_\infty^\circlearrowright,\, \delta_P)$
is naturally a subcomplex of
 $(D_{0\infty}^\circlearrowright, \, \delta_0)$.
 Let
 %$(\Omega_{\mathsf D\rar P^+}\langle P_\infty\rangle, \delta_P)$,
 $$
 F_p:= \{  \mbox{span} \langle f\in P_\infty^\circlearrowright \rangle:\
 |f|_2 - |f|_1\ \leq p\},
 $$
 be the induced filtration, and let $\{E^P_r,d_r^P\}_{r\geq 0}$ be the
 associated spectral sequence. Then $E_1^P=H(E_0^P,d_0^P)$ is a
 subcomplex of $E_1$.

\bip

\no
The main point is that, modulo the action of finite groups,
the spectral sequence $\{E_r,d_r\}_{r\geq 1}$ is isomorphic
to the tensor product of spectral sequences of the form
 $\{E^P_r,d_r^P\}_{r\geq 1}$ and  $\{E^Q_r,d_r^Q\}_{r\geq 1}$.
 By assumption, the latter converge to
$P^\circlearrowright$ and
$Q^\circlearrowright$ respectively.
Which implies the result.

\bip

(ii) The argument is exactly the same as in (i) except the very last paragraph:
the spectral sequences of the form
 $\{E^P_r,d_r^P\}_{r\geq 1}$ and  $\{E^Q_r,d_r^Q\}_{r\geq 1}$ converge, respectively,
to $H(P_\infty^\circlearrowright)$ and to
 $H(Q_\infty^\circlearrowright)$ (rather than to $P^\circlearrowright$ and
$Q^\circlearrowright$).
 \hfill $\Box$

\bip

\no{\bf 3.11. Operadic wheeled extension}. Let $D_\lambda$ be a dioperad and $D_{\lambda\infty}$ its minimal
resolution. Let $D_{\lambda\infty}^\looparrowright$ be a dg subprop of $D_{\lambda\infty}^\circlearrowright$
spanned by graphs with at most operadic wheels (see \S 3.9.1). Similarly one defines
a subprop, $D_\lambda^\looparrowright$, of $D_\lambda^\circlearrowright$.

\bip

\no{\bf 3.11.1. Theorem}. {\em For any Koszul operads $P$ and $Q$ and any Koszul substitution law
$\lambda$,

\bip
(i)  the natural morphism of graph complexes,
$$
(D_{\lambda\infty}^\looparrowright, \, \delta_\lambda)
\lon D_\lambda^\looparrowright,
$$
is a quasi-isomorphism if and only if the operads $P$ and $Q$ are
both stably Koszul.

\bip

(ii) there is, in general, an isomorphism of $\bS$-bimodules,
$$
H(D_{\lambda\infty}^\looparrowright,\, \delta_\lambda)=
H(D_{0\infty}^\looparrowright,\, \delta_0)=
\frac{H(P^\circlearrowright_\infty)*H(Q^\circlearrowright_\infty)^\dag}{I_0}.
$$
where  $H(P^\circlearrowright_\infty)$ and $H(Q^\circlearrowright_\infty)$ are
cohomologies of the wheeled completions of the minimal resolutions
of the operads $P$ and $Q$.  }

\bip

\Proof Use spectral sequence of a filtration, $\{F_p\}$, defined similar to the one introduced
in the  proof of Theorem 3.10. We omit full details
as they are analogous to \S 3.10.
\hfill $\Box$

\bip

In the next section we apply some of the above results to compute cohomology of several
concrete graph complexes with wheels.

\bip

%%%%%%%%%%%%%%%%%%%%%%%%%% 3 %%%%%%%%%%%%%%%%%%%%%%%%%%%%%%%
%%%%%%%%%%%%%%%%%%%%%%%%%%%%%%%%%%%%%%%%%%%%%%%%%%%%%%%%%
\begin{center}
\bf \S 4.  Wheeled Poisson structures and other examples
\end{center}

\bip

%%%%%%%%%%%%%%%%%%%%%%%%%%%%%%%%%%%%%%%%%%%%%%%%%%%%%%%%%%%%%%%%%%%
\no{\bf 4.1. Wheeled operad of strongly homotopy Lie algebras.}
Let $({\Lie_\infty}, \delta)$ be the minimal resolution of the operad,
$\Lie$, of Lie algebras.  It can be identified with the subcomplex of
$(\Liebi_\infty, \delta)$ spanned by connected trees built on degree one
$(1,n)$-corollas, $n\geq 2$,
$$
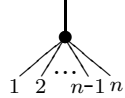
\begin{xy}
 <0mm,0mm>*{\bullet};<0mm,0mm>*{}**@{},
 <0mm,0mm>*{};<0mm,5mm>*{}**@{-},
   %<0mm,0mm>*{};<0mm,5.5mm>*{^{y}}**@{},
   %
<0mm,0mm>*{\bullet};<0mm,0mm>*{}**@{},
 <0mm,0mm>*{};<-6mm,-5mm>*{}**@{-},
 <0mm,0mm>*{};<-3.1mm,-5mm>*{}**@{-},
 <0mm,0mm>*{};<0mm,-4.6mm>*{...}**@{},
 %<0mm,0mm>*{};<1.0mm,-5mm>*{}**@{-},
 <0mm,0mm>*{};<3.1mm,-5mm>*{}**@{-},
 <0mm,0mm>*{};<6mm,-5mm>*{}**@{-},
   <0mm,0mm>*{};<-6.7mm,-6.4mm>*{_1}**@{},
   <0mm,0mm>*{};<-3.2mm,-6.4mm>*{_2}**@{},
   <0mm,0mm>*{};<3.1mm,-6.4mm>*{_{n\mbox{-}1}}**@{},
   <0mm,0mm>*{};<6.9mm,-6.4mm>*{_{n}}**@{},
 \end{xy}
$$
with the differential given by
\Beqrn
\delta\ \begin{xy}
 <0mm,0mm>*{\bullet};<0mm,0mm>*{}**@{},
 <0mm,0mm>*{};<0mm,5mm>*{}**@{-},
   %<0mm,0mm>*{};<0mm,5.5mm>*{^{y}}**@{},
   %
<0mm,0mm>*{\bullet};<0mm,0mm>*{}**@{},
 <0mm,0mm>*{};<-6mm,-5mm>*{}**@{-},
 <0mm,0mm>*{};<-3.1mm,-5mm>*{}**@{-},
 <0mm,0mm>*{};<0mm,-4.6mm>*{...}**@{},
 %<0mm,0mm>*{};<1.0mm,-5mm>*{}**@{-},
 <0mm,0mm>*{};<3.1mm,-5mm>*{}**@{-},
 <0mm,0mm>*{};<6mm,-5mm>*{}**@{-},
   <0mm,0mm>*{};<-6.7mm,-6.4mm>*{_1}**@{},
   <0mm,0mm>*{};<-3.2mm,-6.4mm>*{_2}**@{},
   <0mm,0mm>*{};<3.1mm,-6.4mm>*{_{n\mbox{-}1}}**@{},
   <0mm,0mm>*{};<6.9mm,-6.4mm>*{_{n}}**@{},
 \end{xy}
 &=&
 %
 %                                2nd TERM
 %
 \sum_{[n]=J_1\sqcup J_2\atop {\atop
 {|J_1|\geq 2, |J_2|\geq 1}}
}\hspace{0mm}
\ \
\begin{xy}
 <0mm,0mm>*{\bullet};<0mm,0mm>*{}**@{},
 <0mm,0mm>*{};<0mm,5mm>*{}**@{-},
   %<0mm,0mm>*{};<0mm,5.5mm>*{^{y}}**@{},
   %
<0mm,0mm>*{\bullet};<0mm,0mm>*{}**@{},
 <0mm,0mm>*{};<-8.6mm,-6mm>*{}**@{-},
<-8.6mm,-6mm>*{\bullet};<0mm,0mm>*{}**@{},
 <-8.6mm,-6mm>*{};<-12.6mm,-11mm>*{}**@{-},
 <-8.6mm,-6mm>*{};<-5.6mm,-11mm>*{}**@{-},
 <-8.6mm,-6mm>*{};<-10.6mm,-11mm>*{}**@{-},
  <-8.6mm,-6mm>*{};<-8mm,-10.5mm>*{...}**@{},
 <0mm,0mm>*{};<-9mm,-12.5mm>*{\underbrace{\ \ \ \ \ \ \
      }}**@{},
      <0mm,0mm>*{};<-8mm,-15.6mm>*{^{J_1}}**@{},
 <0mm,0mm>*{};<2mm,-6.4mm>*{\underbrace{\ \ \ \ \ \ \ \ \ \
      }}**@{},
    <0mm,0mm>*{};<2.6mm,-9.5mm>*{^{J_2}}**@{},
 <0mm,0mm>*{};<-3.5mm,-5mm>*{}**@{-},
 <0mm,0mm>*{};<-0mm,-4.6mm>*{...}**@{},
 <0mm,0mm>*{};<3.6mm,-5mm>*{}**@{-},
 %<0mm,0mm>*{};<3.5mm,-5mm>*{}**@{-},
 %<0mm,0mm>*{};<5.4mm,-4.6mm>*{...}**@{},
 <0mm,0mm>*{};<7mm,-5mm>*{}**@{-},
   %<0mm,0mm>*{};<7.6mm,-6.4mm>*{_{x}}**@{},
 \end{xy}
 \\
\Eeqrn
Let $\Lie_\infty^\circlearrowright$ and $\Lie^\circlearrowright$ and wheeled extensions of
$\Lie_\infty$ and, respectively, $\Lie$ (see \S 3.5 for precise definitions).

\bip

\no
%%%%%%%%%%%%%%%%%%%%%%%%%
{\bf 4.1.1. Theorem}. {\em The operad, $\Lie$, of Lie algebras
is stably Koszul, i.e.\ $\sH(\Lie_\infty^\circlearrowright)= {\sf
Lie}^\circlearrowright$}.

\bip

\Proof  We shall show that the natural morphism of dg props,
$$
({\Lie_\infty^\circlearrowright}, \delta) \lon
({\Lie}^\circlearrowright, 0)
$$
is a quasi-isomorphism.
 Consider a surjection of graph complexes (cf.\ Sect.\ 3.4),
$$
u: ({\Lie_\infty^+}, \delta) \lon ({\Lie_\infty^\circlearrowright},
\delta)
$$
where ${\Lie_\infty^+}$ is the marked extension of
${\Lie_\infty^\circlearrowright}$, i.e.\ the one in which one cyclic
edge in every wheel is marked. This surjection respects the
filtrations,
$$
 F_p{\Lie_\infty^+}:= \left\{\mbox{span}\langle f\in {\Lie_\infty^+}\rangle:\
 \mbox{total number of cyclic vertices in}\ f
   \geq p\right\},
 $$
$$
 F_p{\Lie_\infty^\circlearrowright}:= \left\{\mbox{span}\langle f\in \Lie_\infty^\circlearrowright\rangle:\
 \mbox{total number of cyclic vertices in}\ f
   \geq p\right\},
 $$
and hence induces a morphism of the associated 0-th terms of the spectral sequences,
$$
u_0: ({\mathsf E_0^+}, \p_0) \lon ({\mathsf E_0^\circlearrowright},
\p_0).
$$
The point is that the (pro-)cyclic group acting on  $({\mathsf
E_0^+}, \p_0)$ by shifting the marked edge one step further
along orientation commutes with the differential $\p_0$ so that
$u_0$ is nothing but the projection to the coinvariants with respect
to this action. As we work over a field of characteristic 0
coinvariants can be identified with invariants in $({\mathsf E_0^+},
\p_0)$. Hence we get, by Maschke theorem,
$$
\sH({\mathsf E_0^\circlearrowright}, \p_0) = \mbox{cyclic invariants
in}\ \sH({\mathsf E_0^+}, \p_0).
$$
%(In fact we shall see below that this spectral sequence degenerates and
%$H({\mathsf E_0^\circlearrowright}, \delta_0)=H({\mathsf L_\infty^\circlearrowright}, \delta)$.)

\bip

\no
The next step is to compute the cohomology of the complex $({\mathsf
E_0^+}, \p_0)$. Consider its  filtration,
$$
 \cF_p:= \left\{\mbox{span}\langle f\in {\mathsf E^+_0}\rangle:\
 \Ba{l}
 \mbox{total number of non-cyclic input
  }\\
  \mbox{edges at cyclic vertices in}\ f\Ea
   \leq p\right\},
 $$
and let $\{\cE_r, \delta_r\}_{r\geq 0}$ be the associated spectral
sequence. We shall show below that the latter degenerates at the
second term (so that $\cE_2\simeq \sH({\mathsf E_0^+}, \p_0)$). The
differential $\delta_0$ in $\cE_0$ is given by its values on the
vertices as follows: \Bi
\item[(i)]  on every non-cyclic vertex one has $\delta_0=\delta$, the differential in
$\Lie_\infty$;
\item[(ii)] on every cyclic vertex  $\delta_0=0$.
\Ei
Hence the complex
 $(\cE_0, \delta_0)$ is isomorphic to the {direct} sum of tensor products of
complexes $( {\Lie_\infty} ,\delta)$. By
 K\"{u}nneth theorem, we get,
$$
\cE_1= V_1/h(V_2),
$$
where
\Bi
\item[-] $V_1$ is the subspace of ${\Lie_\infty^+}$
consisting of all those graphs whose every non-cyclic vertex is
$
\begin{xy}
 <0mm,0.66mm>*{};<0mm,3mm>*{}**@{-},
 <0.39mm,-0.39mm>*{};<2.2mm,-2.2mm>*{}**@{-},
 <-0.35mm,-0.35mm>*{};<-2.2mm,-2.2mm>*{}**@{-},
 <0mm,0mm>*{\bullet};<0mm,0mm>*{}**@{},
   %<0mm,0.66mm>*{};<0mm,3.4mm>*{^1}**@{},
   %<0.39mm,-0.39mm>*{};<2.9mm,-4mm>*{^2}**@{},
   %<-0.35mm,-0.35mm>*{};<-2.8mm,-4mm>*{^1}**@{},
\end{xy} ;
$
\item[-] $V_2$ is the subspace of  ${\Lie_\infty^+}$
 whose
every non-cyclic vertex is either
$
\begin{xy}
 <0mm,0.66mm>*{};<0mm,3mm>*{}**@{-},
 <0.39mm,-0.39mm>*{};<2.2mm,-2.2mm>*{}**@{-},
 <-0.35mm,-0.35mm>*{};<-2.2mm,-2.2mm>*{}**@{-},
 <0mm,0mm>*{\bullet};<0mm,0mm>*{}**@{},
   %<0mm,0.66mm>*{};<0mm,3.4mm>*{^1}**@{},
   %<0.39mm,-0.39mm>*{};<2.9mm,-4mm>*{^2}**@{},
   %<-0.35mm,-0.35mm>*{};<-2.8mm,-4mm>*{^1}**@{},
\end{xy}
$
 or
$
\begin{xy}
 <0mm,0.66mm>*{};<0mm,3mm>*{}**@{-},
 <0.39mm,-0.39mm>*{};<2.2mm,-2.2mm>*{}**@{-},
 <-0.35mm,-0.35mm>*{};<-2.2mm,-2.2mm>*{}**@{-},
 <0mm,0mm>*{\bullet};<0mm,0mm>*{}**@{},
 <0mm,-0.66mm>*{};<0mm,-2.2mm>*{}**@{-},
   %<0mm,0.66mm>*{};<0mm,3.4mm>*{^1}**@{},
   %<0.39mm,-0.39mm>*{};<2.9mm,-4mm>*{^2}**@{},
   %<-0.35mm,-0.35mm>*{};<-2.8mm,-4mm>*{^1}**@{},
\end{xy}
$
with the number of vertices of the latter type $\geq 1$
;
\item[-] the map $h:V_2\rar V_1$ is given on non-cyclic vertices by
$$
h\left(
\begin{xy}
 <0mm,0.66mm>*{};<0mm,3mm>*{}**@{-},
 <0.39mm,-0.39mm>*{};<2.2mm,-2.2mm>*{}**@{-},
 <-0.35mm,-0.35mm>*{};<-2.2mm,-2.2mm>*{}**@{-},
 <0mm,0mm>*{\bullet};<0mm,0mm>*{}**@{},
\end{xy}
\right) =
\begin{xy}
 <0mm,0.66mm>*{};<0mm,3mm>*{}**@{-},
 <0.39mm,-0.39mm>*{};<2.2mm,-2.2mm>*{}**@{-},
 <-0.35mm,-0.35mm>*{};<-2.2mm,-2.2mm>*{}**@{-},
 <0mm,0mm>*{\bullet};<0mm,0mm>*{}**@{},
   %<0mm,0.66mm>*{};<0mm,3.4mm>*{^1}**@{},
   %<0.39mm,-0.39mm>*{};<2.9mm,-4mm>*{^2}**@{},
   %<-0.35mm,-0.35mm>*{};<-2.8mm,-4mm>*{^1}**@{},
\end{xy} \ \
, \ \ \ \ \ \ \ \
h\left(
\begin{xy}
 <0mm,0.66mm>*{};<0mm,3mm>*{}**@{-},
 <0.39mm,-0.39mm>*{};<2.2mm,-2.2mm>*{}**@{-},
 <-0.35mm,-0.35mm>*{};<-2.2mm,-2.2mm>*{}**@{-},
 <0mm,0mm>*{\bullet};<0mm,0mm>*{}**@{},
 <0mm,-0.66mm>*{};<0mm,-2.2mm>*{}**@{-},
 <0mm,0.66mm>*{};<0mm,-4mm>*{^2}**@{},
   <0.39mm,-0.39mm>*{};<2.9mm,-4mm>*{^3}**@{},
   <-0.35mm,-0.35mm>*{};<-2.8mm,-4mm>*{^1}**@{},
\end{xy}
\right) =
 \begin{xy}
 <0mm,0mm>*{\bullet};<0mm,0mm>*{}**@{},
 <0mm,0.69mm>*{};<0mm,3.0mm>*{}**@{-},
 <0.39mm,-0.39mm>*{};<2.4mm,-2.4mm>*{}**@{-},
 <-0.35mm,-0.35mm>*{};<-1.9mm,-1.9mm>*{}**@{-},
 <-2.4mm,-2.4mm>*{\bullet};<-2.4mm,-2.4mm>*{}**@{},
 <-2.0mm,-2.8mm>*{};<0mm,-4.9mm>*{}**@{-},
 <-2.8mm,-2.9mm>*{};<-4.7mm,-4.9mm>*{}**@{-},
    <0.39mm,-0.39mm>*{};<3.3mm,-4.0mm>*{^3}**@{},
    <-2.0mm,-2.8mm>*{};<0.5mm,-6.7mm>*{^2}**@{},
    <-2.8mm,-2.9mm>*{};<-5.2mm,-6.7mm>*{^1}**@{},
 \end{xy}
\ + \
 \begin{xy}
 <0mm,0mm>*{\bullet};<0mm,0mm>*{}**@{},
 <0mm,0.69mm>*{};<0mm,3.0mm>*{}**@{-},
 <0.39mm,-0.39mm>*{};<2.4mm,-2.4mm>*{}**@{-},
 <-0.35mm,-0.35mm>*{};<-1.9mm,-1.9mm>*{}**@{-},
 <-2.4mm,-2.4mm>*{\bullet};<-2.4mm,-2.4mm>*{}**@{},
 <-2.0mm,-2.8mm>*{};<0mm,-4.9mm>*{}**@{-},
 <-2.8mm,-2.9mm>*{};<-4.7mm,-4.9mm>*{}**@{-},
    <0.39mm,-0.39mm>*{};<3.3mm,-4.0mm>*{^2}**@{},
    <-2.0mm,-2.8mm>*{};<0.5mm,-6.7mm>*{^1}**@{},
    <-2.8mm,-2.9mm>*{};<-5.2mm,-6.7mm>*{^3}**@{},
 \end{xy}
\ + \
 \begin{xy}
 <0mm,0mm>*{\bullet};<0mm,0mm>*{}**@{},
 <0mm,0.69mm>*{};<0mm,3.0mm>*{}**@{-},
 <0.39mm,-0.39mm>*{};<2.4mm,-2.4mm>*{}**@{-},
 <-0.35mm,-0.35mm>*{};<-1.9mm,-1.9mm>*{}**@{-},
 <-2.4mm,-2.4mm>*{\bullet};<-2.4mm,-2.4mm>*{}**@{},
 <-2.0mm,-2.8mm>*{};<0mm,-4.9mm>*{}**@{-},
 <-2.8mm,-2.9mm>*{};<-4.7mm,-4.9mm>*{}**@{-},
    <0.39mm,-0.39mm>*{};<3.3mm,-4.0mm>*{^1}**@{},
    <-2.0mm,-2.8mm>*{};<0.5mm,-6.7mm>*{^3}**@{},
    <-2.8mm,-2.9mm>*{};<-5.2mm,-6.7mm>*{^2}**@{},
 \end{xy}
 $$
and  on all  cyclic vertices $h$ is set to be the identity.
\Ei
 The differential $\delta_1$ in $\cE_1$ is given by its values on
vertices as follows: \Bi
\item[(i)]  on every non-cyclic vertex one has $\delta_1=0$;
\item[(ii)] on every cyclic $(1,n+1)$-vertex with cyclic half-edges
denoted by $x$ and $y$, one has
\Beqrn
\delta_1\ \begin{xy}
 <0mm,0mm>*{\bullet};<0mm,0mm>*{}**@{},
 <0mm,0mm>*{};<0mm,5mm>*{}**@{-},
   <0mm,0mm>*{};<0mm,5.5mm>*{^{y}}**@{},
<0mm,0mm>*{\bullet};<0mm,0mm>*{}**@{},
 <0mm,0mm>*{};<-6mm,-5mm>*{}**@{-},
 <0mm,0mm>*{};<-3.1mm,-5mm>*{}**@{-},
 <0mm,0mm>*{};<0mm,-4.6mm>*{...}**@{},
 %<0mm,0mm>*{};<1.0mm,-5mm>*{}**@{-},
 <0mm,0mm>*{};<3.1mm,-5mm>*{}**@{-},
 <0mm,0mm>*{};<6mm,-5mm>*{}**@{-},
   <0mm,0mm>*{};<-6.7mm,-6.4mm>*{_1}**@{},
   <0mm,0mm>*{};<-3.2mm,-6.4mm>*{_2}**@{},
   <0mm,0mm>*{};<3.5mm,-6.4mm>*{_{n}}**@{},
   <0mm,0mm>*{};<6.9mm,-6.4mm>*{_{x}}**@{},
 \end{xy}
 &=&
 %
 %                                2nd TERM
 %
 \sum_{[n]=J_1\sqcup J_2\atop {\atop
 {|J_1|=2, |J_2|\geq 0}}
}\hspace{0mm}
\ \
\begin{xy}
 <0mm,0mm>*{\bullet};<0mm,0mm>*{}**@{},
 <0mm,0mm>*{};<0mm,5mm>*{}**@{-},
   <0mm,0mm>*{};<0mm,5.5mm>*{^{y}}**@{},
<0mm,0mm>*{\bullet};<0mm,0mm>*{}**@{},
 <0mm,0mm>*{};<-8.6mm,-6mm>*{}**@{-},
<-8.6mm,-6mm>*{\bullet};<0mm,0mm>*{}**@{},
 <-8.6mm,-6mm>*{};<-10.6mm,-9mm>*{}**@{-},
 <-8.6mm,-6mm>*{};<-6.6mm,-9mm>*{}**@{-},
 <0mm,0mm>*{};<-8.6mm,-10mm>*{\underbrace{
      }}**@{},
      <0mm,0mm>*{};<-8.6mm,-13.6mm>*{^{J_1}}**@{},
 <0mm,0mm>*{};<0mm,-6.4mm>*{\underbrace{\ \ \ \ \ \ \
      }}**@{},
    <0mm,0mm>*{};<0mm,-9.5mm>*{^{J_2}}**@{},
 <0mm,0mm>*{};<-3.5mm,-5mm>*{}**@{-},
 <0mm,0mm>*{};<-0mm,-4.6mm>*{...}**@{},
 <0mm,0mm>*{};<3.6mm,-5mm>*{}**@{-},
 <0mm,0mm>*{};<7mm,-5mm>*{}**@{-},
   <0mm,0mm>*{};<7.6mm,-6.4mm>*{_{x}}**@{},
 \end{xy} .
 \\
\Eeqrn
\Ei
 To compute the cohomology of $(\cE_1,\delta_1)$ let us step aside
and compute the cohomology of the minimal resolution,
$({\Lie_\infty}, \delta)$ (which we, of course, already know to be
equal to  ${\Lie}$), in a slightly unusual way:
$$
F^{\Lie}_p := \left\{\mbox{span}\langle f \in { \Lie_\infty}\rangle:
\mbox{number of edges attached to the root vertex of}\ f \leq
p\right\}
$$
is clearly a filtration of the  complex $({\Lie_\infty}, \delta)$.
Let  $\{E^{\Lie}_r,d_r^{\Lie}\}_{r\geq 0}$ be the associated
spectral sequence. The cohomology classes of
$E_1^{\Lie}=\sH(E^{\Lie}_0,d_0^{\Lie})$ resemble elements of
$\cE_1$: they are trees whose root vertex may have any number of
edges while all other vertices are binary, $
\begin{xy}
 <0mm,0.66mm>*{};<0mm,3mm>*{}**@{-},
 <0.39mm,-0.39mm>*{};<2.2mm,-2.2mm>*{}**@{-},
 <-0.35mm,-0.35mm>*{};<-2.2mm,-2.2mm>*{}**@{-},
 <0mm,0mm>*{\bullet};<0mm,0mm>*{}**@{},
   %<0mm,0.66mm>*{};<0mm,3.4mm>*{^1}**@{},
   %<0.39mm,-0.39mm>*{};<2.9mm,-4mm>*{^2}**@{},
   %<-0.35mm,-0.35mm>*{};<-2.8mm,-4mm>*{^1}**@{},
\end{xy}
$. The differential $d_1^{\Lie}$ is non-trivial only on the root
vertex on which it is given by,
\Beqrn d_1^{\Lie}\
\begin{xy}
 <0mm,0mm>*{\bullet};<0mm,0mm>*{}**@{},
 <0mm,0mm>*{};<0mm,5mm>*{}**@{-},
   %<0mm,0mm>*{};<0mm,5.5mm>*{^{y}}**@{},
   %
<0mm,0mm>*{\bullet};<0mm,0mm>*{}**@{},
 <0mm,0mm>*{};<-6mm,-5mm>*{}**@{-},
 <0mm,0mm>*{};<-3.1mm,-5mm>*{}**@{-},
 <0mm,0mm>*{};<0mm,-4.6mm>*{...}**@{},
 %<0mm,0mm>*{};<1.0mm,-5mm>*{}**@{-},
 <0mm,0mm>*{};<3.1mm,-5mm>*{}**@{-},
 <0mm,0mm>*{};<6mm,-5mm>*{}**@{-},
   <0mm,0mm>*{};<-6.7mm,-6.4mm>*{_1}**@{},
   <0mm,0mm>*{};<-3.2mm,-6.4mm>*{_2}**@{},
   <0mm,0mm>*{};<3.1mm,-6.4mm>*{_{n\mbox{-}1}}**@{},
   <0mm,0mm>*{};<6.9mm,-6.4mm>*{_{n}}**@{},
 \end{xy}
 &=&
 %
 %                                2nd TERM
 %
 \sum_{[n]=J_1\sqcup J_2\atop {\atop
 {|J_1|=2, |J_2|\geq 1}}
}\hspace{0mm}
\ \
\begin{xy}
 <0mm,0mm>*{\bullet};<0mm,0mm>*{}**@{},
 <0mm,0mm>*{};<0mm,5mm>*{}**@{-},
   <0mm,0mm>*{};<0mm,5.5mm>*{^{y}}**@{},
<0mm,0mm>*{\bullet};<0mm,0mm>*{}**@{},
 <0mm,0mm>*{};<-8.6mm,-6mm>*{}**@{-},
<-8.6mm,-6mm>*{\bullet};<0mm,0mm>*{}**@{},
 <-8.6mm,-6mm>*{};<-10.6mm,-9mm>*{}**@{-},
 <-8.6mm,-6mm>*{};<-6.6mm,-9mm>*{}**@{-},
 <0mm,0mm>*{};<-8.6mm,-10mm>*{\underbrace{
      }}**@{},
      <0mm,0mm>*{};<-8.6mm,-13.6mm>*{^{J_1}}**@{},
 <0mm,0mm>*{};<0mm,-6.4mm>*{\underbrace{\ \ \ \ \ \ \
      }}**@{},
    <0mm,0mm>*{};<0mm,-9.5mm>*{^{J_2}}**@{},
 <0mm,0mm>*{};<-3.5mm,-5mm>*{}**@{-},
 <0mm,0mm>*{};<-0mm,-4.6mm>*{...}**@{},
 <0mm,0mm>*{};<3.6mm,-5mm>*{}**@{-},
 %<0mm,0mm>*{};<7mm,-5mm>*{}**@{-},
  % <0mm,0mm>*{};<7.6mm,-6.4mm>*{_{x}}**@{},
 \end{xy} .
 \\
\Eeqrn
The cohomology of $(E_1^{\Lie}, d_1^{\Lie})$ is equal
to the operad of Lie algebras. The differential $d_1^{\Lie}$
is identical to
 the differential $\delta_1$ above
except for the term  corresponding to $|J_2|=0$. Thus let us
define another complex, $(E_1^{\Lie+}, d_1^{\Lie+})$, by adding
to $E_1^{\Lie}$ trees whose root vertex is a degree
$-1$ corolla
$\begin{xy}
 <0mm,-0.55mm>*{};<0mm,-2.5mm>*{}**@{-},
 <0mm,0mm>*{};<0mm,2.5mm>*{}**@{-},
 <0mm,0mm>*{\bullet};<0mm,0mm>*{}**@{},
 \end{xy}$ while all other vertices are binary
$
\begin{xy}
 <0mm,0.66mm>*{};<0mm,3mm>*{}**@{-},
 <0.39mm,-0.39mm>*{};<2.2mm,-2.2mm>*{}**@{-},
 <-0.35mm,-0.35mm>*{};<-2.2mm,-2.2mm>*{}**@{-},
 <0mm,0mm>*{\bullet};<0mm,0mm>*{}**@{},
   %<0mm,0.66mm>*{};<0mm,3.4mm>*{^1}**@{},
   %<0.39mm,-0.39mm>*{};<2.9mm,-4mm>*{^2}**@{},
   %<-0.35mm,-0.35mm>*{};<-2.8mm,-4mm>*{^1}**@{},
\end{xy}
$.
The differential  $d_1^{\Lie+}$ is defined on root
 $(1,n)$-corollas with $n\geq 2$ by formally the same formula as for
 $d^{\Lie}_1$ except that the summation range is extended to include
 the term with $|I_1|=0$. We also set
$d^{\Lie+}_1
\begin{xy}
 <0mm,-0.55mm>*{};<0mm,-2.5mm>*{}**@{-},
 <0mm,0mm>*{};<0mm,2.5mm>*{}**@{-},
 <0mm,0mm>*{\bullet};<0mm,0mm>*{}**@{},
 \end{xy}=0$.

\bip

\no
{\em Claim. The cohomology of the complex $(E_1^{\Lie +}, d_1^{\Lie +})$ is a
 one dimensional vector space spanned by $\begin{xy}
 <0mm,-0.55mm>*{};<0mm,-2.5mm>*{}**@{-},
 <0mm,0mm>*{};<0mm,2.5mm>*{}**@{-},
 <0mm,0mm>*{\bullet};<0mm,0mm>*{}**@{},
 \end{xy}$\, .}

\bip

\no
{\em Proof of the claim}. Consider the 2-step filtration, $F_0\subset F_1$ of the complex
$(E_1^{\Lie +}, d_1^{\Lie +})$ by the number of $\begin{xy}
 <0mm,-0.55mm>*{};<0mm,-2.5mm>*{}**@{-},
 <0mm,0mm>*{};<0mm,2.5mm>*{}**@{-},
 <0mm,0mm>*{\bullet};<0mm,0mm>*{}**@{},
 \end{xy}$ . The zero-th term of the associated spectral sequence is isomorphic
to the direct sum of the complexes,
$$
(E_1^{\Lie}, d_1^{\Lie})\oplus (E_1^{\Lie}[1], d_1^{\Lie})
\oplus (\mbox{span}\langle \begin{xy}
 <0mm,-0.55mm>*{};<0mm,-2.5mm>*{}**@{-},
 <0mm,0mm>*{};<0mm,2.5mm>*{}**@{-},
 <0mm,0mm>*{\bullet};<0mm,0mm>*{}**@{},
 \end{xy}\rangle, 0)
$$
so that the next term of the spectral sequence is
$$
{\Lie}\oplus {\Lie}[1] \oplus \langle \begin{xy}
 <0mm,-0.55mm>*{};<0mm,-2.5mm>*{}**@{-},
 <0mm,0mm>*{};<0mm,2.5mm>*{}**@{-},
 <0mm,0mm>*{\bullet};<0mm,0mm>*{}**@{},
 \end{xy}\rangle
$$
with the differential being zero on ${\Lie}[1] \oplus \langle \begin{xy}
 <0mm,-0.55mm>*{};<0mm,-2.5mm>*{}**@{-},
 <0mm,0mm>*{};<0mm,2.5mm>*{}**@{-},
 <0mm,0mm>*{\bullet};<0mm,0mm>*{}**@{},
 \end{xy}\rangle$ and the the natural isomorphism,
$$
{\Lie}\lon {\Lie}[1]
$$
on the remaining summand. Hence the claim follows.

 \bip

\no
The point of the above Claim is that the graph complex  $(\cE_1, \delta_1)$ is isomorphic
to the tensor product of a trivial complex with complexes of the form
 $(E_1^{\Lie +}, d_1^{\Lie +})$.
  Which immediately implies  that
  $\cE_2=\cE_\infty\simeq \sH(E_0^+, \p_0)$ is the direct sum of $\Lie$ and the vector space
spanned by marked wheels
  of the type,
$$
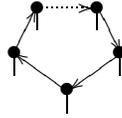
\begin{xy}
<0mm,0mm>*{\bullet};<0mm,0mm>*{}**@{},
<-7mm,5mm>*{\bullet};<0mm,0mm>*{}**@{},
<-4mm,11mm>*{\bullet};<0mm,0mm>*{}**@{},
<7mm,5mm>*{\bullet};<0mm,0mm>*{}**@{},
<4mm,11mm>*{\bullet};<0mm,0mm>*{}**@{},
 <0mm,0mm>*{};<0mm,-3mm>*{}**@{-},
 <-7mm,5mm>*{};<-7mm,2mm>*{}**@{-},
  <-4mm,11mm>*{};<-4mm,8mm>*{}**@{-},
   <7mm,5mm>*{};<7mm,2mm>*{}**@{-},
    <4mm,11mm>*{};<4mm,8mm>*{}**@{-},
{(0,0)*{}\ar@{->} (-7,5)*{}}
\hspace{-7mm}\mbox{\xy (-7,5)*{}\ar@{->} (-4,11)*{}\endxy}
\hspace{-2mm}\mbox{\xy (-4,11)*{}\ar@{.>} (4,11)*{}\endxy}
\hspace{-0.5mm}\mbox{\xy (4,11)*{}\ar@{->} (7,5)*{}\endxy}
\hspace{-8.3mm}\mbox{\xy (7,5)*{}\ar@{->} (0,0)*{}\endxy}
 \end{xy}
$$
whose every vertex is cyclic. Hence the cohomology group
$\sH({\sE_0^\circlearrowright},\p_0)=\sE_1^\circlearrowright$ we started with is equal to
the direct sum of  $\Lie$ and the space, $\mathsf Z$, spanned by unmarked wheels
  of the type,
$$
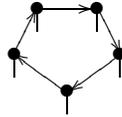
\begin{xy}
<0mm,0mm>*{\bullet};<0mm,0mm>*{}**@{},
<-7mm,5mm>*{\bullet};<0mm,0mm>*{}**@{},
<-4mm,11mm>*{\bullet};<0mm,0mm>*{}**@{},
<7mm,5mm>*{\bullet};<0mm,0mm>*{}**@{},
<4mm,11mm>*{\bullet};<0mm,0mm>*{}**@{},
 <0mm,0mm>*{};<0mm,-3mm>*{}**@{-},
 <-7mm,5mm>*{};<-7mm,2mm>*{}**@{-},
  <-4mm,11mm>*{};<-4mm,8mm>*{}**@{-},
   <7mm,5mm>*{};<7mm,2mm>*{}**@{-},
    <4mm,11mm>*{};<4mm,8mm>*{}**@{-},
{(0,0)*{}\ar@{->} (-7,5)*{}}
\hspace{-7mm}\mbox{\xy (-7,5)*{}\ar@{->} (-4,11)*{}\endxy}
\hspace{-2mm}\mbox{\xy (-4,11)*{}\ar@{->} (4,11)*{}\endxy}
\hspace{-0.5mm}\mbox{\xy (4,11)*{}\ar@{->} (7,5)*{}\endxy}
\hspace{-8.3mm}\mbox{\xy (7,5)*{}\ar@{->} (0,0)*{}\endxy}
 \end{xy}
$$
whose every vertex is cyclic. As every vertex is binary, the induced differential $\p_1$ on
$\sE_1^\circlearrowright$ vanishes, the spectral sequence by the number of cyclic vertices we began with
 degenerates,
and we conclude that this direct sum,  $\Lie \oplus Z$, is isomorphic to the required
cohomology group $\sH({\Lie_\infty^\circlearrowright},d)$.

\bip

\no
Finally one checks using Jacobi identities that every element of  ${\Lie}^\circlearrowright$
containing a wheel can be {\em uniquely}\, represented as a linear combination of graphs from
$\mathsf Z$ implying
$$
{\Lie}^\circlearrowright \simeq {\Lie \oplus {\mathsf Z}} \simeq
\sH({\Lie_\infty^\circlearrowright},d)
$$
and completing the proof.\hfill $\Box$

\bip

%%%%%%%%%%%%%%%%%%%%%%%%%

\no{\bf 4.2. Wheeeled prop of polyvector fields.}
 Let $\Liebi$ be the prop of Lie 1-bialgebras
and $(\Liebi_\infty,\delta)$ its minimal resolution (see \S 2.6.3).
We denote their wheeled extensions by  $\Liebi^\circlearrowright$
and $(\Liebi_\infty^\circlearrowright,\delta)$ respectively (see \S 3.5), and their operadic
wheeled extensions  by  $\Liebi^\looparrowright$
and $(\Liebi_\infty^\looparrowright,\delta)$ (see \S 3.11).
By Theorems 3.11.1 and 4.1.1, we have

\bip

\no{\bf 4.2.1. Proposition.} {\em  The natural epimorphism of dg props,
$$
(\Liebi_\infty^\looparrowright,\delta)\lon (\Liebi^\looparrowright,0)
$$
is a quasi-isomorphism.}
\bip

\noindent We shall study next a subcomplex ({\em not}\, a subprop!) of the complex
 $(\Liebi_\infty^\circlearrowright,\delta)$ which is spanned  by directed
graphs with at most one wheel, i.e.\ with at most one closed path which begins
and ends at the same vertex. We denote this subcomplex by $\Liebi_\infty^\circ$.
Similarly we define a subspace $\Liebi^\circ \subset\Liebi^\circlearrowright$  spanned
by equivalence classes of graphs with at most one wheel.

\bip

%%%%%%%%%%%%%%%%%%%%%%%%%
\no{\bf 4.2.2. Theorem}. $\sH(\Liebi^\circ_\infty,\, \delta)= \Liebi^\circ$.

\bip

\Proof
 {\bf (a)} Consider a two step filtration, $F_0\subset F_1:={\Liebi}^\circ_\infty $
of the complex $({\Liebi}^\circ_\infty, \delta)$,
with $F_0:= {\Liebi}_\infty$ being the subspace spanned by graphs with no wheels.
We shall show below that the cohomology of the associated direct sum complex,
$$
F_0\bigoplus \frac{F_1}{F_0},
$$
is equal to $\Liebi \bigoplus \frac{\Liebi^\circ}{\Liebi}$. In fact, the equality
$\sH(F_0)=\Liebi$ is obvious so that it is enough to show below that the cohomology
of the complex,
$\sC:= \frac{F_1}{F_0}$, is equal to $\frac{\Liebi^\circ}{\Liebi}$.

\bip

 \no  {\bf (b)} Consider a filtration of the complex $(\sC, \delta)$,
$$
 \sF_p\sC:= \left\{\mbox{span}\langle f\in \sC \rangle:\ \mbox{number of cyclic vertices in}\ f
 \geq p\right\},
 $$
and a similar filtration,
$$
 \sF_p\sC^+:= \left\{\mbox{span}\langle f\in \sC^+ \rangle:\ \mbox{number of cyclic vertices in}\ f
 \geq p\right\},
 $$
of the marked version of $\sC$. Let $\{\sE^\circlearrowright_r, \p_r\}_{r\geq 0}$ and $\{\sE_r^+, \p_r\}_{r\geq 0}$
be the associated spectral sequences. There is a natural surjection of {\em  complexes},
$$
u_0: ({\mathsf E_0^+}, \p_0) \lon ({\mathsf E_0^\circlearrowright},
\p_0).
$$
It is easy to see that
the differential $\p_0$ in $\sE_0^+$ commutes with the action of the (pro-)cyclic group on  $({\mathsf
E_0^+}, \delta_0)$ by shifting the marked edge one step further
along orientation so that
$u_0$ is nothing but the projection to the coinvariants with respect
to this action. As we work over a field of characteristic 0,  we get by Maschke theorem,
$$
\sH({\mathsf E_0^\circlearrowright}, \p_0) = \mbox{cyclic invariants
in}\ \sH({\mathsf E_0^+}, \p_0).
$$
so that we can work from now on with the complex $({\mathsf
E_0^+}, \p_0)$. Consider a filtration of the latter,
$$
 \cF_p:= \left\{\mbox{span}\langle f\in {\mathsf E^+_0}\rangle:\
 \Ba{l}
 \mbox{total number of non-cyclic input
  }\\
  \mbox{edges at cyclic vertices in}\ f\Ea
   \leq p\right\},
 $$
and let $\{\cE_r, d_r\}_{r\geq 0}$ be the associated spectral
sequence. The
differential $\delta_0$ in $\cE_0$ is given by its values on the
vertices as follows: \Bi
\item[(i)]  on every non-cyclic vertex one has $d_0=\delta_\lambda$, the differential in
$\Liebi_\infty$;
\item[(ii)] on every cyclic vertex  $d_0=0$.
\Ei
Hence the complex
 $(\cE_0, d_0)$ is isomorphic to the {direct} sum of tensor products of
complexes $( {\Lie_\infty} ,\delta)$ with trivial complexes. By
 K\"{u}nneth theorem, we conclude that $\cE_1=\sH(\cE_0, d_0)$ can be identified with the quotient of the subspace in $\sC$
spanned by graphs whose
every non-cyclic vertex is ternary, e.g.\ either
$\begin{xy}
 <0mm,-0.55mm>*{};<0mm,-2.5mm>*{}**@{-},
 <0.5mm,0.5mm>*{};<2.2mm,2.2mm>*{}**@{-},
 <-0.48mm,0.48mm>*{};<-2.2mm,2.2mm>*{}**@{-},
 <0mm,0mm>*{\bullet};<0mm,0mm>*{}**@{},
   %<0mm,-0.55mm>*{};<0mm,-3.8mm>*{_1}**@{},
   %<0.5mm,0.5mm>*{};<2.7mm,2.8mm>*{^2}**@{},
   %<-0.48mm,0.48mm>*{};<-2.7mm,2.8mm>*{^1}**@{},
 \end{xy}
 $ or
$
\begin{xy}
 <0mm,0.66mm>*{};<0mm,3mm>*{}**@{-},
 <0.39mm,-0.39mm>*{};<2.2mm,-2.2mm>*{}**@{-},
 <-0.35mm,-0.35mm>*{};<-2.2mm,-2.2mm>*{}**@{-},
 <0mm,0mm>*{\bullet};<0mm,0mm>*{}**@{},
   %<0mm,0.66mm>*{};<0mm,3.4mm>*{^1}**@{},
   %<0.39mm,-0.39mm>*{};<2.9mm,-4mm>*{^2}**@{},
   %<-0.35mm,-0.35mm>*{};<-2.8mm,-4mm>*{^1}**@{},
\end{xy}
$\ ,
with respect to the equivalence relation generated by the following equations among non-cyclic vertices,

\Bi
\item[($\star$)]\ \ \ \ \ \
$
 \begin{xy}
 <0mm,0mm>*{\bullet};<0mm,0mm>*{}**@{},
 <0mm,-0.49mm>*{};<0mm,-3.0mm>*{}**@{-},
 <0.49mm,0.49mm>*{};<1.9mm,1.9mm>*{}**@{-},
 <-0.5mm,0.5mm>*{};<-1.9mm,1.9mm>*{}**@{-},
 <-2.3mm,2.3mm>*{\bullet};<-2.3mm,2.3mm>*{}**@{},
 <-1.8mm,2.8mm>*{};<0mm,4.9mm>*{}**@{-},
 <-2.8mm,2.9mm>*{};<-4.6mm,4.9mm>*{}**@{-},
   <0.49mm,0.49mm>*{};<2.7mm,2.3mm>*{^3}**@{},
   <-1.8mm,2.8mm>*{};<0.4mm,5.3mm>*{^2}**@{},
   <-2.8mm,2.9mm>*{};<-5.1mm,5.3mm>*{^1}**@{},
 \end{xy}
\ + \
\begin{xy}
 <0mm,0mm>*{\bullet};<0mm,0mm>*{}**@{},
 <0mm,-0.49mm>*{};<0mm,-3.0mm>*{}**@{-},
 <0.49mm,0.49mm>*{};<1.9mm,1.9mm>*{}**@{-},
 <-0.5mm,0.5mm>*{};<-1.9mm,1.9mm>*{}**@{-},
 <-2.3mm,2.3mm>*{\bullet};<-2.3mm,2.3mm>*{}**@{},
 <-1.8mm,2.8mm>*{};<0mm,4.9mm>*{}**@{-},
 <-2.8mm,2.9mm>*{};<-4.6mm,4.9mm>*{}**@{-},
   <0.49mm,0.49mm>*{};<2.7mm,2.3mm>*{^2}**@{},
   <-1.8mm,2.8mm>*{};<0.4mm,5.3mm>*{^1}**@{},
   <-2.8mm,2.9mm>*{};<-5.1mm,5.3mm>*{^3}**@{},
 \end{xy}
\ + \
\begin{xy}
 <0mm,0mm>*{\bullet};<0mm,0mm>*{}**@{},
 <0mm,-0.49mm>*{};<0mm,-3.0mm>*{}**@{-},
 <0.49mm,0.49mm>*{};<1.9mm,1.9mm>*{}**@{-},
 <-0.5mm,0.5mm>*{};<-1.9mm,1.9mm>*{}**@{-},
 <-2.3mm,2.3mm>*{\bullet};<-2.3mm,2.3mm>*{}**@{},
 <-1.8mm,2.8mm>*{};<0mm,4.9mm>*{}**@{-},
 <-2.8mm,2.9mm>*{};<-4.6mm,4.9mm>*{}**@{-},
   <0.49mm,0.49mm>*{};<2.7mm,2.3mm>*{^1}**@{},
   <-1.8mm,2.8mm>*{};<0.4mm,5.3mm>*{^3}**@{},
   <-2.8mm,2.9mm>*{};<-5.1mm,5.3mm>*{^2}**@{},
 \end{xy}
=0
\ , \ \ \ \
 \begin{xy}
 <0mm,0mm>*{\bullet};<0mm,0mm>*{}**@{},
 <0mm,0.69mm>*{};<0mm,3.0mm>*{}**@{-},
 <0.39mm,-0.39mm>*{};<2.4mm,-2.4mm>*{}**@{-},
 <-0.35mm,-0.35mm>*{};<-1.9mm,-1.9mm>*{}**@{-},
 <-2.4mm,-2.4mm>*{\bullet};<-2.4mm,-2.4mm>*{}**@{},
 <-2.0mm,-2.8mm>*{};<0mm,-4.9mm>*{}**@{-},
 <-2.8mm,-2.9mm>*{};<-4.7mm,-4.9mm>*{}**@{-},
    <0.39mm,-0.39mm>*{};<3.3mm,-4.0mm>*{^3}**@{},
    <-2.0mm,-2.8mm>*{};<0.5mm,-6.7mm>*{^2}**@{},
    <-2.8mm,-2.9mm>*{};<-5.2mm,-6.7mm>*{^1}**@{},
 \end{xy}
\ + \
 \begin{xy}
 <0mm,0mm>*{\bullet};<0mm,0mm>*{}**@{},
 <0mm,0.69mm>*{};<0mm,3.0mm>*{}**@{-},
 <0.39mm,-0.39mm>*{};<2.4mm,-2.4mm>*{}**@{-},
 <-0.35mm,-0.35mm>*{};<-1.9mm,-1.9mm>*{}**@{-},
 <-2.4mm,-2.4mm>*{\bullet};<-2.4mm,-2.4mm>*{}**@{},
 <-2.0mm,-2.8mm>*{};<0mm,-4.9mm>*{}**@{-},
 <-2.8mm,-2.9mm>*{};<-4.7mm,-4.9mm>*{}**@{-},
    <0.39mm,-0.39mm>*{};<3.3mm,-4.0mm>*{^2}**@{},
    <-2.0mm,-2.8mm>*{};<0.5mm,-6.7mm>*{^1}**@{},
    <-2.8mm,-2.9mm>*{};<-5.2mm,-6.7mm>*{^3}**@{},
 \end{xy}
\ + \
 \begin{xy}
 <0mm,0mm>*{\bullet};<0mm,0mm>*{}**@{},
 <0mm,0.69mm>*{};<0mm,3.0mm>*{}**@{-},
 <0.39mm,-0.39mm>*{};<2.4mm,-2.4mm>*{}**@{-},
 <-0.35mm,-0.35mm>*{};<-1.9mm,-1.9mm>*{}**@{-},
 <-2.4mm,-2.4mm>*{\bullet};<-2.4mm,-2.4mm>*{}**@{},
 <-2.0mm,-2.8mm>*{};<0mm,-4.9mm>*{}**@{-},
 <-2.8mm,-2.9mm>*{};<-4.7mm,-4.9mm>*{}**@{-},
    <0.39mm,-0.39mm>*{};<3.3mm,-4.0mm>*{^1}**@{},
    <-2.0mm,-2.8mm>*{};<0.5mm,-6.7mm>*{^3}**@{},
    <-2.8mm,-2.9mm>*{};<-5.2mm,-6.7mm>*{^2}**@{},
 \end{xy}
=0\ , \ \ \ \
\begin{xy}
 <0mm,2.47mm>*{};<0mm,-0.5mm>*{}**@{-},
 <0.5mm,3.5mm>*{};<2.2mm,5.2mm>*{}**@{-},
 <-0.48mm,3.48mm>*{};<-2.2mm,5.2mm>*{}**@{-},
 <0mm,3mm>*{\bullet};<0mm,3mm>*{}**@{},
  <0mm,-0.8mm>*{\bullet};<0mm,-0.8mm>*{}**@{},
<0mm,-0.8mm>*{};<-2.2mm,-3.5mm>*{}**@{-},
 <0mm,-0.8mm>*{};<2.2mm,-3.5mm>*{}**@{-},
     <0.5mm,3.5mm>*{};<2.8mm,5.7mm>*{^2}**@{},
     <-0.48mm,3.48mm>*{};<-2.8mm,5.7mm>*{^1}**@{},
   <0mm,-0.8mm>*{};<-2.7mm,-5.2mm>*{^1}**@{},
   <0mm,-0.8mm>*{};<2.7mm,-5.2mm>*{^2}**@{},
\end{xy} = 0\ ,
$
\Ei
The differential $d_1$ in $\cE_1$ is non-zero only on cyclic vertices,
\Beqrn
d_1\ \begin{xy}
 <0mm,0mm>*{\bullet};<0mm,0mm>*{}**@{},
 <0mm,0mm>*{};<-8mm,5mm>*{}**@{-},
 <0mm,0mm>*{};<-4.5mm,5mm>*{}**@{-},
 <0mm,0mm>*{};<0mm,4.6mm>*{...}**@{},
 <0mm,0mm>*{};<3.5mm,5mm>*{}**@{-},
 <0mm,0mm>*{};<8mm,5mm>*{}**@{.},
   <0mm,0mm>*{};<-8.5mm,5.5mm>*{^1}**@{},
   <0mm,0mm>*{};<-5mm,5.5mm>*{^2}**@{},
   <0mm,0mm>*{};<3.9mm,5.5mm>*{^{m}}**@{},
<0mm,0mm>*{\bullet};<0mm,0mm>*{}**@{},
 <0mm,0mm>*{};<-8mm,-5mm>*{}**@{-},
 <0mm,0mm>*{};<-4.5mm,-5mm>*{}**@{-},
 <0mm,0mm>*{};<0mm,-4.6mm>*{...}**@{},
 %<0mm,0mm>*{};<1.0mm,-5mm>*{}**@{-},
 <0mm,0mm>*{};<3.5mm,-5mm>*{}**@{-},
 <0mm,0mm>*{};<8mm,-5mm>*{}**@{.},
   <0mm,0mm>*{};<-8.5mm,-6.4mm>*{_1}**@{},
   <0mm,0mm>*{};<-5mm,-6.4mm>*{_2}**@{},
   <0mm,0mm>*{};<3.9mm,-6.4mm>*{_{n}}**@{},
 \end{xy}
 &=&
 \sum_{[m]=I_1\sqcup I_2\atop {\atop
 {|I_1|\geq 0, |I_2|= 2}}
}
(-1)^{\sigma(I_1\sqcup I_2)+1}\ \
\begin{xy}
 <0mm,0mm>*{\bullet};<0mm,0mm>*{}**@{},
 <0mm,0mm>*{};<-8mm,5mm>*{}**@{-},
 <0mm,0mm>*{};<-3mm,5mm>*{}**@{-},
 <0mm,0mm>*{};<-1mm,5mm>*{}**@{-},
 <0mm,0mm>*{};<-4.5mm,4.6mm>*{...}**@{},
 <0mm,0mm>*{};<3.5mm,6mm>*{}**@{-},
<3.5mm,6mm>*{\bullet};<0mm,0mm>*{}**@{},
<3.5mm,6mm>*{};<1mm,9.5mm>*{}**@{-},
%<3.5mm,6mm>*{};<1mm,11mm>*{}**@{-},
<3.5mm,6mm>*{};<6mm,9.5mm>*{}**@{-},
%<3.5mm,6mm>*{};<4mm,10mm>*{...}**@{},
 %
 <0mm,0mm>*{};<3.1mm,11mm>*{\overbrace{
      }}**@{},
      <0mm,0mm>*{};<3.1mm,13.6mm>*{^{I_2}}**@{},
 <0mm,0mm>*{};<-5.4mm,6.7mm>*{\overbrace{
      }}**@{},
    <0mm,0mm>*{};<-5.4mm,9mm>*{^{I_1}}**@{},
 %
 %<0mm,0mm>*{};<3.5mm,5mm>*{}**@{.},
 %<0mm,0mm>*{};<5.4mm,4.6mm>*{...}**@{},
 <0mm,0mm>*{};<8mm,5mm>*{}**@{.},
<0mm,0mm>*{\bullet};<0mm,0mm>*{}**@{},
 <0mm,0mm>*{};<-8mm,-5mm>*{}**@{-},
 <0mm,0mm>*{};<-4.5mm,-5mm>*{}**@{-},
 <0mm,0mm>*{};<0mm,-4.6mm>*{...}**@{},
 %<0mm,0mm>*{};<1.0mm,-5mm>*{}**@{-},
 <0mm,0mm>*{};<3.5mm,-5mm>*{}**@{-},
 <0mm,0mm>*{};<8mm,-5mm>*{}**@{.},
   <0mm,0mm>*{};<-8.5mm,-6.4mm>*{_1}**@{},
   <0mm,0mm>*{};<-5mm,-6.4mm>*{_2}**@{},
   <0mm,0mm>*{};<3.9mm,-6.4mm>*{_{n}}**@{},
 \end{xy}
  \ +\
 \sum_{[n]=J_1\sqcup J_2\atop {\atop
 {|J_1|=2, |J_2|= 2}}
}\hspace{0mm}
\ \
\begin{xy}
 <0mm,0mm>*{\bullet};<0mm,0mm>*{}**@{},
 <0mm,0mm>*{};<-8mm,-5mm>*{}**@{-},
 <0mm,0mm>*{};<-3mm,-5mm>*{}**@{-},
 <0mm,0mm>*{};<-1mm,-5mm>*{}**@{-},
 <0mm,0mm>*{};<-4.5mm,-4.6mm>*{...}**@{},
 <0mm,0mm>*{};<3.5mm,-6mm>*{}**@{-},
<3.5mm,-6mm>*{\bullet};<0mm,0mm>*{}**@{},
<3.5mm,-6mm>*{};<1mm,-9.5mm>*{}**@{-},
%<3.5mm,-6mm>*{};<1mm,-11mm>*{}**@{-},
<3.5mm,-6mm>*{};<6mm,-9.5mm>*{}**@{-},
%<3.5mm,-6mm>*{};<4mm,-10mm>*{...}**@{},
 %
 <0mm,0mm>*{};<3.1mm,-11mm>*{\underbrace{\ \
      }}**@{},
      <0mm,0mm>*{};<3.1mm,-14mm>*{_{J_2}}**@{},
 <0mm,0mm>*{};<-5.4mm,-6.7mm>*{\underbrace{
      }}**@{},
    <0mm,0mm>*{};<-5.4mm,-9.9mm>*{_{J_1}}**@{},
 %
 %<0mm,0mm>*{};<3.5mm,5mm>*{}**@{.},
 %<0mm,0mm>*{};<5.4mm,4.6mm>*{...}**@{},
 <0mm,0mm>*{};<8mm,-5mm>*{}**@{.},
<0mm,0mm>*{\bullet};<0mm,0mm>*{}**@{},
 <0mm,0mm>*{};<-8mm,5mm>*{}**@{-},
 <0mm,0mm>*{};<-4.5mm,5mm>*{}**@{-},
 <0mm,0mm>*{};<0mm,4.6mm>*{...}**@{},
 %<0mm,0mm>*{};<1.0mm,5mm>*{}**@{-},
 <0mm,0mm>*{};<3.5mm,5mm>*{}**@{-},
 <0mm,0mm>*{};<8mm,5mm>*{}**@{.},
   <0mm,0mm>*{};<-8.5mm,6.4mm>*{_1}**@{},
   <0mm,0mm>*{};<-5mm,6.4mm>*{_2}**@{},
   <0mm,0mm>*{};<3.9mm,6.4mm>*{_{n}}**@{},
 \end{xy}
 \\
&&\\
%
%                                  3rd TERM
%
  && + \
  \sum_{[m]=I_1\sqcup I_2\atop {[n]=J_1\sqcup J_2\atop
 {|I_1|\geq 0, |I_2|=1 \atop
 |J_1|\geq 0, |J_2|=1}}
}\hspace{0mm}  (-1)^{\sigma(I_1\sqcup I_2)}
\ \
\begin{xy}
 <0mm,0mm>*{\bullet};<0mm,0mm>*{}**@{},
 <0mm,0mm>*{};<-8mm,5mm>*{}**@{-},
 <0mm,0mm>*{};<-4.5mm,5mm>*{}**@{-},
 <0mm,0mm>*{};<-1.5mm,4.6mm>*{...}**@{},
 <0mm,0mm>*{};<1.0mm,5mm>*{}**@{-},
 <0mm,0mm>*{};<4.5mm,5mm>*{}**@{.},
<0mm,0mm>*{};<13.9mm,6mm>*{}**@{-},
<13.9mm,6mm>*{\bullet};<0mm,0mm>*{}**@{},
<13.9mm,6mm>*{};<13.9mm,9.5mm>*{}**@{-},
<13.9mm,6mm>*{};<15.9mm,3mm>*{}**@{-},
<13.9mm,6mm>*{};<13.9mm,11mm>*{^{I_2}}**@{},
<13.9mm,6mm>*{};<16.5mm,1mm>*{_{J_2}}**@{},
 <0mm,0mm>*{};<-3mm,6.4mm>*{\overbrace{\ \ \ \ \ \ \
      }}**@{},
      <0mm,0mm>*{};<-3mm,9mm>*{^{I_1}}**@{},
 <0mm,0mm>*{};<-3.0mm,-6.4mm>*{\underbrace{\ \ \ \ \ \ \
      }}**@{},
    <0mm,0mm>*{};<-3mm,-9.5mm>*{^{J_1}}**@{},
<0mm,0mm>*{\bullet};<0mm,0mm>*{}**@{},
 <0mm,0mm>*{};<-8mm,-5mm>*{}**@{-},
 <0mm,0mm>*{};<-4.5mm,-5mm>*{}**@{-},
 <0mm,0mm>*{};<-1.5mm,-4.6mm>*{...}**@{},
 <0mm,0mm>*{};<1.0mm,-5mm>*{}**@{-},
 <0mm,0mm>*{};<4.5mm,-5mm>*{}**@{.},
 \end{xy}
\\
&&\\
%
%                     4-th TERM
%
 %
  && +
  \
\sum_{[m]=I_1\sqcup I_2\atop {[n]=J_1\sqcup J_2\atop
 {|I_1|=1, |I_2|\geq 0 \atop
 |J_1|=1, |J_2|\geq 0}}
}\hspace{0mm}  (-1)^{\sigma(I_1\sqcup I_2) +m}
\ \
\begin{xy}
 <0mm,0mm>*{\bullet};<0mm,0mm>*{}**@{},
 <0mm,0mm>*{};<-4.5mm,5mm>*{}**@{-},
 <0mm,0mm>*{};<-1.5mm,4.6mm>*{...}**@{},
 <0mm,0mm>*{};<1.0mm,5mm>*{}**@{-},
 <0mm,0mm>*{};<5.5mm,5mm>*{}**@{.},
 %<0mm,0mm>*{};<5.4mm,4.6mm>*{...}**@{},
 %<0mm,0mm>*{};<8mm,5mm>*{}**@{.},
   %
<0mm,0mm>*{\bullet};<0mm,0mm>*{}**@{},
 <0mm,0mm>*{};<-8.1mm,-5.5mm>*{}**@{-},
<-8.6mm,-6mm>*{\circ};<0mm,0mm>*{}**@{},
 <-9mm,-5.5mm>*{};<-10.9mm,-2.7mm>*{}**@{-},
 <-8.6mm,-6.5mm>*{};<-8.6mm,-9mm>*{}**@{-},
 <-9mm,-5.5mm>*{};<-10.9mm,-1.6mm>*{^{I_1}}**@{},
 <-8.6mm,-6.5mm>*{};<-8.6mm,-11mm>*{_{J_1}}**@{},
 <0mm,0mm>*{};<-1.9mm,6.4mm>*{\overbrace{
      }}**@{},
      <0mm,0mm>*{};<-1.7mm,9.4mm>*{_{I_2}}**@{},
 <0mm,0mm>*{};<-1.9mm,-6.4mm>*{\underbrace{
      }}**@{},
    <0mm,0mm>*{};<-1.7mm,-9.5mm>*{^{J_2}}**@{},
 <0mm,0mm>*{};<-4.5mm,-5mm>*{}**@{-},
 <0mm,0mm>*{};<-1.5mm,-4.6mm>*{...}**@{},
 <0mm,0mm>*{};<1.0mm,-5mm>*{}**@{-},
 <0mm,0mm>*{};<4.5mm,-5mm>*{}**@{.},
 \end{xy}
 \\
\Eeqrn
where cyclic half-edges (here and below) are dashed. Then $\cE_1$ can be interpreted as a bicomplex,
$(\cE_1=\bigoplus_{m,n}\cE_1^{m,n}, d_1=\p + \bar{\p})$, with,  say, $m$ counting the number
of vertices attached to cyclic vertices in ``operadic" way (as in the first two summands above), and $n$ counting
the number of vertices attached to cyclic vertices in ``non-operadic" way (as in the last two summands in the above
formula). Note that the assumption that there is only {\em one}\, wheel in $\sC$ is vital for this splitting of the differential
$d_1$ to have sense. The differential $\p$ (respectively, $\bar{\p}$) is equal to the first
(respectively, last) two summands in $d_1$.

\bip

\no
Using Claim in the proof of Theorem~4.1.1 it is not hard check that
$\sH(\cE_1, \p)$ is isomorphic to the quotient of the subspace of $\sC$ spanned by graphs whose

---
every non-cyclic vertex is ternary, e.g.\  either
$\begin{xy}
 <0mm,-0.55mm>*{};<0mm,-2.5mm>*{}**@{-},
 <0.5mm,0.5mm>*{};<2.2mm,2.2mm>*{}**@{-},
 <-0.48mm,0.48mm>*{};<-2.2mm,2.2mm>*{}**@{-},
 <0mm,0mm>*{\bullet};<0mm,0mm>*{}**@{},
   %<0mm,-0.55mm>*{};<0mm,-3.8mm>*{_1}**@{},
   %<0.5mm,0.5mm>*{};<2.7mm,2.8mm>*{^2}**@{},
   %<-0.48mm,0.48mm>*{};<-2.7mm,2.8mm>*{^1}**@{},
 \end{xy}
 $ or
$
\begin{xy}
 <0mm,0.66mm>*{};<0mm,3mm>*{}**@{-},
 <0.39mm,-0.39mm>*{};<2.2mm,-2.2mm>*{}**@{-},
 <-0.35mm,-0.35mm>*{};<-2.2mm,-2.2mm>*{}**@{-},
 <0mm,0mm>*{\bullet};<0mm,0mm>*{}**@{},
   %<0mm,0.66mm>*{};<0mm,3.4mm>*{^1}**@{},
   %<0.39mm,-0.39mm>*{};<2.9mm,-4mm>*{^2}**@{},
   %<-0.35mm,-0.35mm>*{};<-2.8mm,-4mm>*{^1}**@{},
\end{xy}
$\ ;

--- every cyclic vertex is either
$\begin{xy}
 <0mm,0mm>*{};<0mm,-2.5mm>*{}**@{.},
 <0.5mm,0.5mm>*{};<2.2mm,2.2mm>*{}**@{.},
 <-0.48mm,0.48mm>*{};<-2.2mm,2.2mm>*{}**@{-},
 <0mm,0mm>*{\bullet};<0mm,0mm>*{}**@{},
 \end{xy}
 $, or
$
\begin{xy}
 <0mm,0mm>*{};<0mm,3mm>*{}**@{.},
 <0.39mm,-0.39mm>*{};<2.2mm,-2.2mm>*{}**@{.},
 <-0.35mm,-0.35mm>*{};<-2.2mm,-2.2mm>*{}**@{-},
 <0mm,0mm>*{\bullet};<0mm,0mm>*{}**@{},
\end{xy}
$\ ,
or
$
\begin{xy}
 <0mm,0mm>*{};<2.2mm,2.2mm>*{}**@{.},
 <0mm,0mm>*{};<-2.2mm,2.2mm>*{}**@{-},
 <0mm,0mm>*{};<2.2mm,-2.2mm>*{}**@{.},
 <0mm,0mm>*{};<-2.2mm,-2.2mm>*{}**@{-},
 <0mm,0mm>*{\bullet};<0mm,0mm>*{}**@{},
\end{xy}
$\ ,

\noindent with respect to the equivalence relation generated by equations  $(\star)$ and the following ones,
\Bi
\item[($\star\star$)] \ \ \ \ \ \
$
 \begin{xy}
 <0mm,0mm>*{\bullet};<0mm,0mm>*{}**@{},
 <0mm,0mm>*{};<0mm,-3.0mm>*{}**@{.},
 <0mm,0mm>*{};<2mm,2mm>*{}**@{.},
 <-0.5mm,0.5mm>*{};<-1.9mm,1.9mm>*{}**@{-},
 <-2.3mm,2.3mm>*{\bullet};<-2.3mm,2.3mm>*{}**@{},
 <-1.8mm,2.8mm>*{};<0mm,4.9mm>*{}**@{-},
 <-2.8mm,2.9mm>*{};<-4.6mm,4.9mm>*{}**@{-},
 \end{xy}=0\ , \ \ \
 \begin{xy}
 <0mm,0mm>*{\bullet};<0mm,0mm>*{}**@{},
 <0mm,0mm>*{};<0mm,3.0mm>*{}**@{.},
 <0mm,0mm>*{};<2.4mm,-2.4mm>*{}**@{.},
 <0mm,0mm>*{};<-2mm,-2mm>*{}**@{-},
 <-2.4mm,-2.4mm>*{\bullet};<-2.4mm,-2.4mm>*{}**@{},
 <-2.0mm,-2.8mm>*{};<0mm,-4.9mm>*{}**@{-},
 <-2.8mm,-2.9mm>*{};<-4.7mm,-4.9mm>*{}**@{-},
 \end{xy} =0\ , \ \ \
 \begin{xy}
 <0mm,0mm>*{\bullet};<0mm,0mm>*{}**@{},
 <0mm,0mm>*{};<2mm,-2.0mm>*{}**@{.},
<0mm,0mm>*{};<-2mm,-2.0mm>*{}**@{-},
 <0mm,0mm>*{};<2mm,2mm>*{}**@{.},
 <-0.5mm,0.5mm>*{};<-1.9mm,1.9mm>*{}**@{-},
 <-2.3mm,2.3mm>*{\bullet};<-2.3mm,2.3mm>*{}**@{},
 <-1.8mm,2.8mm>*{};<0mm,4.9mm>*{}**@{-},
 <-2.8mm,2.9mm>*{};<-4.6mm,4.9mm>*{}**@{-},
 \end{xy}=0\ , \ \ \
 \begin{xy}
 <0mm,0mm>*{\bullet};<0mm,0mm>*{}**@{},
 <0mm,0mm>*{};<2mm,2.0mm>*{}**@{.},
<0mm,0mm>*{};<-2mm,2.0mm>*{}**@{-},
 <0mm,0mm>*{};<2.4mm,-2.4mm>*{}**@{.},
 <0mm,0mm>*{};<-2mm,-2mm>*{}**@{-},
 <-2.4mm,-2.4mm>*{\bullet};<-2.4mm,-2.4mm>*{}**@{},
 <-2.0mm,-2.8mm>*{};<0mm,-4.9mm>*{}**@{-},
 <-2.8mm,-2.9mm>*{};<-4.7mm,-4.9mm>*{}**@{-},
 \end{xy} =0\ .
$
\Ei
The differential $\bar{\p}$ is non-zero only on cyclic vertices of the type
$
\begin{xy}
 <0mm,0mm>*{};<2.2mm,2.2mm>*{}**@{.},
 <0mm,0mm>*{};<-2.2mm,2.2mm>*{}**@{-},
 <0mm,0mm>*{};<2.2mm,-2.2mm>*{}**@{.},
 <0mm,0mm>*{};<-2.2mm,-2.2mm>*{}**@{-},
 <0mm,0mm>*{\bullet};<0mm,0mm>*{}**@{},
\end{xy}
$\ ,
on which it is given by
$$
\bar{\p} \,
\begin{xy}
 <0mm,0mm>*{};<2.2mm,2.2mm>*{}**@{.},
 <0mm,0mm>*{};<-2.2mm,2.2mm>*{}**@{-},
 <0mm,0mm>*{};<2.2mm,-2.2mm>*{}**@{.},
 <0mm,0mm>*{};<-2.2mm,-2.2mm>*{}**@{-},
 <0mm,0mm>*{\bullet};<0mm,0mm>*{}**@{},
\end{xy}\
=\ - \
\begin{xy}
 <0mm,0mm>*{\bullet};<0mm,0mm>*{}**@{},
 <0mm,0mm>*{};<0mm,-3.0mm>*{}**@{.},
 <0mm,0mm>*{};<2mm,2mm>*{}**@{.},
 <-0.5mm,0.5mm>*{};<-1.9mm,1.9mm>*{}**@{-},
 <-2.3mm,2.3mm>*{\bullet};
 <-2.3mm,2.3mm>*{};<-2.3mm,5mm>*{}**@{-},
 <-2.3mm,2.3mm>*{};<-4.6mm,0mm>*{}**@{-},
 \end{xy}
\ \
-
 \begin{xy}
 <0mm,0mm>*{\bullet};<0mm,0mm>*{}**@{},
 <0mm,0mm>*{};<0mm,3.0mm>*{}**@{.},
 <0mm,0mm>*{};<2.4mm,-2.4mm>*{}**@{.},
 <0mm,0mm>*{};<-2mm,-2mm>*{}**@{-},
 <-2.4mm,-2.4mm>*{\bullet};<-2.4mm,-2.4mm>*{}**@{},
 <-2.4mm,-2.4mm>*{};<-2.4mm,-5mm>*{}**@{-},
 <-2.4mm,-2.4mm>*{};<-4.9mm,0mm>*{}**@{-},
 \end{xy}\ \
$$
Hence $\sA:=\sH(\sH(\cE_1,
\p), \bar{\p})$ can be identified with
 the quotient of the subspace of $\sC$ spanned by graphs whose

---
every non-cyclic vertex is ternary, e.g.\ either
$\begin{xy}
 <0mm,-0.55mm>*{};<0mm,-2.5mm>*{}**@{-},
 <0.5mm,0.5mm>*{};<2.2mm,2.2mm>*{}**@{-},
 <-0.48mm,0.48mm>*{};<-2.2mm,2.2mm>*{}**@{-},
 <0mm,0mm>*{\bullet};<0mm,0mm>*{}**@{},
   %<0mm,-0.55mm>*{};<0mm,-3.8mm>*{_1}**@{},
   %<0.5mm,0.5mm>*{};<2.7mm,2.8mm>*{^2}**@{},
   %<-0.48mm,0.48mm>*{};<-2.7mm,2.8mm>*{^1}**@{},
 \end{xy}
 $ or
$
\begin{xy}
 <0mm,0.66mm>*{};<0mm,3mm>*{}**@{-},
 <0.39mm,-0.39mm>*{};<2.2mm,-2.2mm>*{}**@{-},
 <-0.35mm,-0.35mm>*{};<-2.2mm,-2.2mm>*{}**@{-},
 <0mm,0mm>*{\bullet};<0mm,0mm>*{}**@{},
   %<0mm,0.66mm>*{};<0mm,3.4mm>*{^1}**@{},
   %<0.39mm,-0.39mm>*{};<2.9mm,-4mm>*{^2}**@{},
   %<-0.35mm,-0.35mm>*{};<-2.8mm,-4mm>*{^1}**@{},
\end{xy}
$\ ;

--- every cyclic vertex is also ternary, e.g.\ either
$\begin{xy}
 <0mm,0mm>*{};<0mm,-2.5mm>*{}**@{.},
 <0.5mm,0.5mm>*{};<2.2mm,2.2mm>*{}**@{.},
 <-0.48mm,0.48mm>*{};<-2.2mm,2.2mm>*{}**@{-},
 <0mm,0mm>*{\bullet};<0mm,0mm>*{}**@{},
 \end{xy}
 $ or
$
\begin{xy}
 <0mm,0mm>*{};<0mm,3mm>*{}**@{.},
 <0.39mm,-0.39mm>*{};<2.2mm,-2.2mm>*{}**@{.},
 <-0.35mm,-0.35mm>*{};<-2.2mm,-2.2mm>*{}**@{-},
 <0mm,0mm>*{\bullet};<0mm,0mm>*{}**@{},
\end{xy}
$

\noindent with respect to the equivalence relation generated by equations  $(\star)$, $(\star\star)$
 and, say, the following one,
$$
 \begin{xy}
 <0mm,0mm>*{\bullet};<0mm,0mm>*{}**@{},
 <0mm,0mm>*{};<0mm,-3.0mm>*{}**@{.},
 <0mm,0mm>*{};<2mm,2mm>*{}**@{.},
 <-0.5mm,0.5mm>*{};<-1.9mm,1.9mm>*{}**@{-},
 <-2.3mm,2.3mm>*{\bullet};<-2.3mm,2.3mm>*{}**@{},
 <-2.3mm,2.3mm>*{};<-2.3mm,5mm>*{}**@{-},
 <-2.3mm,2.3mm>*{};<-4.6mm,0mm>*{}**@{-},
 \end{xy}=0\ .
$$
As all vertices are ternary, all higher differentials in our spectral sequences vanish, and we conclude that
$$
\sH(\sC, \delta)\simeq \sA,
$$
which proves the Theorem.
\hfill $\Box$

\bip

\no
{\bf 4.2.3. Remark.} As an independent check of the above arguments one can show using relations $R_1-R_3$ in
\S 2.6.2 that
every element of
$\frac{\Liebi^\circ}{\Liebi}$ can indeed be uniquely represented as a linear combinations
of graphs from the space $\sA$.

\bip

\no{\bf 4.2.4. Remark.} Proposition 4.2.1 and Theorem 4.2.2 can not be extended to the full wheeled
prop $\Liebi_\infty^\circlearrowright$, i.e. the natural surjection,
$$
\pi: (\Liebi_\infty^\circlearrowright, \delta) \lon (\Liebi^\circlearrowright, 0),
$$
is {\em not}\, a quasi-isomorphism. For example, the graph

$$
 \begin{xy}
<-5mm,5mm>*{\bullet};
<-5mm,5mm>*{};<-5mm,8mm>*{}**@{-},
<-5mm,5mm>*{};<-7mm,4mm>*{}**@{-},
 <0mm,0mm>*{\bullet};
<0mm,0mm>*{};<-5mm,-5mm>*{}**@{-},
 <0mm,0mm>*{};<-5mm,5mm>*{}**@{-},
<0mm,0mm>*{};<1.5mm,1.5mm>*{}**@{-},
 <0mm,0mm>*{};<1.5mm,-1.5mm>*{}**@{-};
<-5mm,-5mm>*{\bullet};
<-5mm,-5mm>*{};<-7mm,-2mm>*{}**@{-},
<-5mm,-5mm>*{};<-5mm,-8mm>*{}**@{-},
   \ar@{->}@(ul,dl) (-5.0,8.0)*{};(-7.0,4.0)*{},
   \ar@{->}@(ur,dr) (1.5,1.5)*{};(1.5,-1.5)*{},
   \ar@{->}@(ul,dl) (-7.0,-2.0)*{};(-5.0,-8.0)*{},
\end{xy}
\ - \
 \begin{xy}
<0mm,4mm>*{\bullet};
<0mm,4mm>*{};<0mm,8mm>*{}**@{-},
<0mm,4mm>*{};<3mm,2mm>*{}**@{-},
 <0mm,-5mm>*{\bullet};
<0mm,-5mm>*{};<0mm,4mm>*{}**@{-},
<0mm,-5mm>*{};<2mm,-3mm>*{}**@{-},
<0mm,-5mm>*{};<0mm,-8mm>*{}**@{-},
<-5mm,-5mm>*{};<0mm,4mm>*{}**@{-};
<-5mm,-5mm>*{\bullet};
<-5mm,-5mm>*{};<-7mm,-2mm>*{}**@{-},
<-5mm,-5mm>*{};<-5mm,-8mm>*{}**@{-},
   \ar@{->}@(ur,dr) (0,8.0)*{};(3.0,2.0)*{},
   \ar@{->}@(ur,dr) (2.0,-3.0)*{};(0.0,-8.0)*{},
   \ar@{->}@(ul,dl) (-7.0,-2.0)*{};(-5.0,-8.0)*{},
\end{xy}
\ + \
 \begin{xy}
<0mm,-4mm>*{\bullet};
<0mm,-4mm>*{};<0mm,-8mm>*{}**@{-},
<0mm,-4mm>*{};<3mm,-2mm>*{}**@{-},
 <0mm,5mm>*{\bullet};
<0mm,5mm>*{};<0mm,-4mm>*{}**@{-},
<0mm,5mm>*{};<2mm,3mm>*{}**@{-},
<0mm,5mm>*{};<0mm,8mm>*{}**@{-},
<-5mm,5mm>*{};<0mm,-4mm>*{}**@{-};
<-5mm,5mm>*{\bullet};
<-5mm,5mm>*{};<-7mm,2mm>*{}**@{-},
<-5mm,5mm>*{};<-5mm,8mm>*{}**@{-},
   \ar@{->}@(ur,dr) (3.0,-2.0)*{};(0.0,-8.0)*{},
   \ar@{->}@(ur,dr) (0.0,8.0)*{};(2.0,3.0)*{},
   \ar@{->}@(ul,dl) (-5.0,8.0)*{};(-7.0,2.0)*{},
\end{xy}
$$
represents a non-trivial cohomology class in $\sH^1(\Liebi_\infty^\circlearrowright, \delta)$.

%\sip
%
%It is shown in \cite{Medef} that there exists a dg free prop, $[\Liebi^\circlearrowright]_\infty$,
%which fits a commutative diagram,
%\[
% \xymatrix{
%[\Liebi^\circlearrowright]_\infty \ar[dr]^q \ar[r]^{s} & \Liebi_\infty^\circlearrowright
%  \ar[d]^{\pi} \\
% &
% \Liebi^\circlearrowright
% }
%\]
%with $s$ being an epimorphism of dg props, and $q$ a quasi-isomorphism. Moreover, every choice
%of such a {\em minimal prop resolution}\, of $\Liebi^\circlearrowright$ gives rise \cite{Medef} to a universal quantization $L_\infty$ morphism from
%the Shouten Lie algebra of polyvector fields on a formal manifold $M$
%to the Hochschild dg Lie algebra of polydifferential operators on the structure sheaf, $\f_M$, of $M$.
%This morphism heavily depends on the choice of $[\Liebi^\circlearrowright]_\infty$, and, in view of \cite{Ta}, one
%may draw a conclusion  that the Grothendieck-Teichm\"uller group
%acts on the set of minimal resolutions of the wheeled prop,
%$\Liebi^\circlearrowright$, of Lie 1-bialgebras.

\bip

\no
{\bf 4.3. Wheeled prop of Lie bialgebras.}
Let ${\mathsf L\mathsf i\mathsf e\mathsf B}$ be the prop of Lie bialgebras
which is generated by the dioperad very similar to $\Liebi$
except that both generating Lie and coLie operations are in degree
zero. This dioperad is again Koszul with Koszul substitution law
so that the analogue of Proposition 4.2.1 holds true for the
operadic wheelification, ${\mathsf L\mathsf i\mathsf e\mathsf B}_\infty^\circlearrowright$.
In fact, the analogue of Theorem~4.2.2 holds true for
${\mathsf L\mathsf i\mathsf e\mathsf B}_\infty^\circ$.

\bip

 %%%%%%%%%%%%%%%%%%%%%%%%%%%%%%%%%%%%%%%%%%%%
\no{\bf 4.4. Prop of infinitesimal bialgebras}.
Let $\mathsf I\mathsf B$ be the dioperad of infinitesimal bialgebras \cite{Ag}
which can be represented as a quotient,
$$
 {\mathsf I\mathsf B} = \frac{ {\mathsf D}\langle E \rangle}{{\sf Ideal} <R>}\ ,
$$
of the free prop generated by the following $\bS$-bimodule $E$,
\Bi
\item all $E(m,n)$ vanish except $E(2,1)$ and $E(1,2)$,
\item
$
E(2,1):= k[\bS_2]\ot {\bf 1_1} = span\left(
\begin{xy}
 <0mm,-0.55mm>*{};<0mm,-2.5mm>*{}**@{-},
 <0.5mm,0.5mm>*{};<2.2mm,2.2mm>*{}**@{-},
 <-0.48mm,0.48mm>*{};<-2.2mm,2.2mm>*{}**@{-},
 <0mm,0mm>*{\bullet};<0mm,0mm>*{}**@{},
   <0mm,-0.55mm>*{};<0mm,-3.8mm>*{_1}**@{},
   <0.5mm,0.5mm>*{};<2.7mm,2.8mm>*{^2}**@{},
   <-0.48mm,0.48mm>*{};<-2.7mm,2.8mm>*{^1}**@{},
 \end{xy}
 \ , \
 \begin{xy}
 <0mm,-0.55mm>*{};<0mm,-2.5mm>*{}**@{-},
 <0.5mm,0.5mm>*{};<2.2mm,2.2mm>*{}**@{-},
  <-0.48mm,0.48mm>*{};<-2.2mm,2.2mm>*{}**@{-},
 <0mm,0mm>*{\bullet};<0mm,0mm>*{}**@{},
   <0mm,-0.55mm>*{};<0mm,-3.8mm>*{_1}**@{},
   <0.5mm,0.5mm>*{};<2.7mm,2.8mm>*{^1}**@{},
   <-0.48mm,0.48mm>*{};<-2.7mm,2.8mm>*{^2}**@{},
 \end{xy}
 \right)
 $,
\item
$E(1,2):= {\bf 1_1}\ot k[\bS_2]=span\left(
\begin{xy}
 <0mm,0.66mm>*{};<0mm,3mm>*{}**@{-},
 <0.39mm,-0.39mm>*{};<2.2mm,-2.2mm>*{}**@{-},
 <-0.35mm,-0.35mm>*{};<-2.2mm,-2.2mm>*{}**@{-},
 <0mm,0mm>*{\bullet};<0mm,0mm>*{}**@{},
   <0mm,0.66mm>*{};<0mm,3.4mm>*{^1}**@{},
   <0.39mm,-0.39mm>*{};<2.9mm,-4mm>*{^2}**@{},
   <-0.35mm,-0.35mm>*{};<-2.8mm,-4mm>*{^1}**@{},
\end{xy}
\ , \
 \begin{xy}
 <0mm,0.66mm>*{};<0mm,3mm>*{}**@{-},
 <0.39mm,-0.39mm>*{};<2.2mm,-2.2mm>*{}**@{-},
 <-0.35mm,-0.35mm>*{};<-2.2mm,-2.2mm>*{}**@{-},
 <0mm,0mm>*{\bullet};<0mm,0mm>*{}**@{},
   <0mm,0.66mm>*{};<0mm,3.4mm>*{^1}**@{},
   <0.39mm,-0.39mm>*{};<2.9mm,-4mm>*{^1}**@{},
   <-0.35mm,-0.35mm>*{};<-2.8mm,-4mm>*{^2}**@{},
\end{xy}\right)
$,
\Ei
modulo the ideal  generated by the associativity conditions
for $\begin{xy}
 <0mm,0.66mm>*{};<0mm,3mm>*{}**@{-},
 <0.39mm,-0.39mm>*{};<2.2mm,-2.2mm>*{}**@{-},
 <-0.35mm,-0.35mm>*{};<-2.2mm,-2.2mm>*{}**@{-},
 <0mm,0mm>*{\bullet};<0mm,0mm>*{}**@{},
   %<0mm,0.66mm>*{};<0mm,3.4mm>*{^1}**@{},
   %<0.39mm,-0.39mm>*{};<2.9mm,-4mm>*{^2}**@{},
   %<-0.35mm,-0.35mm>*{};<-2.8mm,-4mm>*{^1}**@{},
\end{xy}$, co-associativity conditions
for
$
\begin{xy}
 <0mm,-0.55mm>*{};<0mm,-2.5mm>*{}**@{-},
 <0.5mm,0.5mm>*{};<2.2mm,2.2mm>*{}**@{-},
 <-0.48mm,0.48mm>*{};<-2.2mm,2.2mm>*{}**@{-},
 <0mm,0mm>*{\bullet};<0mm,0mm>*{}**@{},
   %<0mm,-0.55mm>*{};<0mm,-3.8mm>*{_1}**@{},
   %<0.5mm,0.5mm>*{};<2.7mm,2.8mm>*{^2}**@{},
   %<-0.48mm,0.48mm>*{};<-2.7mm,2.8mm>*{^1}**@{},
 \end{xy}
$,
%and the following compatibility relations,

%%%%%%%%%%%%%%%%%%%%% %% Lie[1]Bi %%%%%%%%%%%%%%%
$$
 \begin{xy}
 <0mm,2.47mm>*{};<0mm,-0.5mm>*{}**@{-},
 <0.5mm,3.5mm>*{};<2.2mm,5.2mm>*{}**@{-},
 <-0.48mm,3.48mm>*{};<-2.2mm,5.2mm>*{}**@{-},
 <0mm,3mm>*{\bullet};<0mm,3mm>*{}**@{},
  <0mm,-0.8mm>*{\bullet};<0mm,-0.8mm>*{}**@{},
<0mm,-0.8mm>*{};<-2.2mm,-3.5mm>*{}**@{-},
 <0mm,-0.8mm>*{};<2.2mm,-3.5mm>*{}**@{-},
\end{xy}
\  - \
\begin{xy}
 <0mm,-1.3mm>*{};<0mm,-3.5mm>*{}**@{-},
 <0.38mm,-0.2mm>*{};<2.2mm,2.2mm>*{}**@{-},
 <-0.38mm,-0.2mm>*{};<-2.2mm,2.2mm>*{}**@{-},
<0mm,-0.8mm>*{\bullet};<0mm,0.8mm>*{}**@{},
 <2.4mm,2.4mm>*{\bullet};<2.4mm,2.4mm>*{}**@{},
 <2.5mm,2.3mm>*{};<4.4mm,-0.8mm>*{}**@{-},
 <2.4mm,2.5mm>*{};<2.4mm,5.2mm>*{}**@{-},
    \end{xy}
\  - \
\begin{xy}
 <0mm,-1.3mm>*{};<0mm,-3.5mm>*{}**@{-},
 <0.38mm,-0.2mm>*{};<2.2mm,2.2mm>*{}**@{-},
 <-0.38mm,-0.2mm>*{};<-2.2mm,2.2mm>*{}**@{-},
<0mm,-0.8mm>*{\bullet};<0mm,0.8mm>*{}**@{},
 <-2.4mm,2.4mm>*{\bullet};<2.4mm,2.4mm>*{}**@{},
 <-2.5mm,2.3mm>*{};<-4.4mm,-0.8mm>*{}**@{-},
 <-2.4mm,2.5mm>*{};<-2.4mm,5.2mm>*{}**@{-},
    \end{xy} \ =\ 0.
$$
This is a Koszul dioperad with a Koszul substitution law.
Its minimal prop resolution,
$({\mathsf I\mathsf B}_\infty,\delta)$  is a dg prop freely generated by the $\bS$-bimodule
 $\mathsf E=\{\mathsf E(m,n)\}_{m,n\geq 1, m+n\geq 3}$, with
\[
\mathsf E(m,n):= k[\bS_m]\ot k[\bS_n][3-m-n]=\mbox{span}\left\langle
\begin{xy}
 <0mm,0mm>*{\bullet};<0mm,0mm>*{}**@{},
 <0mm,0mm>*{};<-8mm,5mm>*{}**@{-},
 <0mm,0mm>*{};<-4.5mm,5mm>*{}**@{-},
 <0mm,0mm>*{};<-1mm,5mm>*{\ldots}**@{},
 <0mm,0mm>*{};<4.5mm,5mm>*{}**@{-},
 <0mm,0mm>*{};<8mm,5mm>*{}**@{-},
   <0mm,0mm>*{};<-8.5mm,5.5mm>*{^1}**@{},
   <0mm,0mm>*{};<-5mm,5.5mm>*{^2}**@{},
   <0mm,0mm>*{};<4.5mm,5.5mm>*{^{m\hspace{-0.5mm}-\hspace{-0.5mm}1}}**@{},
   <0mm,0mm>*{};<9.0mm,5.5mm>*{^m}**@{},
 <0mm,0mm>*{};<-8mm,-5mm>*{}**@{-},
 <0mm,0mm>*{};<-4.5mm,-5mm>*{}**@{-},
 <0mm,0mm>*{};<-1mm,-5mm>*{\ldots}**@{},
 <0mm,0mm>*{};<4.5mm,-5mm>*{}**@{-},
 <0mm,0mm>*{};<8mm,-5mm>*{}**@{-},
   <0mm,0mm>*{};<-8.5mm,-6.9mm>*{^1}**@{},
   <0mm,0mm>*{};<-5mm,-6.9mm>*{^2}**@{},
   <0mm,0mm>*{};<4.5mm,-6.9mm>*{^{n\hspace{-0.5mm}-\hspace{-0.5mm}1}}**@{},
   <0mm,0mm>*{};<9.0mm,-6.9mm>*{^n}**@{},
 \end{xy}
\right\rangle.
\]
By Claim~3.5.4, the analogue of Proposition~4.2.1 does {\em not}\, hold true for ${\mathsf I\mathsf B}_\infty^\looparrowright$.
Moreover, it is not hard to check that the graph
$$
 \begin{xy}
<0mm,0mm>*{\bullet};
<0mm,0mm>*{};<3mm,3mm>*{}**@{-},
<0mm,0mm>*{};<-3mm,3mm>*{}**@{-},
<0mm,0mm>*{};<3mm,-3mm>*{}**@{-},
<0mm,0mm>*{};<-3mm,-3mm>*{}**@{-},
   \ar@{->}@(ur,dr) (-3.0,3.0)*{};(3.0,-3.0)*{},
\end{xy}
\ - \
 \begin{xy}
<0mm,0mm>*{\bullet};
<0mm,0mm>*{};<3mm,3mm>*{}**@{-},
<0mm,0mm>*{};<-3mm,3mm>*{}**@{-},
<0mm,0mm>*{};<3mm,-3mm>*{}**@{-},
<0mm,0mm>*{};<-3mm,-3mm>*{}**@{-},
   \ar@{->}@(ul,dl) (3.0,3.0)*{};(-3.0,-3.0)*{},
\end{xy}
$$
represents a non-trivial cohomology class in $\sH^{-1}({\mathsf I\mathsf B}_\infty^\circlearrowright)$.
Thus neither the analogue of Theorem~4.2.2 holds true for ${\mathsf I\mathsf B}_\infty^\circ$ nor
the natural surjection, ${\mathsf I\mathsf B}_\infty^\circlearrowright\rar {\mathsf I\mathsf B}^\circlearrowright$,
is a quasi-isomorphism. This example is of interest because the wheeled dg prop  ${\mathsf I\mathsf B}_\infty^\circlearrowright$
controls the cohomology of a directed version of Kontsevich's ribbon graph complex.

\bip

%%%%%%%%%%%%%%%%%%%%%%%%%%%%%%%%%%%%%%%%%%%%%%%%%%%%%%%%%%%%%%%%%
 \no{\bf 4.5. Wheeled quasi-minimal resolutions.}
Let $\sP$ be a graded prop with zero differential admitting a minimal resolution,
$$
\pi: (\sP_\infty =\PROP \langle \sE\rangle, \delta) \rar (\sP, 0).
$$
We shall use in the following discussion of this pair of props, $\sP_\infty$ and $\sP$, a so called
{\em Tate-Jozefak}\, grading\footnote{The Tate-Josefak grading of props $\Liebi_\infty$ and $\Liebi$, for example,
 assigns to  generating $(m,n)$ corollas degree $3-m-n$.}
 which, by definition, assignes degree zero to {\em all}\, generators
of $\sP$ and hence make $\sP_\infty$ into a non-positively graded differential prop, $\sP_\infty=\bigoplus_{i\leq 0} \sP_\infty^i$,
with cohomology concentrated in degree zero, $\sH^0(\sP_\infty, \delta)=\sP$.
Both props $\sP_\infty$ and $\sP$ admit
canonically wheeled extensions,
\Beqrn
\sP_\infty^\circlearrowright &:=& \bigoplus_{G\in \fG^\circlearrowright} G\langle \sE\rangle\\
\sP^\circlearrowright &:=& \sH^0(\sP_\infty^\circlearrowright, \delta).
\Eeqrn
However, the natural extension of the epimorphism $\pi$,
$$
\pi^\circlearrowright: (\sP_\infty^\circlearrowright, \delta) \rar (\sP, 0).
$$
fails in general to stay a quasi-isomorphism.

\bip

\no
Note that the dg prop, $(\sP_\infty^\circlearrowright, \delta)$, defined above is a {\em free}\, prop
$$
\sP_\infty^\circlearrowright := \bigoplus_{G\in \fG^\uparrow} G\langle \sE^\circlearrowright\rangle\\
$$
on the $\bS$-bimodule, $\sE^\circlearrowright=\{\sE(m,n)\}_{m,n\geq 0}$,
$$
\sE^\circlearrowright(m,n):= \bigoplus_{G\in \fG^\circlearrowright_{indec}(m,n)} G\langle \sE\rangle
$$
generated by indecomposable (with respect to graftings and disjoint unions) decorated wheeled graphs. Note that the induced differential
is {\em not}\, quadratic with respect to the generating set $\sE^\circlearrowright$.

\bip

%%%%%%%%%%%%%%%%%%%%%%%%%%%%%%%%%%%%%%%%%%%%%%%%%%%%%%%%%%%%%%%%%%%%%%%%%%%%%%%%%%%
 \no{\bf 4.5.1. Theorem-definition.} {\em There exists a dg free prop,
$([\sP^\circlearrowright]_\infty, \delta)$, which fits into a
commutative diagram of morphisms of  props,
\[
 \xymatrix{
[\sP^\circlearrowright]_\infty \ar[dr]^{qis} \ar[r]^{\al} & \sP^\circlearrowright_\infty
  \ar[d]^{\pi^\circlearrowright} \\
 &
 \sP^\circlearrowright
 }
\]
where $\al$ is an epimorphism of (nondifferential) props $qis$ a quasi-isomorphism of dg props.
The prop $[\sP^\circlearrowright]_\infty$ is called
a {\em quasi-minimal resolution} of $\sP^\circlearrowright$.
}

\bip

\Proof Let $s_1: \sH^{-1}(\sP_\infty^\circlearrowright)\rar \sP_\infty^\circlearrowright$ be any representation
of degree $-1$ cohomology classes (if there are any) as cycles. Set
$\sE_1:= \sH^{-1}(\sP_\infty^\circlearrowright)[1]$ and define a differential graded prop,
$$
\sQ_1:= \PROP \langle \sE^\circlearrowright \oplus \sE_1 \rangle
$$
with the differential $\delta$ extended to new generators as $s_1[1]$. By construction, $\sH^0(\sQ_1)=\sP^\circlearrowright$,
and $\sH^{-1}(\sQ_1)=0$.

\bip

\no
 Let $s_2: \sH^{-2}(\sQ_1)\rar \sQ_1$ be any representation
of degree $-2$ cohomology classes (if there are any) as cycles. Set
$\sE_2:= \sH^{-2}(\sQ_1)[1]$ and define a differential graded prop,
$$
\sQ_2:= \PROP \langle \sE^\circlearrowright \oplus \sE_1 \oplus \sE_2  \rangle
$$
with the differential $\delta$ extended to new generators as $s_2[1]$. By construction, $\sH^0(\sQ_2)=\sP^\circlearrowright$,
and $\sH^{-1}(\sQ_2)=\sH^{-2}(\sQ_2)=0$.

\bip

\no
Continuing by induction we construct a dg free prop, $[\sP^\circlearrowright]_\infty:=\lim_{n\rar\infty}\sQ_n=
\PROP \langle \sE^\circlearrowright \oplus \sE_1 \oplus \sE_2 \oplus \sE_3\oplus \ldots \rangle$
with all the cohomology concentrated in Tate-Jozefak degree $0$ and equal to $\sP^\circlearrowright$.
\hfill $\Box$

\bip

\no {\bf 4.5.2. Example.}
The prop $[\Ass^\circlearrowright]_\infty$ has been explicitly described in \cite{MMS}:
this is a dg free prop
$[\Ass^\circlearrowright]_\infty$  generated by
 planar (1,n)-corollas in degree $2-n$,
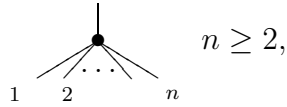
$$
\begin{xy}
 <0mm,0mm>*{\bullet};<0mm,0mm>*{}**@{},
 <0mm,0mm>*{};<-8mm,-5mm>*{}**@{-},
 <0mm,0mm>*{};<-4.5mm,-5mm>*{}**@{-},
 <0mm,0mm>*{};<0mm,-4mm>*{\ldots}**@{},
 <0mm,0mm>*{};<4.5mm,-5mm>*{}**@{-},
 <0mm,0mm>*{};<8mm,-5mm>*{}**@{-},
   <0mm,0mm>*{};<-11mm,-7.9mm>*{^{1}}**@{},
   <0mm,0mm>*{};<-4mm,-7.9mm>*{^{2}}**@{},
   %<0mm,0mm>*{};<4.5mm,5.5mm>*{^{n\hspace{-0.5mm}-\hspace{-0.5mm}1}}**@{},
   <0mm,0mm>*{};<10.0mm,-7.9mm>*{^{n}}**@{},
 <0mm,0mm>*{};<0mm,5mm>*{}**@{-},
 \end{xy} \ \ n\geq 2,
$$
{and} planar $(0,m+n)$-corollas in degree $-m-n$,
$$
\begin{xy}
 <0mm,-0.5mm>*{\blacktriangledown};<0mm,0mm>*{}**@{},
 <0mm,0mm>*{};<-16mm,-5mm>*{}**@{-},
 <0mm,0mm>*{};<-12mm,-5mm>*{}**@{-},
 <0mm,0mm>*{};<-3.5mm,-5mm>*{}**@{-},
 <0mm,0mm>*{};<-6mm,-5mm>*{...}**@{},
   <0mm,0mm>*{};<-21mm,-7.9mm>*{^{1}}**@{},
   <0mm,0mm>*{};<-14mm,-7.9mm>*{^{2}}**@{},
   <0mm,0mm>*{};<-4mm,-7.9mm>*{^{m}}**@{},
 <0mm,0mm>*{};<16mm,-5mm>*{}**@{-},
 <0mm,0mm>*{};<12mm,-5mm>*{}**@{-},
 <0mm,0mm>*{};<3.5mm,-5mm>*{}**@{-},
 <0mm,0mm>*{};<6.6mm,-5mm>*{...}**@{},
   <0mm,0mm>*{};<22mm,-7.9mm>*{^{m\hspace{-0.5mm}+\hspace{-0.5mm}n}}**@{},
   %<0mm,0mm>*{};<14mm,-7.9mm>*{^{2}}**@{},
   <0mm,0mm>*{};<6mm,-7.9mm>*{^{m\hspace{-0.5mm}+\hspace{-0.5mm}1}}**@{},
 \end{xy}, \ m,n\geq 1,
$$
having the cyclic skew-symmetry
\[
\hskip 3em
\begin{xy}
 <0mm,-0.5mm>*{\blacktriangledown};<0mm,0mm>*{}**@{},
 <0mm,0mm>*{};<-16mm,-5mm>*{}**@{-},
 <0mm,0mm>*{};<-12mm,-5mm>*{}**@{-},
 <0mm,0mm>*{};<-3.5mm,-5mm>*{}**@{-},
 <0mm,0mm>*{};<-6mm,-5mm>*{...}**@{},
   <0mm,0mm>*{};<-19mm,-7.9mm>*{^{1}}**@{},
   <0mm,0mm>*{};<-13mm,-7.9mm>*{^{2}}**@{},
   <0mm,0mm>*{};<-4mm,-7.9mm>*{^{m}}**@{},
 <0mm,0mm>*{};<16mm,-5mm>*{}**@{-},
 <0mm,0mm>*{};<12mm,-5mm>*{}**@{-},
 <0mm,0mm>*{};<3.5mm,-5mm>*{}**@{-},
 <0mm,0mm>*{};<6.6mm,-5mm>*{...}**@{},
   <0mm,0mm>*{};<19mm,-7.9mm>*{^{m\hspace{-0.5mm}+\hspace{-0.5mm}n}}**@{},
   %<0mm,0mm>*{};<14mm,-7.9mm>*{^{\sigma(2)}}**@{},
   <0mm,0mm>*{};<5mm,-7.9mm>*{^{m\hspace{-0.5mm}+\hspace{-0.5mm}1}}**@{},
 \end{xy}
\hskip -.5em = (-1)^{\sgn(\zeta)}\hspace{-13mm}
\begin{xy}
 <0mm,-0.5mm>*{\blacktriangledown};<0mm,0mm>*{}**@{},
 <0mm,0mm>*{};<-16mm,-5mm>*{}**@{-},
 <0mm,0mm>*{};<-12mm,-5mm>*{}**@{-},
 <0mm,0mm>*{};<-3.5mm,-5mm>*{}**@{-},
 <0mm,0mm>*{};<-6mm,-5mm>*{...}**@{},
   <0mm,0mm>*{};<-21mm,-7.9mm>*{^{\zeta(1)}}**@{},
   <0mm,0mm>*{};<-14mm,-7.9mm>*{^{\zeta(2)}}**@{},
   <0mm,0mm>*{};<-4mm,-7.9mm>*{^{\zeta(m)}}**@{},
 <0mm,0mm>*{};<16mm,-5mm>*{}**@{-},
 <0mm,0mm>*{};<12mm,-5mm>*{}**@{-},
 <0mm,0mm>*{};<3.5mm,-5mm>*{}**@{-},
 <0mm,0mm>*{};<6.6mm,-5mm>*{...}**@{},
   <0mm,0mm>*{};<20mm,-7.9mm>*{^{m\hspace{-0.5mm}+\hspace{-0.5mm}n}}**@{},
   %<0mm,0mm>*{};<14mm,-7.9mm>*{^{\sigma(2)}}**@{},
   <0mm,0mm>*{};<6mm,-7.9mm>*{^{m\hspace{-0.5mm}+\hspace{-0.5mm}1}}**@{},
 \end{xy}
\hskip -.5em
 =(-1)^{\sgn(\xi)}\hspace{-11mm}
\begin{xy}
 <0mm,-0.5mm>*{\blacktriangledown};<0mm,0mm>*{}**@{},
 <0mm,0mm>*{};<-16mm,-5mm>*{}**@{-},
 <0mm,0mm>*{};<-12mm,-5mm>*{}**@{-},
 <0mm,0mm>*{};<-3.5mm,-5mm>*{}**@{-},
 <0mm,0mm>*{};<-6mm,-5mm>*{...}**@{},
   <0mm,0mm>*{};<-21mm,-7.9mm>*{^{1}}**@{},
   <0mm,0mm>*{};<-14mm,-7.9mm>*{^{2}}**@{},
   <0mm,0mm>*{};<-4mm,-7.9mm>*{^{m}}**@{},
 <0mm,0mm>*{};<16mm,-5mm>*{}**@{-},
 <0mm,0mm>*{};<12mm,-5mm>*{}**@{-},
 <0mm,0mm>*{};<3.5mm,-5mm>*{}**@{-},
 <0mm,0mm>*{};<6.6mm,-5mm>*{...}**@{},
   <0mm,0mm>*{};<20mm,-7.9mm>*{^{\xi(m\hspace{-0.5mm}+\hspace{-0.5mm}n)}}**@{},
   %<0mm,0mm>*{};<14mm,-7.9mm>*{^{\sigma(2)}}**@{},
   <0mm,0mm>*{};<6mm,-7.9mm>*{^{\xi(m\hspace{-0.5mm}+\hspace{-0.5mm}1)}}**@{},
 \end{xy}
\]
with respect to the cyclic permutations $\zeta=(12\ldots m)$ and
$\xi=((m+1)(m+2)\ldots(m+n))$.
 The differential is given on
generators as
\begin{eqnarray*}
\label{ass1}
\delta
\begin{xy}
 <0mm,0mm>*{\bullet};<0mm,0mm>*{}**@{},
 <0mm,0mm>*{};<-8mm,-5mm>*{}**@{-},
 <0mm,0mm>*{};<-4.5mm,-5mm>*{}**@{-},
 <0mm,0mm>*{};<0mm,-4mm>*{\ldots}**@{},
 <0mm,0mm>*{};<4.5mm,-5mm>*{}**@{-},
 <0mm,0mm>*{};<8mm,-5mm>*{}**@{-},
   <0mm,0mm>*{};<-9mm,-7.9mm>*{^1}**@{},
   <0mm,0mm>*{};<-5mm,-7.9mm>*{^2}**@{},
   %<0mm,0mm>*{};<4.5mm,5.5mm>*{^{n\hspace{-0.5mm}-\hspace{-0.5mm}1}}**@{},
   <0mm,0mm>*{};<10.0mm,-7.9mm>*{^n}**@{},
 <0mm,0mm>*{};<0mm,5mm>*{}**@{-},
 \end{xy}
&=&\sum_{k=0}^{n-2}\sum_{l=2}^{n-k} (-1)^{k+l(n-k-l)+1}
\begin{xy}
<0mm,0mm>*{\bullet}, <0mm,5mm>*{}**@{-}, <4mm,-7mm>*{^{1\  \dots \
k\qquad\ \ k+l+1\ \ \dots\ \  n }}, <-14mm,-5mm>*{}**@{-},
<-6mm,-5mm>*{}**@{-}, <20mm,-5mm>*{}**@{-}, <8mm,-5mm>*{}**@{-},
<0mm,-5mm>*{}**@{-}, <0mm,-5mm>*{\bullet}; <-5mm,-10mm>*{}**@{-},
<-2mm,-10mm>*{}**@{-}, <2mm,-10mm>*{}**@{-}, <5mm,-10mm>*{}**@{-},
<0mm,-12mm>*{_{k+1\ \, \dots\ \, k+l}},
\end{xy},
\\
\label{ass2}
%druha formule
\hskip 2em \delta \hskip -1em
\begin{xy}
 <0mm,-0.5mm>*{\blacktriangledown};<0mm,0mm>*{}**@{},
 <0mm,0mm>*{};<-16mm,-5mm>*{}**@{-},
 <0mm,0mm>*{};<-12mm,-5mm>*{}**@{-},
 <0mm,0mm>*{};<-3.5mm,-5mm>*{}**@{-},
 <0mm,0mm>*{};<-6mm,-5mm>*{...}**@{},
   <0mm,0mm>*{};<-19mm,-7.9mm>*{^{1}}**@{},
   <0mm,0mm>*{};<-13mm,-7.9mm>*{^{2}}**@{},
   <0mm,0mm>*{};<-4mm,-7.9mm>*{^{m}}**@{},
 <0mm,0mm>*{};<16mm,-5mm>*{}**@{-},
 <0mm,0mm>*{};<12mm,-5mm>*{}**@{-},
 <0mm,0mm>*{};<3.5mm,-5mm>*{}**@{-},
 <0mm,0mm>*{};<6.6mm,-5mm>*{...}**@{},
   <0mm,0mm>*{};<19mm,-7.9mm>*{^{m\hspace{-0.5mm}+\hspace{-0.5mm}n}}**@{},
   %<0mm,0mm>*{};<14mm,-7.9mm>*{^{\sigma(2)}}**@{},
   <0mm,0mm>*{};<5mm,-7.9mm>*{^{m\hspace{-0.5mm}+\hspace{-0.5mm}1}}**@{},
 \end{xy}
\hskip -1.2em & = &
\sum_{i=0}^{m-1}\left((-1)^{m+1}\zeta\right)^i\sum_{j=1}^{n-1}
\left( (-1)^{n+1}\xi\right)^j \hskip -.2em
\left(\rule{0em}{2.3em}\right.
\begin{xy}
 <0mm,0mm>*{\bullet};<0mm,0mm>*{}**@{},
 <0mm,0mm>*{};<-11mm,-5mm>*{}**@{-},
 <0mm,0mm>*{};<-8.5mm,-5mm>*{}**@{-},
 <0mm,0mm>*{};<-5mm,-5mm>*{...}**@{},
 <0mm,0mm>*{};<-2.5mm,-5mm>*{}**@{-},
   <0mm,0mm>*{};<-13mm,-7.9mm>*{^1}**@{},
   <0mm,0mm>*{};<-9mm,-7.9mm>*{^2}**@{},
   <0mm,0mm>*{};<-3mm,-7.9mm>*{^m}**@{},
<0mm,0mm>*{};<11mm,-5mm>*{}**@{-},
 <0mm,0mm>*{};<8.5mm,-5mm>*{}**@{-},
 <0mm,0mm>*{};<5.5mm,-5mm>*{...}**@{},
 <0mm,0mm>*{};<2.5mm,-5mm>*{}**@{-},
   <0mm,0mm>*{};<14mm,-7.9mm>*{^{m\hspace{-0.5mm}+\hspace{-0.5mm}n}}**@{},
   <0mm,0mm>*{};<3mm,-7.9mm>*{^{m\hspace{-0.5mm}+\hspace{-0.5mm}1}}**@{},
 <0mm,-10mm>*{};<0mm,5mm>*{}**@{-},
(0,5)*{}
   \ar@{->}@(ur,dr) (0,-10)*{}
 \end{xy}
% \right.
 \\
 &&
\hskip -10em + \sum_{k=2}^{m} (-1)^{k(m+n)} \hskip -1em
\begin{xy}
 <0mm,-0.5mm>*{\blacktriangledown};<0mm,0mm>*{}**@{},
 <0mm,0mm>*{};<-16mm,-5mm>*{}**@{-},
 <0mm,0mm>*{};<-11mm,-5mm>*{}**@{-},
 <0mm,0mm>*{};<-3.5mm,-5mm>*{}**@{-},
 <0mm,0mm>*{};<-6mm,-5mm>*{...}**@{},
   <0mm,0mm>*{};<-10mm,-7.9mm>*{^{k\hspace{-0.5mm}+\hspace{-0.5mm}1}}**@{},
   <0mm,0mm>*{};<-4mm,-7.9mm>*{^{m}}**@{},
 <0mm,0mm>*{};<16mm,-5mm>*{}**@{-},
 <0mm,0mm>*{};<12mm,-5mm>*{}**@{-},
 <0mm,0mm>*{};<3.5mm,-5mm>*{}**@{-},
 <0mm,0mm>*{};<6.6mm,-5mm>*{...}**@{},
   <0mm,0mm>*{};<19mm,-7.9mm>*{^{m\hspace{-0.5mm}+\hspace{-0.5mm}n}}**@{},
   <0mm,0mm>*{};<5mm,-7.9mm>*{^{m\hspace{-0.5mm}+\hspace{-0.5mm}1}}**@{},
 <-16mm,-5.5mm>*{\bullet};<0mm,0mm>*{}**@{},
 <-16mm,-5.5mm>*{};<-20mm,-11mm>*{}**@{-},
 <-16mm,-5.5mm>*{};<-12mm,-11mm>*{}**@{-},
 <-16mm,-5.5mm>*{};<-18mm,-11mm>*{}**@{-},
 <-16mm,-5.5mm>*{};<-14mm,-11mm>*{}**@{-},
 <-16mm,-13mm>*{_{1\ \, \dots\ \ k}},
 \end{xy}
\nonumber
%\\
%&&
%\left.
+ \sum_{k=2}^{n-2} (-1)^{m+k+nk+1}
\begin{xy}
<0mm,-0.5mm>*{\blacktriangledown};<0mm,0mm>*{}**@{},
 <0mm,0mm>*{};<-16mm,-5mm>*{}**@{-},
 <0mm,0mm>*{};<-12mm,-5mm>*{}**@{-},
 <0mm,0mm>*{};<-3.5mm,-5mm>*{}**@{-},
 <0mm,0mm>*{};<-6mm,-5mm>*{...}**@{},
   <0mm,0mm>*{};<-18mm,-7.9mm>*{^{1}}**@{},
   <0mm,0mm>*{};<-13mm,-7.9mm>*{^{2}}**@{},
   <0mm,0mm>*{};<-4mm,-7.9mm>*{^{m}}**@{},
 <0mm,0mm>*{};<17mm,-5mm>*{}**@{-},
 <0mm,0mm>*{};<7mm,-5mm>*{}**@{-},
 <0mm,0mm>*{};<3.5mm,-5mm>*{}**@{-},
 <0mm,0mm>*{};<11.6mm,-5mm>*{...}**@{},
   <0mm,0mm>*{};<22mm,-7.9mm>*{^{m\hspace{-0.5mm}+\hspace{-0.5mm}n}}**@{},
   <0mm,0mm>*{};<11mm,-7.9mm>*{^{m\hspace{-0.5mm}+\hspace{-0.5mm}k+\hspace{-0.5mm}1}}**@{},
 <3.5mm,-5.5mm>*{\bullet};<0mm,0mm>*{}**@{},
 <3.5mm,-5.5mm>*{};<-0.5mm,-11mm>*{}**@{-},
 <3.5mm,-5.5mm>*{};<1.5mm,-11mm>*{}**@{-},
 <3.5mm,-5.5mm>*{};<5.5mm,-11mm>*{}**@{-},
 <3.5mm,-5.5mm>*{};<7.5mm,-11mm>*{}**@{-},
 <3.5mm,-13mm>*{_{m+1\, \dots\ m+k}},
 \end{xy}\left. \rule{0em}{2.3em}\right).
\nonumber
\end{eqnarray*}

\bip

\no {\bf 4.6. Wheeled Poisson structures.} By Theorem~4.5.1, there exists a natural extension,
$[\Liebi^\circlearrowright]_\infty$, of the dg prop
$\Liebi_\infty^\circlearrowright$ which provides us with a quasi-minimal prop resolution of $\Liebi^\circlearrowright$.

\bip

\no
It follows from Proposition~3.4.1 that representations of the dg prop $\Liebi_\infty^\circlearrowright$ in a
finite-dimensional dg space $V$ are the same as Poisson structures on the formal manifold $V$. This fact promts
us to make the following

\bip
%%%%%%%%%%%%%%%%%%%%%%%%%%%%%%%%%%%%%%%%%%%%%%

\no {\bf Definition.}
 Representations of the dg prop $[\Liebi^\circlearrowright]_\infty$ in a finite-dimensional dg space
$V$ are called {\em wheeled Poisson structures}\, on the formal graded manifold $M$.

\bip

\no
As generators of $[\Liebi^\circlearrowright]_\infty$ contain the generators of $\Liebi^\circlearrowright_\infty$,
a wheeled Poisson structure includes a degree one polyvector field $\ga\in \wedge^\bullet \cT_V$ satisfying
the Poisson equations, $[\ga, \ga]=0$. However, there are other generators of $[\Liebi^\circlearrowright]_\infty]$, hence
other {\em new}\, tensor fields on $V$ enter, in general, the content list of a wheeled Poisson structure.
One can get some insight into that content list from the above example of the prop
$[\Ass^\circlearrowright]_\infty$ where new fields
 are cyclically invariant
``functions" on the noncommutative manifold $V$ --- representations of the type $(0,m+n)$ corollas for all
$n\geq 1, m\geq 1$.

\bip

\no
The { new} fields of a wheeled Poisson structure
satisfy systems of differential equations which involve {\em traces}\, of tensors formed from polyvector fields
and their partial derivatives. Again the above example \S 4.6.5 (more precisely the value of the differential
on $(0,m+n)$-corollas) gives some intuition into the possible structure of the equations but not that much:
it follows from Theorem~4.2.2 that the terms in that equations which involve traces of polyvector fields $\ga$ can
 be neither linear in $\ga$ nor contain only wheels of genus one (i.e. ``wheeled" genus must be at least 2).

\bip

\no
Geometric meaning of a wheeled Poisson structure is not clear to the author at present. One can only be sure that this
notion is
a {\em canonical}\, generalization of ordinary Poisson structure in {\em finite}\, dimensions.
 Clearer picture might emerge from computation of the
cohomology of the complex $(\Liebi_\infty^\circlearrowright, \delta)$ which is a highly non-trivial problem comparable
in complexity with the problem of computing the homology of the directed versions of famous Kontsevich's
graph complexes \cite{Kon}.

\bip

%%%%%%%%%%%%%%%%%%%%%%%%%%%%%%%%%%%%%%%%%%%%%%%%%%%%%%%%%%%%%%%%%%%%%%%
%%%%%%%%%%%%%%%%%%%%%%%%%%%%%%%%%%%%%%%%%%%%%%%%%%%%%%%%%
\begin{center}
\bf \S 5. Deformation quantization via dg props
\end{center}

\bip

\no{\bf 5.1. Reminder.} Recall that in \S 2.6 we introduced the dg
prop of polyvector fields, ${\Liebi}_\infty$, and then in \S 4 we
studied its wheeled completion $\Liebi_\infty^\circlearrowright$ and
proved that the latter can be further  extended into a dg free prop,
$[\Liebi^\circlearrowright]_\infty$, which fits  the commutative
diagram,
\[
 \xymatrix{
[\Liebi^\circlearrowright]_\infty \ar[dr]^{qis} \ar[r]^{\al} &
\Liebi^\circlearrowright_\infty
  \ar[d]^{\pi^\circlearrowright} \\
 &
 \Liebi^\circlearrowright
 }
\]
with $\al$ being an epimorphism of non-differential props and $qis$ a
quasi-isomorphism of differential props. Representations of $[\Liebi^\circlearrowright]_\infty$ in a dg vector space
$V$ are called wheeled Poisson structures (see \S 4.6).

\bip

 In \S 2.7 we introduced the dg prop,
$\DefQ$, of ``star products".

\bip

\no
% Let us denote by
%$\widehat{[\Liebi^\circlearrowright]}_\infty$ the completion of
%$[\Liebi^\circlearrowright]_\infty$ with respect to the number of
%vertices.
The main purpose of this section is  to prove Theorem
5.2 which says that wheeled Poisson structures
can be deformation quantized.

\bip

\no{\bf 5.2. Theorem.} {\em There exists a morphism of dg
props,
\[
\hat{Q}:\ \ \DefQ \lon {[\Liebi^\circlearrowright]}_\infty,
\]
such that
\[
\pi_{1}\circ\al\circ \hat{Q}\left( \xy
%\begin{xy}
 <0mm,0mm>*{\mbox{$\xy *=<20mm,3mm>\txt{}*\frm{-}\endxy$}};<0mm,0mm>*{}**@{},
  <-10mm,1.5mm>*{};<-12mm,7mm>*{}**@{-},
  <-10mm,1.5mm>*{};<-11mm,7mm>*{}**@{-},
  <-10mm,1.5mm>*{};<-9.5mm,6mm>*{}**@{-},
  <-10mm,1.5mm>*{};<-8mm,7mm>*{}**@{-},
 <-10mm,1.5mm>*{};<-9.5mm,6.6mm>*{.\hspace{-0.4mm}.\hspace{-0.4mm}.}**@{},
 <0mm,0mm>*{};<-6.5mm,3.6mm>*{.\hspace{-0.1mm}.\hspace{-0.1mm}.}**@{},
  <-3mm,1.5mm>*{};<-5mm,7mm>*{}**@{-},
  <-3mm,1.5mm>*{};<-4mm,7mm>*{}**@{-},
  <-3mm,1.5mm>*{};<-2.5mm,6mm>*{}**@{-},
  <-3mm,1.5mm>*{};<-1mm,7mm>*{}**@{-},
 <-3mm,1.5mm>*{};<-2.5mm,6.6mm>*{.\hspace{-0.4mm}.\hspace{-0.4mm}.}**@{},
  <2mm,1.5mm>*{};<0mm,7mm>*{}**@{-},
  <2mm,1.5mm>*{};<1mm,7mm>*{}**@{-},
  <2mm,1.5mm>*{};<2.5mm,6mm>*{}**@{-},
  <2mm,1.5mm>*{};<4mm,7mm>*{}**@{-},
 <2mm,1.5mm>*{};<2.5mm,6.6mm>*{.\hspace{-0.4mm}.\hspace{-0.4mm}.}**@{},
 <0mm,0mm>*{};<6mm,3.6mm>*{.\hspace{-0.1mm}.\hspace{-0.1mm}.}**@{},
<10mm,1.5mm>*{};<8mm,7mm>*{}**@{-},
  <10mm,1.5mm>*{};<9mm,7mm>*{}**@{-},
  <10mm,1.5mm>*{};<10.5mm,6mm>*{}**@{-},
  <10mm,1.5mm>*{};<12mm,7mm>*{}**@{-},
 <10mm,1.5mm>*{};<10.5mm,6.6mm>*{.\hspace{-0.4mm}.\hspace{-0.4mm}.}**@{},
 <-10mm,-1.5mm>*{};<-12mm,-6mm>*{}**@{-},
 <-7mm,-1.5mm>*{};<-8mm,-6mm>*{}**@{-},
 <-4mm,-1.5mm>*{};<-4.5mm,-6mm>*{}**@{-},
 <0mm,0mm>*{};<0mm,-4.6mm>*{.\hspace{0.1mm}.\hspace{0.1mm}.}**@{},
<10mm,-1.5mm>*{};<12mm,-6mm>*{}**@{-},
 <7mm,-1.5mm>*{};<8mm,-6mm>*{}**@{-},
  <4mm,-1.5mm>*{};<4.5mm,-6mm>*{}**@{-},
  %%%%
<0mm,0mm>*{};<-9.5mm,8.2mm>*{^{I_{ 1}}}**@{},
<0mm,0mm>*{};<-3mm,8.2mm>*{^{I_{ i}}}**@{},
<0mm,0mm>*{};<2mm,8.2mm>*{^{I_{ i+1}}}**@{},
<0mm,0mm>*{};<10mm,8.2mm>*{^{I_{ k}}}**@{},
 %%%
<0mm,0mm>*{};<-12mm,-7.4mm>*{_1}**@{},
<0mm,0mm>*{};<-8mm,-7.4mm>*{_2}**@{},
<0mm,0mm>*{};<-4mm,-7.4mm>*{_3}**@{},
<0mm,0mm>*{};<6mm,-7.4mm>*{\ldots}**@{},
%<0mm,0mm>*{};<8.5mm,-7.4mm>*{_{n-1}}**@{},
<0mm,0mm>*{};<12.5mm,-7.4mm>*{_n}**@{},
\endxy
\right) = \left\{\Ba{cl}
 \begin{xy}
 <0mm,0mm>*{\bullet};<0mm,0mm>*{}**@{},
 <0mm,0mm>*{};<-8mm,5mm>*{}**@{-},
 <0mm,0mm>*{};<-4.5mm,5mm>*{}**@{-},
 <0mm,0mm>*{};<-1mm,5mm>*{\ldots}**@{},
 <0mm,0mm>*{};<4.5mm,5mm>*{}**@{-},
 <0mm,0mm>*{};<8mm,5mm>*{}**@{-},
   <0mm,0mm>*{};<-8.5mm,5.5mm>*{^1}**@{},
   <0mm,0mm>*{};<-5mm,5.5mm>*{^2}**@{},
   <0mm,0mm>*{};<4.5mm,5.5mm>*{^{k\hspace{-0.5mm}-\hspace{-0.5mm}1}}**@{},
   <0mm,0mm>*{};<9.0mm,5.5mm>*{^k}**@{},
 <0mm,0mm>*{};<-8mm,-5mm>*{}**@{-},
 <0mm,0mm>*{};<-4.5mm,-5mm>*{}**@{-},
 <0mm,0mm>*{};<-1mm,-5mm>*{\ldots}**@{},
 <0mm,0mm>*{};<4.5mm,-5mm>*{}**@{-},
 <0mm,0mm>*{};<8mm,-5mm>*{}**@{-},
   <0mm,0mm>*{};<-8.5mm,-6.9mm>*{^1}**@{},
   <0mm,0mm>*{};<-5mm,-6.9mm>*{^2}**@{},
   <0mm,0mm>*{};<4.5mm,-6.9mm>*{^{n\hspace{-0.5mm}-\hspace{-0.5mm}1}}**@{},
   <0mm,0mm>*{};<9.0mm,-6.9mm>*{^n}**@{},
 \end{xy}
&
\mbox{for}\  |I_1|=\ldots=|I_k|=1 \\
0 & \mbox{otherwise}. \Ea \right.
\]
where $\pi_1$ is the projection to to the subspace of $\Liebi_\infty^\circlearrowright$
spanned by one-vertex graphs. }

\sip

%%%%%%%%%%%%%%%%%%%%%%%%%%%%%%%%%%%%%%%%%%%%%%%%%%%%%%%%%%%%%%
\noindent{\bf Proof}. We prove it below in three movements:

\sip

\noindent{\sc Step I:} we construct a morphism of dg props, $q:
(\DefQh, \delta) \lon (\Liebi^\circlearrowright,0)$.

\noindent{\sc Step II:} we prove  cofibrancy of $(\DefQh, \delta)$
and then use this property to show that $q$ can be lifted to a
morphism $Q$ making the diagram
\[
 \xymatrix{
  & [\Liebi^\circlearrowright]_\infty
  \ar[d]^{qis} \\
 \DefQh\ar[ur]^Q\ar[r]_{q} &
 \Liebi^\circlearrowright
 }
\]
commutative;

\noindent{\sc Step III:} we set $\hat{Q}$ to be the composition,
$$
\DefQ \stackrel{\chi_\hbar}{\lon}
 \DefQ[[\hbar]] \stackrel{Q}{\lon}
[\Liebi^\circlearrowright]_\infty[[\hbar]]
$$
at $\hbar=1$ (which makes sense as a morphism between completed
props).

\bip

\noindent{\bf Step I.} To construct  a morphism $q: (\DefQh, \delta)
\lon (\Liebi^\circlearrowright,0)$ is the same as to deformation
quantize an arbitrary finite-dimensional Lie 1-bialgebra, that is, a
pair $(\nu, \xi)$, consisting of a  linear Poisson
structure $\nu$ and a quadratic homological vector field $\xi$
such that $[\xi, \nu]_S=0$.

\sip

As a first approximation to $q$ we discuss first a well-known
Poincare-Birkhoff-Witt quantization of linear Poisson structures,
which, in the language of props, translates into existence of a
morphism,
 $$
\PBW : (\DefQh, \delta) \lon (\CoLie,0),
$$
 where $\CoLie$ is the prop
 of coLie algebras.  The morphism $q$ we construct below in  II.2 will make the diagram
\[
 \xymatrix{
  & \Liebi^\circlearrowright
  \ar[d]^{s} \\
 \DefQh\ar[ur]^q\ar[r]_{\PBW} &
 \CoLie^\circlearrowright
 }
\]
commutative, where $s$ is a natural surjection,
$$
\Ba{rccc}
s: & \Liebi & \lon & \CoLie \\
&    (\nu,\xi)& \lon & \nu, \Ea
$$
``forgetting" the quadratic homological vector field $\xi$.

\bip

%%%%%%%%%%%%%%%%%%%%%%%%%%%%%%%%%%%%%%
\noindent{\bf I.1. Poincare-Birkhoff-Witt quantization}. A representation of
$\CoLie$ in a vector space $V$ is the same as a linear Poisson
structure on $V$ viewed as a graded manifold.
This Poisson structure makes the graded commutative algebra
$\f_{V}[[\hbar]]:=\widehat{\odot^\bullet}V^*\ot \K[[\hbar]]$ into a  Lie algebra with respect to the
Poisson brackets, $\{\ ,\ \}$. Clearly, $V^*[[\hbar]]\subset
\f_{V}[[\hbar]]$ is a Lie subalgebra consisting of linear functions.
Let ${\cU}_\hbar$ be the associated universal enveloping
algebra defined as a quotient,
$$
\cU_\hbar:= \widehat{\ot^\bullet} V^*[[\hbar]]/ {\mathcal I},
$$
where the ideal $\mathcal I$ is generated by all expressions of the form
$$
X_1\ot X_2 - (-1)^{|X_1||X_2|}X_2\ot X_1 - \hbar\{X_1, X_2\}, \ \
X_i\in V^*.
$$

To construct a morphism of props $\PBW$ is the same as to
deformation quantize an arbitrary
 linear Poisson structure $\nu$. Which is a well-known trick:
the Poincare-Birkhoff-Witt theorem says that the natural
morphism,
$$
s: \f_{V}[[\hbar]] \lon
\cU_\hbar,
$$
is an isomorphism of vector spaces so that  one can quantize $\nu$
by the formula,
$$
f \star_{\hbar} g := s^{-1}(s(f)\circ s(g)), \ \ \
\forall f,g\in \f_V,
$$
where $\circ$ is the product in ${\cU}(V_\hbar)$. This proves the
existence of
 the required prop
morphism $\PBW$.

\bip

%%%%%%%%%%%%%%%%%%%%%%%%%%%%%
\noindent{\bf I.2. Deformation quantization of $\Liebi$ algebras}. A
representation of $\Liebi$ in a graded vector space $V$ is the same
as a degree 0 linear Poisson structure,
$$
\nu= \hbar\sum_{\al,\be,\ga} \Phi^{\al\be}_\ga x^\ga\frac{\p}{\p
x^\al}\wedge \frac{\p}{\p x^\be},
$$
together with a degree 1 quadratic vector field,
$$
\xi= \hbar\sum_{\al,\be,\ga} C^\ga_{\al\be}
x^{\al}x^{\be}\frac{\p}{\p x^\ga},
$$
 on $V$ satisfying $[\xi,\xi]_S=0$ and $Lie_\xi \nu=0$. Here $\{x^a\}$ are arbitrary linear coordinates on $V$
and  $Lie_\xi$ stands for the Lie
derivative along the vector field $\xi$. The association,
$$
\Phi^{\al\be}_\ga \simeq \begin{xy}
 <0mm,-0.55mm>*{};<0mm,-2.5mm>*{}**@{-},
 <0.5mm,0.5mm>*{};<2.2mm,2.2mm>*{}**@{-},
 <-0.48mm,0.48mm>*{};<-2.2mm,2.2mm>*{}**@{-},
 <0mm,0mm>*{\circ};<0mm,0mm>*{}**@{},
   <0mm,-0.55mm>*{};<0mm,-3.8mm>*{_\ga}**@{},
   <0.5mm,0.5mm>*{};<2.7mm,2.8mm>*{^\be}**@{},
   <-0.48mm,0.48mm>*{};<-2.7mm,2.8mm>*{^\al}**@{},
 \end{xy}
\ \ , \ \
 C^\ga_{\al\be}  \simeq
\begin{xy}
 <0mm,0.66mm>*{};<0mm,3mm>*{}**@{-},
 <0.39mm,-0.39mm>*{};<2.2mm,-2.2mm>*{}**@{-},
 <-0.35mm,-0.35mm>*{};<-2.2mm,-2.2mm>*{}**@{-},
 <0mm,0mm>*{\bullet};<0mm,0mm>*{}**@{},
   <0mm,0.66mm>*{};<0mm,3.4mm>*{^\ga}**@{},
   <0.39mm,-0.39mm>*{};<2.9mm,-4mm>*{^\be}**@{},
   <-0.35mm,-0.35mm>*{};<-2.8mm,-4mm>*{^\al}**@{},
\end{xy}
$$
translates the equations,
$$
[\nu,\nu]_S=0, \ \ [\xi,\xi]_S=0, \ \ \ Lie_\xi\nu=0,
$$
precisely into the graph relations $R_1$-$R_3$ in \S 2.6.2.

\bip

\no
To prove existence of a  morphism  $q: \DefQh \rar
\Liebi^\circlearrowright$ is the same as  to deformation quantize
such a pair $(\nu,\xi)$, that is,  to construct from $(\nu,\xi)$ a
degree 2 function $\Ga_0\in
\Hom_2(\K,\f_V)[[\hbar]]=\f_V[2][[\hbar]]$, a differential operator
$\Gamma_1\in \Hom_1(\f_V,\f_V)[[\hbar]]$,  and a bi-differential
operator $\Gamma_2\in \Hom_0(\f_V^{\ot 2},\f_V)[[\hbar]]$ such that
the following equations are satisfied\footnote{There can {\em not}\,
be non-vanishing terms $\Gamma_k\in \Hom_{2-k}(\f_V^{\ot
k},\f_V)[[\hbar]]$ with $k\geq 3$ for degree reasons.},

$$
\Ga_1\Ga_0=0, \ \ \ \ \ \ \ \ \ \ \ \  \ \ \ \Ga_1^2
+[\Ga_0,\Ga_2]_{\rm H}=0,
$$
$$
 d_{\rm H}\Ga_1 +
[\Ga_1, \Ga_2]_{\rm H}=0,  \ \ \ \ \ \  d_{\rm H}\Ga_2+ \frac{1}{2}
[\Ga_2, \Ga_2]_{\rm H}=0.
$$
We can solve the last equation by setting $\Ga_2$ to be related to
$\Phi_\ga^{\al\be}$ via the $\PBW$ quantization, i.e.\ we choose the
star product,
$$
\star_\hbar = \mbox{usual product of functions} + \Ga_2
$$
to be given as before, $f \star_{\hbar} g :=
\sigma^{-1}(\sigma(f)\cdot \sigma(g))$, $\forall f,g\in \f_V$.

\bip

\no
Our next task is to find a degree 1 differential operator $\Ga_1$ such
that $ d_{\rm H}\Ga_1 +  [\Ga_1, \Ga_2]_{\rm H}=0$ which is
equivalent to saying that $\Ga_1$ is a derivation of the star
product,
$$
\Ga_1(f \star_{\hbar} g )= (\Ga_1 f )\star_{\hbar} g + (-1)^{|f|} f
\star_{\hbar}( \Ga_1 g ), \ \ \ \ \ \forall f,g\in \f_V.
$$

\bip

\no
Consider first a derivation of the tensor algebra  ${\ot^\bullet
V}[[\hbar]]$ given on the generators, $\{t^\ga\}$, by
$$
\hat{\xi}(t^\ga):=  \hbar\sum_{\al,\be} C^\ga_{\al\be} x^{\al}\ot
x^{\be}.
$$
It is straightforward to check using equations $[\xi,\xi]_S=0$ and
$Lie_\xi\nu=0$
 that
$$
\hat{\xi}\left(x^\al\ot x^\be - (-1)^{|\al||\be|}x^\be\ot x^\al -
\hbar\sum_\ga \Phi^{\al\be}_\ga x^\ga\right)=0 \bmod {\cI}
$$
so that $\hat{\xi}$ descends to a derivation of the star product.
Hence setting $\Ga_1=\hat{\xi}$ we solve the equation $d_{\rm
H}\Ga_1 +  [\Ga_1, \Ga_2]_{\rm H}=0$. However, $\hat{\xi}^2\neq 0$,
which is enough to check on a generator $t^\al$,
$$
\hat{\xi}^2(x^\al)= - \sum \frac{1}{3}
C^\al_{\be\ga}C^\be_{\mu\nu}\Phi^{\nu\ga}_\tau \Phi^{\mu\tau}_\var
x^\var,
$$
or, in terms of graphs,
$$
-\frac{1}{3}\ \
\begin{xy}
 <0mm,-1.3mm>*{};<0mm,-3.5mm>*{}**@{-},
 <0.4mm,0.0mm>*{};<2.4mm,2.1mm>*{}**@{-},
 <-0.38mm,-0.2mm>*{};<-2.8mm,2.5mm>*{}**@{-},
<0mm,-0.8mm>*{\circ};<0mm,0.8mm>*{}**@{},
 <2.96mm,2.4mm>*{\circ};<2.45mm,2.35mm>*{}**@{},
 <2.4mm,2.8mm>*{};<0mm,5mm>*{}**@{-},
 <3.35mm,2.9mm>*{};<5.5mm,5mm>*{}**@{-},
     <0mm,-1.3mm>*{};<0mm,-5.3mm>*{_\var}**@{},
     %<2.5mm,2.3mm>*{};<5.1mm,-2.6mm>*{^1}**@{},
%
<-2.8mm,2.5mm>*{};<0mm,5mm>*{\bullet}**@{},
<-2.8mm,2.5mm>*{};<0mm,5mm>*{}**@{-},
<2.96mm,5mm>*{};<2.96mm,7.5mm>*{\bullet}**@{},
<0.2mm,5.1mm>*{};<2.8mm,7.5mm>*{}**@{-},
<0.2mm,5.1mm>*{};<2.8mm,7.5mm>*{}**@{-},
<5.5mm,5mm>*{};<2.8mm,7.5mm>*{}**@{-},
<2.9mm,7.5mm>*{};<2.9mm,10.5mm>*{}**@{-},
<2.9mm,7.5mm>*{};<2.9mm,11.1mm>*{^\al}**@{},
    \end{xy} .
$$
The data, $(\Ga_2, \Ga_1=\hat{\xi})$, are given by directed graphs
without cycles  (in particular, these data make sense for {\em
infinite}\, dimensional representations of the prop  $\Liebi$).
However, to solve the next equation, $\Ga_1^2 + [\Ga_0,\Ga_2]_{\rm
H}=0$,
 one has to construct from the generators of
  $\Liebi$ a non-vanishing graph, $\Ga_0$, with no output legs which is impossible to do
 without using
 graphs with oriented wheels. This is the reason why in Theorem~5.2 we attempt to construct
a morphism into a {\em wheeled extension}\, of the prop $\Liebi_\infty$:   there already does {\em not}\, exists a
morphism,
  $\DefQ^\hbar \rar \Liebi$, into the original unwheeled version of the prop of
  Lie 1-bialgebras which satisfies quasi-classical limit condition that terms linear in $\hbar$ are $\nu+\xi$.

\bip

\no
 In fact, one has to modify the naive choice,
 $\Ga_1=\hat{\xi}$, to get a solution of  $\Ga_1^2 + [\Ga_0,\Ga_2]_{\rm H}=0$.
Consider first a derivation of the tensor algebra  ${\ot^\bullet
V}[[\hbar]]$ given on the generators, $\{x^\ga\}$, by
$$
\breve{\xi}(t^\ga):= %\frac{1}{24}
a \hbar^3 \Theta^\ga %C^\ga_{\al\be}\Phi^{\al\mu}_\nu \Phi^{\be\nu}_\mu
 +
  \hbar\sum_{\al,\be} C^\ga_{\al\be}
x^{\al}\ot x^{\be},
$$
where
$$
\Theta^\ga:= \sum_{\al,\be,\mu,\nu} C^\ga_{\al\be}\Phi^{\al\mu}_\nu
\Phi^{\be\nu}_\mu = \ \
\begin{xy}
 <0.4mm,0.0mm>*{};<2.4mm,2.1mm>*{}**@{-},
 <-0.38mm,-0.2mm>*{};<-2.8mm,2.5mm>*{}**@{-},
<0mm,-0.8mm>*{\circ};<0mm,0.8mm>*{}**@{},
 <2.96mm,2.4mm>*{\circ};<2.45mm,2.35mm>*{}**@{},
 <2.4mm,2.8mm>*{};<0mm,5mm>*{}**@{-},
<-2.8mm,2.5mm>*{};<0mm,5mm>*{\bullet}**@{},
<-2.8mm,2.5mm>*{};<0mm,5mm>*{}**@{-},
<0mm,5mm>*{};<0mm,8.8mm>*{}**@{-},
<2.9mm,7.5mm>*{};<0mm,9.4mm>*{^\ga}**@{},
(3.2,3)*{}
   \ar@{->}@(ur,dr) (0.1,-1.4)*{}
    \end{xy} .
$$
and $a$ is a constant.

\bip

%%%%%%%%%%%%%%%%%%%%%%%%%%%%%%
\noindent{\bf I.2.1. Lemma}. $\Phi^{\al\be}_\ga \Theta^\ga=0$, {\em i.e.},
$
\begin{xy}
 <0.4mm,0.0mm>*{};<2.4mm,2.1mm>*{}**@{-},
 <-0.38mm,-0.2mm>*{};<-2.8mm,2.5mm>*{}**@{-},
<0mm,-0.8mm>*{\circ};<0mm,0.8mm>*{}**@{},
 <2.96mm,2.4mm>*{\circ};<2.45mm,2.35mm>*{}**@{},
 <2.4mm,2.8mm>*{};<0mm,5mm>*{}**@{-},
<-2.8mm,2.5mm>*{};<0mm,5mm>*{\bullet}**@{},
<-2.8mm,2.5mm>*{};<0mm,5mm>*{}**@{-},
<0mm,5mm>*{};<0mm,8.6mm>*{}**@{-},
<0mm,9mm>*{\circ};<0mm,0mm>*{}**@{},
<-0.5mm,9.6mm>*{};<-3.1mm,11.6mm>*{}**@{-},
<0.5mm,9.6mm>*{};<3.1mm,11.6mm>*{}**@{-},
<2.9mm,7.5mm>*{};<-3.1mm,12.4mm>*{^\al}**@{},
<2.9mm,7.5mm>*{};<3.1mm,12.4mm>*{^\be}**@{},
(3.2,3)*{}
   \ar@{->}@(ur,dr) (0.1,-1.4)*{}
    \end{xy} \ =\ 0.
$

\sip

\noindent\Proof Using relations
$$
%%%%%%%%%%%%%%%%%%%%% %% Lie[1]Bi %%%%%%%%%%%%%%%
 \begin{xy}
 <0mm,2.47mm>*{};<0mm,-0.5mm>*{}**@{-},
 <0.5mm,3.5mm>*{};<2.2mm,5.2mm>*{}**@{-},
 <-0.48mm,3.48mm>*{};<-2.2mm,5.2mm>*{}**@{-},
 <0mm,3mm>*{\circ};<0mm,3mm>*{}**@{},
  <0mm,-0.8mm>*{\bullet};<0mm,-0.8mm>*{}**@{},
<0mm,-0.8mm>*{};<-2.2mm,-3.5mm>*{}**@{-},
 <0mm,-0.8mm>*{};<2.2mm,-3.5mm>*{}**@{-},
     <0.5mm,3.5mm>*{};<2.8mm,5.9mm>*{^\be}**@{},
     <-0.48mm,3.48mm>*{};<-2.8mm,5.7mm>*{^\al}**@{},
   <0mm,-0.8mm>*{};<-2.7mm,-5.2mm>*{^1}**@{},
   <0mm,-0.8mm>*{};<2.7mm,-5.2mm>*{^2}**@{},
\end{xy}
\  = \
\begin{xy}
 <0mm,-1.3mm>*{};<0mm,-3.5mm>*{}**@{-},
 <0.38mm,-0.2mm>*{};<2.2mm,2.2mm>*{}**@{-},
 <-0.38mm,-0.2mm>*{};<-2.2mm,2.2mm>*{}**@{-},
<0mm,-0.8mm>*{\circ};<0mm,0.8mm>*{}**@{},
% <-2.25mm,2.2mm>*{};<-2.2mm,5.2mm>*{}**@{-},
 <2.4mm,2.4mm>*{\bullet};<2.4mm,2.4mm>*{}**@{},
 <2.5mm,2.3mm>*{};<4.4mm,-0.8mm>*{}**@{-},
% <4.4mm,-0.8mm>*{};<4.4mm,-3.5mm>*{}**@{-},
 <2.4mm,2.5mm>*{};<2.4mm,5.2mm>*{}**@{-},
     <0mm,-1.3mm>*{};<0mm,-5.3mm>*{^1}**@{},
     <2.5mm,2.3mm>*{};<5.1mm,-2.6mm>*{^2}**@{},
    <2.4mm,2.5mm>*{};<2.4mm,5.9mm>*{^\be}**@{},
    <-0.38mm,-0.2mm>*{};<-2.8mm,2.5mm>*{^\al}**@{},
    \end{xy}
\  - \
\begin{xy}
 <0mm,-1.3mm>*{};<0mm,-3.5mm>*{}**@{-},
 <0.38mm,-0.2mm>*{};<2.2mm,2.2mm>*{}**@{-},
 <-0.38mm,-0.2mm>*{};<-2.2mm,2.2mm>*{}**@{-},
<0mm,-0.8mm>*{\circ};<0mm,0.8mm>*{}**@{},
% <-2.25mm,2.2mm>*{};<-2.2mm,5.2mm>*{}**@{-},
 <2.4mm,2.4mm>*{\bullet};<2.4mm,2.4mm>*{}**@{},
 <2.5mm,2.3mm>*{};<4.4mm,-0.8mm>*{}**@{-},
% <4.4mm,-0.8mm>*{};<4.4mm,-3.5mm>*{}**@{-},
 <2.4mm,2.5mm>*{};<2.4mm,5.2mm>*{}**@{-},
     <0mm,-1.3mm>*{};<0mm,-5.3mm>*{^1}**@{},
     <2.5mm,2.3mm>*{};<5.1mm,-2.6mm>*{^2}**@{},
    <2.4mm,2.5mm>*{};<2.4mm,5.7mm>*{^\al}**@{},
    <-0.38mm,-0.2mm>*{};<-2.8mm,2.8mm>*{^\be}**@{},
    \end{xy}
\  + \
\begin{xy}
 <0mm,-1.3mm>*{};<0mm,-3.5mm>*{}**@{-},
 <0.38mm,-0.2mm>*{};<2.2mm,2.2mm>*{}**@{-},
 <-0.38mm,-0.2mm>*{};<-2.2mm,2.2mm>*{}**@{-},
<0mm,-0.8mm>*{\circ};<0mm,0.8mm>*{}**@{},
% <-2.25mm,2.2mm>*{};<-2.2mm,5.2mm>*{}**@{-},
 <2.4mm,2.4mm>*{\bullet};<2.4mm,2.4mm>*{}**@{},
 <2.5mm,2.3mm>*{};<4.4mm,-0.8mm>*{}**@{-},
% <4.4mm,-0.8mm>*{};<4.4mm,-3.5mm>*{}**@{-},
 <2.4mm,2.5mm>*{};<2.4mm,5.2mm>*{}**@{-},
     <0mm,-1.3mm>*{};<0mm,-5.3mm>*{^2}**@{},
     <2.5mm,2.3mm>*{};<5.1mm,-2.6mm>*{^1}**@{},
    <2.4mm,2.5mm>*{};<2.4mm,5.9mm>*{^\be}**@{},
    <-0.38mm,-0.2mm>*{};<-2.8mm,2.5mm>*{^\al}**@{},
    \end{xy}
\  - \
\begin{xy}
 <0mm,-1.3mm>*{};<0mm,-3.5mm>*{}**@{-},
 <0.38mm,-0.2mm>*{};<2.2mm,2.2mm>*{}**@{-},
 <-0.38mm,-0.2mm>*{};<-2.2mm,2.2mm>*{}**@{-},
<0mm,-0.8mm>*{\circ};<0mm,0.8mm>*{}**@{},
% <-2.25mm,2.2mm>*{};<-2.2mm,5.2mm>*{}**@{-},
 <2.4mm,2.4mm>*{\bullet};<2.4mm,2.4mm>*{}**@{},
 <2.5mm,2.3mm>*{};<4.4mm,-0.8mm>*{}**@{-},
% <4.4mm,-0.8mm>*{};<4.4mm,-3.5mm>*{}**@{-},
 <2.4mm,2.5mm>*{};<2.4mm,5.2mm>*{}**@{-},
     <0mm,-1.3mm>*{};<0mm,-5.3mm>*{^2}**@{},
     <2.5mm,2.3mm>*{};<5.1mm,-2.6mm>*{^1}**@{},
    <2.4mm,2.5mm>*{};<2.4mm,5.7mm>*{^\al}**@{},
    <-0.38mm,-0.2mm>*{};<-2.8mm,2.8mm>*{^\be}**@{},
    \end{xy},
$$
we first obtain,
$$
\begin{xy}
 <0.4mm,0.0mm>*{};<2.4mm,2.1mm>*{}**@{-},
 <-0.38mm,-0.2mm>*{};<-2.8mm,2.5mm>*{}**@{-},
<0mm,-0.8mm>*{\circ};<0mm,0.8mm>*{}**@{},
 <2.96mm,2.4mm>*{\circ};<2.45mm,2.35mm>*{}**@{},
 <2.4mm,2.8mm>*{};<0mm,5mm>*{}**@{-},
<-2.8mm,2.5mm>*{};<0mm,5mm>*{\bullet}**@{},
<-2.8mm,2.5mm>*{};<0mm,5mm>*{}**@{-},
<0mm,5mm>*{};<0mm,8.6mm>*{}**@{-},
<0mm,9mm>*{\circ};<0mm,0mm>*{}**@{},
<-0.5mm,9.6mm>*{};<-3.1mm,11.6mm>*{}**@{-},
<0.5mm,9.6mm>*{};<3.1mm,11.6mm>*{}**@{-},
<2.9mm,7.5mm>*{};<-3.1mm,12.4mm>*{^\al}**@{},
<2.9mm,7.5mm>*{};<3.1mm,12.4mm>*{^\be}**@{},
(3.2,3)*{}
   \ar@{->}@(ur,dr) (0.1,-1.4)*{}
    \end{xy}
    \ \ =\ \ 2\
\begin{xy}
 <0.4mm,0.0mm>*{};<2.4mm,2.1mm>*{}**@{-},
<0mm,-0.8mm>*{\circ};<0mm,0.8mm>*{}**@{},
 <2.96mm,2.4mm>*{\circ};<2.45mm,2.35mm>*{}**@{},
 <2.4mm,2.8mm>*{};<0mm,5mm>*{}**@{-},
<-2.8mm,2.5mm>*{};<0mm,5mm>*{\bullet}**@{},
<-2.4mm,2.8mm>*{};<0mm,5mm>*{}**@{-},
<0mm,5mm>*{};<0mm,8.6mm>*{}**@{-},
<-2.96mm,2.4mm>*{\circ};<0mm,0mm>*{}**@{},
 <-0.4mm,0.0mm>*{};<-2.4mm,2.1mm>*{}**@{-},
<-3.2mm,3mm>*{};<-5.1mm,5mm>*{}**@{-},
<2.9mm,7.5mm>*{};<-5.1mm,5.6mm>*{^\al}**@{},
<2.9mm,7.5mm>*{};<0mm,9.3mm>*{^\be}**@{},
(3.2,3)*{}
   \ar@{->}@(ur,dr) (0.1,-1.4)*{}
\end{xy}
    \ - \
2\
\begin{xy}
 <0.4mm,0.0mm>*{};<2.4mm,2.1mm>*{}**@{-},
<0mm,-0.8mm>*{\circ};<0mm,0.8mm>*{}**@{},
 <2.96mm,2.4mm>*{\circ};<2.45mm,2.35mm>*{}**@{},
 <2.4mm,2.8mm>*{};<0mm,5mm>*{}**@{-},
<-2.8mm,2.5mm>*{};<0mm,5mm>*{\bullet}**@{},
<-2.4mm,2.8mm>*{};<0mm,5mm>*{}**@{-},
<0mm,5mm>*{};<0mm,8.6mm>*{}**@{-},
<-2.96mm,2.4mm>*{\circ};<0mm,0mm>*{}**@{},
 <-0.4mm,0.0mm>*{};<-2.4mm,2.1mm>*{}**@{-},
<-3.2mm,3mm>*{};<-5.1mm,5mm>*{}**@{-},
<2.9mm,7.5mm>*{};<-5.1mm,5.6mm>*{^\be}**@{},
<2.9mm,7.5mm>*{};<0mm,9.3mm>*{^\al}**@{},
(3.2,3)*{}
   \ar@{->}@(ur,dr) (0.1,-1.4)*{}
\end{xy} .
$$
Next, using relations,
$$
\begin{xy}
 <0mm,0mm>*{\circ};<0mm,0mm>*{}**@{},
 <0mm,-0.49mm>*{};<0mm,-3.0mm>*{}**@{-},
 <0.49mm,0.49mm>*{};<1.9mm,1.9mm>*{}**@{-},
 <-0.5mm,0.5mm>*{};<-1.9mm,1.9mm>*{}**@{-},
 <-2.3mm,2.3mm>*{\circ};<-2.3mm,2.3mm>*{}**@{},
 <-1.8mm,2.8mm>*{};<0mm,4.9mm>*{}**@{-},
 <-2.8mm,2.9mm>*{};<-4.6mm,4.9mm>*{}**@{-},
   <0.49mm,0.49mm>*{};<2.7mm,2.3mm>*{^3}**@{},
   <-1.8mm,2.8mm>*{};<0.4mm,5.3mm>*{^2}**@{},
   <-2.8mm,2.9mm>*{};<-5.1mm,5.3mm>*{^1}**@{},
 \end{xy}
\ + \
\begin{xy}
 <0mm,0mm>*{\circ};<0mm,0mm>*{}**@{},
 <0mm,-0.49mm>*{};<0mm,-3.0mm>*{}**@{-},
 <0.49mm,0.49mm>*{};<1.9mm,1.9mm>*{}**@{-},
 <-0.5mm,0.5mm>*{};<-1.9mm,1.9mm>*{}**@{-},
 <-2.3mm,2.3mm>*{\circ};<-2.3mm,2.3mm>*{}**@{},
 <-1.8mm,2.8mm>*{};<0mm,4.9mm>*{}**@{-},
 <-2.8mm,2.9mm>*{};<-4.6mm,4.9mm>*{}**@{-},
   <0.49mm,0.49mm>*{};<2.7mm,2.3mm>*{^2}**@{},
   <-1.8mm,2.8mm>*{};<0.4mm,5.3mm>*{^1}**@{},
   <-2.8mm,2.9mm>*{};<-5.1mm,5.3mm>*{^3}**@{},
 \end{xy}
\ + \
\begin{xy}
 <0mm,0mm>*{\circ};<0mm,0mm>*{}**@{},
 <0mm,-0.49mm>*{};<0mm,-3.0mm>*{}**@{-},
 <0.49mm,0.49mm>*{};<1.9mm,1.9mm>*{}**@{-},
 <-0.5mm,0.5mm>*{};<-1.9mm,1.9mm>*{}**@{-},
 <-2.3mm,2.3mm>*{\circ};<-2.3mm,2.3mm>*{}**@{},
 <-1.8mm,2.8mm>*{};<0mm,4.9mm>*{}**@{-},
 <-2.8mm,2.9mm>*{};<-4.6mm,4.9mm>*{}**@{-},
   <0.49mm,0.49mm>*{};<2.7mm,2.3mm>*{^1}**@{},
   <-1.8mm,2.8mm>*{};<0.4mm,5.3mm>*{^3}**@{},
   <-2.8mm,2.9mm>*{};<-5.1mm,5.3mm>*{^2}**@{},
 \end{xy}
 \ = \ 0,
$$
we finally obtain the desired result,
$$
\hspace{48mm}
\begin{xy}
 <0.4mm,0.0mm>*{};<2.4mm,2.1mm>*{}**@{-},
<0mm,-0.8mm>*{\circ};<0mm,0.8mm>*{}**@{},
 <2.96mm,2.4mm>*{\circ};<2.45mm,2.35mm>*{}**@{},
 <2.4mm,2.8mm>*{};<0mm,5mm>*{}**@{-},
<-2.8mm,2.5mm>*{};<0mm,5mm>*{\bullet}**@{},
<-2.4mm,2.8mm>*{};<0mm,5mm>*{}**@{-},
<0mm,5mm>*{};<0mm,8.6mm>*{}**@{-},
<-2.96mm,2.4mm>*{\circ};<0mm,0mm>*{}**@{},
 <-0.4mm,0.0mm>*{};<-2.4mm,2.1mm>*{}**@{-},
<-3.2mm,3mm>*{};<-5.1mm,5mm>*{}**@{-},
<2.9mm,7.5mm>*{};<-5.1mm,5.6mm>*{^\al}**@{},
<2.9mm,7.5mm>*{};<0mm,9.3mm>*{^\be}**@{},
(3.2,3)*{}
   \ar@{->}@(ur,dr) (0.1,-1.4)*{}
\end{xy}
\ \ = \ \
\begin{xy}
 <0.4mm,0.0mm>*{};<2.4mm,2.1mm>*{}**@{-},
 <-0.38mm,-0.2mm>*{};<-2.8mm,2.5mm>*{}**@{-},
<0mm,-0.8mm>*{\circ};<0mm,0.8mm>*{}**@{},
 <2.96mm,2.4mm>*{\circ};<2.45mm,2.35mm>*{}**@{},
 <2.4mm,2.8mm>*{};<0mm,5mm>*{}**@{-},
<3.4mm,2.9mm>*{};<5.1mm,4.6mm>*{}**@{-},
<-2.8mm,2.5mm>*{};<5.5mm,5mm>*{\circ}**@{},
<5.1mm,5.7mm>*{};<3.1mm,7.7mm>*{}**@{-},
<2.9mm,7.5mm>*{};<3.1mm,8.1mm>*{^\al}**@{},
<-2.8mm,2.5mm>*{};<0mm,5mm>*{\bullet}**@{},
<-2.8mm,2.5mm>*{};<0mm,5mm>*{}**@{-},
<0mm,5mm>*{};<0mm,8.8mm>*{}**@{-},
<2.9mm,7.5mm>*{};<0mm,9.6mm>*{^\be}**@{},
(5.2,5)*{}
   \ar@{->}@(ur,dr) (0.1,-1.4)*{}
\end{xy}
\ - \
\begin{xy}
 <0.4mm,0.0mm>*{};<2.4mm,2.1mm>*{}**@{-},
 <-0.38mm,-0.2mm>*{};<-2.8mm,2.5mm>*{}**@{-},
<0mm,-0.8mm>*{\circ};<0mm,0.8mm>*{}**@{},
 <2.96mm,2.4mm>*{\circ};<2.45mm,2.35mm>*{}**@{},
 <2.4mm,2.8mm>*{};<0mm,5mm>*{}**@{-},
<3.4mm,2.9mm>*{};<5.1mm,4.6mm>*{}**@{-},
<-2.8mm,2.5mm>*{};<-3.2mm,-3.4mm>*{\circ}**@{},
<-2.7mm,-2.9mm>*{};<-0.5mm,-0.9mm>*{}**@{-},
<-3.7mm,-2.7mm>*{};<-5.9mm,0mm>*{}**@{-},
<2.9mm,7.5mm>*{};<-5.9mm,0.5mm>*{^\al}**@{},
<-2.8mm,2.5mm>*{};<0mm,5mm>*{\bullet}**@{},
<-2.8mm,2.5mm>*{};<0mm,5mm>*{}**@{-},
<0mm,5mm>*{};<0mm,8.8mm>*{}**@{-},
<2.9mm,7.5mm>*{};<0mm,9.6mm>*{^\be}**@{},
(3.2,3)*{}
   \ar@{->}@(ur,dr) (-3.2,-3.7)*{}
\end{xy}
\ =\ 0. \hspace{37mm} \Box
$$

\bip

%%%%%%%%%%%%%%%%%%%%%%%%%%%%%%%%%%%%%%%%%%%%%%%%%%%%%%%%%%
\noindent{\bf I.2.2. Corollary-definition}. {\em For any constant $a$ the
operator $\check{\xi}$ descends to a derivation of the star product
$\star_\hbar$ which we denote from now on by $\Ga_1$}.
\vspace{-6mm}

\Beqrn
\Proof \hspace{10mm} \breve{\xi}\left(x^\al\ot x^\be -
(-1)^{|\al||\be|}x^\be\ot x^\al - \hbar\sum_\ga \Phi^{\al\be}_\ga
x^\ga\right)&=& - a\hbar \sum_\ga \Phi_\ga^{\al\be} \Theta^\ga \bmod {\cI}\\
&=& 0  \bmod {\cI}.\hspace{30mm} \Box \Eeqrn

\bip

\noindent Let us next define a linear function,
$
\Gamma_0:= b \hbar^3\sum_\ga \Xi_\ga x^\ga
$,
where $b$ is a constant and
$$
\Xi_\ga:=\sum_{\al,\be,\mu,\nu} C^\mu_{\al\be}C^{\be}_{\mu\nu}
\Phi^{\nu\al}_\ga = \
\Ba{c}
 \begin{xy}
 <0mm,0mm>*{\bullet};<0mm,0mm>*{}**@{},
 %<0mm,0.69mm>*{};<0mm,3.0mm>*{}**@{-},
 <0.39mm,-0.39mm>*{};<2.4mm,-2.4mm>*{}**@{-},
 <-0.35mm,-0.35mm>*{};<-1.9mm,-1.9mm>*{}**@{-},
 <-2.4mm,-2.4mm>*{\bullet};<-2.4mm,-2.4mm>*{}**@{},
 <-2.0mm,-2.8mm>*{};<-0.4mm,-4.5mm>*{}**@{-},
 %<-2.8mm,-2.9mm>*{};<-4.7mm,-4.9mm>*{}**@{-},
 %
 <2.4mm,-2.4mm>*{};<0.4mm,-4.5mm>*{}**@{-},
  <0mm,-5.1mm>*{\circ};<0mm,0mm>*{}**@{},
  <0mm,-5.5mm>*{};<0mm,-7.7mm>*{}**@{-},
    <0.39mm,-0.39mm>*{};<0mm,-9.1mm>*{_\ga}**@{},
(0,0)*{}
   \ar@{->}@(ul,dl) (-2.4,-2.4)*{}
 \end{xy}
 \\
 \ \
\Ea
$$

%%%%%%%%%%%%%%%%%%%%%%%%%%
\noindent{\bf I.2.3. Lemma}. $C_{\al\be}^\ga \Xi_\ga =0$, i.e.\
$
\Ba{c}
 \begin{xy}
 <0mm,0mm>*{\bullet};<0mm,0mm>*{}**@{},
 %<0mm,0.69mm>*{};<0mm,3.0mm>*{}**@{-},
 <0.39mm,-0.39mm>*{};<2.4mm,-2.4mm>*{}**@{-},
 <-0.35mm,-0.35mm>*{};<-1.9mm,-1.9mm>*{}**@{-},
 <-2.4mm,-2.4mm>*{\bullet};<-2.4mm,-2.4mm>*{}**@{},
 <-2.0mm,-2.8mm>*{};<-0.4mm,-4.5mm>*{}**@{-},
 %<-2.8mm,-2.9mm>*{};<-4.7mm,-4.9mm>*{}**@{-},
 %
 <2.4mm,-2.4mm>*{};<0.4mm,-4.5mm>*{}**@{-},
  <0mm,-5.1mm>*{\circ};<0mm,0mm>*{}**@{},
  <0mm,-5.5mm>*{};<0mm,-7.7mm>*{}**@{-},
<0mm,-8.2mm>*{\bullet};<0mm,0mm>*{}**@{},
 <0mm,-8.2mm>*{};<-2.4mm,-10.1mm>*{}**@{-},
 <0mm,-8.2mm>*{};<2.4mm,-10.1mm>*{}**@{-},
<0.39mm,-0.39mm>*{};<2.5mm,-11.7mm>*{_\be}**@{},
<0.39mm,-0.39mm>*{};<-2.5mm,-11.2mm>*{_\al}**@{},
(0,0)*{}
   \ar@{->}@(ul,dl) (-2.4,-2.4)*{}
 \end{xy}
 \\
 \ \
\Ea \ = \ 0.
$

\noindent Proof is very similar to the proof of Lemma~I.2.1. We omit
details.

\bip

%%%%%%%%%%%%%%%%%%%%%%%%%%
\noindent{\bf I.2.4. Corollary}. $\Ga_1\Ga_0=0$.

\bip

Therefore, for any Lie 1-bialgebra, $(C_{\al\be}^\ga,
\Phi^{\mu\nu}_\var)$, we constructed differential operators,
$\Ga_0\in  \Hom_2(\R,\f_V)[[\hbar]]$, $\Gamma_1\in
\Hom_1(\f_V,\f_V)[[\hbar]]$ and $\Gamma_2\in
\Hom_0(\f_V^{\ot 2},\f_V)[[\hbar]]$ such that the equations,
$\Ga_1\Ga_0=0$ , $d_{\rm H}\Ga_1 + [\Ga_1, \Ga_2]_{\rm H}=0$ and
$d_{\rm H}\Ga_2+ \frac{1}{2} [\Ga_2, \Ga_2]_{\rm H}=0$, are
satisfied. It remains to check whether or not we can adjust the free
parameters $a$ and $b$ in such a way that the last equation,
$$
\Ga_1^2 +[\Ga_0,\Ga_2]_{\rm H}=0,
$$
holds. As the l.h.s.\ of this equation is obviously a derivation of
the star product, it is enough to check the latter only on
generators $x^\al$,
$$
\Ga_1^2(x^\al) + \Ga_0\star_\hbar x^\al  - x^\al\star_\hbar \Ga_0=0.
$$
We have
\Beqrn
\hbar^{-4}(\Ga_0\star_\hbar x^\al  - x^\al\star_\hbar
\Ga_0) &=& b\sum\Xi_\ga\Phi^{\ga\al}_\be
x^\be\\
&\simeq& b \Ba{c}
 \begin{xy}
 <0mm,0mm>*{\bullet};<0mm,0mm>*{}**@{},
 %<0mm,0.69mm>*{};<0mm,3.0mm>*{}**@{-},
 <0.39mm,-0.39mm>*{};<2.4mm,-2.4mm>*{}**@{-},
 <-0.35mm,-0.35mm>*{};<-1.9mm,-1.9mm>*{}**@{-},
 <-2.4mm,-2.4mm>*{\bullet};<-2.4mm,-2.4mm>*{}**@{},
 <-2.0mm,-2.8mm>*{};<-0.4mm,-4.5mm>*{}**@{-},
 %<-2.8mm,-2.9mm>*{};<-4.7mm,-4.9mm>*{}**@{-},
 %
 <2.4mm,-2.4mm>*{};<0.4mm,-4.5mm>*{}**@{-},
  <0mm,-5.1mm>*{\circ};<0mm,0mm>*{}**@{},
  <0.3mm,-5.5mm>*{};<2.0mm,-7.5mm>*{}**@{-},
<2.4mm,-8.3mm>*{\circ};<0mm,0mm>*{}**@{},
<2.4mm,-8.8mm>*{};<2.4mm,-11.2mm>*{}**@{-},
<2.9mm,-7.9mm>*{};<5.1mm,-6mm>*{}**@{-},
<0.39mm,-0.39mm>*{};<5.5mm,-5.5mm>*{^\al}**@{},
<0.39mm,-0.39mm>*{};<2.4mm,-12.9mm>*{_\be}**@{},
(0,0)*{}
   \ar@{->}@(ul,dl) (-2.4,-2.4)*{}
 \end{xy}
 \Ea \
= \ 2b\
\begin{xy}
 <0.4mm,0.0mm>*{};<2.4mm,2.1mm>*{}**@{-},
<0mm,-0.8mm>*{\circ};<0mm,0.8mm>*{}**@{},
 <2.96mm,2.4mm>*{\circ};<2.45mm,2.35mm>*{}**@{},
 <2.4mm,2.8mm>*{};<0mm,5mm>*{}**@{-},
<-2.8mm,2.5mm>*{};<0mm,5mm>*{\bullet}**@{},
<-2.4mm,2.8mm>*{};<0mm,5mm>*{}**@{-},
<0mm,-1.2mm>*{};<0mm,-4mm>*{}**@{-},
<-2.96mm,2.4mm>*{\bullet};<0mm,0mm>*{}**@{},
 <-0.4mm,0.0mm>*{};<-2.4mm,2.1mm>*{}**@{-},
<3.4mm,3mm>*{};<5.3mm,5mm>*{}**@{-},
<0mm,0mm>*{};<0mm,-6.5mm>*{^\be}**@{},
<0mm,0mm>*{};<5.2mm,5.3mm>*{^\al}**@{},
(0,5)*{}
   \ar@{->}@(ul,dl) (-2.9,2.4)*{}
\end{xy}
\Eeqrn
We also have,
\Beqrn
\hbar^{-4}\Ga_1^2(x^\al)&=&
-\frac{1}{3}\ \
\begin{xy}
 <0mm,-1.3mm>*{};<0mm,-3.5mm>*{}**@{-},
 <0.4mm,0.0mm>*{};<2.4mm,2.1mm>*{}**@{-},
 <-0.38mm,-0.2mm>*{};<-2.8mm,2.5mm>*{}**@{-},
<0mm,-0.8mm>*{\circ};<0mm,0.8mm>*{}**@{},
 <2.96mm,2.4mm>*{\circ};<2.45mm,2.35mm>*{}**@{},
 <2.4mm,2.8mm>*{};<0mm,5mm>*{}**@{-},
 <3.35mm,2.9mm>*{};<5.5mm,5mm>*{}**@{-},
     <0mm,-1.3mm>*{};<0mm,-5.4mm>*{_\be}**@{},
     %<2.5mm,2.3mm>*{};<5.1mm,-2.6mm>*{^1}**@{},
%
<-2.8mm,2.5mm>*{};<0mm,5mm>*{\bullet}**@{},
<-2.8mm,2.5mm>*{};<0mm,5mm>*{}**@{-},
<2.96mm,5mm>*{};<2.96mm,7.5mm>*{\bullet}**@{},
<0.2mm,5.1mm>*{};<2.8mm,7.5mm>*{}**@{-},
<0.2mm,5.1mm>*{};<2.8mm,7.5mm>*{}**@{-},
<5.5mm,5mm>*{};<2.8mm,7.5mm>*{}**@{-},
<2.9mm,7.5mm>*{};<2.9mm,10.5mm>*{}**@{-},
<2.9mm,7.5mm>*{};<2.9mm,11.1mm>*{^\al}**@{},
    \end{xy}
\ \ - \ 2a\ \,
\begin{xy}
 <0.4mm,0.0mm>*{};<2.4mm,2.1mm>*{}**@{-},
 <-0.38mm,-0.2mm>*{};<-2.8mm,2.5mm>*{}**@{-},
<0mm,-0.8mm>*{\circ};<0mm,0.8mm>*{}**@{},
 <2.96mm,2.4mm>*{\circ};<2.45mm,2.35mm>*{}**@{},
 <2.4mm,2.8mm>*{};<0mm,5mm>*{}**@{-},
<-2.8mm,2.5mm>*{};<0mm,5mm>*{\bullet}**@{},
<-2.8mm,2.5mm>*{};<0mm,5mm>*{}**@{-},
<0mm,5mm>*{};<-2.96mm,8mm>*{}**@{-},
<0mm,0mm>*{};<-2.96mm,11.2mm>*{^\al}**@{},
<0mm,0mm>*{};<-6.1mm,2.8mm>*{^\be}**@{},
 <-2.96mm,8mm>*{\bullet};<0mm,0mm>*{}**@{},
<-2.96mm,8mm>*{};<-2.96mm,11mm>*{}**@{-},
<-2.96mm,8mm>*{};<-6.1mm,5mm>*{}**@{-},
(3.2,3)*{}
   \ar@{->}@(ur,dr) (0.1,-1.4)*{}
    \end{xy}
    \
=\ -\frac{1}{3}\ \
\begin{xy}
 <0mm,-1.3mm>*{};<0mm,-3.5mm>*{}**@{-},
 <0.4mm,0.0mm>*{};<2.4mm,2.1mm>*{}**@{-},
 <-0.38mm,-0.2mm>*{};<-2.8mm,2.5mm>*{}**@{-},
<0mm,-0.8mm>*{\circ};<0mm,0.8mm>*{}**@{},
 <2.96mm,2.4mm>*{\circ};<2.45mm,2.35mm>*{}**@{},
 <2.4mm,2.8mm>*{};<0mm,5mm>*{}**@{-},
 <3.35mm,2.9mm>*{};<5.5mm,5mm>*{}**@{-},
     <0mm,-1.3mm>*{};<0mm,-5.4mm>*{_\be}**@{},
     %<2.5mm,2.3mm>*{};<5.1mm,-2.6mm>*{^1}**@{},
%
<-2.8mm,2.5mm>*{};<0mm,5mm>*{\bullet}**@{},
<-2.8mm,2.5mm>*{};<0mm,5mm>*{}**@{-},
<2.96mm,5mm>*{};<2.96mm,7.5mm>*{\bullet}**@{},
<0.2mm,5.1mm>*{};<2.8mm,7.5mm>*{}**@{-},
<0.2mm,5.1mm>*{};<2.8mm,7.5mm>*{}**@{-},
<5.5mm,5mm>*{};<2.8mm,7.5mm>*{}**@{-},
<2.9mm,7.5mm>*{};<2.9mm,10.5mm>*{}**@{-},
<2.9mm,7.5mm>*{};<2.9mm,11.1mm>*{^\al}**@{},
    \end{xy}
\ \ + \ 4a\ \,
\begin{xy}
 <0.4mm,0.0mm>*{};<2.4mm,2.1mm>*{}**@{-},
<0mm,-0.8mm>*{\circ};<0mm,0.8mm>*{}**@{},
 <2.96mm,2.4mm>*{\circ};<2.45mm,2.35mm>*{}**@{},
 <2.4mm,2.8mm>*{};<0mm,5mm>*{}**@{-},
<-2.8mm,2.5mm>*{};<0mm,5mm>*{\bullet}**@{},
<-2.4mm,2.8mm>*{};<0mm,5mm>*{}**@{-},
<0mm,5mm>*{};<0mm,8.6mm>*{}**@{-},
<-2.96mm,2.4mm>*{\bullet};<0mm,0mm>*{}**@{},
 <-0.4mm,0.0mm>*{};<-2.4mm,2.1mm>*{}**@{-},
<-3.2mm,2.4mm>*{};<-5.7mm,0mm>*{}**@{-},
<2.9mm,7.5mm>*{};<-6.1mm,-2.2mm>*{^\be}**@{},
<2.9mm,7.5mm>*{};<0mm,9.3mm>*{^\al}**@{},
(3.2,3)*{}
   \ar@{->}@(ur,dr) (0.1,-1.4)*{}
\end{xy} \ \ .
    \\
\Eeqrn
Consider a graph, \ \ $
\begin{xy}
 <0mm,2.47mm>*{};<0mm,-0.5mm>*{}**@{-},
 <0.5mm,3.5mm>*{};<2.2mm,5.2mm>*{}**@{-},
 <-0.48mm,3.48mm>*{};<-2.2mm,5.2mm>*{}**@{-},
 <0mm,3mm>*{\circ};<0mm,3mm>*{}**@{},
  <0mm,-0.8mm>*{\bullet};<0mm,-0.8mm>*{}**@{},
<0mm,-0.8mm>*{};<-2.2mm,-3.5mm>*{}**@{-},
   <0mm,0mm>*{};<-2.5mm,9.7mm>*{^\al}**@{},
 <-2.5mm,5.7mm>*{\bullet};<0mm,0mm>*{}**@{},
<-2.5mm,5.7mm>*{};<-2.5mm,9.4mm>*{}**@{-},
<-2.5mm,5.7mm>*{};<-5mm,3mm>*{}**@{-},
<-5mm,3mm>*{};<-5mm,-0.8mm>*{}**@{-},
 <-2.5mm,-4.2mm>*{\circ};<0mm,3mm>*{}**@{},
 <-2.8mm,-3.6mm>*{};<-5mm,-0.8mm>*{}**@{-},
 <-2.5mm,-4.6mm>*{};<-2.5mm,-7.3mm>*{}**@{-},
  <0mm,0mm>*{};<-2.5mm,-8.9mm>*{_\be}**@{},
%   %
 (0.4,3.6)*{}
   \ar@{->}@(ur,dr) (0.1,-0.6)*{}
\end{xy}
$.
Replacing its  lower two vertices together with their half-edges by
the following linear combination of four graphs,
$$
\begin{xy}
 <0mm,-1.3mm>*{};<0mm,-3.5mm>*{}**@{-},
 <0.38mm,-0.2mm>*{};<2.2mm,2.2mm>*{}**@{-},
 <-0.38mm,-0.2mm>*{};<-2.2mm,2.2mm>*{}**@{-},
<0mm,-0.8mm>*{\circ};<0mm,0.8mm>*{}**@{},
% <-2.25mm,2.2mm>*{};<-2.2mm,5.2mm>*{}**@{-},
 <2.4mm,2.4mm>*{\bullet};<2.4mm,2.4mm>*{}**@{},
 <2.5mm,2.3mm>*{};<4.4mm,-0.8mm>*{}**@{-},
% <4.4mm,-0.8mm>*{};<4.4mm,-3.5mm>*{}**@{-},
 <2.4mm,2.5mm>*{};<2.4mm,5.2mm>*{}**@{-},
     <0mm,-1.3mm>*{};<0mm,-5.5mm>*{^\be}**@{},
     %<2.5mm,2.3mm>*{};<5.1mm,-2.6mm>*{^2}**@{},
    <2.4mm,2.5mm>*{};<2.4mm,5.7mm>*{^2}**@{},
    <-0.38mm,-0.2mm>*{};<-2.8mm,2.8mm>*{^1}**@{},
    \end{xy}
\ = \ -\
 \begin{xy}
 <0mm,2.47mm>*{};<0mm,-0.5mm>*{}**@{-},
 <0.5mm,3.5mm>*{};<2.2mm,5.2mm>*{}**@{-},
 <-0.48mm,3.48mm>*{};<-2.2mm,5.2mm>*{}**@{-},
 <0mm,3mm>*{\circ};<0mm,3mm>*{}**@{},
  <0mm,-0.8mm>*{\bullet};<0mm,-0.8mm>*{}**@{},
<0mm,-0.8mm>*{};<-2.2mm,-3.5mm>*{}**@{-},
 <0mm,-0.8mm>*{};<2.2mm,-3.5mm>*{}**@{-},
     <0.5mm,3.5mm>*{};<2.8mm,5.9mm>*{^1}**@{},
     <-0.48mm,3.48mm>*{};<-2.8mm,5.7mm>*{^2}**@{},
   <0mm,-0.8mm>*{};<-2.7mm,-5.5mm>*{^\be}**@{},
   %<0mm,-0.8mm>*{};<2.7mm,-5.2mm>*{}**@{},
\end{xy}
\  + \
\begin{xy}
 <0mm,-1.3mm>*{};<0mm,-3.5mm>*{}**@{-},
 <0.38mm,-0.2mm>*{};<2.2mm,2.2mm>*{}**@{-},
 <-0.38mm,-0.2mm>*{};<-2.2mm,2.2mm>*{}**@{-},
<0mm,-0.8mm>*{\circ};<0mm,0.8mm>*{}**@{},
% <-2.25mm,2.2mm>*{};<-2.2mm,5.2mm>*{}**@{-},
 <2.4mm,2.4mm>*{\bullet};<2.4mm,2.4mm>*{}**@{},
 <2.5mm,2.3mm>*{};<4.4mm,-0.8mm>*{}**@{-},
% <4.4mm,-0.8mm>*{};<4.4mm,-3.5mm>*{}**@{-},
 <2.4mm,2.5mm>*{};<2.4mm,5.2mm>*{}**@{-},
     <0mm,-1.3mm>*{};<0mm,-5.5mm>*{^\be}**@{},
     %<2.5mm,2.3mm>*{};<5.1mm,-2.6mm>*{^2}**@{},
    <2.4mm,2.5mm>*{};<2.4mm,5.9mm>*{^1}**@{},
    <-0.38mm,-0.2mm>*{};<-2.8mm,2.5mm>*{^2}**@{},
    \end{xy}
\  + \
\begin{xy}
 <0mm,-1.3mm>*{};<0mm,-3.5mm>*{}**@{-},
 <0.38mm,-0.2mm>*{};<2.2mm,2.2mm>*{}**@{-},
 <-0.38mm,-0.2mm>*{};<-2.2mm,2.2mm>*{}**@{-},
<0mm,-0.8mm>*{\circ};<0mm,0.8mm>*{}**@{},
% <-2.25mm,2.2mm>*{};<-2.2mm,5.2mm>*{}**@{-},
 <2.4mm,2.4mm>*{\bullet};<2.4mm,2.4mm>*{}**@{},
 <2.5mm,2.3mm>*{};<4.4mm,-0.8mm>*{}**@{-},
% <4.4mm,-0.8mm>*{};<4.4mm,-3.5mm>*{}**@{-},
 <2.4mm,2.5mm>*{};<2.4mm,5.2mm>*{}**@{-},
     %<0mm,-1.3mm>*{};<0mm,-5.3mm>*{^2}**@{},
     <2.5mm,2.3mm>*{};<5.1mm,-2.8mm>*{^\be}**@{},
    <2.4mm,2.5mm>*{};<2.4mm,5.9mm>*{^1}**@{},
    <-0.38mm,-0.2mm>*{};<-2.8mm,2.5mm>*{^2}**@{},
    \end{xy}
\  - \
\begin{xy}
 <0mm,-1.3mm>*{};<0mm,-3.5mm>*{}**@{-},
 <0.38mm,-0.2mm>*{};<2.2mm,2.2mm>*{}**@{-},
 <-0.38mm,-0.2mm>*{};<-2.2mm,2.2mm>*{}**@{-},
<0mm,-0.8mm>*{\circ};<0mm,0.8mm>*{}**@{},
% <-2.25mm,2.2mm>*{};<-2.2mm,5.2mm>*{}**@{-},
 <2.4mm,2.4mm>*{\bullet};<2.4mm,2.4mm>*{}**@{},
 <2.5mm,2.3mm>*{};<4.4mm,-0.8mm>*{}**@{-},
% <4.4mm,-0.8mm>*{};<4.4mm,-3.5mm>*{}**@{-},
 <2.4mm,2.5mm>*{};<2.4mm,5.2mm>*{}**@{-},
     %<0mm,-1.3mm>*{};<0mm,-5.3mm>*{^2}**@{},
     <2.5mm,2.3mm>*{};<5.1mm,-2.6mm>*{^\be}**@{},
    <2.4mm,2.5mm>*{};<2.4mm,5.7mm>*{^2}**@{},
    <-0.38mm,-0.2mm>*{};<-2.8mm,2.8mm>*{^1}**@{},
    \end{xy},
$$
one gets, after a cancelation of two terms, the identity,
$$
\begin{xy}
 <0mm,2.47mm>*{};<0mm,-0.5mm>*{}**@{-},
 <0.5mm,3.5mm>*{};<2.2mm,5.2mm>*{}**@{-},
 <-0.48mm,3.48mm>*{};<-2.2mm,5.2mm>*{}**@{-},
 <0mm,3mm>*{\circ};<0mm,3mm>*{}**@{},
  <0mm,-0.8mm>*{\bullet};<0mm,-0.8mm>*{}**@{},
<0mm,-0.8mm>*{};<-2.2mm,-3.5mm>*{}**@{-},
   <0mm,0mm>*{};<-2.5mm,9.7mm>*{^\al}**@{},
 <-2.5mm,5.7mm>*{\bullet};<0mm,0mm>*{}**@{},
<-2.5mm,5.7mm>*{};<-2.5mm,9.4mm>*{}**@{-},
<-2.5mm,5.7mm>*{};<-5mm,3mm>*{}**@{-},
<-5mm,3mm>*{};<-5mm,-0.8mm>*{}**@{-},
 <-2.5mm,-4.2mm>*{\circ};<0mm,3mm>*{}**@{},
 <-2.8mm,-3.6mm>*{};<-5mm,-0.8mm>*{}**@{-},
 <-2.5mm,-4.6mm>*{};<-2.5mm,-7.3mm>*{}**@{-},
  <0mm,0mm>*{};<-2.5mm,-8.9mm>*{_\be}**@{},
%   %
 (0.4,3.6)*{}
   \ar@{->}@(ur,dr) (0.1,-0.6)*{}
\end{xy}
\ = \ \
\begin{xy}
 <0mm,-1.3mm>*{};<0mm,-3.5mm>*{}**@{-},
 <0.4mm,0.0mm>*{};<2.4mm,2.1mm>*{}**@{-},
 <-0.38mm,-0.2mm>*{};<-2.8mm,2.5mm>*{}**@{-},
<0mm,-0.8mm>*{\circ};<0mm,0.8mm>*{}**@{},
 <2.96mm,2.4mm>*{\circ};<2.45mm,2.35mm>*{}**@{},
 <2.4mm,2.8mm>*{};<0mm,5mm>*{}**@{-},
 <3.35mm,2.9mm>*{};<5.5mm,5mm>*{}**@{-},
     <0mm,-1.3mm>*{};<0mm,-5.4mm>*{_\be}**@{},
     %<2.5mm,2.3mm>*{};<5.1mm,-2.6mm>*{^1}**@{},
%
<-2.8mm,2.5mm>*{};<0mm,5mm>*{\bullet}**@{},
<-2.8mm,2.5mm>*{};<0mm,5mm>*{}**@{-},
<2.96mm,5mm>*{};<2.96mm,7.5mm>*{\bullet}**@{},
<0.2mm,5.1mm>*{};<2.8mm,7.5mm>*{}**@{-},
<0.2mm,5.1mm>*{};<2.8mm,7.5mm>*{}**@{-},
<5.5mm,5mm>*{};<2.8mm,7.5mm>*{}**@{-},
<2.9mm,7.5mm>*{};<2.9mm,10.5mm>*{}**@{-},
<2.9mm,7.5mm>*{};<2.9mm,11.1mm>*{^\al}**@{},
    \end{xy}
    \ - \
\begin{xy}
 <0.4mm,0.0mm>*{};<2.4mm,2.1mm>*{}**@{-},
<0mm,-0.8mm>*{\circ};<0mm,0.8mm>*{}**@{},
 <2.96mm,2.4mm>*{\circ};<2.45mm,2.35mm>*{}**@{},
 <2.4mm,2.8mm>*{};<0mm,5mm>*{}**@{-},
<-2.8mm,2.5mm>*{};<0mm,5mm>*{\bullet}**@{},
<-2.4mm,2.8mm>*{};<0mm,5mm>*{}**@{-},
<0mm,5mm>*{};<0mm,8.6mm>*{}**@{-},
<-2.96mm,2.4mm>*{\bullet};<0mm,0mm>*{}**@{},
 <-0.4mm,0.0mm>*{};<-2.4mm,2.1mm>*{}**@{-},
<-3.2mm,2.4mm>*{};<-5.7mm,0mm>*{}**@{-},
<2.9mm,7.5mm>*{};<-6.1mm,-2.2mm>*{^\be}**@{},
<2.9mm,7.5mm>*{};<0mm,9.3mm>*{^\al}**@{},
(3.2,3)*{}
   \ar@{->}@(ur,dr) (0.1,-1.4)*{}
\end{xy}
$$
Playing a similar trick with upper two vertices one gets another
identity,
$$
\begin{xy}
 <0mm,2.47mm>*{};<0mm,-0.5mm>*{}**@{-},
 <0.5mm,3.5mm>*{};<2.2mm,5.2mm>*{}**@{-},
 <-0.48mm,3.48mm>*{};<-2.2mm,5.2mm>*{}**@{-},
 <0mm,3mm>*{\circ};<0mm,3mm>*{}**@{},
  <0mm,-0.8mm>*{\bullet};<0mm,-0.8mm>*{}**@{},
<0mm,-0.8mm>*{};<-2.2mm,-3.5mm>*{}**@{-},
   <0mm,0mm>*{};<-2.5mm,9.7mm>*{^\al}**@{},
 <-2.5mm,5.7mm>*{\bullet};<0mm,0mm>*{}**@{},
<-2.5mm,5.7mm>*{};<-2.5mm,9.4mm>*{}**@{-},
<-2.5mm,5.7mm>*{};<-5mm,3mm>*{}**@{-},
<-5mm,3mm>*{};<-5mm,-0.8mm>*{}**@{-},
 <-2.5mm,-4.2mm>*{\circ};<0mm,3mm>*{}**@{},
 <-2.8mm,-3.6mm>*{};<-5mm,-0.8mm>*{}**@{-},
 <-2.5mm,-4.6mm>*{};<-2.5mm,-7.3mm>*{}**@{-},
  <0mm,0mm>*{};<-2.5mm,-8.9mm>*{_\be}**@{},
%   %
 (0.4,3.6)*{}
   \ar@{->}@(ur,dr) (0.1,-0.6)*{}
\end{xy}
\ = \ \ - \
\begin{xy}
 <0mm,-1.3mm>*{};<0mm,-3.5mm>*{}**@{-},
 <0.4mm,0.0mm>*{};<2.4mm,2.1mm>*{}**@{-},
 <-0.38mm,-0.2mm>*{};<-2.8mm,2.5mm>*{}**@{-},
<0mm,-0.8mm>*{\circ};<0mm,0.8mm>*{}**@{},
 <2.96mm,2.4mm>*{\circ};<2.45mm,2.35mm>*{}**@{},
 <2.4mm,2.8mm>*{};<0mm,5mm>*{}**@{-},
 <3.35mm,2.9mm>*{};<5.5mm,5mm>*{}**@{-},
     <0mm,-1.3mm>*{};<0mm,-5.4mm>*{_\be}**@{},
     %<2.5mm,2.3mm>*{};<5.1mm,-2.6mm>*{^1}**@{},
%
<-2.8mm,2.5mm>*{};<0mm,5mm>*{\bullet}**@{},
<-2.8mm,2.5mm>*{};<0mm,5mm>*{}**@{-},
<2.96mm,5mm>*{};<2.96mm,7.5mm>*{\bullet}**@{},
<0.2mm,5.1mm>*{};<2.8mm,7.5mm>*{}**@{-},
<0.2mm,5.1mm>*{};<2.8mm,7.5mm>*{}**@{-},
<5.5mm,5mm>*{};<2.8mm,7.5mm>*{}**@{-},
<2.9mm,7.5mm>*{};<2.9mm,10.5mm>*{}**@{-},
<2.9mm,7.5mm>*{};<2.9mm,11.1mm>*{^\al}**@{},
    \end{xy}
    \ - \
\begin{xy}
 <0.4mm,0.0mm>*{};<2.4mm,2.1mm>*{}**@{-},
<0mm,-0.8mm>*{\circ};<0mm,0.8mm>*{}**@{},
 <2.96mm,2.4mm>*{\circ};<2.45mm,2.35mm>*{}**@{},
 <2.4mm,2.8mm>*{};<0mm,5mm>*{}**@{-},
<-2.8mm,2.5mm>*{};<0mm,5mm>*{\bullet}**@{},
<-2.4mm,2.8mm>*{};<0mm,5mm>*{}**@{-},
<0mm,-1.2mm>*{};<0mm,-4mm>*{}**@{-},
<-2.96mm,2.4mm>*{\bullet};<0mm,0mm>*{}**@{},
 <-0.4mm,0.0mm>*{};<-2.4mm,2.1mm>*{}**@{-},
<3.4mm,3mm>*{};<5.3mm,5mm>*{}**@{-},
<0mm,0mm>*{};<0mm,-6.5mm>*{^\be}**@{},
<0mm,0mm>*{};<5.2mm,5.3mm>*{^\al}**@{},
(0,5)*{}
   \ar@{->}@(ul,dl) (-2.9,2.4)*{}
\end{xy}
$$
These both identities imply
$$
\begin{xy}
 <0mm,-1.3mm>*{};<0mm,-3.5mm>*{}**@{-},
 <0.4mm,0.0mm>*{};<2.4mm,2.1mm>*{}**@{-},
 <-0.38mm,-0.2mm>*{};<-2.8mm,2.5mm>*{}**@{-},
<0mm,-0.8mm>*{\circ};<0mm,0.8mm>*{}**@{},
 <2.96mm,2.4mm>*{\circ};<2.45mm,2.35mm>*{}**@{},
 <2.4mm,2.8mm>*{};<0mm,5mm>*{}**@{-},
 <3.35mm,2.9mm>*{};<5.5mm,5mm>*{}**@{-},
     <0mm,-1.3mm>*{};<0mm,-5.4mm>*{_\be}**@{},
     %<2.5mm,2.3mm>*{};<5.1mm,-2.6mm>*{^1}**@{},
%
<-2.8mm,2.5mm>*{};<0mm,5mm>*{\bullet}**@{},
<-2.8mm,2.5mm>*{};<0mm,5mm>*{}**@{-},
<2.96mm,5mm>*{};<2.96mm,7.5mm>*{\bullet}**@{},
<0.2mm,5.1mm>*{};<2.8mm,7.5mm>*{}**@{-},
<0.2mm,5.1mm>*{};<2.8mm,7.5mm>*{}**@{-},
<5.5mm,5mm>*{};<2.8mm,7.5mm>*{}**@{-},
<2.9mm,7.5mm>*{};<2.9mm,10.5mm>*{}**@{-},
<2.9mm,7.5mm>*{};<2.9mm,11.1mm>*{^\al}**@{},
    \end{xy}
\ = \ \frac{1}{2}\ \,
\begin{xy}
 <0.4mm,0.0mm>*{};<2.4mm,2.1mm>*{}**@{-},
<0mm,-0.8mm>*{\circ};<0mm,0.8mm>*{}**@{},
 <2.96mm,2.4mm>*{\circ};<2.45mm,2.35mm>*{}**@{},
 <2.4mm,2.8mm>*{};<0mm,5mm>*{}**@{-},
<-2.8mm,2.5mm>*{};<0mm,5mm>*{\bullet}**@{},
<-2.4mm,2.8mm>*{};<0mm,5mm>*{}**@{-},
<0mm,5mm>*{};<0mm,8.6mm>*{}**@{-},
<-2.96mm,2.4mm>*{\bullet};<0mm,0mm>*{}**@{},
 <-0.4mm,0.0mm>*{};<-2.4mm,2.1mm>*{}**@{-},
<-3.2mm,2.4mm>*{};<-5.7mm,0mm>*{}**@{-},
<2.9mm,7.5mm>*{};<-6.1mm,-2.2mm>*{^\be}**@{},
<2.9mm,7.5mm>*{};<0mm,9.3mm>*{^\al}**@{},
(3.2,3)*{}
   \ar@{->}@(ur,dr) (0.1,-1.4)*{}
\end{xy}
\ - \ \frac{1}{2}
\begin{xy}
 <0.4mm,0.0mm>*{};<2.4mm,2.1mm>*{}**@{-},
<0mm,-0.8mm>*{\circ};<0mm,0.8mm>*{}**@{},
 <2.96mm,2.4mm>*{\circ};<2.45mm,2.35mm>*{}**@{},
 <2.4mm,2.8mm>*{};<0mm,5mm>*{}**@{-},
<-2.8mm,2.5mm>*{};<0mm,5mm>*{\bullet}**@{},
<-2.4mm,2.8mm>*{};<0mm,5mm>*{}**@{-},
<0mm,-1.2mm>*{};<0mm,-4mm>*{}**@{-},
<-2.96mm,2.4mm>*{\bullet};<0mm,0mm>*{}**@{},
 <-0.4mm,0.0mm>*{};<-2.4mm,2.1mm>*{}**@{-},
<3.4mm,3mm>*{};<5.3mm,5mm>*{}**@{-},
<0mm,0mm>*{};<0mm,-6.5mm>*{^\be}**@{},
<0mm,0mm>*{};<5.2mm,5.3mm>*{^\al}**@{},
(0,5)*{}
   \ar@{->}@(ul,dl) (-2.9,2.4)*{}
\end{xy}
$$
which in turn implies, \Beqrn \hbar^{-4}\left(\Ga_1^2(x^\al) +
\Ga_0\star_\hbar x^\al  - x^\al\star_\hbar \Ga_0\right) &=&
-\frac{1}{3}\ \
\begin{xy}
 <0mm,-1.3mm>*{};<0mm,-3.5mm>*{}**@{-},
 <0.4mm,0.0mm>*{};<2.4mm,2.1mm>*{}**@{-},
 <-0.38mm,-0.2mm>*{};<-2.8mm,2.5mm>*{}**@{-},
<0mm,-0.8mm>*{\circ};<0mm,0.8mm>*{}**@{},
 <2.96mm,2.4mm>*{\circ};<2.45mm,2.35mm>*{}**@{},
 <2.4mm,2.8mm>*{};<0mm,5mm>*{}**@{-},
 <3.35mm,2.9mm>*{};<5.5mm,5mm>*{}**@{-},
     <0mm,-1.3mm>*{};<0mm,-5.4mm>*{_\be}**@{},
     %<2.5mm,2.3mm>*{};<5.1mm,-2.6mm>*{^1}**@{},
%
<-2.8mm,2.5mm>*{};<0mm,5mm>*{\bullet}**@{},
<-2.8mm,2.5mm>*{};<0mm,5mm>*{}**@{-},
<2.96mm,5mm>*{};<2.96mm,7.5mm>*{\bullet}**@{},
<0.2mm,5.1mm>*{};<2.8mm,7.5mm>*{}**@{-},
<0.2mm,5.1mm>*{};<2.8mm,7.5mm>*{}**@{-},
<5.5mm,5mm>*{};<2.8mm,7.5mm>*{}**@{-},
<2.9mm,7.5mm>*{};<2.9mm,10.5mm>*{}**@{-},
<2.9mm,7.5mm>*{};<2.9mm,11.1mm>*{^\al}**@{},
    \end{xy}
\ \ + \ 4a\ \,
\begin{xy}
 <0.4mm,0.0mm>*{};<2.4mm,2.1mm>*{}**@{-},
<0mm,-0.8mm>*{\circ};<0mm,0.8mm>*{}**@{},
 <2.96mm,2.4mm>*{\circ};<2.45mm,2.35mm>*{}**@{},
 <2.4mm,2.8mm>*{};<0mm,5mm>*{}**@{-},
<-2.8mm,2.5mm>*{};<0mm,5mm>*{\bullet}**@{},
<-2.4mm,2.8mm>*{};<0mm,5mm>*{}**@{-},
<0mm,5mm>*{};<0mm,8.6mm>*{}**@{-},
<-2.96mm,2.4mm>*{\bullet};<0mm,0mm>*{}**@{},
 <-0.4mm,0.0mm>*{};<-2.4mm,2.1mm>*{}**@{-},
<-3.2mm,2.4mm>*{};<-5.7mm,0mm>*{}**@{-},
<2.9mm,7.5mm>*{};<-6.1mm,-2.2mm>*{^\be}**@{},
<2.9mm,7.5mm>*{};<0mm,9.3mm>*{^\al}**@{},
(3.2,3)*{}
   \ar@{->}@(ur,dr) (0.1,-1.4)*{}
\end{xy}
\ + \ 2b\
\begin{xy}
 <0.4mm,0.0mm>*{};<2.4mm,2.1mm>*{}**@{-},
<0mm,-0.8mm>*{\circ};<0mm,0.8mm>*{}**@{},
 <2.96mm,2.4mm>*{\circ};<2.45mm,2.35mm>*{}**@{},
 <2.4mm,2.8mm>*{};<0mm,5mm>*{}**@{-},
<-2.8mm,2.5mm>*{};<0mm,5mm>*{\bullet}**@{},
<-2.4mm,2.8mm>*{};<0mm,5mm>*{}**@{-},
<0mm,-1.2mm>*{};<0mm,-4mm>*{}**@{-},
<-2.96mm,2.4mm>*{\bullet};<0mm,0mm>*{}**@{},
 <-0.4mm,0.0mm>*{};<-2.4mm,2.1mm>*{}**@{-},
<3.4mm,3mm>*{};<5.3mm,5mm>*{}**@{-},
<0mm,0mm>*{};<0mm,-6.5mm>*{^\be}**@{},
<0mm,0mm>*{};<5.2mm,5.3mm>*{^\al}**@{},
(0,5)*{}
   \ar@{->}@(ul,dl) (-2.9,2.4)*{}
\end{xy}
 \\
 &=&
\ (4a-\frac{1}{6})\ \,
\begin{xy}
 <0.4mm,0.0mm>*{};<2.4mm,2.1mm>*{}**@{-},
<0mm,-0.8mm>*{\circ};<0mm,0.8mm>*{}**@{},
 <2.96mm,2.4mm>*{\circ};<2.45mm,2.35mm>*{}**@{},
 <2.4mm,2.8mm>*{};<0mm,5mm>*{}**@{-},
<-2.8mm,2.5mm>*{};<0mm,5mm>*{\bullet}**@{},
<-2.4mm,2.8mm>*{};<0mm,5mm>*{}**@{-},
<0mm,5mm>*{};<0mm,8.6mm>*{}**@{-},
<-2.96mm,2.4mm>*{\bullet};<0mm,0mm>*{}**@{},
 <-0.4mm,0.0mm>*{};<-2.4mm,2.1mm>*{}**@{-},
<-3.2mm,2.4mm>*{};<-5.7mm,0mm>*{}**@{-},
<2.9mm,7.5mm>*{};<-6.1mm,-2.2mm>*{^\be}**@{},
<2.9mm,7.5mm>*{};<0mm,9.3mm>*{^\al}**@{},
(3.2,3)*{}
   \ar@{->}@(ur,dr) (0.1,-1.4)*{}
\end{xy}
\ + \ (2b+\frac{1}{6})\
\begin{xy}
 <0.4mm,0.0mm>*{};<2.4mm,2.1mm>*{}**@{-},
<0mm,-0.8mm>*{\circ};<0mm,0.8mm>*{}**@{},
 <2.96mm,2.4mm>*{\circ};<2.45mm,2.35mm>*{}**@{},
 <2.4mm,2.8mm>*{};<0mm,5mm>*{}**@{-},
<-2.8mm,2.5mm>*{};<0mm,5mm>*{\bullet}**@{},
<-2.4mm,2.8mm>*{};<0mm,5mm>*{}**@{-},
<0mm,-1.2mm>*{};<0mm,-4mm>*{}**@{-},
<-2.96mm,2.4mm>*{\bullet};<0mm,0mm>*{}**@{},
 <-0.4mm,0.0mm>*{};<-2.4mm,2.1mm>*{}**@{-},
<3.4mm,3mm>*{};<5.3mm,5mm>*{}**@{-},
<0mm,0mm>*{};<0mm,-6.5mm>*{^\be}**@{},
<0mm,0mm>*{};<5.2mm,5.3mm>*{^\al}**@{},
(0,5)*{}
   \ar@{->}@(ul,dl) (-2.9,2.4)*{}
\end{xy}
\Eeqrn
Hence setting $a=1/24$ and $b=-1/12$, i.e.\ adding to the PBW
star product, $\Ga_2$, the operators,
$$
\Ga_0=-\hbar^3\frac{1}{12} \sum_{\al,\be,\ga,\mu,\nu}
C^\mu_{\al\be}C^{\be}_{\mu\nu} \Phi^{\al\nu}_\ga x^\ga
=-\hbar^3\frac{1}{12} \sum_{\ga} \Ba{c}
 \begin{xy}
 <0mm,0mm>*{\bullet};<0mm,0mm>*{}**@{},
 %<0mm,0.69mm>*{};<0mm,3.0mm>*{}**@{-},
 <0.39mm,-0.39mm>*{};<2.4mm,-2.4mm>*{}**@{-},
 <-0.35mm,-0.35mm>*{};<-1.9mm,-1.9mm>*{}**@{-},
 <-2.4mm,-2.4mm>*{\bullet};<-2.4mm,-2.4mm>*{}**@{},
 <-2.0mm,-2.8mm>*{};<-0.4mm,-4.5mm>*{}**@{-},
 %<-2.8mm,-2.9mm>*{};<-4.7mm,-4.9mm>*{}**@{-},
 %
 <2.4mm,-2.4mm>*{};<0.4mm,-4.5mm>*{}**@{-},
  <0mm,-5.1mm>*{\circ};<0mm,0mm>*{}**@{},
  <0mm,-5.5mm>*{};<0mm,-7.7mm>*{}**@{-},
    <0.39mm,-0.39mm>*{};<0mm,-9.1mm>*{_\ga}**@{},
(0,0)*{}
   \ar@{->}@(ul,dl) (-2.4,-2.4)*{}
 \end{xy}
\Ea \ x^\ga
$$
$$
\Ga_1(x^\ga)=  \hbar\sum_{\al,\be} C^\ga_{\al\be} x^{\al}x^{\be} +
\frac{1}{24} \hbar^3
 \sum_{\al,\be,\mu,\nu} C^\ga_{\al\be}\Phi^{\al\mu}_\nu \Phi^{\be\nu}_\mu=
  \hbar\sum_{\al,\be}
  \begin{xy}
 <0mm,0.66mm>*{};<0mm,3mm>*{}**@{-},
 <0.39mm,-0.39mm>*{};<2.2mm,-2.2mm>*{}**@{-},
 <-0.35mm,-0.35mm>*{};<-2.2mm,-2.2mm>*{}**@{-},
 <0mm,0mm>*{\bullet};<0mm,0mm>*{}**@{},
   <0mm,0.66mm>*{};<0mm,3.4mm>*{^\ga}**@{},
   <0.39mm,-0.39mm>*{};<2.9mm,-4mm>*{^\be}**@{},
   <-0.35mm,-0.35mm>*{};<-2.8mm,-4mm>*{^\al}**@{},
\end{xy}\ x^{\al}x^{\be} +    \frac{1}{24} \hbar^3
\begin{xy}
 <0.4mm,0.0mm>*{};<2.4mm,2.1mm>*{}**@{-},
 <-0.38mm,-0.2mm>*{};<-2.8mm,2.5mm>*{}**@{-},
<0mm,-0.8mm>*{\circ};<0mm,0.8mm>*{}**@{},
 <2.96mm,2.4mm>*{\circ};<2.45mm,2.35mm>*{}**@{},
 <2.4mm,2.8mm>*{};<0mm,5mm>*{}**@{-},
<-2.8mm,2.5mm>*{};<0mm,5mm>*{\bullet}**@{},
<-2.8mm,2.5mm>*{};<0mm,5mm>*{}**@{-},
<0mm,5mm>*{};<0mm,8.8mm>*{}**@{-},
<2.9mm,7.5mm>*{};<0mm,9.4mm>*{^\ga}**@{},
(3.2,3)*{}
   \ar@{->}@(ur,dr) (0.1,-1.4)*{}
    \end{xy} .
$$
we complete  deformation quantization of Lie 1-bialgebras.
Thus we proved the following,

\bip

%%%%%%%%%%%%%%%%%%%%%%%%%%%
\noindent{\bf I.2.5. Proposition}. {\em There exists a morphism of
dg props,\, $q:\DefQh \rar \Liebi^\circlearrowright$, making the
diagram
\[
 \xymatrix{
  & \Liebi^\circlearrowright
  \ar[d]^{s} \\
 \DefQh\ar[ur]^q\ar[r]_{\PBW} &
 \CoLie^\circlearrowright
 }
\]
commutative.}

\sip

\noindent Hence Step 1 is done.

\bip

%%%%%%%%%%%%%%%%%%%%%%%%%%%%%%%%%%%%%%%%%%%%%%%%%%%%%%%%%%%%%%%%%%%%%%
\noindent{\bf Step II.} First we explain iterative procedure behind
our construction of a morphism $Q$ fitting the commutative diagram
\[
 \xymatrix{
  & [\Liebi^\circlearrowright]_\infty
  \ar[d]^{qis} \\
 \DefQh\ar[ur]^Q\ar[r]_{q} &
 \Liebi^\circlearrowright
 }
\]
and then illustrate it with concrete examples. Define $E_s$ to be
zero for negative $s$ and, for $s\geq 0$,
 \Beqrn
E_s &:=& {\rm span}\left\{ \xy
%\begin{xy}
 <0mm,0mm>*{\mbox{$\xy *=<20mm,3mm>\txt{\em a}*\frm{-}\endxy$}};<0mm,0mm>*{}**@{},
  <-10mm,1.5mm>*{};<-12mm,7mm>*{}**@{-},
  <-10mm,1.5mm>*{};<-11mm,7mm>*{}**@{-},
  <-10mm,1.5mm>*{};<-9.5mm,6mm>*{}**@{-},
  <-10mm,1.5mm>*{};<-8mm,7mm>*{}**@{-},
 <-10mm,1.5mm>*{};<-9.5mm,6.6mm>*{.\hspace{-0.4mm}.\hspace{-0.4mm}.}**@{},
 <0mm,0mm>*{};<-6.5mm,3.6mm>*{.\hspace{-0.1mm}.\hspace{-0.1mm}.}**@{},
  <-3mm,1.5mm>*{};<-5mm,7mm>*{}**@{-},
  <-3mm,1.5mm>*{};<-4mm,7mm>*{}**@{-},
  <-3mm,1.5mm>*{};<-2.5mm,6mm>*{}**@{-},
  <-3mm,1.5mm>*{};<-1mm,7mm>*{}**@{-},
 <-3mm,1.5mm>*{};<-2.5mm,6.6mm>*{.\hspace{-0.4mm}.\hspace{-0.4mm}.}**@{},
  <2mm,1.5mm>*{};<0mm,7mm>*{}**@{-},
  <2mm,1.5mm>*{};<1mm,7mm>*{}**@{-},
  <2mm,1.5mm>*{};<2.5mm,6mm>*{}**@{-},
  <2mm,1.5mm>*{};<4mm,7mm>*{}**@{-},
 <2mm,1.5mm>*{};<2.5mm,6.6mm>*{.\hspace{-0.4mm}.\hspace{-0.4mm}.}**@{},
 <0mm,0mm>*{};<6mm,3.6mm>*{.\hspace{-0.1mm}.\hspace{-0.1mm}.}**@{},
<10mm,1.5mm>*{};<8mm,7mm>*{}**@{-},
  <10mm,1.5mm>*{};<9mm,7mm>*{}**@{-},
  <10mm,1.5mm>*{};<10.5mm,6mm>*{}**@{-},
  <10mm,1.5mm>*{};<12mm,7mm>*{}**@{-},
 <10mm,1.5mm>*{};<10.5mm,6.6mm>*{.\hspace{-0.4mm}.\hspace{-0.4mm}.}**@{},
 <-10mm,-1.5mm>*{};<-12mm,-6mm>*{}**@{-},
 <-7mm,-1.5mm>*{};<-8mm,-6mm>*{}**@{-},
 <-4mm,-1.5mm>*{};<-4.5mm,-6mm>*{}**@{-},
 <0mm,0mm>*{};<0mm,-4.6mm>*{.\hspace{0.1mm}.\hspace{0.1mm}.}**@{},
<10mm,-1.5mm>*{};<12mm,-6mm>*{}**@{-},
 <7mm,-1.5mm>*{};<8mm,-6mm>*{}**@{-},
  <4mm,-1.5mm>*{};<4.5mm,-6mm>*{}**@{-},
  %%%%
<0mm,0mm>*{};<-9.5mm,8.2mm>*{^{I_{ 1}}}**@{},
<0mm,0mm>*{};<-3mm,8.2mm>*{^{I_{ i}}}**@{},
<0mm,0mm>*{};<2mm,8.2mm>*{^{I_{ i+1}}}**@{},
<0mm,0mm>*{};<10mm,8.2mm>*{^{I_{ k}}}**@{},
 %%%
<0mm,0mm>*{};<-12mm,-7.4mm>*{_1}**@{},
<0mm,0mm>*{};<-8mm,-7.4mm>*{_2}**@{},
<0mm,0mm>*{};<-4mm,-7.4mm>*{_3}**@{},
<0mm,0mm>*{};<6mm,-7.4mm>*{\ldots}**@{},
%<0mm,0mm>*{};<8.5mm,-7.4mm>*{_{n-1}}**@{},
<0mm,0mm>*{};<12.5mm,-7.4mm>*{_n}**@{},
\endxy
\in \DefQh \right\}_{2a + k-2=s\atop k\geq 0, a\geq 1, n\geq 0}
\Eeqrn For example,
\[
E_{0} =
 {\rm span}
\left\{ \xy
 <0mm,0mm>*{\mbox{$\xy *=<5mm,3mm>\txt{\em 1}*\frm{-}\endxy$}};<0mm,0mm>*{}**@{},
 <-2.5mm,-1.6mm>*{};<-4mm,-6mm>*{}**@{-},
 <-1.6mm,-1.6mm>*{};<-2mm,-6mm>*{}**@{-},
 <1.6mm,-1.6mm>*{};<2.5mm,-6mm>*{}**@{-},
 <2.5mm,-1.6mm>*{};<4mm,-6mm>*{}**@{-},
 <0mm,0mm>*{};<0mm,-4.6mm>*{.\hspace{0mm}.\hspace{0mm}.}**@{},
<0mm,0mm>*{};<-4.1mm,-7.4mm>*{_1}**@{},
<0mm,0mm>*{};<-2mm,-7.4mm>*{_2}**@{},
%<0mm,0mm>*{};<-4mm,-7.4mm>*{_3}**@{},
<0mm,0mm>*{};<1.7mm,-7.4mm>*{...}**@{},
%<0mm,0mm>*{};<8.5mm,-7.4mm>*{_{n-1}}**@{},
<0mm,0mm>*{};<4.6mm,-7.4mm>*{_n}**@{},
\endxy
\right\} \ \ \ \ E_{1} =
 {\rm span}
\left\{ \xy
 <0mm,0mm>*{\mbox{$\xy *=<5mm,3mm>\txt{\em 1}*\frm{-}\endxy$}};<0mm,0mm>*{}**@{},
 <-2.5mm,-1.6mm>*{};<-4mm,-6mm>*{}**@{-},
 <-1.6mm,-1.6mm>*{};<-2mm,-6mm>*{}**@{-},
 <1.6mm,-1.6mm>*{};<2.5mm,-6mm>*{}**@{-},
 <2.5mm,-1.6mm>*{};<4mm,-6mm>*{}**@{-},
 <0mm,0mm>*{};<0mm,-4.6mm>*{.\hspace{0mm}.\hspace{0mm}.}**@{},
<0mm,0mm>*{};<-4.1mm,-7.4mm>*{_1}**@{},
<0mm,0mm>*{};<-2mm,-7.4mm>*{_2}**@{},
%<0mm,0mm>*{};<-4mm,-7.4mm>*{_3}**@{},
<0mm,0mm>*{};<1.7mm,-7.4mm>*{...}**@{},
%<0mm,0mm>*{};<8.5mm,-7.4mm>*{_{n-1}}**@{},
<0mm,0mm>*{};<4.6mm,-7.4mm>*{_n}**@{},
<0mm,1.5mm>*{};<2mm,7mm>*{}**@{-},
  <0mm,1.5mm>*{};<1mm,6mm>*{}**@{-},
  <0mm,1.5mm>*{};<-0.5mm,7mm>*{}**@{-},
  <0mm,1.5mm>*{};<-2mm,7mm>*{}**@{-},
 <0mm,1.5mm>*{};<0.7mm,6.6mm>*{.\hspace{-0.4mm}.\hspace{-0.4mm}.}**@{},
\endxy
\right\} \ \ \ \ \ E_{2} =
 {\rm span}
\left\{ \xy
 <0mm,0mm>*{\mbox{$\xy *=<5mm,3mm>\txt{\em 2}*\frm{-}\endxy$}};<0mm,0mm>*{}**@{},
 <-2.5mm,-1.6mm>*{};<-4mm,-6mm>*{}**@{-},
 <-1.6mm,-1.6mm>*{};<-2mm,-6mm>*{}**@{-},
 <1.6mm,-1.6mm>*{};<2.5mm,-6mm>*{}**@{-},
 <2.5mm,-1.6mm>*{};<4mm,-6mm>*{}**@{-},
 <0mm,0mm>*{};<0mm,-4.6mm>*{.\hspace{0mm}.\hspace{0mm}.}**@{},
<0mm,0mm>*{};<-4.1mm,-7.4mm>*{_1}**@{},
<0mm,0mm>*{};<-2mm,-7.4mm>*{_2}**@{},
%<0mm,0mm>*{};<-4mm,-7.4mm>*{_3}**@{},
<0mm,0mm>*{};<1.7mm,-7.4mm>*{...}**@{},
%<0mm,0mm>*{};<8.5mm,-7.4mm>*{_{n-1}}**@{},
<0mm,0mm>*{};<4.6mm,-7.4mm>*{_n}**@{},
\endxy
\ \ , \ \ \xy
 <0mm,0mm>*{\mbox{$\xy *=<5mm,3mm>\txt{\em 1}*\frm{-}\endxy$}};<0mm,0mm>*{}**@{},
 <-2.5mm,-1.6mm>*{};<-4mm,-6mm>*{}**@{-},
 <-1.6mm,-1.6mm>*{};<-2mm,-6mm>*{}**@{-},
 <1.6mm,-1.6mm>*{};<2.5mm,-6mm>*{}**@{-},
 <2.5mm,-1.6mm>*{};<4mm,-6mm>*{}**@{-},
 <0mm,0mm>*{};<0mm,-4.6mm>*{.\hspace{0mm}.\hspace{0mm}.}**@{},
<0mm,0mm>*{};<-4.1mm,-7.4mm>*{_1}**@{},
<0mm,0mm>*{};<-2mm,-7.4mm>*{_2}**@{},
%<0mm,0mm>*{};<-4mm,-7.4mm>*{_3}**@{},
<0mm,0mm>*{};<1.7mm,-7.4mm>*{...}**@{},
%<0mm,0mm>*{};<8.5mm,-7.4mm>*{_{n-1}}**@{},
<0mm,0mm>*{};<4.6mm,-7.4mm>*{_n}**@{},
 <-2.5mm,1.5mm>*{};<-5.5mm,7mm>*{}**@{-},
  <-2.5mm,1.5mm>*{};<-4.5mm,7mm>*{}**@{-},
  <-2.5mm,1.5mm>*{};<-3mm,6mm>*{}**@{-},
  <-2.5mm,1.5mm>*{};<-1.5mm,7mm>*{}**@{-},
 <-2.5mm,1.5mm>*{};<-3mm,6.6mm>*{.\hspace{-0.4mm}.\hspace{-0.4mm}.}**@{},
 <2.5mm,1.5mm>*{};<5.5mm,7mm>*{}**@{-},
  <2.5mm,1.5mm>*{};<4.5mm,7mm>*{}**@{-},
  <2.5mm,1.5mm>*{};<3mm,6mm>*{}**@{-},
  <2.5mm,1.5mm>*{};<1.5mm,7mm>*{}**@{-},
 <2.5mm,1.5mm>*{};<3mm,6.6mm>*{.\hspace{-0.4mm}.\hspace{-0.4mm}.}**@{},
\endxy
\right\}.
\]
Let  $\DefQh_s\subset \DefQh$
 be the free prop generated by $\oplus_{i=0}^{s} E_i$. Thereby we get an increasing filtration,
$0\subset \DefQh_0\subset \ldots\subset  \DefQh_s\subset
\DefQh_{s+1} \ldots$
  with
 \[
\lim_{s\rar\infty} \DefQh_s = \DefQh.
\]
A straightforward inspection of the formula for differential
$\delta$ in \S 2.9 implies
$$
\delta E_{s+1}\subset \DefQh_s,
$$
i.e.\ that the dg prop $(\DefQh, \delta)$ has a ``cell" structure
analogous to that of CW complex. This simple but crucial for our
purposes observation permits us to apply the well-known in algebraic
topology Whitehead lifting trick and
 construct a morphism $Q: \DefQh \rar
[\Liebi^\circlearrowright]_\infty$ by defining its values on the
generators $E_s$ via an induction on $s$.
 We begin the induction by specifying its values on
$E_0\oplus E_1\oplus E_2$ as follows
\[
Q\left( \xy
 <0mm,0mm>*{\mbox{$\xy *=<5mm,3mm>\txt{\em 1}*\frm{-}\endxy$}};<0mm,0mm>*{}**@{},
 <-2.5mm,-1.6mm>*{};<-4mm,-6mm>*{}**@{-},
 <-1.6mm,-1.6mm>*{};<-2mm,-6mm>*{}**@{-},
 <1.6mm,-1.6mm>*{};<2.5mm,-6mm>*{}**@{-},
 <2.5mm,-1.6mm>*{};<4mm,-6mm>*{}**@{-},
 <0mm,0mm>*{};<0mm,-4.6mm>*{.\hspace{0mm}.\hspace{0mm}.}**@{},
<0mm,0mm>*{};<-4.1mm,-7.4mm>*{_1}**@{},
<0mm,0mm>*{};<-2mm,-7.4mm>*{_2}**@{},
%<0mm,0mm>*{};<-4mm,-7.4mm>*{_3}**@{},
<0mm,0mm>*{};<1.7mm,-7.4mm>*{...}**@{},
%<0mm,0mm>*{};<8.5mm,-7.4mm>*{_{n-1}}**@{},
<0mm,0mm>*{};<4.6mm,-7.4mm>*{_n}**@{},
\endxy
\right)=0, \ \ \ \ \ \ \ Q\left( \xy
 <0mm,0mm>*{\mbox{$\xy *=<5mm,3mm>\txt{\em 1}*\frm{-}\endxy$}};<0mm,0mm>*{}**@{},
 <-2.5mm,-1.6mm>*{};<-4mm,-6mm>*{}**@{-},
 <-1.6mm,-1.6mm>*{};<-2mm,-6mm>*{}**@{-},
 <1.6mm,-1.6mm>*{};<2.5mm,-6mm>*{}**@{-},
 <2.5mm,-1.6mm>*{};<4mm,-6mm>*{}**@{-},
 <0mm,0mm>*{};<0mm,-4.6mm>*{.\hspace{0mm}.\hspace{0mm}.}**@{},
<0mm,0mm>*{};<-4.1mm,-7.4mm>*{_1}**@{},
<0mm,0mm>*{};<-2mm,-7.4mm>*{_2}**@{},
%<0mm,0mm>*{};<-4mm,-7.4mm>*{_3}**@{},
<0mm,0mm>*{};<1.7mm,-7.4mm>*{...}**@{},
%<0mm,0mm>*{};<8.5mm,-7.4mm>*{_{n-1}}**@{},
<0mm,0mm>*{};<4.6mm,-7.4mm>*{_n}**@{},
<0mm,1.5mm>*{};<2mm,7mm>*{}**@{-},
  <0mm,1.5mm>*{};<1mm,6mm>*{}**@{-},
  <0mm,1.5mm>*{};<-0.5mm,7mm>*{}**@{-},
  <0mm,1.5mm>*{};<-2mm,7mm>*{}**@{-},
 <0mm,1.5mm>*{};<0.7mm,6.6mm>*{.\hspace{-0.4mm}.\hspace{-0.4mm}.}**@{},
 <0mm,1.5mm>*{};<0.7mm,8.6mm>*{_{I_1}}**@{},
\endxy
\right)= \left\{\Ba{rr}
\begin{xy}
 <0mm,0.66mm>*{};<0mm,3mm>*{}**@{-},
 <0.39mm,-0.39mm>*{};<2.2mm,-2.2mm>*{}**@{-},
 <-0.35mm,-0.35mm>*{};<-2.2mm,-2.2mm>*{}**@{-},
 <0mm,0mm>*{\bullet};<0mm,0mm>*{}**@{},
\end{xy} & \mbox{for}\ |I_1|=1, n=2 \vspace{3mm}\\
0 & \mbox{otherwise.} \Ea \right.
\]

\[
Q\left( \xy
 <0mm,0mm>*{\mbox{$\xy *=<5mm,3mm>\txt{\em 2}*\frm{-}\endxy$}};<0mm,0mm>*{}**@{},
 <-2.5mm,-1.6mm>*{};<-4mm,-6mm>*{}**@{-},
 <-1.6mm,-1.6mm>*{};<-2mm,-6mm>*{}**@{-},
 <1.6mm,-1.6mm>*{};<2.5mm,-6mm>*{}**@{-},
 <2.5mm,-1.6mm>*{};<4mm,-6mm>*{}**@{-},
 <0mm,0mm>*{};<0mm,-4.6mm>*{.\hspace{0mm}.\hspace{0mm}.}**@{},
<0mm,0mm>*{};<-4.1mm,-7.4mm>*{_1}**@{},
<0mm,0mm>*{};<-2mm,-7.4mm>*{_2}**@{},
%<0mm,0mm>*{};<-4mm,-7.4mm>*{_3}**@{},
<0mm,0mm>*{};<1.7mm,-7.4mm>*{...}**@{},
%<0mm,0mm>*{};<8.5mm,-7.4mm>*{_{n-1}}**@{},
<0mm,0mm>*{};<4.6mm,-7.4mm>*{_n}**@{},
\endxy
\right)=0, \ \ \ \ \ \ \ Q\left( \xy
 <0mm,0mm>*{\mbox{$\xy *=<5mm,3mm>\txt{\em 1}*\frm{-}\endxy$}};<0mm,0mm>*{}**@{},
 <-2.5mm,-1.6mm>*{};<-4mm,-6mm>*{}**@{-},
 <-1.6mm,-1.6mm>*{};<-2mm,-6mm>*{}**@{-},
 <1.6mm,-1.6mm>*{};<2.5mm,-6mm>*{}**@{-},
 <2.5mm,-1.6mm>*{};<4mm,-6mm>*{}**@{-},
 <0mm,0mm>*{};<0mm,-4.6mm>*{.\hspace{0mm}.\hspace{0mm}.}**@{},
<0mm,0mm>*{};<-4.1mm,-7.4mm>*{_1}**@{},
<0mm,0mm>*{};<-2mm,-7.4mm>*{_2}**@{},
%<0mm,0mm>*{};<-4mm,-7.4mm>*{_3}**@{},
<0mm,0mm>*{};<1.7mm,-7.4mm>*{...}**@{},
%<0mm,0mm>*{};<8.5mm,-7.4mm>*{_{n-1}}**@{},
<0mm,0mm>*{};<4.6mm,-7.4mm>*{_n}**@{},
 <-2.5mm,1.5mm>*{};<-5.5mm,7mm>*{}**@{-},
  <-2.5mm,1.5mm>*{};<-4.5mm,7mm>*{}**@{-},
  <-2.5mm,1.5mm>*{};<-3mm,6mm>*{}**@{-},
  <-2.5mm,1.5mm>*{};<-1.5mm,7mm>*{}**@{-},
 <-2.5mm,1.5mm>*{};<-3mm,6.6mm>*{.\hspace{-0.4mm}.\hspace{-0.4mm}.}**@{},
 <0mm,0mm>*{};<-3mm,8.6mm>*{_{I_1}}**@{},
 <2.5mm,1.5mm>*{};<5.5mm,7mm>*{}**@{-},
  <2.5mm,1.5mm>*{};<4.5mm,7mm>*{}**@{-},
  <2.5mm,1.5mm>*{};<3mm,6mm>*{}**@{-},
  <2.5mm,1.5mm>*{};<1.5mm,7mm>*{}**@{-},
 <2.5mm,1.5mm>*{};<3mm,6.6mm>*{.\hspace{-0.4mm}.\hspace{-0.4mm}.}**@{},
<0mm,0mm>*{};<3mm,8.6mm>*{_{I_2}}**@{},
\endxy
\right)= \left\{\Ba{rr}
\begin{xy}
 <0mm,-0.55mm>*{};<0mm,-2.5mm>*{}**@{-},
 <0.5mm,0.5mm>*{};<2.2mm,2.2mm>*{}**@{-},
 <-0.48mm,0.48mm>*{};<-2.2mm,2.2mm>*{}**@{-},
 <0mm,0mm>*{\circ};<0mm,0mm>*{}**@{},
 \end{xy}
 & \mbox{for}\ |I_1|=1, |I_2|=1,n=1 \vspace{3mm}\\
0 & \mbox{otherwise.} \Ea \right.
\]
Note that this choice respects both commutativity of the diagram and
the condition on $\pi_1\circ Q$ in the Proposition. Thus we
constructed a morphism of dg props,
$$
Q_3: \DefQh_3 \lon [\Liebi^\circlearrowright]_\infty
$$
satisfying the required conditions.

\bip

\no
Assume now that a morphism
$$
Q_s: \DefQh_s \lon [\Liebi^\circlearrowright]_\infty
$$
satisfying the required conditions is already constructed for some
$s\geq 3$. We want to extend $Q_s$ to a morphism of dg props,
$$
\Ba{rccc}
Q_{s+1}:&  \DefQh_{s+1}&\lon& [\Liebi^\circlearrowright]_\infty\\
& e & \lon & Q_{s+1}(e) \Ea
$$
such that  ${\it qis}\circ Q_{s+1}=q$ and the condition on
$\pi_1\circ Q_{s+1}$ is fulfilled. Let $e'$ be a lift of $q(e)$
along the surjection $qis$. Then
$$
Q_s(\delta e) - \delta e'
$$
is a cycle in  $[\Liebi^\circlearrowright]_\infty$ which projects
under $qis$ to zero. As $qis$ is a homology isomorphism, this
element is exact,
$$
Q_s(\delta e) - \delta e' = \delta e'',
$$
for some $e''\in [\Liebi^\circlearrowright]_\infty$. We set
$Q_{s+1}(e):=e'+e''$ completing thereby the inductive construction
of $Q$ as a morphism of dg props. The condition $qis\circ Q=q$ is
automatically satisfied. Another condition on $\pi_1\circ\al\circ Q$
is also satisfied because the Hochschild brackets, $[\ga,\ga]_H$, of
a degree 1 polyvector field viewed as a polydifferential operator
contain the Schouten brackets, $[\ga,\ga]_S$, as one of the
irreducible (in the $\bS$-bimodule sense) summands. Using this fact
it is easy to check by induction on $2a+k-2$ that
$$
\pi_{1}\circ\al\circ Q\left( \xy
%\begin{xy}
 <0mm,0mm>*{\mbox{$\xy *=<20mm,3mm>\txt{\em a}*\frm{-}\endxy$}};<0mm,0mm>*{}**@{},
  <-10mm,1.5mm>*{};<-12mm,7mm>*{}**@{-},
  <-10mm,1.5mm>*{};<-11mm,7mm>*{}**@{-},
  <-10mm,1.5mm>*{};<-9.5mm,6mm>*{}**@{-},
  <-10mm,1.5mm>*{};<-8mm,7mm>*{}**@{-},
 <-10mm,1.5mm>*{};<-9.5mm,6.6mm>*{.\hspace{-0.4mm}.\hspace{-0.4mm}.}**@{},
 <0mm,0mm>*{};<-6.5mm,3.6mm>*{.\hspace{-0.1mm}.\hspace{-0.1mm}.}**@{},
  <-3mm,1.5mm>*{};<-5mm,7mm>*{}**@{-},
  <-3mm,1.5mm>*{};<-4mm,7mm>*{}**@{-},
  <-3mm,1.5mm>*{};<-2.5mm,6mm>*{}**@{-},
  <-3mm,1.5mm>*{};<-1mm,7mm>*{}**@{-},
 <-3mm,1.5mm>*{};<-2.5mm,6.6mm>*{.\hspace{-0.4mm}.\hspace{-0.4mm}.}**@{},
  <2mm,1.5mm>*{};<0mm,7mm>*{}**@{-},
  <2mm,1.5mm>*{};<1mm,7mm>*{}**@{-},
  <2mm,1.5mm>*{};<2.5mm,6mm>*{}**@{-},
  <2mm,1.5mm>*{};<4mm,7mm>*{}**@{-},
 <2mm,1.5mm>*{};<2.5mm,6.6mm>*{.\hspace{-0.4mm}.\hspace{-0.4mm}.}**@{},
 <0mm,0mm>*{};<6mm,3.6mm>*{.\hspace{-0.1mm}.\hspace{-0.1mm}.}**@{},
<10mm,1.5mm>*{};<8mm,7mm>*{}**@{-},
  <10mm,1.5mm>*{};<9mm,7mm>*{}**@{-},
  <10mm,1.5mm>*{};<10.5mm,6mm>*{}**@{-},
  <10mm,1.5mm>*{};<12mm,7mm>*{}**@{-},
 <10mm,1.5mm>*{};<10.5mm,6.6mm>*{.\hspace{-0.4mm}.\hspace{-0.4mm}.}**@{},
 <-10mm,-1.5mm>*{};<-12mm,-6mm>*{}**@{-},
 <-7mm,-1.5mm>*{};<-8mm,-6mm>*{}**@{-},
 <-4mm,-1.5mm>*{};<-4.5mm,-6mm>*{}**@{-},
 <0mm,0mm>*{};<0mm,-4.6mm>*{.\hspace{0.1mm}.\hspace{0.1mm}.}**@{},
<10mm,-1.5mm>*{};<12mm,-6mm>*{}**@{-},
 <7mm,-1.5mm>*{};<8mm,-6mm>*{}**@{-},
  <4mm,-1.5mm>*{};<4.5mm,-6mm>*{}**@{-},
  %%%%
<0mm,0mm>*{};<-9.5mm,8.2mm>*{^{I_{ 1}}}**@{},
<0mm,0mm>*{};<-3mm,8.2mm>*{^{I_{ i}}}**@{},
<0mm,0mm>*{};<2mm,8.2mm>*{^{I_{ i+1}}}**@{},
<0mm,0mm>*{};<10mm,8.2mm>*{^{I_{ k}}}**@{},
 %%%
<0mm,0mm>*{};<-12mm,-7.4mm>*{_1}**@{},
<0mm,0mm>*{};<-8mm,-7.4mm>*{_2}**@{},
<0mm,0mm>*{};<-4mm,-7.4mm>*{_3}**@{},
<0mm,0mm>*{};<6mm,-7.4mm>*{\ldots}**@{},
%<0mm,0mm>*{};<8.5mm,-7.4mm>*{_{n-1}}**@{},
<0mm,0mm>*{};<12.5mm,-7.4mm>*{_n}**@{},
\endxy
\right) = \left\{\Ba{cl}
 \begin{xy}
 <0mm,0mm>*{\bullet};<0mm,0mm>*{}**@{},
 <0mm,0mm>*{};<-8mm,5mm>*{}**@{-},
 <0mm,0mm>*{};<-4.5mm,5mm>*{}**@{-},
 <0mm,0mm>*{};<-1mm,5mm>*{\ldots}**@{},
 <0mm,0mm>*{};<4.5mm,5mm>*{}**@{-},
 <0mm,0mm>*{};<8mm,5mm>*{}**@{-},
   <0mm,0mm>*{};<-8.5mm,5.5mm>*{^1}**@{},
   <0mm,0mm>*{};<-5mm,5.5mm>*{^2}**@{},
   <0mm,0mm>*{};<4.5mm,5.5mm>*{^{k\hspace{-0.5mm}-\hspace{-0.5mm}1}}**@{},
   <0mm,0mm>*{};<9.0mm,5.5mm>*{^k}**@{},
 <0mm,0mm>*{};<-8mm,-5mm>*{}**@{-},
 <0mm,0mm>*{};<-4.5mm,-5mm>*{}**@{-},
 <0mm,0mm>*{};<-1mm,-5mm>*{\ldots}**@{},
 <0mm,0mm>*{};<4.5mm,-5mm>*{}**@{-},
 <0mm,0mm>*{};<8mm,-5mm>*{}**@{-},
   <0mm,0mm>*{};<-8.5mm,-6.9mm>*{^1}**@{},
   <0mm,0mm>*{};<-5mm,-6.9mm>*{^2}**@{},
   <0mm,0mm>*{};<4.5mm,-6.9mm>*{^{n\hspace{-0.5mm}-\hspace{-0.5mm}1}}**@{},
   <0mm,0mm>*{};<9.0mm,-6.9mm>*{^n}**@{},
 \end{xy}
&
\mbox{for}\ a=n+k-2\ \mbox{and}\ |I_1|=\ldots=|I_k|=1 \\
0 & \mbox{otherwise}, \Ea \right.
$$
 where $\pi_{1}:
{\Liebi}^\circlearrowright_\infty\rar
({\Liebi}_\infty^\circlearrowright)_{1}$ is the projection to the
subspace, $ ({\Liebi}_\infty^\circlearrowright)_{1}$, consisting of
graphs with precisely $1$ internal vertex.

\bip

\noindent{\bf II.2. Corollary.} (i) If $a_s\in E_s$ is such that
$Q_{s-1}(\delta a_s)=0$ then we can set $Q_s(a_s):= qis^{-1}\circ
q(a_s)$, an arbitrary lifting of $q(a_s)\in
\Liebi^\circlearrowright$ to $\Liebi^\circlearrowright_\infty \subset [\Liebi^\circlearrowright]_\infty$.

(ii) If $a_s\in E_s$ is such that $Q_{s-1}(\delta a_s)=0$  and
$q(a_s)=0$, then we can set $Q(a_s)=0$.

\bip

\noindent{\bf II.3. Iteration.} Iteration (with respect
to the ''weight" parameter $s=2a+k-2$)
formula has the form,
$$
Q\left( \xy
%\begin{xy}
 <0mm,0mm>*{\mbox{$\xy *=<20mm,3mm>\txt{\em a}*\frm{-}\endxy$}};<0mm,0mm>*{}**@{},
  <-10mm,1.5mm>*{};<-12mm,7mm>*{}**@{-},
  <-10mm,1.5mm>*{};<-11mm,7mm>*{}**@{-},
  <-10mm,1.5mm>*{};<-9.5mm,6mm>*{}**@{-},
  <-10mm,1.5mm>*{};<-8mm,7mm>*{}**@{-},
 <-10mm,1.5mm>*{};<-9.5mm,6.6mm>*{.\hspace{-0.4mm}.\hspace{-0.4mm}.}**@{},
 <0mm,0mm>*{};<-6.5mm,3.6mm>*{.\hspace{-0.1mm}.\hspace{-0.1mm}.}**@{},
  <-3mm,1.5mm>*{};<-5mm,7mm>*{}**@{-},
  <-3mm,1.5mm>*{};<-4mm,7mm>*{}**@{-},
  <-3mm,1.5mm>*{};<-2.5mm,6mm>*{}**@{-},
  <-3mm,1.5mm>*{};<-1mm,7mm>*{}**@{-},
 <-3mm,1.5mm>*{};<-2.5mm,6.6mm>*{.\hspace{-0.4mm}.\hspace{-0.4mm}.}**@{},
  <2mm,1.5mm>*{};<0mm,7mm>*{}**@{-},
  <2mm,1.5mm>*{};<1mm,7mm>*{}**@{-},
  <2mm,1.5mm>*{};<2.5mm,6mm>*{}**@{-},
  <2mm,1.5mm>*{};<4mm,7mm>*{}**@{-},
 <2mm,1.5mm>*{};<2.5mm,6.6mm>*{.\hspace{-0.4mm}.\hspace{-0.4mm}.}**@{},
 <0mm,0mm>*{};<6mm,3.6mm>*{.\hspace{-0.1mm}.\hspace{-0.1mm}.}**@{},
<10mm,1.5mm>*{};<8mm,7mm>*{}**@{-},
  <10mm,1.5mm>*{};<9mm,7mm>*{}**@{-},
  <10mm,1.5mm>*{};<10.5mm,6mm>*{}**@{-},
  <10mm,1.5mm>*{};<12mm,7mm>*{}**@{-},
 <10mm,1.5mm>*{};<10.5mm,6.6mm>*{.\hspace{-0.4mm}.\hspace{-0.4mm}.}**@{},
 <-10mm,-1.5mm>*{};<-12mm,-6mm>*{}**@{-},
 <-7mm,-1.5mm>*{};<-8mm,-6mm>*{}**@{-},
 <-4mm,-1.5mm>*{};<-4.5mm,-6mm>*{}**@{-},
 <0mm,0mm>*{};<0mm,-4.6mm>*{.\hspace{0.1mm}.\hspace{0.1mm}.}**@{},
<10mm,-1.5mm>*{};<12mm,-6mm>*{}**@{-},
 <7mm,-1.5mm>*{};<8mm,-6mm>*{}**@{-},
  <4mm,-1.5mm>*{};<4.5mm,-6mm>*{}**@{-},
  %%%%
<0mm,0mm>*{};<-9.5mm,8.2mm>*{^{I_{ 1}}}**@{},
<0mm,0mm>*{};<-3mm,8.2mm>*{^{I_{ i}}}**@{},
<0mm,0mm>*{};<2mm,8.2mm>*{^{I_{ i+1}}}**@{},
<0mm,0mm>*{};<10mm,8.2mm>*{^{I_{ k}}}**@{},
 %%%
<0mm,0mm>*{};<-12mm,-7.4mm>*{_1}**@{},
<0mm,0mm>*{};<-8mm,-7.4mm>*{_2}**@{},
<0mm,0mm>*{};<-4mm,-7.4mm>*{_3}**@{},
<0mm,0mm>*{};<6mm,-7.4mm>*{\ldots}**@{},
%<0mm,0mm>*{};<8.5mm,-7.4mm>*{_{n-1}}**@{},
<0mm,0mm>*{};<12.5mm,-7.4mm>*{_n}**@{},
\endxy
\right)= e'+ e''
$$
where $e'$ is an arbitrary lift of
$
q\left( \xy
%\begin{xy}
 <0mm,0mm>*{\mbox{$\xy *=<20mm,3mm>\txt{\em a}*\frm{-}\endxy$}};<0mm,0mm>*{}**@{},
  <-10mm,1.5mm>*{};<-12mm,7mm>*{}**@{-},
  <-10mm,1.5mm>*{};<-11mm,7mm>*{}**@{-},
  <-10mm,1.5mm>*{};<-9.5mm,6mm>*{}**@{-},
  <-10mm,1.5mm>*{};<-8mm,7mm>*{}**@{-},
 <-10mm,1.5mm>*{};<-9.5mm,6.6mm>*{.\hspace{-0.4mm}.\hspace{-0.4mm}.}**@{},
 <0mm,0mm>*{};<-6.5mm,3.6mm>*{.\hspace{-0.1mm}.\hspace{-0.1mm}.}**@{},
  <-3mm,1.5mm>*{};<-5mm,7mm>*{}**@{-},
  <-3mm,1.5mm>*{};<-4mm,7mm>*{}**@{-},
  <-3mm,1.5mm>*{};<-2.5mm,6mm>*{}**@{-},
  <-3mm,1.5mm>*{};<-1mm,7mm>*{}**@{-},
 <-3mm,1.5mm>*{};<-2.5mm,6.6mm>*{.\hspace{-0.4mm}.\hspace{-0.4mm}.}**@{},
  <2mm,1.5mm>*{};<0mm,7mm>*{}**@{-},
  <2mm,1.5mm>*{};<1mm,7mm>*{}**@{-},
  <2mm,1.5mm>*{};<2.5mm,6mm>*{}**@{-},
  <2mm,1.5mm>*{};<4mm,7mm>*{}**@{-},
 <2mm,1.5mm>*{};<2.5mm,6.6mm>*{.\hspace{-0.4mm}.\hspace{-0.4mm}.}**@{},
 <0mm,0mm>*{};<6mm,3.6mm>*{.\hspace{-0.1mm}.\hspace{-0.1mm}.}**@{},
<10mm,1.5mm>*{};<8mm,7mm>*{}**@{-},
  <10mm,1.5mm>*{};<9mm,7mm>*{}**@{-},
  <10mm,1.5mm>*{};<10.5mm,6mm>*{}**@{-},
  <10mm,1.5mm>*{};<12mm,7mm>*{}**@{-},
 <10mm,1.5mm>*{};<10.5mm,6.6mm>*{.\hspace{-0.4mm}.\hspace{-0.4mm}.}**@{},
 <-10mm,-1.5mm>*{};<-12mm,-6mm>*{}**@{-},
 <-7mm,-1.5mm>*{};<-8mm,-6mm>*{}**@{-},
 <-4mm,-1.5mm>*{};<-4.5mm,-6mm>*{}**@{-},
 <0mm,0mm>*{};<0mm,-4.6mm>*{.\hspace{0.1mm}.\hspace{0.1mm}.}**@{},
<10mm,-1.5mm>*{};<12mm,-6mm>*{}**@{-},
 <7mm,-1.5mm>*{};<8mm,-6mm>*{}**@{-},
  <4mm,-1.5mm>*{};<4.5mm,-6mm>*{}**@{-},
  %%%%
<0mm,0mm>*{};<-9.5mm,8.2mm>*{^{I_{ 1}}}**@{},
<0mm,0mm>*{};<-3mm,8.2mm>*{^{I_{ i}}}**@{},
<0mm,0mm>*{};<2mm,8.2mm>*{^{I_{ i+1}}}**@{},
<0mm,0mm>*{};<10mm,8.2mm>*{^{I_{ k}}}**@{},
 %%%
<0mm,0mm>*{};<-12mm,-7.4mm>*{_1}**@{},
<0mm,0mm>*{};<-8mm,-7.4mm>*{_2}**@{},
<0mm,0mm>*{};<-4mm,-7.4mm>*{_3}**@{},
<0mm,0mm>*{};<6mm,-7.4mm>*{\ldots}**@{},
%<0mm,0mm>*{};<8.5mm,-7.4mm>*{_{n-1}}**@{},
<0mm,0mm>*{};<12.5mm,-7.4mm>*{_n}**@{},
\endxy
\right)\in \Liebi^\circlearrowright
$
to a cycle in $\Liebi^\circlearrowright_\infty \subset [\Liebi^\circlearrowright]_\infty$, and $e''$ is a
solution of the following equation in
 $ [\Liebi^\circlearrowright]_\infty$,
\Beqrn \delta e''= Q\left(\delta \xy
%\begin{xy}
 <0mm,0mm>*{\mbox{$\xy *=<20mm,3mm>\txt{\em a}*\frm{-}\endxy$}};<0mm,0mm>*{}**@{},
  <-10mm,1.5mm>*{};<-12mm,7mm>*{}**@{-},
  <-10mm,1.5mm>*{};<-11mm,7mm>*{}**@{-},
  <-10mm,1.5mm>*{};<-9.5mm,6mm>*{}**@{-},
  <-10mm,1.5mm>*{};<-8mm,7mm>*{}**@{-},
 <-10mm,1.5mm>*{};<-9.5mm,6.6mm>*{.\hspace{-0.4mm}.\hspace{-0.4mm}.}**@{},
 <0mm,0mm>*{};<-6.5mm,3.6mm>*{.\hspace{-0.1mm}.\hspace{-0.1mm}.}**@{},
  <-3mm,1.5mm>*{};<-5mm,7mm>*{}**@{-},
  <-3mm,1.5mm>*{};<-4mm,7mm>*{}**@{-},
  <-3mm,1.5mm>*{};<-2.5mm,6mm>*{}**@{-},
  <-3mm,1.5mm>*{};<-1mm,7mm>*{}**@{-},
 <-3mm,1.5mm>*{};<-2.5mm,6.6mm>*{.\hspace{-0.4mm}.\hspace{-0.4mm}.}**@{},
  <2mm,1.5mm>*{};<0mm,7mm>*{}**@{-},
  <2mm,1.5mm>*{};<1mm,7mm>*{}**@{-},
  <2mm,1.5mm>*{};<2.5mm,6mm>*{}**@{-},
  <2mm,1.5mm>*{};<4mm,7mm>*{}**@{-},
 <2mm,1.5mm>*{};<2.5mm,6.6mm>*{.\hspace{-0.4mm}.\hspace{-0.4mm}.}**@{},
 <0mm,0mm>*{};<6mm,3.6mm>*{.\hspace{-0.1mm}.\hspace{-0.1mm}.}**@{},
<10mm,1.5mm>*{};<8mm,7mm>*{}**@{-},
  <10mm,1.5mm>*{};<9mm,7mm>*{}**@{-},
  <10mm,1.5mm>*{};<10.5mm,6mm>*{}**@{-},
  <10mm,1.5mm>*{};<12mm,7mm>*{}**@{-},
 <10mm,1.5mm>*{};<10.5mm,6.6mm>*{.\hspace{-0.4mm}.\hspace{-0.4mm}.}**@{},
 <-10mm,-1.5mm>*{};<-12mm,-6mm>*{}**@{-},
 <-7mm,-1.5mm>*{};<-8mm,-6mm>*{}**@{-},
 <-4mm,-1.5mm>*{};<-4.5mm,-6mm>*{}**@{-},
 <0mm,0mm>*{};<0mm,-4.6mm>*{.\hspace{0.1mm}.\hspace{0.1mm}.}**@{},
<10mm,-1.5mm>*{};<12mm,-6mm>*{}**@{-},
 <7mm,-1.5mm>*{};<8mm,-6mm>*{}**@{-},
  <4mm,-1.5mm>*{};<4.5mm,-6mm>*{}**@{-},
  %%%%
<0mm,0mm>*{};<-9.5mm,8.2mm>*{^{I_{ 1}}}**@{},
<0mm,0mm>*{};<-3mm,8.2mm>*{^{I_{ i}}}**@{},
<0mm,0mm>*{};<2mm,8.2mm>*{^{I_{ i+1}}}**@{},
<0mm,0mm>*{};<10mm,8.2mm>*{^{I_{ k}}}**@{},
 %%%
<0mm,0mm>*{};<-12mm,-7.4mm>*{_1}**@{},
<0mm,0mm>*{};<-8mm,-7.4mm>*{_2}**@{},
<0mm,0mm>*{};<-4mm,-7.4mm>*{_3}**@{},
<0mm,0mm>*{};<6mm,-7.4mm>*{\ldots}**@{},
%<0mm,0mm>*{};<8.5mm,-7.4mm>*{_{n-1}}**@{},
<0mm,0mm>*{};<12.5mm,-7.4mm>*{_n}**@{},
\endxy
\right)
&=&\sum_{i=1}^k(-1)^{i+1}
%%%%%%%%%%%%%%%%%%%%%%%%%%%%%%%%%%%%%%%%%%%%%%%%%%%%
Q\left( \xy
%\begin{xy}
 <0mm,0mm>*
{\mbox{$\xy *=<20mm,3mm>\txt{\em
a}*\frm{-}\endxy$}};<0mm,0mm>*{}**@{},
  <-10mm,1.5mm>*{};<-12mm,7mm>*{}**@{-},
  <-10mm,1.5mm>*{};<-11mm,7mm>*{}**@{-},
  <-10mm,1.5mm>*{};<-9.5mm,6mm>*{}**@{-},
  <-10mm,1.5mm>*{};<-8mm,7mm>*{}**@{-},
 <-10mm,1.5mm>*{};<-9.5mm,6.6mm>*{.\hspace{-0.4mm}.\hspace{-0.4mm}.}**@{},
 <0mm,0mm>*{};<-5.5mm,3.6mm>*{.\hspace{-0.1mm}.\hspace{-0.1mm}.}**@{},
%
  %<0mm,1.5mm>*{};<-5mm,7mm>*{}**@{-},
  <0mm,1.5mm>*{};<-4mm,7mm>*{}**@{-},
  <0mm,1.5mm>*{};<-2.0mm,6mm>*{}**@{-},
  <0mm,1.5mm>*{};<-1mm,7mm>*{}**@{-},
 <0mm,1.5mm>*{};<-2.3mm,6.6mm>*{.\hspace{-0.4mm}.\hspace{-0.4mm}.}**@{},
%
  %<0mm,1.5mm>*{};<0mm,7mm>*{}**@{-},
  <0mm,1.5mm>*{};<1mm,7mm>*{}**@{-},
  <0mm,1.5mm>*{};<2.0mm,6mm>*{}**@{-},
  <0mm,1.5mm>*{};<4mm,7mm>*{}**@{-},
 <0mm,1.5mm>*{};<2.3mm,6.6mm>*{.\hspace{-0.4mm}.\hspace{-0.4mm}.}**@{},
 <0mm,0mm>*{};<6mm,3.6mm>*{.\hspace{-0.1mm}.\hspace{-0.1mm}.}**@{},
<10mm,1.5mm>*{};<8mm,7mm>*{}**@{-},
  <10mm,1.5mm>*{};<9mm,7mm>*{}**@{-},
  <10mm,1.5mm>*{};<10.5mm,6mm>*{}**@{-},
  <10mm,1.5mm>*{};<12mm,7mm>*{}**@{-},
 <10mm,1.5mm>*{};<10.5mm,6.6mm>*{.\hspace{-0.4mm}.\hspace{-0.4mm}.}**@{},
 <-10mm,-1.5mm>*{};<-12mm,-6mm>*{}**@{-},
 <-7mm,-1.5mm>*{};<-8mm,-6mm>*{}**@{-},
 <-4mm,-1.5mm>*{};<-4.5mm,-6mm>*{}**@{-},
 <0mm,0mm>*{};<0mm,-4.6mm>*{.\hspace{0.1mm}.\hspace{0.1mm}.}**@{},
<10mm,-1.5mm>*{};<12mm,-6mm>*{}**@{-},
 <7mm,-1.5mm>*{};<8mm,-6mm>*{}**@{-},
  <4mm,-1.5mm>*{};<4.5mm,-6mm>*{}**@{-},
  %%%%
<0mm,0mm>*{};<-9.5mm,8.2mm>*{^{I_{ 1}}}**@{},
<0mm,0mm>*{};<-2.5mm,8.2mm>*{{^{I_{ i}\sqcup}}}**@{},
<0mm,0mm>*{};<2.7mm,8.2mm>*{^{I_{ i+1}}}**@{},
<0mm,0mm>*{};<10mm,8.2mm>*{^{I_{ k}}}**@{},
 %%%
<0mm,0mm>*{};<-12mm,-7.4mm>*{_1}**@{},
<0mm,0mm>*{};<-8mm,-7.4mm>*{_2}**@{},
<0mm,0mm>*{};<-4mm,-7.4mm>*{_3}**@{},
<0mm,0mm>*{};<6mm,-7.4mm>*{\ldots}**@{},
%<0mm,0mm>*{};<8.5mm,-7.4mm>*{_{n-1}}**@{},
<0mm,0mm>*{};<12.5mm,-7.4mm>*{_n}**@{},
\endxy
\right) %%
\\
&& +\ \sum_{b+c=a\atop b,c\geq 1}\sum_{p+q=k+1\atop p\geq 1,q\geq
0}\sum_{i=0}^{p-1} \sum_{  {I_{i+1}=I_{i+1}'\sqcup I''_{i+1}\atop
.......................} \atop I_{i+q}=I_{i+q}'\sqcup I''_{i+q}}
\sum_{[n]=J_1\sqcup J_2}\sum_{s\geq 0}(-1)^{(p+1)q + i(q-1)}\\
&&
%%%%%%%%%%%%%%%%%%%%%%%%%%%%%%%%%%%%%%%%%%%%%%%%%%%%%%%%%%%
\frac{1}{s!}\ \ Q(\xy
%\begin{xy}
 <19mm,0mm>*{\mbox{$\xy *=<58mm,3mm>\txt{\em b}*\frm{-}\endxy$}};
 <0mm,0mm>*{}**@{},
  <-10mm,1.5mm>*{};<-12mm,7mm>*{}**@{-},
  <-10mm,1.5mm>*{};<-11mm,7mm>*{}**@{-},
  <-10mm,1.5mm>*{};<-9.5mm,6mm>*{}**@{-},
  <-10mm,1.5mm>*{};<-8mm,7mm>*{}**@{-},
 <-10mm,1.5mm>*{};<-9.5mm,6.6mm>*{.\hspace{-0.4mm}.\hspace{-0.4mm}.}**@{},
 <0mm,0mm>*{};<-6.5mm,3.6mm>*{.\hspace{-0.1mm}.\hspace{-0.1mm}.}**@{},
  <-3mm,1.5mm>*{};<-5mm,7mm>*{}**@{-},
  <-3mm,1.5mm>*{};<-4mm,7mm>*{}**@{-},
  <-3mm,1.5mm>*{};<-2.5mm,6mm>*{}**@{-},
  <-3mm,1.5mm>*{};<-1mm,7mm>*{}**@{-},
 <-3mm,1.5mm>*{};<-2.5mm,6.6mm>*{.\hspace{-0.4mm}.\hspace{-0.4mm}.}**@{},
%%%%%
  <10mm,1.5mm>*{};<0mm,7mm>*{}**@{-},
  <10mm,1.5mm>*{};<4mm,7mm>*{}**@{-},
<10mm,1.5mm>*{};<7.3mm,5.9mm>*{.\hspace{-0.0mm}.\hspace{-0.0mm}.}**@{},
   <10mm,1.5mm>*{};<3.8mm,6.0mm>*{}**@{-},
 <10mm,1.5mm>*{};<2.5mm,6.6mm>*{.\hspace{-0.4mm}.\hspace{-0.4mm}.}**@{},
%
% <0mm,0mm>*{};<6mm,3.6mm>*{.\hspace{-0.1mm}.\hspace{-0.1mm}.}**@{},
%
<10mm,1.5mm>*{};<9mm,7mm>*{}**@{-},
  <10mm,1.5mm>*{};<10.5mm,6mm>*{}**@{-},
  <10mm,1.5mm>*{};<12mm,7mm>*{}**@{-},
 <10mm,1.5mm>*{};<10.5mm,6.6mm>*{.\hspace{-0.4mm}.\hspace{-0.4mm}.}**@{},
 %%%%%%%%%%%%%%%%%%%%%%%%%%%%%%%%%%%%%%%
%%%% lower legs of 1st vertex %%%%%%%
 <-10mm,-1.5mm>*{};<-12mm,-6mm>*{}**@{-},
 <-7mm,-1.5mm>*{};<-8mm,-6mm>*{}**@{-},
 <-4mm,-1.5mm>*{};<-4.5mm,-6mm>*{}**@{-},
 <-1mm,-1.5mm>*{};<-1.1mm,-6mm>*{}**@{-},
 <2mm,-1.5mm>*{};<2.0mm,-6mm>*{}**@{-},
 <0mm,0mm>*{};<20mm,-4.6mm>*{.\hspace{2mm}.\hspace{2mm}.\hspace{2mm}
 .\hspace{2mm}.\hspace{2mm}.\hspace{2mm}.\hspace{2mm}
 .\hspace{2mm}.\hspace{2mm}}**@{},
<48mm,-1.5mm>*{};<50mm,-6mm>*{}**@{-},
 <45mm,-1.5mm>*{};<46mm,-6mm>*{}**@{-},
  <42mm,-1.5mm>*{};<42.5mm,-6mm>*{}**@{-},
  <39mm,-1.5mm>*{};<39.2mm,-6mm>*{}**@{-},
   <36mm,-1.5mm>*{};<36mm,-6mm>*{}**@{-},
 <20mm,-1.5mm>*{};<20.0mm,-8mm>*{\underbrace{\hspace{66mm}}}**@{},
 <20mm,-1.5mm>*{};<20.0mm,-11mm>*{_{J_1}}**@{},
  %%%%
<0mm,0mm>*{};<-9.5mm,8.4mm>*{^{I_{ 1}}}**@{},
<0mm,0mm>*{};<-3mm,8.4mm>*{^{I_{ i}}}**@{},
<0mm,0mm>*{};<2mm,8.6mm>*{^{I_{ i+1}'}}**@{},
<0mm,0mm>*{};<10.5mm,8.6mm>*{^{I_{ i+q}'}}**@{},
%%%% edges connecting two vertices %%%%%%%%%%
<10mm,1.5mm>*{};<18mm,12mm>*{}**@{-},
<10mm,1.5mm>*{};<20.0mm,12mm>*{}**@{-},
<10mm,1.5mm>*{};<25mm,12mm>*{}**@{-},
<10mm,1.5mm>*{};<18.7mm,8.6mm>*{.\hspace{-0.4mm}.\hspace{-0.4mm}.}**@{},
<10mm,1.5mm>*{};<19.6mm,9.0mm>*{^s}**@{},
%%%%%%%%%%%%%%%%%%%%%%%%%%%%%%%%%%%%%%%%%%%%%%
%%%% second vertex %%%%%%%%%%%%%%%%%%%%%%%%%%
%%%%%%%%%%%%%%%%%%%%%%%%%%%%%%%%%%%%%%%%%%%%%
<25mm,13.75mm>*{\mbox{$\xy *=<14mm,3mm>\txt{\em c}*\frm{-}\endxy$}};
<0mm,0mm>*{}**@{},
 <18mm,15mm>*{};<16mm,20.5mm>*{}**@{-},
 <18mm,15mm>*{};<17mm,20.5mm>*{}**@{-},
 <18mm,15mm>*{};<18.5mm,19.6mm>*{}**@{-},
 <18mm,15mm>*{};<20mm,20.5mm>*{}**@{-},
 <18mm,15mm>*{};<18.6mm,20.3mm>*{.\hspace{-0.4mm}.\hspace{-0.4mm}.}**@{},
<18mm,15mm>*{};<14mm,14mm>*{Q(}**@{},
<18mm,15mm>*{};<34mm,14mm>*{)}**@{},
<22mm,15mm>*{};<25.5mm,17.7mm>*{\cdots}**@{},
 <32mm,15.2mm>*{};<30mm,20.5mm>*{}**@{-},
 <32mm,15.2mm>*{};<31mm,20.5mm>*{}**@{-},
 <32mm,15.2mm>*{};<32.5mm,19.6mm>*{}**@{-},
 <32mm,15mm>*{};<34mm,20.5mm>*{}**@{-},
 <32mm,15mm>*{};<32.3mm,20.3mm>*{.\hspace{-0.4mm}.\hspace{-0.4mm}.}**@{},
%
  %%%%
<0mm,0mm>*{};<18mm,22.6mm>*{^{I_{ i+1}''}}**@{},
<0mm,0mm>*{};<32.5mm,22.6mm>*{^{I_{ i+q}''}}**@{},
  %%%%
%%%% lower legs of 2nd vertex %%%%%%
 <26mm,12mm>*{};<25mm,9mm>*{}**@{-},
 <27mm,12mm>*{};<26.8mm,9mm>*{}**@{-},
 <29mm,12mm>*{};<29.2mm,10mm>*{.\hspace{-0.1mm}.\hspace{-0.1mm}.}**@{},
 <32mm,12mm>*{};<33mm,9mm>*{}**@{-},
  <31mm,12mm>*{};<31.5mm,9mm>*{}**@{-},
 <29mm,12mm>*{};<29mm,8mm>*{\underbrace{\ \ \ \ \ \ \ \  }}**@{},
 <29mm,12mm>*{};<29mm,5.3mm>*{_{J_2}}**@{},
%%%%%%%%%%%%%%%%%% 1st vertex contd %%%%%%%%%%%
<38mm,1.5mm>*{};<36mm,7mm>*{}**@{-},
<38mm,1.5mm>*{};<37mm,7mm>*{}**@{-},
<38mm,1.5mm>*{};<38.5mm,6mm>*{}**@{-},
<38mm,1.5mm>*{};<40mm,7mm>*{}**@{-},
<38mm,1.5mm>*{};<38.5mm,6.6mm>*{.\hspace{-0.4mm}.\hspace{-0.4mm}.}**@{},
<38mm,1.5mm>*{};<43mm,4mm>*{.\hspace{-0.0mm}.\hspace{-0.0mm}.}**@{},
<48mm,1.5mm>*{};<46mm,7mm>*{}**@{-},
<48mm,1.5mm>*{};<47mm,7mm>*{}**@{-},
<48mm,1.5mm>*{};<48.5mm,6mm>*{}**@{-},
<48mm,1.5mm>*{};<50mm,7mm>*{}**@{-},
<48mm,1.5mm>*{};<48.5mm,6.6mm>*{.\hspace{-0.4mm}.\hspace{-0.4mm}.}**@{},
<0mm,0mm>*{};<40.3mm,8.6mm>*{^{I_{ i+q+1}}}**@{},
<0mm,0mm>*{};<48.5mm,8.6mm>*{^{I_{ k}}}**@{},
  %%%%
\endxy
) \hspace{29mm} (\bigstar) \Eeqrn Note that the r.h.s.\ of the above
equation contains values of $Q$ on $[a',k',n']$-corollas with the
weight $2a'+k'-2< 2a+k-2$ which are, by induction assumption, are
already known.

\bip

%%%%%%%%%%%%%%%%%%%%%%%%%%%%%%%%%%%%%%%%%%%%%%%%%%%%%%%
%\noindent{\bf II.3.1. Remark.} The sum over $s$ in equation $(\bigstar)$ is
%always {\em finite}: due to the skew symmetry of output legs and
%symmetry of input legs of the generating corollas in
%$\LieBi_\infty$, the number of internal edges in an element $G\in
%\Liebi_\infty^\circlearrowright$ with $l$ vertices can {\em not}\,
%be larger than  $2l$ (no two vertices can be connected by two
%internal edges with the same orientation). In fact, by Proposition
%III.1(ii) below, we can replace in $(\bigstar)$ infinite sum $\sum_{s\geq
%0}$ with finite sum $\sum_{s=0}^{2a}$.

%\bip

%%%%%%%%%%%%%%%%%%%%%%%%%%%%%%%%%%%%%%%%%%%%%%%%%%%%
\noindent{\bf II.4. How it works.} We shall illustrate the above construction
of $Q$ in a few examples. \sip

\noindent{\bf Iteration level  s=3}: We have
\[
E_3=\mbox{\sf span} \left\{ \xy
 <0mm,0mm>*{\mbox{$\xy *=<5mm,3mm>\txt{\em 2}*\frm{-}\endxy$}};<0mm,0mm>*{}**@{},
 <-2.5mm,-1.6mm>*{};<-4mm,-6mm>*{}**@{-},
 <-1.6mm,-1.6mm>*{};<-2mm,-6mm>*{}**@{-},
 <1.6mm,-1.6mm>*{};<2.5mm,-6mm>*{}**@{-},
 <2.5mm,-1.6mm>*{};<4mm,-6mm>*{}**@{-},
 <0mm,0mm>*{};<0mm,-4.6mm>*{.\hspace{0mm}.\hspace{0mm}.}**@{},
<0mm,0mm>*{};<-4.1mm,-7.4mm>*{_1}**@{},
<0mm,0mm>*{};<-2mm,-7.4mm>*{_2}**@{},
%<0mm,0mm>*{};<-4mm,-7.4mm>*{_3}**@{},
<0mm,0mm>*{};<1.7mm,-7.4mm>*{...}**@{},
%<0mm,0mm>*{};<8.5mm,-7.4mm>*{_{n-1}}**@{},
<0mm,0mm>*{};<4.6mm,-7.4mm>*{_n}**@{},
<0mm,1.5mm>*{};<2mm,7mm>*{}**@{-},
  <0mm,1.5mm>*{};<1mm,6mm>*{}**@{-},
  <0mm,1.5mm>*{};<-0.5mm,7mm>*{}**@{-},
  <0mm,1.5mm>*{};<-2mm,7mm>*{}**@{-},
 <0mm,1.5mm>*{};<0.7mm,6.6mm>*{.\hspace{-0.4mm}.\hspace{-0.4mm}.}**@{},
\endxy
\ \ , \ \ \xy
 <0mm,0mm>*{\mbox{$\xy *=<10mm,3mm>\txt{\em 1}*\frm{-}\endxy$}};<0mm,0mm>*{}**@{},
 <-2.5mm,-1.6mm>*{};<-4mm,-6mm>*{}**@{-},
 <-1.6mm,-1.6mm>*{};<-2mm,-6mm>*{}**@{-},
 <1.6mm,-1.6mm>*{};<2.5mm,-6mm>*{}**@{-},
 <2.5mm,-1.6mm>*{};<4mm,-6mm>*{}**@{-},
 <0mm,0mm>*{};<0mm,-4.6mm>*{.\hspace{0mm}.\hspace{0mm}.}**@{},
<0mm,0mm>*{};<-4.1mm,-7.4mm>*{_1}**@{},
<0mm,0mm>*{};<-2mm,-7.4mm>*{_2}**@{},
%<0mm,0mm>*{};<-4mm,-7.4mm>*{_3}**@{},
<0mm,0mm>*{};<1.7mm,-7.4mm>*{...}**@{},
%<0mm,0mm>*{};<8.5mm,-7.4mm>*{_{n-1}}**@{},
<0mm,0mm>*{};<4.6mm,-7.4mm>*{_n}**@{},
<0mm,1.5mm>*{};<2mm,7mm>*{}**@{-},
  <0mm,1.5mm>*{};<1mm,6mm>*{}**@{-},
  <0mm,1.5mm>*{};<-0.5mm,7mm>*{}**@{-},
  <0mm,1.5mm>*{};<-2mm,7mm>*{}**@{-},
 <0mm,1.5mm>*{};<0.7mm,6.6mm>*{.\hspace{-0.4mm}.\hspace{-0.4mm}.}**@{},
 <-5mm,1.5mm>*{};<-8.5mm,7mm>*{}**@{-},
  <-5mm,1.5mm>*{};<-7.5mm,7mm>*{}**@{-},
  <-5mm,1.5mm>*{};<-6mm,6mm>*{}**@{-},
  <-5mm,1.5mm>*{};<-4.5mm,7mm>*{}**@{-},
 <-3.5mm,1.5mm>*{};<-6mm,6.6mm>*{.\hspace{-0.4mm}.\hspace{-0.4mm}.}**@{},
 <5mm,1.5mm>*{};<8.5mm,7mm>*{}**@{-},
  <5mm,1.5mm>*{};<7.5mm,7mm>*{}**@{-},
  <5mm,1.5mm>*{};<6mm,6mm>*{}**@{-},
  <5mm,1.5mm>*{};<4.5mm,7mm>*{}**@{-},
 <3.5mm,1.5mm>*{};<6mm,6.6mm>*{.\hspace{-0.4mm}.\hspace{-0.4mm}.}**@{},
\endxy
\right\}
\]
The equation $(\bigstar)$  for the first generator takes the form,
\Beqrn \delta e''= Q\left(\delta \xy
 <0mm,0mm>*{\mbox{$\xy *=<5mm,3mm>\txt{\em 2}*\frm{-}\endxy$}};<0mm,0mm>*{}**@{},
 <-2.5mm,-1.6mm>*{};<-4mm,-6mm>*{}**@{-},
 <-1.6mm,-1.6mm>*{};<-2mm,-6mm>*{}**@{-},
 <1.6mm,-1.6mm>*{};<2.5mm,-6mm>*{}**@{-},
 <2.5mm,-1.6mm>*{};<4mm,-6mm>*{}**@{-},
 <0mm,0mm>*{};<0mm,-4.6mm>*{.\hspace{0mm}.\hspace{0mm}.}**@{},
<0mm,0mm>*{};<-4.1mm,-7.4mm>*{_1}**@{},
<0mm,0mm>*{};<-2mm,-7.4mm>*{_2}**@{},
%<0mm,0mm>*{};<-4mm,-7.4mm>*{_3}**@{},
<0mm,0mm>*{};<1.7mm,-7.4mm>*{...}**@{},
%<0mm,0mm>*{};<8.5mm,-7.4mm>*{_{n-1}}**@{},
<0mm,0mm>*{};<4.6mm,-7.4mm>*{_n}**@{},
<0mm,1.5mm>*{};<2mm,7mm>*{}**@{-},
  <0mm,1.5mm>*{};<1mm,6mm>*{}**@{-},
  <0mm,1.5mm>*{};<-0.5mm,7mm>*{}**@{-},
  <0mm,1.5mm>*{};<-2mm,7mm>*{}**@{-},
 <0mm,1.5mm>*{};<0.7mm,6.6mm>*{.\hspace{-0.4mm}.\hspace{-0.4mm}.}**@{},
<0mm,1.5mm>*{};<0mm,8.6mm>*{_{I_1}}**@{},
\endxy
\right) &=& \sum_{I_{1}=I_{1}'\sqcup I''_{1} \atop [n]=J_1\sqcup
J_2}\sum_{s\geq 0}
%%%%%%%%%%%%%%%%%%%%%%%%%%%%%%%%%%%%%%%%%%%%%%%%%%%%%%%%%%%
\frac{1}{s!}\ \ \xy
%\begin{xy}
 <10mm,0mm>*{\mbox{$\xy *=<5mm,3mm>\txt{\em 1}*\frm{-}\endxy$}};
<0mm,0mm>*{};<4mm,0mm>*{Q(}**@{},
<0mm,0mm>*{};<15.5mm,0mm>*{)}**@{},
<10mm,1.5mm>*{};<8mm,7mm>*{}**@{-},
  <10mm,1.5mm>*{};<10.5mm,6mm>*{}**@{-},
  <10mm,1.5mm>*{};<12mm,7mm>*{}**@{-},
 <10mm,1.5mm>*{};<10.5mm,6.6mm>*{.\hspace{-0.4mm}.\hspace{-0.4mm}.}**@{},
 %%%%%%%%%%%%%%%%%%%%%%%%%%%%%%%%%%%%%%%
%%%% lower legs of 1st vertex %%%%%%%
 <7.5mm,-1.5mm>*{};<6mm,-6mm>*{}**@{-},
 <9mm,-1.5mm>*{};<8mm,-6mm>*{}**@{-},
 <12.5mm,-1.5mm>*{};<14mm,-6mm>*{}**@{-},
 <11mm,-1.5mm>*{};<12mm,-6mm>*{}**@{-},
 <10mm,-1.5mm>*{};<10mm,-6mm>*{...}**@{},
 <20mm,-1.5mm>*{};<10.0mm,-8mm>*{\underbrace{\hspace{10mm}}}**@{},
 <20mm,-1.5mm>*{};<10.0mm,-11mm>*{_{J_1}}**@{},
  %%%%
<0mm,0mm>*{};<10.5mm,8.6mm>*{^{I_{ 1}'}}**@{},
%%%% edges connecting two vertices %%%%%%%%%%
<10mm,1.5mm>*{};<18mm,12mm>*{}**@{-},
<10mm,1.5mm>*{};<20.0mm,12mm>*{}**@{-},
<10mm,1.5mm>*{};<25mm,12mm>*{}**@{-},
<10mm,1.5mm>*{};<18.7mm,8.6mm>*{.\hspace{-0.4mm}.\hspace{-0.4mm}.}**@{},
<10mm,1.5mm>*{};<19.6mm,9.0mm>*{^s}**@{},
%%%%%%%%%%%%%%%%%%%%%%%%%%%%%%%%%%%%%%%%%%%%%%
%%%% second vertex %%%%%%%%%%%%%%%%%%%%%%%%%%
%%%%%%%%%%%%%%%%%%%%%%%%%%%%%%%%%%%%%%%%%%%%%
<25mm,13.75mm>*{\mbox{$\xy *=<14mm,3mm>\txt{\em 1}*\frm{-}\endxy$}};
<0mm,0mm>*{}**@{},
 <25mm,15mm>*{};<23mm,20.5mm>*{}**@{-},
 <25mm,15mm>*{};<24mm,20.5mm>*{}**@{-},
 <25mm,15mm>*{};<25.5mm,19.6mm>*{}**@{-},
 <25mm,15mm>*{};<27mm,20.5mm>*{}**@{-},
 <25mm,15mm>*{};<25.6mm,20.3mm>*{.\hspace{-0.4mm}.\hspace{-0.4mm}.}**@{},
<18mm,15mm>*{};<15mm,14mm>*{Q(}**@{},
<18mm,15mm>*{};<34mm,14mm>*{)}**@{},
%
  %%%%
<0mm,0mm>*{};<25mm,22.6mm>*{^{I_{ 1}''}}**@{},
  %%%%
%%%% lower legs of 2nd vertex %%%%%%
 <26mm,12mm>*{};<25mm,9mm>*{}**@{-},
 <27mm,12mm>*{};<26.8mm,9mm>*{}**@{-},
 <29mm,12mm>*{};<29.2mm,10mm>*{.\hspace{-0.1mm}.\hspace{-0.1mm}.}**@{},
 <32mm,12mm>*{};<33mm,9mm>*{}**@{-},
  <31mm,12mm>*{};<31.5mm,9mm>*{}**@{-},
 <29mm,12mm>*{};<29mm,8mm>*{\underbrace{\ \ \ \ \ \ \ \  }}**@{},
 <29mm,12mm>*{};<29mm,5.3mm>*{_{J_2}}**@{},
%%%%%%%%%%%%%%%%%% 1st vertex contd %%%%%%%%%%%
\endxy
\\
&=&\left\{ \Ba{rr}
\begin{xy}
 <0mm,0mm>*{\bullet};<0mm,0mm>*{}**@{},
 <0mm,0.69mm>*{};<0mm,3.0mm>*{}**@{-},
 <0.39mm,-0.39mm>*{};<2.4mm,-2.4mm>*{}**@{-},
 <-0.35mm,-0.35mm>*{};<-1.9mm,-1.9mm>*{}**@{-},
 <-2.4mm,-2.4mm>*{\bullet};<-2.4mm,-2.4mm>*{}**@{},
 <-2.0mm,-2.8mm>*{};<0mm,-4.9mm>*{}**@{-},
 <-2.8mm,-2.9mm>*{};<-4.7mm,-4.9mm>*{}**@{-},
    <0.39mm,-0.39mm>*{};<3.3mm,-4.0mm>*{^3}**@{},
    <-2.0mm,-2.8mm>*{};<0.5mm,-6.7mm>*{^2}**@{},
    <-2.8mm,-2.9mm>*{};<-5.2mm,-6.7mm>*{^1}**@{},
 \end{xy}
+
 \begin{xy}
 <0mm,0mm>*{\bullet};<0mm,0mm>*{}**@{},
 <0mm,0.69mm>*{};<0mm,3.0mm>*{}**@{-},
 <0.39mm,-0.39mm>*{};<2.4mm,-2.4mm>*{}**@{-},
 <-0.35mm,-0.35mm>*{};<-1.9mm,-1.9mm>*{}**@{-},
 <-2.4mm,-2.4mm>*{\bullet};<-2.4mm,-2.4mm>*{}**@{},
 <-2.0mm,-2.8mm>*{};<0mm,-4.9mm>*{}**@{-},
 <-2.8mm,-2.9mm>*{};<-4.7mm,-4.9mm>*{}**@{-},
    <0.39mm,-0.39mm>*{};<3.3mm,-4.0mm>*{^2}**@{},
    <-2.0mm,-2.8mm>*{};<0.5mm,-6.7mm>*{^1}**@{},
    <-2.8mm,-2.9mm>*{};<-5.2mm,-6.7mm>*{^3}**@{},
 \end{xy}
+
 \begin{xy}
 <0mm,0mm>*{\bullet};<0mm,0mm>*{}**@{},
 <0mm,0.69mm>*{};<0mm,3.0mm>*{}**@{-},
 <0.39mm,-0.39mm>*{};<2.4mm,-2.4mm>*{}**@{-},
 <-0.35mm,-0.35mm>*{};<-1.9mm,-1.9mm>*{}**@{-},
 <-2.4mm,-2.4mm>*{\bullet};<-2.4mm,-2.4mm>*{}**@{},
 <-2.0mm,-2.8mm>*{};<0mm,-4.9mm>*{}**@{-},
 <-2.8mm,-2.9mm>*{};<-4.7mm,-4.9mm>*{}**@{-},
    <0.39mm,-0.39mm>*{};<3.3mm,-4.0mm>*{^1}**@{},
    <-2.0mm,-2.8mm>*{};<0.5mm,-6.7mm>*{^3}**@{},
    <-2.8mm,-2.9mm>*{};<-5.2mm,-6.7mm>*{^2}**@{},
 \end{xy} & \mbox{for}\ |I_1|=1,n=3 \\
0 & \mbox{otherwise.} \Ea \right. \Eeqrn implying
$$
e''= \left\{ \Ba{rr}
\begin{xy}
 <0mm,0mm>*{\bullet};<0mm,0mm>*{}**@{},
 <0mm,0.69mm>*{};<0mm,4.0mm>*{}**@{-},
 <0mm,0mm>*{};<4mm,-4mm>*{}**@{-},
 <0mm,0mm>*{};<-4mm,-4mm>*{}**@{-},
 <0mm,0mm>*{};<0mm,-4mm>*{}**@{-},
   <0.39mm,-0.39mm>*{};<4.9mm,-6.4mm>*{^3}**@{},
    <-2.0mm,-2.8mm>*{};<0mm,-6.4mm>*{^2}**@{},
    <-2.8mm,-2.9mm>*{};<-4.9mm,-6.4mm>*{^1}**@{},
 \end{xy}
 & \mbox{for}\ |I_1|=1,n=3 \\
0 & \mbox{otherwise.} \Ea \right.
$$
As
$$
q\left( \xy
 <0mm,0mm>*{\mbox{$\xy *=<5mm,3mm>\txt{\em 2}*\frm{-}\endxy$}};<0mm,0mm>*{}**@{},
 <-2.5mm,-1.6mm>*{};<-4mm,-6mm>*{}**@{-},
 <-1.6mm,-1.6mm>*{};<-2mm,-6mm>*{}**@{-},
 <1.6mm,-1.6mm>*{};<2.5mm,-6mm>*{}**@{-},
 <2.5mm,-1.6mm>*{};<4mm,-6mm>*{}**@{-},
 <0mm,0mm>*{};<0mm,-4.6mm>*{.\hspace{0mm}.\hspace{0mm}.}**@{},
<0mm,0mm>*{};<-4.1mm,-7.4mm>*{_1}**@{},
<0mm,0mm>*{};<-2mm,-7.4mm>*{_2}**@{},
%<0mm,0mm>*{};<-4mm,-7.4mm>*{_3}**@{},
<0mm,0mm>*{};<1.7mm,-7.4mm>*{...}**@{},
%<0mm,0mm>*{};<8.5mm,-7.4mm>*{_{n-1}}**@{},
<0mm,0mm>*{};<4.6mm,-7.4mm>*{_n}**@{},
<0mm,1.5mm>*{};<2mm,7mm>*{}**@{-},
  <0mm,1.5mm>*{};<1mm,6mm>*{}**@{-},
  <0mm,1.5mm>*{};<-0.5mm,7mm>*{}**@{-},
  <0mm,1.5mm>*{};<-2mm,7mm>*{}**@{-},
 <0mm,1.5mm>*{};<0.7mm,6.6mm>*{.\hspace{-0.4mm}.\hspace{-0.4mm}.}**@{},
<0mm,1.5mm>*{};<0mm,8.6mm>*{_{I_1}}**@{},
\endxy
\right) = \left\{ \Ba{rr} -\begin{xy}
 <0mm,-1.3mm>*{};<0mm,-3.5mm>*{}**@{-},
 <0.38mm,-0.2mm>*{};<2.2mm,2.2mm>*{}**@{-},
 <-0.38mm,-0.2mm>*{};<-2.2mm,2.2mm>*{}**@{-},
<0mm,-0.8mm>*{\circ};<0mm,0.8mm>*{}**@{},
% <-2.25mm,2.2mm>*{};<-2.2mm,5.2mm>*{}**@{-},
 <2.4mm,2.4mm>*{\bullet};<2.4mm,2.4mm>*{}**@{},
 <2.5mm,2.3mm>*{};<4.4mm,-0.8mm>*{}**@{-},
% <4.4mm,-0.8mm>*{};<4.4mm,-3.5mm>*{}**@{-},
 <2.4mm,2.5mm>*{};<2.4mm,5.2mm>*{}**@{-},
     <0mm,-1.3mm>*{};<0mm,-5.3mm>*{^1}**@{},
     <2.5mm,2.3mm>*{};<5.1mm,-2.6mm>*{^2}**@{},
    <2.4mm,2.5mm>*{};<2.4mm,5.7mm>*{^2}**@{},
    <-0.38mm,-0.2mm>*{};<-2.8mm,2.5mm>*{^1}**@{},
    \end{xy}
-
\begin{xy}
 <0mm,-1.3mm>*{};<0mm,-3.5mm>*{}**@{-},
 <0.38mm,-0.2mm>*{};<2.2mm,2.2mm>*{}**@{-},
 <-0.38mm,-0.2mm>*{};<-2.2mm,2.2mm>*{}**@{-},
<0mm,-0.8mm>*{\circ};<0mm,0.8mm>*{}**@{},
% <-2.25mm,2.2mm>*{};<-2.2mm,5.2mm>*{}**@{-},
 <2.4mm,2.4mm>*{\bullet};<2.4mm,2.4mm>*{}**@{},
 <2.5mm,2.3mm>*{};<4.4mm,-0.8mm>*{}**@{-},
% <4.4mm,-0.8mm>*{};<4.4mm,-3.5mm>*{}**@{-},
 <2.4mm,2.5mm>*{};<2.4mm,5.2mm>*{}**@{-},
     <0mm,-1.3mm>*{};<0mm,-5.3mm>*{^1}**@{},
     <2.5mm,2.3mm>*{};<5.1mm,-2.6mm>*{^2}**@{},
    <2.4mm,2.5mm>*{};<2.4mm,5.7mm>*{^1}**@{},
    <-0.38mm,-0.2mm>*{};<-2.8mm,2.5mm>*{^2}**@{},
    \end{xy}
-
\begin{xy}
 <0mm,-1.3mm>*{};<0mm,-3.5mm>*{}**@{-},
 <0.38mm,-0.2mm>*{};<2.2mm,2.2mm>*{}**@{-},
 <-0.38mm,-0.2mm>*{};<-2.2mm,2.2mm>*{}**@{-},
<0mm,-0.8mm>*{\circ};<0mm,0.8mm>*{}**@{},
% <-2.25mm,2.2mm>*{};<-2.2mm,5.2mm>*{}**@{-},
 <2.4mm,2.4mm>*{\bullet};<2.4mm,2.4mm>*{}**@{},
 <2.5mm,2.3mm>*{};<4.4mm,-0.8mm>*{}**@{-},
% <4.4mm,-0.8mm>*{};<4.4mm,-3.5mm>*{}**@{-},
 <2.4mm,2.5mm>*{};<2.4mm,5.2mm>*{}**@{-},
     <0mm,-1.3mm>*{};<0mm,-5.3mm>*{^2}**@{},
     <2.5mm,2.3mm>*{};<5.1mm,-2.6mm>*{^1}**@{},
    <2.4mm,2.5mm>*{};<2.4mm,5.7mm>*{^2}**@{},
    <-0.38mm,-0.2mm>*{};<-2.8mm,2.5mm>*{^1}**@{},
    \end{xy}
-
\begin{xy}
 <0mm,-1.3mm>*{};<0mm,-3.5mm>*{}**@{-},
 <0.38mm,-0.2mm>*{};<2.2mm,2.2mm>*{}**@{-},
 <-0.38mm,-0.2mm>*{};<-2.2mm,2.2mm>*{}**@{-},
<0mm,-0.8mm>*{\circ};<0mm,0.8mm>*{}**@{},
% <-2.25mm,2.2mm>*{};<-2.2mm,5.2mm>*{}**@{-},
 <2.4mm,2.4mm>*{\bullet};<2.4mm,2.4mm>*{}**@{},
 <2.5mm,2.3mm>*{};<4.4mm,-0.8mm>*{}**@{-},
% <4.4mm,-0.8mm>*{};<4.4mm,-3.5mm>*{}**@{-},
 <2.4mm,2.5mm>*{};<2.4mm,5.2mm>*{}**@{-},
     <0mm,-1.3mm>*{};<0mm,-5.3mm>*{^2}**@{},
     <2.5mm,2.3mm>*{};<5.1mm,-2.6mm>*{^1}**@{},
    <2.4mm,2.5mm>*{};<2.4mm,5.7mm>*{^1}**@{},
    <-0.38mm,-0.2mm>*{};<-2.8mm,2.5mm>*{^2}**@{},
    \end{xy}
 & \mbox{for}\ |I_1|=2,n=2 \\
0 & \mbox{otherwise.} \Ea \right.
$$
and identifying the r.h.s.\ with its natural lift, $e'$, into
$\Liebi_\infty$, we finally obtain
$$
Q\left( \xy
 <0mm,0mm>*{\mbox{$\xy *=<5mm,3mm>\txt{\em 2}*\frm{-}\endxy$}};<0mm,0mm>*{}**@{},
 <-2.5mm,-1.6mm>*{};<-4mm,-6mm>*{}**@{-},
 <-1.6mm,-1.6mm>*{};<-2mm,-6mm>*{}**@{-},
 <1.6mm,-1.6mm>*{};<2.5mm,-6mm>*{}**@{-},
 <2.5mm,-1.6mm>*{};<4mm,-6mm>*{}**@{-},
 <0mm,0mm>*{};<0mm,-4.6mm>*{.\hspace{0mm}.\hspace{0mm}.}**@{},
<0mm,0mm>*{};<-4.1mm,-7.4mm>*{_1}**@{},
<0mm,0mm>*{};<-2mm,-7.4mm>*{_2}**@{},
%<0mm,0mm>*{};<-4mm,-7.4mm>*{_3}**@{},
<0mm,0mm>*{};<1.7mm,-7.4mm>*{...}**@{},
%<0mm,0mm>*{};<8.5mm,-7.4mm>*{_{n-1}}**@{},
<0mm,0mm>*{};<4.6mm,-7.4mm>*{_n}**@{},
<0mm,1.5mm>*{};<2mm,7mm>*{}**@{-},
  <0mm,1.5mm>*{};<1mm,6mm>*{}**@{-},
  <0mm,1.5mm>*{};<-0.5mm,7mm>*{}**@{-},
  <0mm,1.5mm>*{};<-2mm,7mm>*{}**@{-},
 <0mm,1.5mm>*{};<0.7mm,6.6mm>*{.\hspace{-0.4mm}.\hspace{-0.4mm}.}**@{},
<0mm,1.5mm>*{};<0mm,8.6mm>*{_{I_1}}**@{},
\endxy
\right) = \left\{ \Ba{rr}
\begin{xy}
 <0mm,0mm>*{\bullet};<0mm,0mm>*{}**@{},
 <0mm,0.69mm>*{};<0mm,4.0mm>*{}**@{-},
 <0mm,0mm>*{};<4mm,-4mm>*{}**@{-},
 <0mm,0mm>*{};<-4mm,-4mm>*{}**@{-},
 <0mm,0mm>*{};<0mm,-4mm>*{}**@{-},
   <0.39mm,-0.39mm>*{};<4.9mm,-6.4mm>*{^3}**@{},
    <-2.0mm,-2.8mm>*{};<0mm,-6.4mm>*{^2}**@{},
    <-2.8mm,-2.9mm>*{};<-4.9mm,-6.4mm>*{^1}**@{},
 \end{xy}
 & \mbox{for}\ |I_1|=1,n=3 \\
-\begin{xy}
 <0mm,-1.3mm>*{};<0mm,-3.5mm>*{}**@{-},
 <0.38mm,-0.2mm>*{};<2.2mm,2.2mm>*{}**@{-},
 <-0.38mm,-0.2mm>*{};<-2.2mm,2.2mm>*{}**@{-},
<0mm,-0.8mm>*{\circ};<0mm,0.8mm>*{}**@{},
% <-2.25mm,2.2mm>*{};<-2.2mm,5.2mm>*{}**@{-},
 <2.4mm,2.4mm>*{\bullet};<2.4mm,2.4mm>*{}**@{},
 <2.5mm,2.3mm>*{};<4.4mm,-0.8mm>*{}**@{-},
% <4.4mm,-0.8mm>*{};<4.4mm,-3.5mm>*{}**@{-},
 <2.4mm,2.5mm>*{};<2.4mm,5.2mm>*{}**@{-},
     <0mm,-1.3mm>*{};<0mm,-5.3mm>*{^1}**@{},
     <2.5mm,2.3mm>*{};<5.1mm,-2.6mm>*{^2}**@{},
    <2.4mm,2.5mm>*{};<2.4mm,5.7mm>*{^2}**@{},
    <-0.38mm,-0.2mm>*{};<-2.8mm,2.5mm>*{^1}**@{},
    \end{xy}
-
\begin{xy}
 <0mm,-1.3mm>*{};<0mm,-3.5mm>*{}**@{-},
 <0.38mm,-0.2mm>*{};<2.2mm,2.2mm>*{}**@{-},
 <-0.38mm,-0.2mm>*{};<-2.2mm,2.2mm>*{}**@{-},
<0mm,-0.8mm>*{\circ};<0mm,0.8mm>*{}**@{},
% <-2.25mm,2.2mm>*{};<-2.2mm,5.2mm>*{}**@{-},
 <2.4mm,2.4mm>*{\bullet};<2.4mm,2.4mm>*{}**@{},
 <2.5mm,2.3mm>*{};<4.4mm,-0.8mm>*{}**@{-},
% <4.4mm,-0.8mm>*{};<4.4mm,-3.5mm>*{}**@{-},
 <2.4mm,2.5mm>*{};<2.4mm,5.2mm>*{}**@{-},
     <0mm,-1.3mm>*{};<0mm,-5.3mm>*{^1}**@{},
     <2.5mm,2.3mm>*{};<5.1mm,-2.6mm>*{^2}**@{},
    <2.4mm,2.5mm>*{};<2.4mm,5.7mm>*{^1}**@{},
    <-0.38mm,-0.2mm>*{};<-2.8mm,2.5mm>*{^2}**@{},
    \end{xy}
-
\begin{xy}
 <0mm,-1.3mm>*{};<0mm,-3.5mm>*{}**@{-},
 <0.38mm,-0.2mm>*{};<2.2mm,2.2mm>*{}**@{-},
 <-0.38mm,-0.2mm>*{};<-2.2mm,2.2mm>*{}**@{-},
<0mm,-0.8mm>*{\circ};<0mm,0.8mm>*{}**@{},
% <-2.25mm,2.2mm>*{};<-2.2mm,5.2mm>*{}**@{-},
 <2.4mm,2.4mm>*{\bullet};<2.4mm,2.4mm>*{}**@{},
 <2.5mm,2.3mm>*{};<4.4mm,-0.8mm>*{}**@{-},
% <4.4mm,-0.8mm>*{};<4.4mm,-3.5mm>*{}**@{-},
 <2.4mm,2.5mm>*{};<2.4mm,5.2mm>*{}**@{-},
     <0mm,-1.3mm>*{};<0mm,-5.3mm>*{^2}**@{},
     <2.5mm,2.3mm>*{};<5.1mm,-2.6mm>*{^1}**@{},
    <2.4mm,2.5mm>*{};<2.4mm,5.7mm>*{^2}**@{},
    <-0.38mm,-0.2mm>*{};<-2.8mm,2.5mm>*{^1}**@{},
    \end{xy}
-
\begin{xy}
 <0mm,-1.3mm>*{};<0mm,-3.5mm>*{}**@{-},
 <0.38mm,-0.2mm>*{};<2.2mm,2.2mm>*{}**@{-},
 <-0.38mm,-0.2mm>*{};<-2.2mm,2.2mm>*{}**@{-},
<0mm,-0.8mm>*{\circ};<0mm,0.8mm>*{}**@{},
% <-2.25mm,2.2mm>*{};<-2.2mm,5.2mm>*{}**@{-},
 <2.4mm,2.4mm>*{\bullet};<2.4mm,2.4mm>*{}**@{},
 <2.5mm,2.3mm>*{};<4.4mm,-0.8mm>*{}**@{-},
% <4.4mm,-0.8mm>*{};<4.4mm,-3.5mm>*{}**@{-},
 <2.4mm,2.5mm>*{};<2.4mm,5.2mm>*{}**@{-},
     <0mm,-1.3mm>*{};<0mm,-5.3mm>*{^2}**@{},
     <2.5mm,2.3mm>*{};<5.1mm,-2.6mm>*{^1}**@{},
    <2.4mm,2.5mm>*{};<2.4mm,5.7mm>*{^1}**@{},
    <-0.38mm,-0.2mm>*{};<-2.8mm,2.5mm>*{^2}**@{},
    \end{xy}
 & \mbox{for}\ |I_1|=2,n=2 \\
0 & \mbox{otherwise.} \Ea \right.
$$
By Corollary~II.2(ii), we can set
$$
Q\left( \xy
 <0mm,0mm>*{\mbox{$\xy *=<10mm,3mm>\txt{\em 1}*\frm{-}\endxy$}};<0mm,0mm>*{}**@{},
 <-2.5mm,-1.6mm>*{};<-4mm,-6mm>*{}**@{-},
 <-1.6mm,-1.6mm>*{};<-2mm,-6mm>*{}**@{-},
 <1.6mm,-1.6mm>*{};<2.5mm,-6mm>*{}**@{-},
 <2.5mm,-1.6mm>*{};<4mm,-6mm>*{}**@{-},
 <0mm,0mm>*{};<0mm,-4.6mm>*{.\hspace{0mm}.\hspace{0mm}.}**@{},
<0mm,0mm>*{};<-4.1mm,-7.4mm>*{_1}**@{},
<0mm,0mm>*{};<-2mm,-7.4mm>*{_2}**@{},
%<0mm,0mm>*{};<-4mm,-7.4mm>*{_3}**@{},
<0mm,0mm>*{};<1.7mm,-7.4mm>*{...}**@{},
%<0mm,0mm>*{};<8.5mm,-7.4mm>*{_{n-1}}**@{},
<0mm,0mm>*{};<4.6mm,-7.4mm>*{_n}**@{},
<0mm,1.5mm>*{};<2mm,7mm>*{}**@{-},
  <0mm,1.5mm>*{};<1mm,6mm>*{}**@{-},
  <0mm,1.5mm>*{};<-0.5mm,7mm>*{}**@{-},
  <0mm,1.5mm>*{};<-2mm,7mm>*{}**@{-},
 <0mm,1.5mm>*{};<0.7mm,6.6mm>*{.\hspace{-0.4mm}.\hspace{-0.4mm}.}**@{},
 <-5mm,1.5mm>*{};<-8.5mm,7mm>*{}**@{-},
  <-5mm,1.5mm>*{};<-7.5mm,7mm>*{}**@{-},
  <-5mm,1.5mm>*{};<-6mm,6mm>*{}**@{-},
  <-5mm,1.5mm>*{};<-4.5mm,7mm>*{}**@{-},
 <-3.5mm,1.5mm>*{};<-6mm,6.6mm>*{.\hspace{-0.4mm}.\hspace{-0.4mm}.}**@{},
 <5mm,1.5mm>*{};<8.5mm,7mm>*{}**@{-},
  <5mm,1.5mm>*{};<7.5mm,7mm>*{}**@{-},
  <5mm,1.5mm>*{};<6mm,6mm>*{}**@{-},
  <5mm,1.5mm>*{};<4.5mm,7mm>*{}**@{-},
 <3.5mm,1.5mm>*{};<6mm,6.6mm>*{.\hspace{-0.4mm}.\hspace{-0.4mm}.}**@{},
\endxy
\right)=0
$$
completing thereby the $s=3$ iteration step.

\bip

%%%%%%%%%%%%%%%%%%%%%%%%%%%%%%%%%%%%%%%%%%%%%%%%%%%%% s=4 %%%%%%%%%%%%%%%%%%%%%%%%%%%

\noindent{\bf Iteration level s=4}: We have
\[
E_4=\mbox{\sf span} \left\{ \xy
 <0mm,0mm>*{\mbox{$\xy *=<5mm,3mm>\txt{\em 3}*\frm{-}\endxy$}};<0mm,0mm>*{}**@{},
 <-2.5mm,-1.6mm>*{};<-4mm,-6mm>*{}**@{-},
 <-1.6mm,-1.6mm>*{};<-2mm,-6mm>*{}**@{-},
 <1.6mm,-1.6mm>*{};<2.5mm,-6mm>*{}**@{-},
 <2.5mm,-1.6mm>*{};<4mm,-6mm>*{}**@{-},
 <0mm,0mm>*{};<0mm,-4.6mm>*{.\hspace{0mm}.\hspace{0mm}.}**@{},
<0mm,0mm>*{};<-4.1mm,-7.4mm>*{_1}**@{},
<0mm,0mm>*{};<-2mm,-7.4mm>*{_2}**@{},
%<0mm,0mm>*{};<-4mm,-7.4mm>*{_3}**@{},
<0mm,0mm>*{};<1.7mm,-7.4mm>*{...}**@{},
%<0mm,0mm>*{};<8.5mm,-7.4mm>*{_{n-1}}**@{},
<0mm,0mm>*{};<4.6mm,-7.4mm>*{_n}**@{},
\endxy
\ \ , \ \ \xy
 <0mm,0mm>*{\mbox{$\xy *=<5mm,3mm>\txt{\em 2}*\frm{-}\endxy$}};<0mm,0mm>*{}**@{},
 <-2.5mm,-1.6mm>*{};<-4mm,-6mm>*{}**@{-},
 <-1.6mm,-1.6mm>*{};<-2mm,-6mm>*{}**@{-},
 <1.6mm,-1.6mm>*{};<2.5mm,-6mm>*{}**@{-},
 <2.5mm,-1.6mm>*{};<4mm,-6mm>*{}**@{-},
 <0mm,0mm>*{};<0mm,-4.6mm>*{.\hspace{0mm}.\hspace{0mm}.}**@{},
<0mm,0mm>*{};<-4.1mm,-7.4mm>*{_1}**@{},
<0mm,0mm>*{};<-2mm,-7.4mm>*{_2}**@{},
%<0mm,0mm>*{};<-4mm,-7.4mm>*{_3}**@{},
<0mm,0mm>*{};<1.7mm,-7.4mm>*{...}**@{},
%<0mm,0mm>*{};<8.5mm,-7.4mm>*{_{n-1}}**@{},
<0mm,0mm>*{};<4.6mm,-7.4mm>*{_n}**@{},
 <-2.5mm,1.5mm>*{};<-5.5mm,7mm>*{}**@{-},
  <-2.5mm,1.5mm>*{};<-4.5mm,7mm>*{}**@{-},
  <-2.5mm,1.5mm>*{};<-3mm,6mm>*{}**@{-},
  <-2.5mm,1.5mm>*{};<-1.5mm,7mm>*{}**@{-},
 <-2.5mm,1.5mm>*{};<-3mm,6.6mm>*{.\hspace{-0.4mm}.\hspace{-0.4mm}.}**@{},
 <2.5mm,1.5mm>*{};<5.5mm,7mm>*{}**@{-},
  <2.5mm,1.5mm>*{};<4.5mm,7mm>*{}**@{-},
  <2.5mm,1.5mm>*{};<3mm,6mm>*{}**@{-},
  <2.5mm,1.5mm>*{};<1.5mm,7mm>*{}**@{-},
 <2.5mm,1.5mm>*{};<3mm,6.6mm>*{.\hspace{-0.4mm}.\hspace{-0.4mm}.}**@{},
\endxy
\ \ , \ \ \xy
 <0mm,0mm>*{\mbox{$\xy *=<10mm,3mm>\txt{\em 1}*\frm{-}\endxy$}};<0mm,0mm>*{}**@{},
 <-2.5mm,-1.6mm>*{};<-4mm,-6mm>*{}**@{-},
<-1.6mm,-1.6mm>*{};<-2.5mm,-6mm>*{}**@{-},
 <1.6mm,-1.6mm>*{};<2.5mm,-6mm>*{}**@{-},
 <2.5mm,-1.6mm>*{};<4mm,-6mm>*{}**@{-},
 <0mm,0mm>*{};<0mm,-4.6mm>*{.\hspace{0mm}.\hspace{0mm}.}**@{},
<0mm,0mm>*{};<-4.1mm,-7.4mm>*{_1}**@{},
<0mm,0mm>*{};<-2mm,-7.4mm>*{_2}**@{},
%<0mm,0mm>*{};<-4mm,-7.4mm>*{_3}**@{},
<0mm,0mm>*{};<1.7mm,-7.4mm>*{...}**@{},
%<0mm,0mm>*{};<8.5mm,-7.4mm>*{_{n-1}}**@{},
<0mm,0mm>*{};<4.6mm,-7.4mm>*{_n}**@{},
<2mm,1.5mm>*{};<3.9mm,7mm>*{}**@{-},
  <2mm,1.5mm>*{};<0.9mm,7mm>*{}**@{-},
  <-2mm,1.5mm>*{};<-3.9mm,7mm>*{}**@{-},
  <-2mm,1.5mm>*{};<-0.9mm,7mm>*{}**@{-},
 <0mm,1.5mm>*{};<2.4mm,6.6mm>*{.\hspace{-0.4mm}.\hspace{-0.4mm}.}**@{},
<2mm,1.5mm>*{};<2.4mm,6mm>*{}**@{-},
<0mm,1.5mm>*{};<-2.4mm,6.6mm>*{.\hspace{-0.4mm}.\hspace{-0.4mm}.}**@{},
<-2mm,1.5mm>*{};<-2.4mm,6mm>*{}**@{-},
 <-5mm,1.5mm>*{};<-8.5mm,7mm>*{}**@{-},
  <-5mm,1.5mm>*{};<-7mm,6mm>*{}**@{-},
  <-5mm,1.5mm>*{};<-6mm,7mm>*{}**@{-},
 <3.5mm,1.5mm>*{};<-7mm,6.6mm>*{.\hspace{-0.4mm}.\hspace{-0.4mm}.}**@{},
 <5mm,1.5mm>*{};<8.5mm,7mm>*{}**@{-},
  <5mm,1.5mm>*{};<7mm,6mm>*{}**@{-},
  <5mm,1.5mm>*{};<6mm,7mm>*{}**@{-},
 <3.5mm,1.5mm>*{};<7mm,6.6mm>*{.\hspace{-0.4mm}.\hspace{-0.4mm}.}**@{},
\endxy
\right\}
\]
By Corollary~II.1(i) the values of the morphism $Q$ on corollas of
the first type are completely determined by lifts of the values of
the morphism $q$, i.e.
$$
Q\left( \xy
 <0mm,0mm>*{\mbox{$\xy *=<5mm,3mm>\txt{\em 3}*\frm{-}\endxy$}};<0mm,0mm>*{}**@{},
 <-2.5mm,-1.6mm>*{};<-4mm,-6mm>*{}**@{-},
 <-1.6mm,-1.6mm>*{};<-2mm,-6mm>*{}**@{-},
 <1.6mm,-1.6mm>*{};<2.5mm,-6mm>*{}**@{-},
 <2.5mm,-1.6mm>*{};<4mm,-6mm>*{}**@{-},
 <0mm,0mm>*{};<0mm,-4.6mm>*{.\hspace{0mm}.\hspace{0mm}.}**@{},
<0mm,0mm>*{};<-4.1mm,-7.4mm>*{_1}**@{},
<0mm,0mm>*{};<-2mm,-7.4mm>*{_2}**@{},
%<0mm,0mm>*{};<-4mm,-7.4mm>*{_3}**@{},
<0mm,0mm>*{};<1.7mm,-7.4mm>*{...}**@{},
%<0mm,0mm>*{};<8.5mm,-7.4mm>*{_{n-1}}**@{},
<0mm,0mm>*{};<4.6mm,-7.4mm>*{_n}**@{},
\endxy
\right)= qis^{-1}\circ q\left( \xy
 <0mm,0mm>*{\mbox{$\xy *=<5mm,3mm>\txt{\em 3}*\frm{-}\endxy$}};<0mm,0mm>*{}**@{},
 <-2.5mm,-1.6mm>*{};<-4mm,-6mm>*{}**@{-},
 <-1.6mm,-1.6mm>*{};<-2mm,-6mm>*{}**@{-},
 <1.6mm,-1.6mm>*{};<2.5mm,-6mm>*{}**@{-},
 <2.5mm,-1.6mm>*{};<4mm,-6mm>*{}**@{-},
 <0mm,0mm>*{};<0mm,-4.6mm>*{.\hspace{0mm}.\hspace{0mm}.}**@{},
<0mm,0mm>*{};<-4.1mm,-7.4mm>*{_1}**@{},
<0mm,0mm>*{};<-2mm,-7.4mm>*{_2}**@{},
%<0mm,0mm>*{};<-4mm,-7.4mm>*{_3}**@{},
<0mm,0mm>*{};<1.7mm,-7.4mm>*{...}**@{},
%<0mm,0mm>*{};<8.5mm,-7.4mm>*{_{n-1}}**@{},
<0mm,0mm>*{};<4.6mm,-7.4mm>*{_n}**@{},
\endxy
\right)= \left\{ \Ba{rr} -\frac{1}{12} \Ba{c}
 \begin{xy}
 <0mm,0mm>*{\bullet};<0mm,0mm>*{}**@{},
 %<0mm,0.69mm>*{};<0mm,3.0mm>*{}**@{-},
 <0.39mm,-0.39mm>*{};<2.4mm,-2.4mm>*{}**@{-},
 <-0.35mm,-0.35mm>*{};<-1.9mm,-1.9mm>*{}**@{-},
 <-2.4mm,-2.4mm>*{\bullet};<-2.4mm,-2.4mm>*{}**@{},
 <-2.0mm,-2.8mm>*{};<-0.4mm,-4.5mm>*{}**@{-},
 %<-2.8mm,-2.9mm>*{};<-4.7mm,-4.9mm>*{}**@{-},
 %
 <2.4mm,-2.4mm>*{};<0.4mm,-4.5mm>*{}**@{-},
  <0mm,-5.1mm>*{\circ};<0mm,0mm>*{}**@{},
  <0mm,-5.5mm>*{};<0mm,-7.7mm>*{}**@{-},
(0,0)*{}
   \ar@{->}@(ul,dl) (-2.4,-2.4)*{}
 \end{xy}
\Ea
 & \mbox{for}\ n=1 \\
0 & \mbox{otherwise.} \Ea \right.
$$
It is also not hard to check that
$$
Q\left( \xy
 <0mm,0mm>*{\mbox{$\xy *=<5mm,3mm>\txt{\em 2}*\frm{-}\endxy$}};<0mm,0mm>*{}**@{},
 <-2.5mm,-1.6mm>*{};<-4mm,-6mm>*{}**@{-},
 <-1.6mm,-1.6mm>*{};<-2mm,-6mm>*{}**@{-},
 <1.6mm,-1.6mm>*{};<2.5mm,-6mm>*{}**@{-},
 <2.5mm,-1.6mm>*{};<4mm,-6mm>*{}**@{-},
 <0mm,0mm>*{};<0mm,-4.6mm>*{.\hspace{0mm}.\hspace{0mm}.}**@{},
<0mm,0mm>*{};<-4.1mm,-7.4mm>*{_1}**@{},
<0mm,0mm>*{};<-2mm,-7.4mm>*{_2}**@{},
%<0mm,0mm>*{};<-4mm,-7.4mm>*{_3}**@{},
<0mm,0mm>*{};<1.7mm,-7.4mm>*{...}**@{},
%<0mm,0mm>*{};<8.5mm,-7.4mm>*{_{n-1}}**@{},
<0mm,0mm>*{};<4.6mm,-7.4mm>*{_n}**@{},
 <-2.5mm,1.5mm>*{};<-5.5mm,7mm>*{}**@{-},
  <-2.5mm,1.5mm>*{};<-4.5mm,7mm>*{}**@{-},
  <-2.5mm,1.5mm>*{};<-3mm,6mm>*{}**@{-},
  <-2.5mm,1.5mm>*{};<-1.5mm,7mm>*{}**@{-},
 <-2.5mm,1.5mm>*{};<-3mm,6.6mm>*{.\hspace{-0.4mm}.\hspace{-0.4mm}.}**@{},
 <2.5mm,1.5mm>*{};<5.5mm,7mm>*{}**@{-},
  <2.5mm,1.5mm>*{};<4.5mm,7mm>*{}**@{-},
  <2.5mm,1.5mm>*{};<3mm,6mm>*{}**@{-},
  <2.5mm,1.5mm>*{};<1.5mm,7mm>*{}**@{-},
 <2.5mm,1.5mm>*{};<3mm,6.6mm>*{.\hspace{-0.4mm}.\hspace{-0.4mm}.}**@{},
\endxy
\right)
= \left\{ \Ba{rr} \frac{1}{2} \begin{xy}
 <1mm,0mm>*{\circ};<0mm,0mm>*{}**@{},
<-1mm,0mm>*{\circ};<0mm,0mm>*{}**@{},
 <0.6mm,0.4mm>*{};<-4mm,4mm>*{}**@{-},
<1.6mm,0.4mm>*{};<5mm,4mm>*{}**@{-},
<1mm,-0.4mm>*{};<1mm,-3mm>*{}**@{-},
 <-1.4mm,0.4mm>*{};<-6mm,4mm>*{}**@{-},
<-0.6mm,0.4mm>*{};<3mm,4mm>*{}**@{-},
<-1mm,-0.4mm>*{};<-1mm,-3mm>*{}**@{-},
 <0.6mm,0.4mm>*{};<-3.4mm,5.5mm>*{^{2)}}**@{},
<1.6mm,0.4mm>*{};<6mm,5.5mm>*{^{4)}}**@{},
 <0.6mm,0.4mm>*{};<3mm,5.5mm>*{^{(3}}**@{},
<1.6mm,0.4mm>*{};<-7mm,5.5mm>*{^{(1}}**@{},
\end{xy}
 & \mbox{for}\ |I_1|=2, |I_2|=2, n=2 \vspace{3mm} \\
\frac{1}{6}\begin{xy} <0mm,0mm>*{\circ};<0mm,0mm>*{}**@{},
<3mm,3mm>*{\circ};<0mm,0mm>*{}**@{},
 <0.6mm,0.4mm>*{};<2.5mm,2.5mm>*{}**@{-},
<-0.4mm,0.4mm>*{};<-4mm,6.5mm>*{}**@{-},
<0mm,-0.4mm>*{};<0mm,-3mm>*{}**@{-},
 <3.5mm,3.5mm>*{};<6.5mm,6.5mm>*{}**@{-},
 <2.5mm,3.5mm>*{};<-2mm,6.5mm>*{}**@{-},
 <3.5mm,3.5mm>*{};<7mm,7.4mm>*{^3}**@{-},
 <0.6mm,0.4mm>*{};<-4.6mm,7.4mm>*{^{(1}}**@{},
<1.6mm,0.4mm>*{};<-1.6mm,7.4mm>*{^{2)}}**@{},
\end{xy}
&  \mbox{for}\ |I_1|=2, |I_2|=1, n=1  \vspace{3mm}\\
\frac{1}{6}\begin{xy} <0mm,0mm>*{\circ};<0mm,0mm>*{}**@{},
<-3mm,3mm>*{\circ};<0mm,0mm>*{}**@{},
 <-0.6mm,0.4mm>*{};<-2.5mm,2.5mm>*{}**@{-},
<0.4mm,0.4mm>*{};<4mm,6.5mm>*{}**@{-},
<0mm,-0.4mm>*{};<0mm,-3mm>*{}**@{-},
 <-3.5mm,3.5mm>*{};<-6.5mm,6.5mm>*{}**@{-},
 <-2.5mm,3.5mm>*{};<2mm,6.5mm>*{}**@{-},
 <-3.5mm,3.5mm>*{};<-7mm,7.4mm>*{^1}**@{-},
 <-0.6mm,0.4mm>*{};<4.6mm,7.4mm>*{^{3)}}**@{},
<-1.6mm,0.4mm>*{};<1.6mm,7.4mm>*{^{(2}}**@{},
\end{xy}
 &  \mbox{for}\ |I_1|=1, |I_2|=2, n=1  \vspace{3mm}\\
\\
\begin{xy} <0mm,0mm>*{\circ};
 <0.4mm,0.4mm>*{};<2.5mm,2.5mm>*{}**@{-},
<-0.4mm,0.4mm>*{};<-2.5mm,2.5mm>*{}**@{-},
<-0.4mm,-0.4mm>*{};<-2.5mm,-2.5mm>*{}**@{-},
 <0.4mm,-0.4mm>*{};<2.5mm,-2.5mm>*{}**@{-},
\end{xy} \ \ \ \
&  \mbox{for}\ |I_1|=1, |I_2|=1, n=2  \vspace{3mm}\\
0 & \mbox{otherwise.} \Ea\right.
$$
Here round brackets stand for symmetrization of the labels.

\sip

\no
By Corollary~II.2(ii), we can set
$
Q\left( \xy
 <0mm,0mm>*{\mbox{$\xy *=<10mm,3mm>\txt{\em 1}*\frm{-}\endxy$}};<0mm,0mm>*{}**@{},
 <-2.5mm,-1.6mm>*{};<-4mm,-6mm>*{}**@{-},
<-1.6mm,-1.6mm>*{};<-2.5mm,-6mm>*{}**@{-},
 <1.6mm,-1.6mm>*{};<2.5mm,-6mm>*{}**@{-},
 <2.5mm,-1.6mm>*{};<4mm,-6mm>*{}**@{-},
 <0mm,0mm>*{};<0mm,-4.6mm>*{.\hspace{0mm}.\hspace{0mm}.}**@{},
<0mm,0mm>*{};<-4.1mm,-7.4mm>*{_1}**@{},
<0mm,0mm>*{};<-2mm,-7.4mm>*{_2}**@{},
%<0mm,0mm>*{};<-4mm,-7.4mm>*{_3}**@{},
<0mm,0mm>*{};<1.7mm,-7.4mm>*{...}**@{},
%<0mm,0mm>*{};<8.5mm,-7.4mm>*{_{n-1}}**@{},
<0mm,0mm>*{};<4.6mm,-7.4mm>*{_n}**@{},
<2mm,1.5mm>*{};<3.9mm,7mm>*{}**@{-},
  <2mm,1.5mm>*{};<0.9mm,7mm>*{}**@{-},
  <-2mm,1.5mm>*{};<-3.9mm,7mm>*{}**@{-},
  <-2mm,1.5mm>*{};<-0.9mm,7mm>*{}**@{-},
 <0mm,1.5mm>*{};<2.4mm,6.6mm>*{.\hspace{-0.4mm}.\hspace{-0.4mm}.}**@{},
<2mm,1.5mm>*{};<2.4mm,6mm>*{}**@{-},
<0mm,1.5mm>*{};<-2.4mm,6.6mm>*{.\hspace{-0.4mm}.\hspace{-0.4mm}.}**@{},
<-2mm,1.5mm>*{};<-2.4mm,6mm>*{}**@{-},
 <-5mm,1.5mm>*{};<-8.5mm,7mm>*{}**@{-},
  <-5mm,1.5mm>*{};<-7mm,6mm>*{}**@{-},
  <-5mm,1.5mm>*{};<-6mm,7mm>*{}**@{-},
 <3.5mm,1.5mm>*{};<-7mm,6.6mm>*{.\hspace{-0.4mm}.\hspace{-0.4mm}.}**@{},
 <5mm,1.5mm>*{};<8.5mm,7mm>*{}**@{-},
  <5mm,1.5mm>*{};<7mm,6mm>*{}**@{-},
  <5mm,1.5mm>*{};<6mm,7mm>*{}**@{-},
 <3.5mm,1.5mm>*{};<7mm,6.6mm>*{.\hspace{-0.4mm}.\hspace{-0.4mm}.}**@{},
\endxy
\right)=0
$
completing thereby the $s=4$ iteration step.

\bip

%%%%%%%%%%%%%%%%%%%%%%%%%%%%%%%%%%%%%%%%%%%%%%
\noindent{\bf Iteration level s=6}: as a last but least trivial
example we compute the value of the morphism $Q$ on the generator
$$
\xy
 <0mm,0mm>*{\mbox{$\xy *=<5mm,3mm>\txt{\em 4}*\frm{-}\endxy$}};<0mm,0mm>*{}**@{},
 <-2.5mm,-1.6mm>*{};<-4mm,-6mm>*{}**@{-},
 <-1.6mm,-1.6mm>*{};<-2mm,-6mm>*{}**@{-},
 <1.6mm,-1.6mm>*{};<2.5mm,-6mm>*{}**@{-},
 <2.5mm,-1.6mm>*{};<4mm,-6mm>*{}**@{-},
 <0mm,0mm>*{};<0mm,-4.6mm>*{.\hspace{0mm}.\hspace{0mm}.}**@{},
<0mm,0mm>*{};<-4.1mm,-7.4mm>*{_1}**@{},
<0mm,0mm>*{};<-2mm,-7.4mm>*{_2}**@{},
%<0mm,0mm>*{};<-4mm,-7.4mm>*{_3}**@{},
<0mm,0mm>*{};<1.7mm,-7.4mm>*{...}**@{},
%<0mm,0mm>*{};<8.5mm,-7.4mm>*{_{n-1}}**@{},
<0mm,0mm>*{};<4.6mm,-7.4mm>*{_n}**@{},
\endxy
\in E_6.
$$
As the value of the morphism $q$ on such a generator is zero, we
have
$$
Q\left( \xy
 <0mm,0mm>*{\mbox{$\xy *=<5mm,3mm>\txt{\em 4}*\frm{-}\endxy$}};<0mm,0mm>*{}**@{},
 <-2.5mm,-1.6mm>*{};<-4mm,-6mm>*{}**@{-},
 <-1.6mm,-1.6mm>*{};<-2mm,-6mm>*{}**@{-},
 <1.6mm,-1.6mm>*{};<2.5mm,-6mm>*{}**@{-},
 <2.5mm,-1.6mm>*{};<4mm,-6mm>*{}**@{-},
 <0mm,0mm>*{};<0mm,-4.6mm>*{.\hspace{0mm}.\hspace{0mm}.}**@{},
<0mm,0mm>*{};<-4.1mm,-7.4mm>*{_1}**@{},
<0mm,0mm>*{};<-2mm,-7.4mm>*{_2}**@{},
%<0mm,0mm>*{};<-4mm,-7.4mm>*{_3}**@{},
<0mm,0mm>*{};<1.7mm,-7.4mm>*{...}**@{},
%<0mm,0mm>*{};<8.5mm,-7.4mm>*{_{n-1}}**@{},
<0mm,0mm>*{};<4.6mm,-7.4mm>*{_n}**@{},
\endxy
\right)=e'',
$$
where $e''$ is a solution of the equation in
$[\Liebi^\circlearrowright]_\infty$, \Beqrn \delta e''=Q\left(\delta
\xy
 <0mm,0mm>*{\mbox{$\xy *=<5mm,3mm>\txt{\em 4}*\frm{-}\endxy$}};<0mm,0mm>*{}**@{},
 <-2.5mm,-1.6mm>*{};<-4mm,-6mm>*{}**@{-},
 <-1.6mm,-1.6mm>*{};<-2mm,-6mm>*{}**@{-},
 <1.6mm,-1.6mm>*{};<2.5mm,-6mm>*{}**@{-},
 <2.5mm,-1.6mm>*{};<4mm,-6mm>*{}**@{-},
 <0mm,0mm>*{};<0mm,-4.6mm>*{.\hspace{0mm}.\hspace{0mm}.}**@{},
<0mm,0mm>*{};<-4.1mm,-7.4mm>*{_1}**@{},
<0mm,0mm>*{};<-2mm,-7.4mm>*{_2}**@{},
%<0mm,0mm>*{};<-4mm,-7.4mm>*{_3}**@{},
<0mm,0mm>*{};<1.7mm,-7.4mm>*{...}**@{},
%<0mm,0mm>*{};<8.5mm,-7.4mm>*{_{n-1}}**@{},
<0mm,0mm>*{};<4.6mm,-7.4mm>*{_n}**@{},
\endxy
\right)&=&\sum_{b+c=4}\sum_{s\geq 0}\frac{1}{s!}\ \ \ \ \xy
 <0mm,0mm>*{\hspace{-3mm}Q(\,\mbox{$\xy *=<5mm,3mm>\txt{\em b}*\frm{-}\endxy$}\,)};<0mm,0mm>*{}**@{},
 <-2.5mm,-1.6mm>*{};<-4mm,-6mm>*{}**@{-},
 <-1.6mm,-1.6mm>*{};<-2mm,-6mm>*{}**@{-},
 <1.6mm,-1.6mm>*{};<2.5mm,-6mm>*{}**@{-},
 <2.5mm,-1.6mm>*{};<4mm,-6mm>*{}**@{-},
<0mm,0mm>*{};<3.2mm,4mm>*{_s}**@{},
 <0mm,0mm>*{};<0mm,-4.6mm>*{.\hspace{0mm}.\hspace{0mm}.}**@{},
<0mm,0mm>*{};<-4.1mm,-7.4mm>*{_1}**@{},
<0mm,0mm>*{};<-2mm,-7.4mm>*{_2}**@{},
%<0mm,0mm>*{};<-4mm,-7.4mm>*{_3}**@{},
<0mm,0mm>*{};<1.7mm,-7.4mm>*{...}**@{},
%<0mm,0mm>*{};<8.5mm,-7.4mm>*{_{n-1}}**@{},
<0mm,0mm>*{};<4.6mm,-7.4mm>*{_n}**@{},
<0mm,1.5mm>*{};<2mm,7mm>*{}**@{-},
  <0mm,1.5mm>*{};<1mm,6mm>*{}**@{-},
  <0mm,1.5mm>*{};<-0.5mm,7mm>*{}**@{-},
  <0mm,1.5mm>*{};<-2mm,7mm>*{}**@{-},
 <0mm,1.5mm>*{};<0.7mm,6.6mm>*{.\hspace{-0.4mm}.\hspace{-0.4mm}.}**@{},
<0mm,1.5mm>*{};<0mm,8.6mm>*{\hspace{-3mm}Q(\, \mbox{$\xy
*=<5mm,3mm>\txt{\em c}*\frm{-}\endxy$}\,)}**@{},
\endxy
\vspace{5mm}\\
&=&\left\{ -\frac{1}{12}\hspace{-4mm} \Ba{rr} \Ba{c}
 \begin{xy}
 <0mm,0mm>*{\bullet};<0mm,0mm>*{}**@{},
 %<0mm,0.69mm>*{};<0mm,3.0mm>*{}**@{-},
 <0.39mm,-0.39mm>*{};<2.4mm,-2.4mm>*{}**@{-},
 <-0.35mm,-0.35mm>*{};<-1.9mm,-1.9mm>*{}**@{-},
 <-2.4mm,-2.4mm>*{\bullet};<-2.4mm,-2.4mm>*{}**@{},
 <-2.0mm,-2.8mm>*{};<-0.4mm,-4.5mm>*{}**@{-},
 %<-2.8mm,-2.9mm>*{};<-4.7mm,-4.9mm>*{}**@{-},
 %
 <2.4mm,-2.4mm>*{};<0.4mm,-4.5mm>*{}**@{-},
  <0mm,-5.1mm>*{\circ};<0mm,0mm>*{}**@{},
  <0mm,-5.5mm>*{};<0mm,-7.7mm>*{}**@{-},
<0mm,-8.2mm>*{\bullet};<0mm,0mm>*{}**@{},
 <0mm,-8.2mm>*{};<-2.4mm,-10.1mm>*{}**@{-},
 <0mm,-8.2mm>*{};<2.4mm,-10.1mm>*{}**@{-},
(0,0)*{}
   \ar@{->}@(ul,dl) (-2.4,-2.4)*{}
 \end{xy}
 \\
 \ \
\Ea
 & \mbox{for}\  n=2 \\
0 & \mbox{otherwise.} \Ea \right. \Eeqrn It is not hard to solve the
latter for $e''$ and finally get
$$
Q\left( \xy
 <0mm,0mm>*{\mbox{$\xy *=<5mm,3mm>\txt{\em 4}*\frm{-}\endxy$}};<0mm,0mm>*{}**@{},
 <-2.5mm,-1.6mm>*{};<-4mm,-6mm>*{}**@{-},
 <-1.6mm,-1.6mm>*{};<-2mm,-6mm>*{}**@{-},
 <1.6mm,-1.6mm>*{};<2.5mm,-6mm>*{}**@{-},
 <2.5mm,-1.6mm>*{};<4mm,-6mm>*{}**@{-},
 <0mm,0mm>*{};<0mm,-4.6mm>*{.\hspace{0mm}.\hspace{0mm}.}**@{},
<0mm,0mm>*{};<-4.1mm,-7.4mm>*{_1}**@{},
<0mm,0mm>*{};<-2mm,-7.4mm>*{_2}**@{},
%<0mm,0mm>*{};<-4mm,-7.4mm>*{_3}**@{},
<0mm,0mm>*{};<1.7mm,-7.4mm>*{...}**@{},
%<0mm,0mm>*{};<8.5mm,-7.4mm>*{_{n-1}}**@{},
<0mm,0mm>*{};<4.6mm,-7.4mm>*{_n}**@{},
\endxy
\right)= \left\{ \Ba{rr} -\frac{1}{12}\hspace{-3mm} \Ba{c}
 \begin{xy}
 <0mm,0mm>*{\bullet};<0mm,0mm>*{}**@{},
 %<0mm,0.69mm>*{};<0mm,3.0mm>*{}**@{-},
 <0.39mm,-0.39mm>*{};<2.4mm,-2.4mm>*{}**@{-},
 <-0.35mm,-0.35mm>*{};<-1.9mm,-1.9mm>*{}**@{-},
 <-2.4mm,-2.4mm>*{\bullet};<-2.4mm,-2.4mm>*{}**@{},
 <-2.0mm,-2.8mm>*{};<-0.4mm,-4.5mm>*{}**@{-},
 %<-2.8mm,-2.9mm>*{};<-4.7mm,-4.9mm>*{}**@{-},
 %
 <2.4mm,-2.4mm>*{};<0.4mm,-4.5mm>*{}**@{-},
  <0mm,-5.1mm>*{\circ};<0mm,0mm>*{}**@{},
  <0.4mm,-5.5mm>*{};<2mm,-7.7mm>*{}**@{-},
<-0.4mm,-5.5mm>*{};<-2mm,-7.7mm>*{}**@{-},
  <0.4mm,-5.5mm>*{};<2.9mm,-9.7mm>*{^2}**@{},
 <0.4mm,-5.5mm>*{};<-2.9mm,-9.7mm>*{^1}**@{},
(0,0)*{}
   \ar@{->}@(ul,dl) (-2.4,-2.4)*{}
 \end{xy}
 \\
 \ \
\Ea - \frac{1}{6}\hspace{-2mm} \Ba{c}
 \begin{xy}
 <0mm,0mm>*{\bullet};<0mm,0mm>*{}**@{},
 %<0mm,0.69mm>*{};<0mm,3.0mm>*{}**@{-},
 <0mm,-0mm>*{};<1.8mm,-2.4mm>*{}**@{-},
 <0mm,0mm>*{};<-1.9mm,-1.9mm>*{}**@{-},
 <-2.4mm,-2.4mm>*{\bullet};<-2.4mm,-2.4mm>*{}**@{},
 <-2.0mm,-2.8mm>*{};<-0.4mm,-4.5mm>*{}**@{-},
 %<-2.8mm,-2.9mm>*{};<-4.7mm,-4.9mm>*{}**@{-},
<0mm,0mm>*{};<5.4mm,-2.4mm>*{}**@{-},
 <1.8mm,-2.4mm>*{};<0.4mm,-4.5mm>*{}**@{-},
  <0mm,-5.1mm>*{\circ};<0mm,0mm>*{}**@{},
  <0mm,-5.5mm>*{};<0mm,-7.7mm>*{}**@{-},
 <0.4mm,-5.5mm>*{};<6mm,-5mm>*{^2}**@{},
 <0.4mm,-5.5mm>*{};<0mm,-10.7mm>*{^1}**@{},
(0,0)*{}
   \ar@{->}@(ul,dl) (-2.4,-2.4)*{}
 \end{xy}
 \\
 \ \
\Ea - \frac{1}{6}\hspace{-2mm} \Ba{c}
 \begin{xy}
 <0mm,0mm>*{\bullet};<0mm,0mm>*{}**@{},
 %<0mm,0.69mm>*{};<0mm,3.0mm>*{}**@{-},
 <0mm,-0mm>*{};<1.8mm,-2.4mm>*{}**@{-},
 <0mm,0mm>*{};<-1.9mm,-1.9mm>*{}**@{-},
 <-2.4mm,-2.4mm>*{\bullet};<-2.4mm,-2.4mm>*{}**@{},
 <-2.0mm,-2.8mm>*{};<-0.4mm,-4.5mm>*{}**@{-},
 %<-2.8mm,-2.9mm>*{};<-4.7mm,-4.9mm>*{}**@{-},
<0mm,0mm>*{};<5.4mm,-2.4mm>*{}**@{-},
 <1.8mm,-2.4mm>*{};<0.4mm,-4.5mm>*{}**@{-},
  <0mm,-5.1mm>*{\circ};<0mm,0mm>*{}**@{},
  <0mm,-5.5mm>*{};<0mm,-7.7mm>*{}**@{-},
 <0.4mm,-5.5mm>*{};<6mm,-5mm>*{^1}**@{},
 <0.4mm,-5.5mm>*{};<0mm,-10.7mm>*{^2}**@{},
(0,0)*{}
   \ar@{->}@(ul,dl) (-2.4,-2.4)*{}
 \end{xy}
 \\
 \ \
\Ea
 & \mbox{for}\  n=2 \\
0 & \mbox{otherwise.} \Ea \right.
$$

%%%%%%%%%%%%%%%%%%%%%%%%%%%%%%%%%%%%%%%%%%%%%%%%%%%%%%%%%%%%%%%%%%%%%%%%%%

\noindent{\bf Step III.} One can show by induction on the parameter $s=2a+k-2$ (we omit these details) that the map
$Q$ can be chosen so that the composition  $Q\circ
\chi_{\hbar=1}$ (see \S 2.9.2) makes sense as a morphism into a completion of $[\Liebi^\circlearrowright]_\infty$ with respect to the number of
vertices.  We finally set $\hat{Q}:= Q\circ
\chi_{\hbar=1}$ completing the construction. \hfill  $\Box$

\bip

\no
{\bf 5.3. Remark.} There exists a canonical monomorphism of dg
props, $i:\Liebi^\circlearrowright_\infty\rar
[\Liebi^\circlearrowright]_\infty$. Hence any morphism $\hat{Q}$ which factors through $i$ gives rise
to a universal quantization of ordinary Poisson structures.

\bip

\no
{\bf 5.4. Remark.} The condition on the projection $\pi_1\circ\al\circ \hat{Q}$ in Theorem 5.2
implies high non-triviality of the quantization morphism $\hat{Q}$ in the sense that $\hat{Q}$ involves
{\em all}\, possible jets of the polyvector part of the input wheeled Poisson structure.

%%%%%%%%%%%%%%%%%%%%% cyclic complex %%%%%%%%%%%%%%%%%%%%%%%%%%%%%%%%%%%%

\bip

\bip

\begin{center}
{\bf Appendix}
\end{center}

\bip

%%%%%%%%%%%%%%%%%%%%%%%%%
{\bf A.1. Genus 1 wheels.}
Let $(\PROP\langle E\rangle, \delta)$ be a dg free prop,
and let $(\PROP^\circlearrowright\langle E\rangle, \delta)$ be its wheeled extension.
We assume in this section that the differential $\delta$ preserves the number of
wheels\footnote{This is not that dramatic loss
of generality in the sense there always exits a filtration of
$({\mathsf P}^\circlearrowright\langle E\rangle, \delta)$ by the number of wheels
whose spectral sequence has zero-th term satisfying our condition on the differential.}.
Then it makes sense to define
 a subcomplex, $({\mathsf T}^\circlearrowright\langle E\rangle \subset \PROP^\circlearrowright\langle E\rangle$,
spanned by graphs with precisely one wheel. In this section we use the ideas of cyclic homology to define
a new cyclic bicomplex which computes cohomology of $({\mathsf T}^\circlearrowright\langle E\rangle, \delta)$.

\sip

All the above assumptions are satisfied automatically  if $(\PROP\langle E\rangle, \delta)$
is the free dg prop associated with a free dg operad.

\sip

We denote by  ${\mathsf T}^+\langle E\rangle$ the obvious ``marked wheel" extension of
${\mathsf T}^\circlearrowright\langle E\rangle$ (see \S 3.2).

\bip

{\bf A.2. Abbreviated notations for graphs in %${\mathsf T}^\circlearrowright\langle E\rangle$
 ${\mathsf T}^+\langle E\rangle$.}
The half-edges attached to any
internal vertex of split into, say $m$, ingoing
and, say $n$, outgoing ones. The differential $\delta$ is uniquely determined by its values
on such $(m,n)$-vertices for all possible $m,n\geq 1$. If the vertex is
cyclic, then one of its input half-edges is cyclic and one of its output
half-edges is also cyclic. In this section we show in pictures only
those (half-)edges attached to vertices which are
cyclic (unless otherwise is explicitly stated), so that
\Bi
\item[-]
$\xy
 <0mm,0mm>*{\mbox{$\xy *=<4mm,3mm>
\txt{{{e}}}*\frm{-}\endxy$}};<0mm,0mm>*{}**@{}\endxy$\ \
stands for a non-cyclic $(m,n)$-vertex decorated by an element $e\in E(m,n)$,
\item[-]
 $\xy <0mm,0mm>*{\mbox{$\xy *=<4mm,4mm>
\txt{{{e}}}*\frm{-}\endxy$}};<0mm,0mm>*{}**@{},
  <0mm,2mm>*{};<0mm,4mm>*{}**@{-},
  <0mm,-2mm>*{};<0mm,-4mm>*{}**@{-},\endxy$ \ \
is a decorated cyclic $(m,n)$-vertex with no input or output cyclic
half-edges marked,
\item[-]
 $\xy <0mm,0mm>*{\mbox{$\xy *=<4mm,4mm>
\txt{{{e}}}*\frm{-}\endxy$}};<0mm,0mm>*{}**@{},
  <0mm,2mm>*{};<0mm,4mm>*{}**@{.},
  <0mm,-2mm>*{};<0mm,-4mm>*{}**@{-}\endxy$\ \
is a decorated cyclic $(m,n)$-vertex with the output cyclic half-edge marked,
\item[-]
$\xy <0mm,0mm>*{\mbox{$\xy *=<4mm,4mm>
\txt{{{e}}}*\frm{-}\endxy$}};<0mm,0mm>*{}**@{},
  <0mm,2mm>*{};<0mm,4mm>*{}**@{-},
  <0mm,-2.2mm>*{};<0mm,-4mm>*{}**@{.}\endxy$\ \
is a decorated cyclic $(m,n)$-vertex with the input cyclic half-edge marked.
\Ei
\bip

The differential $\delta$ applied to any vertex of the last three
types can be uniquely decomposed into the sum of the
following three groups of terms,
$$
\delta\ \xy
 <0mm,0mm>*{\mbox{$\xy *=<4mm,4mm>
\txt{{{e}}}*\frm{-}\endxy$}};<0mm,0mm>*{}**@{},
  <0mm,2mm>*{};<0mm,4mm>*{}**@{-},
  <0mm,-2.2mm>*{};<0mm,-4mm>*{}**@{.},\endxy =
  \sum_{\al\in I_1}\
\xy <0mm,0mm>*{\mbox{$\xy *=<6mm,4mm>
\txt{{{e$_{a'}$}}}*\frm{-}\endxy$}};<0mm,0mm>*{}**@{},
  <0mm,2mm>*{};<0mm,4mm>*{}**@{-},
  <0mm,-2.2mm>*{};<0mm,-4mm>*{}**@{.},
<0mm,6.2mm>*{\mbox{$\xy *=<6mm,4mm>
\txt{{{e$_{a''}$}}}*\frm{-}\endxy$}};<0mm,6mm>*{}**@{},
<0mm,8.2mm>*{};<0mm,10.2mm>*{}**@{-},
  \endxy\
  + \
\sum_{a\in I_2}\
  \xy <0mm,0mm>*{\mbox{$\xy *=<6mm,4mm>
\txt{{{e$_{a'}$}}}*\frm{-}\endxy$}};<0mm,0mm>*{}**@{},
   <0mm,2mm>*{};<0mm,4mm>*{}**@{-},
  <2.7mm,2mm>*{};<5mm,4mm>*{}**@{-},
  <0mm,-2.2mm>*{};<0mm,-4mm>*{}**@{.},
<6mm,6.2mm>*{\mbox{$\xy *=<6mm,4mm>
\txt{{{e$_{a''}$}}}*\frm{-}\endxy$}};<0mm,6mm>*{}**@{},
  \endxy
  + \
\sum_{b\in I_3}\
  \xy <0mm,0mm>*{\mbox{$\xy *=<6mm,4mm>
\txt{{{e$_{b'}$}}}*\frm{-}\endxy$}};<0mm,0mm>*{}**@{},
   <0mm,2mm>*{};<0mm,4mm>*{}**@{-},
  <2.7mm,-2.2mm>*{};<5mm,-4.2mm>*{}**@{-},
  <0mm,-2.2mm>*{};<0mm,-4mm>*{}**@{.},
<6mm,-6.2mm>*{\mbox{$\xy *=<6mm,4mm>
\txt{{{e$_{b''}$}}}*\frm{-}\endxy$}};<0mm,6mm>*{}**@{},
  \endxy
$$
where we have shown also non-cyclic {\em internal}\, edges in the last two
groups of terms. The differential $\delta$ applied to\ \ $\xy
 <0mm,0mm>*{\mbox{$\xy *=<4mm,4mm>
\txt{{{e}}}*\frm{-}\endxy$}};<0mm,0mm>*{}**@{},
  <0mm,2mm>*{};<0mm,4mm>*{}**@{-},
  <0mm,-2mm>*{};<0mm,-4mm>*{}**@{-},\endxy$\ \
  and\ $\xy
 <0mm,0mm>*{\mbox{$\xy *=<4mm,4mm>
\txt{{{e}}}*\frm{-}\endxy$}};<0mm,0mm>*{}**@{},
  <0mm,2mm>*{};<0mm,4mm>*{}**@{.},
  <0mm,-2mm>*{};<0mm,-4mm>*{}**@{-},\endxy$\ \
is given by exactly the same formula except for the presence/position
of dashed markings.

\bip

{\bf A.3. New differential in ${\mathsf T}^+\langle E\rangle$.}
 Let us define a new derivation, $b$,
in $\sT^+\langle E\rangle$, as follows:
\Bi
\item[-]
$b\ \xy
 <0mm,0mm>*{\mbox{$\xy *=<4mm,3mm>
\txt{{{e}}}*\frm{-}\endxy$}};<0mm,0mm>*{}**@{}\endxy
\, := \,
\delta\ \xy
 <0mm,0mm>*{\mbox{$\xy *=<4mm,3mm>
\txt{{{e}}}*\frm{-}\endxy$}};<0mm,0mm>*{}**@{}\endxy
$\ ,
\item[-]
 $b\ \xy <0mm,0mm>*{\mbox{$\xy *=<4mm,4mm>
\txt{{{e}}}*\frm{-}\endxy$}};<0mm,0mm>*{}**@{},
  <0mm,2mm>*{};<0mm,4mm>*{}**@{-},
  <0mm,-2mm>*{};<0mm,-4mm>*{}**@{-},\endxy
\ := \ \delta\
\xy <0mm,0mm>*{\mbox{$\xy *=<4mm,4mm>
\txt{{{e}}}*\frm{-}\endxy$}};<0mm,0mm>*{}**@{},
  <0mm,2mm>*{};<0mm,4mm>*{}**@{-},
  <0mm,-2mm>*{};<0mm,-4mm>*{}**@{-},\endxy$\ ,
\item[-]
 $b\ \xy <0mm,0mm>*{\mbox{$\xy *=<4mm,4mm>
\txt{{{e}}}*\frm{-}\endxy$}};<0mm,0mm>*{}**@{},
  <0mm,2mm>*{};<0mm,4mm>*{}**@{.},
  <0mm,-2mm>*{};<0mm,-4mm>*{}**@{-}\endxy
  \ = \
\delta\ \xy <0mm,0mm>*{\mbox{$\xy *=<4mm,4mm>
\txt{{{e}}}*\frm{-}\endxy$}};<0mm,0mm>*{}**@{},
  <0mm,2mm>*{};<0mm,4mm>*{}**@{.},
  <0mm,-2mm>*{};<0mm,-4mm>*{}**@{-}\endxy
  $\ ,
\item[-]
$b\ \xy <0mm,0mm>*{\mbox{$\xy *=<4mm,4mm>
\txt{{{e}}}*\frm{-}\endxy$}};<0mm,0mm>*{}**@{},
  <0mm,2mm>*{};<0mm,4mm>*{}**@{-},
  <0mm,-2.2mm>*{};<0mm,-4mm>*{}**@{.}\endxy
  \ : = \
\delta\ \xy <0mm,0mm>*{\mbox{$\xy *=<4mm,4mm>
\txt{{{e}}}*\frm{-}\endxy$}};<0mm,0mm>*{}**@{},
  <0mm,2mm>*{};<0mm,4mm>*{}**@{-},
  <0mm,-2.2mm>*{};<0mm,-4mm>*{}**@{.}\endxy \ \ + \ \
  \sum_{\al\in I_1}\
\xy <0mm,-0.3mm>*{\mbox{$\xy *=<6mm,4mm>
\txt{{{e$_{a'}$}}}*\frm{-}\endxy$}};<0mm,0mm>*{}**@{},
  <0mm,2mm>*{};<0mm,4mm>*{}**@{.},
  <0mm,-2.4mm>*{};<0mm,-4.2mm>*{}**@{-},
<0mm,6.3mm>*{\mbox{$\xy *=<6mm,4mm>
\txt{{{e$_{a''}$}}}*\frm{-}\endxy$}};<0mm,6mm>*{}**@{},
<0mm,8.2mm>*{};<0mm,10.2mm>*{}**@{-},
  \endxy\
$\ .
\Ei

\bip

%%%%%%%%%%%%%%%%%%%%%%%%%%%%%%%%%%%%%%%%%%%%%%%%%%%%%%%%
{\bf A.3.1. Lemma.} {\em The derivation $b$ satisfies
$b^2=0$, i.e.\ $(\sT^+\langle E\rangle, b)$ is a complex}.

\sip

Proof is a straightforward but tedious calculation based solely on the
relation $\delta^2=0$.

\bip

{\bf A.4. Action of cyclic groups.} The vector space
$\sT^+\langle E\rangle$ is naturally bigraded,
$$
\sT^+\langle E\rangle= \sum_{m\geq 0, n\geq 1}
\sT^+\langle E\rangle_{m, n},
$$
where the summand  $\sT^+\langle E\rangle_{m,n}$ consists of all
graphs with $m$ non-cyclic and $n$ cyclic vertices. Note that
 $\sT^+\langle E\rangle_{m,n}$ is naturally a representation
space of the cyclic group $\Z_n$ whose generator, $t$, moves the mark
to the next cyclic edge along the orientation. Define also the operator,
$N:=1+t+\ldots+ t^n: \sT^+\langle E\rangle_{m,n}\rar \sT^+\langle
E\rangle_{m,n}$, which symmetrizes the marked graphs.

\bip

%%%%%%%%%%%%%%%%%%%%%%%%%%%%%%%%%%%%%%%%%%%%%%%%%%%%%%%
{\bf A.4.1. Lemma.} %{\em One has in  $\PPP\langle E\rangle^1$
$\delta (1-t)=(1-t)b$ {\em and }\,
$N\delta= bN$.

\sip
Proof is a straightforward calculation based on the definition of $b$.

\bip

Following the ideas of the theory of cyclic homology (see, e.g.,
\cite{Lo}) we introduce a 4th quadrant bicomplex,
$$
C_{p,q}:= C_q, \ \ \ \ C_q:= \sum_{m+n=q}
\sT^+\langle E\rangle_{m, n}, \ \ \ \  p\leq 0, q\geq 1,
$$
with the differentials given by the following diagram,
$$
\Ba{ccccccccc}
 && \ldots && \ldots
&& \ldots && \ldots\\
 && \Big\uparrow\vcenter{\rlap{$b$}} &&
\Big\uparrow\vcenter{\rlap{$\delta$}}
&& \Big\uparrow\vcenter{\rlap{$b$}} && \Big\uparrow\vcenter{\rlap{$\delta$}}\\
\ldots &\stackrel{N}{\lon}& C_4 &\stackrel{1-t}{\lon}& C_4
&\stackrel{N}{\lon}& C_4 &\stackrel{1-t}{\lon}& C_4\\
 && \Big\uparrow\vcenter{\rlap{$b$}} &&
\Big\uparrow\vcenter{\rlap{$\delta$}}
&& \Big\uparrow\vcenter{\rlap{$b$}} && \Big\uparrow\vcenter{\rlap{$\delta$}}\\
\ldots &\stackrel{N}{\lon}& C_3 &\stackrel{1-t}{\lon}& C_3
&\stackrel{N}{\lon}& C_3 &\stackrel{1-t}{\lon}& C_3\\
 && \Big\uparrow\vcenter{\rlap{$b$}} &&
\Big\uparrow\vcenter{\rlap{$\delta$}}
&& \Big\uparrow\vcenter{\rlap{$b$}} && \Big\uparrow\vcenter{\rlap{$\delta$}}\\
\ldots &\stackrel{N}{\lon}& C_2 &\stackrel{1-t}{\lon}& C_2
&\stackrel{N}{\lon}& C_2 &\stackrel{1-t}{\lon}& C_2\\
 && \Big\uparrow\vcenter{\rlap{$b$}} &&
\Big\uparrow\vcenter{\rlap{$\delta$}}
&& \Big\uparrow\vcenter{\rlap{$b$}} && \Big\uparrow\vcenter{\rlap{$\delta$}}\\
\ldots &\stackrel{N}{\lon}& C_1 &\stackrel{1-t}{\lon}& C_1
&\stackrel{N}{\lon}& C_1 &\stackrel{1-t}{\lon}& C_1\\
\Ea
$$

%%%%%%%%%%%%%%%%%%%%%%%%%%%%%%%%%%%%%%%%%%%%%%%%%%%%%%%
{\bf A.4.2. Theorem.} {\em The cohomology group of the unmarked
graph complex, $H(\sT^\circlearrowright\langle E \rangle,\delta)$,
% where
%$\sT_1^\circlearrowright\langle E \rangle$ is the subcomplex of $\sT^\circlearrowright\langle E \rangle$
%consisting of graphs with only one wheel,
 is equal to the cohomology of the total complex associated with
 the cyclic bicomplex $C_{\bullet,\bullet}$}.

\bip

\Proof The complex $(\sT^\circlearrowright\langle E \rangle,\delta)$ can be identified
with the cokernel, $C_\bullet/(1-t)$, of the endomorphism $(1-t)$
of the total complex, $C_{\bullet}$, associated with the bicomplex
$C_{\bullet,\bullet}$. As the rows of $C_{\bullet,\bullet}$
are exact \cite{Lo}, the claim follows.
{\hfill $\Box$}

\bip

%%%%%%%%%%%%%%%%%%%%%%%%%%%%%%%%%%%%%%%%%%%%%%%%%%%%%%%
%{\bf 4.4.3. Facts.} {(i)}  ${\mathsf L}^+\langle E \rangle:=(1-t)\sT^+\langle E \rangle$ is a
%subcomplex of  $(\sT^+\langle E \rangle,\delta)$.

%\sip

%\noindent{ (ii)} There is a short exact sequence of complexes,
%$$
%0\lon ({\mathsf L}^+_1\langle E \rangle^1, \delta) \lon (\sT^+_1\langle E \rangle^1,\delta)
% \stackrel{p}{\lon} (\sT^\circlearrowright_1\langle E \rangle,\delta) \lon 0
%$$
%where $p$ is the surjection which forgets the markings.

%\sip

%\noindent{(iii)} $(N\sT^+\langle E \rangle, b)$ is a subcomplex
%of  $(\sT^+\langle E \rangle, b)$ whose cohomology is isomorphic
%to $H(\sT^\circlearrowright\langle E \rangle,\delta)$.

\bip

\bip

 \noindent{\small
{\em Acknowledgement.} It is a pleasure to thank Sergei Shadrin and
Bruno Vallette for helpful discussions.

\bip


\newpage

  \begin{thebibliography}{99}

\bibitem[Ag]{Ag} M.\ Aguiar, {\em Infinitesimal Hopf Algebras}, In: New trends in Hopf algebra theory,
vol.\ 267, Contemporary Math. (2000), 1-29.
\vspace{-1mm}



\bibitem[BFFLS]{BF}  F.\ Bayen, M.\ Flato, C.\ Fronsdal,
A.\ Lichnerowics, and D.\ Sternheimer, {\em Deformation theory and
quantization. I. Deformations of symplectic structures.}
Ann.\ Phys. {\bf 111} (1978), 61-110.
\vspace{-1mm}

\bibitem[Ga]{G} W.L.\ Gan, {\em Koszul duality for dioperads}, Math.\ Res.\ Lett.
{\bf 10} (2003), 109-124.
\vspace{-1mm}


\bibitem[Ko1]{Kon}
M.~Kontsevich.
\newblock {\em Formal (non)commutative symplectic geometry}.
\newblock In: { The Gel'fand
mathematics seminars 1990--1992}, Birkh\"auser,
  1993.
 \vspace{-1mm}

 \bibitem[Ko2]{Ko} M.\ Kontsevich, {\em Deformation quantization
 of Poisson manifolds I},
 math/9709040.
 \vspace{-1mm}

 \bibitem[Ko3]{Ko2} M.\ Kontsevich, {letter to Martin Markl}, November 2002.
 \vspace{-1mm}


\bibitem[Lo]{Lo} J.L.\ Loday, { Cyclic homology}, Springer 1998.
 \vspace{-1mm}



\bibitem[MMS]{MMS} M.\ Markl, S.\ Merkulov and S.\ Shadrin, {Wheeled PROPs and the master
equation}, math.AG/0610683.
 \vspace{-1mm}



 \bibitem[MaVo]{MV} M.\ Markl and A.A.\ Voronov,
{\em PROPped up graph cohomology}, math.QA/0307081.
 \vspace{-1mm}



\bibitem[Mc]{Mc} S.\ McLane, {\em Categorical algebra},
Bull.\ Amer.\ Math.\ Soc. {\bf 71} (1965), 40-106.
\vspace{-1mm}

 %\bibitem[Me1]{Me} S.A.\ Merkulov, {\em Operads, deformation theory and
 %$F$-manifolds
 %}, math.AG/0210478.
 %\vspace{-1mm}


\bibitem[Me1]{Me1} S.A.\ Merkulov, {\em PROP profile of Poisson geometry},
preprint math.DG/0401034, Commun.Math.Phys. {\bf 262} (2006), 117-135.
 \vspace{-1mm}

\bibitem[Me2]{Me2} S.A.\ Merkulov, {\em Nijenhuis infinity and contractible
dg manifolds}, math.AG/0403244, Compositio Mathematica {\bf 141} (2005), 1238-1254.
\vspace{-1mm}


\bibitem[Ta]{Ta} D.E.\ Tamarkin, {\em Another proof of M. Kontsevich formality
theorem}, math.QA/9803025.
\vspace{-1mm}

\bibitem[Va]{V} B.\ Vallette, {\em Dualit$\acute{e}$ de Koszul des props},
PhD thesis, Universit$\acute{e}$ Louis Pasteur, 2003, math.QA/0405057.
English translation:
{\em A Koszul duality for props}, math.AT/0411542.
\vspace{-1mm}



 %\bibitem[St]{St} J.D.\ Stasheff, {\em On the homotopy
  %associativity of $H$-spaces, I \@ II}, Trans.\ Amer.\
  %Math.\ Soc.\ {\bf 108} (1963), 272-292 \& 293-312.
 %\vspace{-1mm}


\bip \bip

\begin{tabular}{l}

Matematiska institutionen,
 Stockholms universitet\\
 10691 Stockholm,
 Sweden\\
sm@@math.su.se
 \end{tabular}



  \end{thebibliography}
  \end{document}